\documentclass[10pt,leqno]{article}

\usepackage[francais]{babel}
\usepackage{amsmath,amssymb,amscd}
\usepackage{a4wide}
\usepackage[latin1]{inputenc}

\usepackage[all]{xy} 
\usepackage[T1]{fontenc}

\title{Le transfert singulier pour la formule des traces \\ de Jacquet-Rallis}
\author{Pierre-Henri Chaudouard et Micha\l\  Zydor}
\date{}

\makeatletter
\@addtoreset{equation}{subsection}         
\makeatother

\newenvironment{paragr}[1][]{\refstepcounter{subsubsection} \noindent \textbf{\thesubsubsection . \ #1}}{\medskip}

\newenvironment{theoreme}{ \medskip\refstepcounter{theo}  \noindent\textbf{Th\'eor\`eme \thetheo}. ---\em}{\em \medskip}
\newenvironment{proposition}{\medskip\refstepcounter{theo}   \noindent\textbf{Proposition \thetheo}. ---\em}{\em\medskip}
\newenvironment{corollaire}{\medskip\refstepcounter{theo}  \noindent\textbf{Corollaire \thetheo}. ---\em}{\em\medskip}

\newenvironment{lemme}{\medskip\refstepcounter{theo}   \noindent\textbf{Lemme \thetheo}. ---\em}{\em\medskip}

\newenvironment{preuve}[1][]{\noindent \textbf{Démonstration.} #1 --- }{\hfill
  \ensuremath{\square} \medskip}

\newenvironment{remarque}{\medskip\refstepcounter{theo}  \noindent\textbf{Remarque \thetheo}. ---}{\medskip}

\DeclareMathOperator{\Orb}{Orb}
\DeclareMathOperator{\rs}{rss}
\DeclareMathOperator{\Grs}{\mathit{G}\text{-}rss}
\DeclareMathOperator{\Urs}{\mathit{U}\text{-}rss}

\DeclareMathOperator{\Lie}{Lie}

\DeclareMathOperator{\vol}{vol}

\DeclareMathOperator{\End}{End}

\DeclareMathOperator{\supp}{supp}

\DeclareMathOperator{\Norm}{Norm}

\DeclareMathOperator{\Gal}{Gal}
\DeclareMathOperator{\Hom}{Hom}

\DeclareMathOperator{\Id}{Id}

\DeclareMathOperator{\Res}{Res}

\DeclareMathOperator{\disc}{disc}
\DeclareMathOperator{\Ker}{Ker}

\DeclareMathOperator{\vect}{vect}

\DeclareMathOperator{\trace}{trace}

\newcommand{\ZZ}{\mathbb{Z}}

\newcommand{\NN}{\mathbb{N}}
\newcommand{\RR}{\mathbb{R}}
\newcommand{\AAA}{\mathbb{A}}
\newcommand{\CC}{\mathbb{C}}

\newcommand{\QQ}{\mathbb{Q}}

\newcommand{\oc}{\mathcal{O}}

\newcommand{\Sc}{\mathcal{S}}

\newcommand{\vc}{\mathcal{V}}

\newcommand{\yc}{\mathcal{Y}}

\newcommand{\fc}{\mathcal{F}}
\newcommand{\pc}{\mathcal{P}}
\newcommand{\hc}{\mathcal{H}}
\newcommand{\nc}{\mathcal{N}}
\newcommand{\Ac}{\mathcal{A}}

\newcommand{\Ic}{\mathcal{I}}

\newcommand{\xc}{\mathcal{X}}

\newcommand{\MM}{\mathbf{M}}

\newcommand{\tggo}{\widetilde{\mathfrak{g}}}
\newcommand{\thgo}{\widetilde{\mathfrak{h}}}
\newcommand{\tsgo}{\widetilde{\mathfrak{s}}}
\newcommand{\tmgo}{\widetilde{\mathfrak{m}}}
\newcommand{\tngo}{\widetilde{\mathfrak{n}}}
\newcommand{\tpgo}{\widetilde{\mathfrak{p}}}

\newcommand{\tgl}{\widetilde{\mathfrak{gl}}}
\newcommand{\tugo}{\widetilde{\mathfrak{u}}}
\newcommand{\gl}{\mathfrak{gl}}
\newcommand{\ggo}{\mathfrak{g}}

\newcommand{\of}{\mathfrak{o}}

\newcommand{\mgo}{\mathfrak{m}}
\newcommand{\ngo}{\mathfrak{n}}
\newcommand{\ago}{\mathfrak{a}}

\newcommand{\pgo}{\mathfrak{p}}

\newcommand{\ugo}{\mathfrak{u}}
\newcommand{\zgo}{\mathfrak{z}}

\newcommand{\Fgo}{\mathfrak{F}}

\newcommand{\al}{\alpha}
\newcommand{\be}{\beta}

\newcommand{\om}{\omega}
\newcommand{\Om}{\Omega}

\newcommand{\la}{\lambda}

\newcommand{\back}{\backslash}

\newcommand{\Cc}{C_c^\infty}

\newcommand{\bg}{\langle}
\newcommand{\bd}{\rangle}

\newcommand{\eps}{\varepsilon}

\renewcommand{\leq}{\leqslant}
\renewcommand{\geq}{\geqslant}

\newcommand{\tlA}{\widetilde{A}}
\newcommand{\tlB}{\widetilde{B}}

\newcommand{\tlG}{\widetilde{G}}
\newcommand{\tlH}{\widetilde{H}}

\newcommand{\tlK}{\widetilde{K}}
\newcommand{\tlL}{\widetilde{L}}
\newcommand{\tlM}{\widetilde{M}}
\newcommand{\tlN}{\widetilde{N}}

\newcommand{\tlP}{\widetilde{P}}
\newcommand{\tlR}{\widetilde{R}}
\newcommand{\tlQ}{\widetilde{Q}}
\newcommand{\tlS}{\widetilde{S}}
\newcommand{\tlT}{\widetilde{T}}
\newcommand{\tlU}{\widetilde{U}}
\newcommand{\tlV}{\widetilde{V}}

\newcommand{\tla}{\widetilde{a}}

\newcommand{\tlg}{\widetilde{g}}
\newcommand{\tlh}{\widetilde{h}}

\newcommand{\tlj}{\widetilde{j}}

\newcommand{\tlx}{\widetilde{x}}

\newcommand{\tl}{\widetilde{0}}

\begin{document}

\newcommand{\upla}{\underline{\rho}} 
\newcommand{\dsl}{\displaystyle \left(}
\newcommand{\rb}{\right)}
\newcommand{\rmF}{F} 
\newcommand{\ogo}{\mathfrak{o}}

\maketitle

\begin{abstract}
La formule des traces relative de Jacquet-Rallis (pour les groupes unitaires ou les groupes linéaires généraux) est une identité entre des périodes des représentations automorphes  et des distributions géométriques. Dans cet article, nous établissons un transfert entre tous  les termes géométriques des groupes unitaires et ceux des groupes linéaires. Nous prouvons en particulier que tous les termes géométriques sont dans l'adhérence faible des intégrales orbitales semi-simples régulières locales. En guise d'illustration, on donne une application à la conjecture de Gan-Gross-Prasad pour les groupes unitaires.
\end{abstract}

\renewcommand{\abstractname}{Abstract}
\begin{abstract}
The relative trace formula of Jacquet-Rallis (for unitary groups  or general linear groups) is an identity between periods of automorphic representations and geometric distributions. In this paper, we prove the  transfer between  all geometric terms of unitary groups and those of linear groups. We also show that all geometric terms are in the weak closure of local regular semi-simple orbital integrals. We mention an application to the Gan-Gross-Prasad conjecture for unitary groups.
\end{abstract}

\tableofcontents

\section{Introduction}

\subsection{Principaux résultats}\label{ssec:ppauxres}

\begin{paragr}[Conjectures de Gan-Gross-Prasad pour les groupes unitaires.] ---\label{S:GGP}
  Soit $E/F$ une extension quadratique de corps de nombres. Soit $n\geq 1$ un entier et $W$ un $F$-espace vectoriel de dimension $n+ 1$. Soit $W=V\oplus F e_0$ une décomposition en sous-espaces. Soit $\Phi$ une forme hermitienne non-dégénérée sur $V\otimes_F E$. Soit $\tilde{\Phi}$ la forme non-dégénérée sur $W\otimes_F E$ qui coïncide sur $V\otimes_F E$ avec $\Phi$ et qui est telle que $\tilde{\Phi}(e_0,e_0)=1$ et  les sous-espaces $V\otimes_F E$  et $Ee_0$ sont orthogonaux. Soit $U=U_\Phi$ et $\tlU=\tlU_{\Phi}$ les groupes unitaires des espaces hermitiens $(V\otimes_F E,\Phi)$ et  $(W\otimes_F E,\tilde{\Phi})$. Le groupe $U$ s'identifie naturellement au sous-groupe de $\tlU$ qui fixe le vecteur $e_0$. Soit $U'=U'_\Phi=U\times \tlU$. Le plongement diagonal identifie $U$ à un sous-groupe de $U'$. Soit $\AAA$ l'anneau des adèles de $F$. Soit $\pi\otimes\tilde{\pi}$ une représentation automorphe cuspidale de $U'(\AAA)$. L'espace de cette représentation est  vu comme un espace de fonctions $f$  sur $U'(F)\back U'(\AAA)$ ; sur celui-ci, on définit la forme linéaire $P_\Phi$, une \og période\fg{}, de la manière  suivante
  \begin{equation}
    \label{eq:periode}
    P_\Phi(f)=\int_{U(F)\back U(\AAA)} f(x)\,dx.
  \end{equation}
Soit $G_E$ et $\tlG_E$ les $F$-groupes obtenus par restriction des scalaires de $E$ à $F$ respectivement de $GL_E(V\otimes_FE)$ et  $GL_E(W\otimes_FE)$. 
Soit $L(\pi\otimes\tilde{\pi},s)$ la fonction $L$ de Rankin-Selberg du changement de base de $\pi\otimes\tilde{\pi}$ à $G_E(\AAA)\times \tlG_E(\AAA)$.  Considérons les deux assertions suivantes :
  \begin{enumerate}
  \item $L(\pi\otimes\tilde{\pi},1/2)\not=0$ ;
  \item Il existe une forme hermitienne $\Phi_1$ et une représentation  automorphe cuspidale  $\pi_1\otimes\tilde{\pi}_1$ sur $U'_{\Phi_1}(\AAA)$ presque équivalente à  $\pi\otimes\tilde{\pi}$ telle que la période $P_{\Phi_1}$ définisse une forme linéaire non nulle sur l'espace de $\pi_1\otimes\tilde{\pi}_1$.
  \end{enumerate}
En presque toute place de $F$, les formes $\Phi$ et $\Phi_1$ sont équivalentes et les groupes associés à $\Phi$ et $\Phi_1$ s'identifient naturellement. L'expression \og presque équivalentes \fg{} signifie que les composantes locales des représentations considérées sont équivalentes en presque toute place.
L'équivalence possible des assertions 1 et 2 est l'objet de la conjecture de Gan-Gross-Prasad (\cite{ggp} conjecture 24.1). Dans \cite{Z1}, Zhang montre que l'équivalence des assertions 1 et 2  vaut pourvu que  $\pi\otimes\tilde{\pi}$ satisfasse certaines hypothèses locales un peu artificielles. Dans cet article, Zhang suit la stratégie élaborée par Jacquet et Rallis (cf. \cite{jacqrall}) qui repose sur la comparaison de formules des traces relatives. Les hypothèses techniques chez W. Zhang sont liées à l'utilisation de variantes simples de ces formules des traces (on fait en général  des hypothèses restrictives en deux places distinctes). Si l'on veut passer outre ces limitations, il faut utiliser des formules des traces relatives générales que Zydor établit dans \cite{leMoi3}.
\end{paragr}

\begin{paragr}[Le côté géométrique de la  formule des traces relative de Jacquet-Rallis pour les groupes linéaires.] --- Comme dans le cas unitaire, on identifie naturellement $G_E$ à un sous-groupe de $\tlG_E$. Soit  $G'_E=G_E\times \tlG_E$. On va s'intéresser à deux sous-$F$-groupes de $G'_E$ et à leur période. Le premier est  $G_E$ plongé diagonalement dans $G'_E$. L'intégration sur $G_E(\AAA)$ donne la période de Rankin-Selberg qui est une forme linéaire sur l'espace d'une représentation automorphe de $G_E'(\AAA)$. Cette période est directement reliée par la théorie des intégrales zêta à la valeur en $1/2$ de la fonction $L$ de Rankin-Selberg. Le second sous-groupe est $G'=GL_F(V)\times GL_F(W)$. Soit $\eta$ le caractère quadratique de $\AAA^\times$ associé à l'extension $E/F$. On note encore $\eta$ le caractère de $G'(\AAA)$ donné par 
$$\eta(g, \tilde{g})=\eta(\det(g))^{n+1}\eta(\det(\tilde g))^{n}.
$$
L'intégration 
$$f\mapsto \int_{G'(F)\back G'(\AAA)} f(g) \eta(g) \, dg
$$
définit la période de Flicker-Rallis qui est une forme linéaire sur l'espace d'une représentation automorphe de $G_E'(\AAA)$.  Elle est intimement reliée à la théorie de la représentation intégrale des fonctions $L$ d'Asai. Son rôle ici est d'éliminer dans le spectre automorphe de $G_E'$ ce qui ne vient pas par changement de base d'un groupe $U'_\Phi$.

La formule des traces relative de Jacquet-Rallis consiste à considérer l'intégrale (\emph{a priori} divergente)
$$I^\eta(f)"="\int_{G_E(F)\back G_E(\AAA)} \int_{G'(F)\back G'(\AAA)} k_f(x,y) \, \eta(y)dxdy
$$
où $k_f$ est le noyau automorphe associé à une fonction test $f$ sur $G'_E(\AAA)$ et à la développer suivant des données spectrales ou géométriques. Le développement spectral  s'exprime, du moins formellement, en terme de périodes de Rankin-Selberg et de Flicker-Rallis.  Dans cet article on ne s'intéressera qu'au côté géométrique. Comme le montre Zydor dans \cite{leMoi3} (cf. aussi théorème \ref{thm:cvgHe}), il y a un procédé de régularisation de l'intégrale ci-dessous. Le côté géométrique est alors une égalité de la forme
$$I^\eta(f)=\sum_{a\in \tilde{\Ac}(F)} I_a^\eta(f),
$$
où l'on introduit le morphisme canonique vers le quotient catégorique
\begin{equation}
  \label{eq:int-canolin}
G'_E\to \tilde{\Ac}=G_E\backslash\back G'_E// G'
\end{equation}
pour les actions naturelles par translation à gauche et à droite des sous-groupes $G_E$ et $G'$. Pour tout $a\in  \tilde{\Ac}$ soit $G'_{E,a}$ la fibre au-dessus de  $a$ du  morphisme \eqref{eq:int-canolin}. Les propriétés les plus importantes des distributions  $I^\eta_a$ sont les suivantes : elles sont invariantes par l'action à gauche de $G_E(\AAA)$, $\eta$-équivariantes pour l'action à droite de $G'(\AAA)$ et enfin leur support est contenu dans $G'_{E,a}(\AAA)$ (cf. théorème \ref{thm:cvgHe}). 

\end{paragr}

\begin{paragr}[Le côté géométrique de la  formule des traces relative de Jacquet-Rallis pour les groupes unitaires.] --- Là encore il s'agit de partir de  l'intégrale (\emph{a priori} divergente)
$$I^\Phi(f)"="\int_{U(F)\back U(\AAA)} \int_{U(F)\back U(\AAA)} k_f(x,y) \, dxdy
$$
où $k_f$ est le noyau automorphe associée à une fonction test $f$ sur $U'(\AAA)$ et de la développer suivant des données spectrales ou géométriques. Du côté spectral, on s'attend à obtenir les périodes \eqref{eq:periode}. Attardons-nous sur le côté géométrique. Tout d'abord, on a un morphisme
\begin{equation}
  \label{eq:int-cano}
  U' =U_\phi'\to \tilde{\Ac}
\end{equation}
bi-invariant par $U$ qui identifie $ \tilde{\Ac}$ au quotient catégorique $U\back\back U'//U$. Le procédé de régularisation de Zydor (cf.  \cite{leMoi3} et théorème \ref{thm:cvgU'}) donne un sens à chaque membre de l'égalité ci-dessus et établit le développement
$$I^\Phi(f)=\sum_{a\in \tilde{\Ac}(F)} I_a^\Phi(f).$$
Les distributions  $I^{\Phi}_a$  sont bi-invariantes et à support inclus dans $U'_a(\AAA)$ où $U'_a$ est la fibre en $a$ du morphisme \eqref{eq:int-cano}.
\end{paragr}

\begin{paragr}[Stratégie de Jacquet-Rallis.] --- Jacquet et Rallis proposent de comparer les côtés géométriques des formules des traces relatives pour $G'_E$ et $U'$. La première observation qui rend plausible cette comparaison est que les distributions du côté géométrique sont indexées par le même ensemble $\tilde{\Ac}(F)$. Cette comparaison doit être possible $a$ par $a$ dans  $\tilde{\Ac}(F)$ pour des fonctions qui vérifient des identités locales d'intégrales orbitales que nous allons expliquer dans le paragraphe suivant. Une fois cette comparaison faite, on a égalité des développements spectraux et des arguments d'analyse fonctionnelle combinée à une analyse fine du développement spectral ramènent la conjecture de Gan-Gross-Prasad au cas (connu) du groupe $G'_E$.  
\end{paragr}

\begin{paragr}[Intégrales orbitales locales.]  --- Soit $\tilde{\Ac}^{\rs}$ l'ouvert \og semi-simple régulier\fg{} de  $\tilde{\Ac}$. Soit $v$ une place de $F$ et $F_v$ le complété de $F$. On peut poser des définitions similaires au cas global relativement à la $F_v$-algèbre quadratique $E_v=E\otimes_F F_v$. On inclut \emph{mutatis mutandis} dans la discussion le cas dégénéré où la place $v$ est scindée dans $E$. On considère une forme hermitienne $\Phi_v$ non dégénérée sur $V\otimes_F E_v$ et on note $U$ et $U'$ les $F_v$-groupes associés. Soit $a\in\tilde{\Ac}^{\rs}(F_v)$ tel que la fibre $U'_a(F_v)$ soit non vide (cette condition détermine d'ailleurs la classe d'équivalence de $\Phi_v$). Dans ce cas, cette fibre est une $U(F_v)\times  U(F_v)$-orbite dont on choisit un élément $y \in U(F_v)$. On  définit alors l'intégrale orbitale  locale
\begin{equation}
  \label{eq:IOP-int-U}
  I_a^{\Phi_v}(f^{\Phi_v})= \int_{U(F_v)}\int_{U(F_v)}f^{\Phi_v}(h_1^{-1} y h_2)\,dh_1dh_2
\end{equation}
pour toute fonction $f^{\Phi_v} \in \Cc(U'(F_v))$.

La situation est un peu différente sur le groupe  $G'_{E}(F_v)$. D'une part pour tout $a\in\tilde{\Ac}^{\rs}(F_v)$, la fibre $G'_{E,a}(F_v)$ n'est jamais vide et elle est toujours formée d'une seule $G_E(F_v)\times G'(F_v)$ -orbite dont on choisit un élément  $x$. Si l'on prend une définition naïve, à cause de la présence du caractère $\eta$, les intégrales orbitales locales sur le groupe $G'_{E}(F_v)$ ne sont pas \emph{a priori} définies sur la base $\tilde{\Ac}^{\rs}(F_v)$ (elles vont dépendre du choix de $x$). C'est pourquoi on fixe un certain facteur de transfert $\Om$ (le choix est \og global\fg{}, cf. §\ref{S:factOmega}) défini sur l'ouvert régulier semi-simple $(G'_{E})^{\rs}(F_v)$ qui vérifie
$$ \Omega(g^{-1} y h ) = \eta(h ) \Omega(y)$$
pour $y\in (G'_{E})^{\rs}(F_v)$,  $ g \in G_E(F_v) $ et $ h\in G'(F_v)$. On pose alors pour tout $f\in \Cc(G_E'(F_v))$
\begin{equation}
  \label{eq:IOP-int-G}
I_a^\Om(f)=  \int_{G_{E}(F_v)}\int_{G'(F_v)}
 f(g^{-1} x h)\Omega(g^{-1}xh)\, dh dg,
\end{equation}
définition qui est évidemment indépendante du choix de $x$.

On dira que les fonctions $f$ et $f^{\Phi_v}$ sont des transferts l'une de l'autre si on a l'égalité 
$$I_a^{\Phi}(f^{\Phi_v})= I_a^\Om(f)$$
pour tout $a\in\tilde{\Ac}^{\rs}(F_v)$ tel que $U'_a(F_v)$ soit non vide.
\end{paragr}

\begin{paragr}[Le résultat principal.] ---
  Revenons au cadre global. Soit $f\in \Cc(G_E'(\AAA))$ et $f^\Phi \in \Cc(U'(\AAA))$.  On dira que $f$ et $f^\Phi$  sont des transferts l'une de l'autre si ces fonctions sont des tenseurs purs (par exemple $f=\otimes_v f_v$ où le produit est pris sur l'ensemble des places de $F$) et si leurs composantes locales le sont (au sens précédent) en toute place. Notons que l'existence de tels couples $(f,f^\Phi)$ est non triviale et repose sur un lemme fondamental démontré par Yun et Gordon (cf. \cite{yun}). Soit $a\in \tilde{\Ac}^{\rs}(F)$ tel que $U_a'(F)\not=\emptyset$. Dans ce cas, les distributions globales $I_a^\eta(f)$ et $I^\Phi_a(f^\Phi)$ s'expriment à l'aide des intégrales orbitales locales définies dans la paragraphe précédent. On déduit du fait que $f$ et $f^\Phi$  sont des transferts l'une de l'autre qu'on  a
$$I_a^\eta(f)=I^\Phi_a(f^\Phi).
$$

Le but de cet article est d'obtenir la généralisation de cet énoncé à toute classe $a\in \tilde{\Ac}(F)$. Pour cela, il nous faut considérer l'ensemble $|\hc|$ des classes de formes hermitiennes non dégénérées sur $V\otimes_F E$ (qu'on identifie à un système de représentants). Voici le principal théorème que nous démontrons dans cet article.

\begin{theoreme}\label{thm:int-main}(cf. théorème \ref{thm:transfertGp})
 Soit $f\in   \Cc(G_E'(\AAA))$ et pour tout $\Phi\in |\hc|$ soit $f^{\Phi}\in  \Cc(U'_\Phi(\AAA))$. On suppose que pour tout $\Phi$  les fonctions $f$ et $f^\Phi$  sont des transferts l'une de l'autre. On a alors pour tout $a\in \tilde{\Ac}(F)$
 \begin{equation}
   \label{eq:identite-tr}
   I_a^\eta(f)=\sum_{\Phi\in|\hc|} I^\Phi_a(f^\Phi).
 \end{equation}
 \end{theoreme}

Notons que le transfert $p$-adique ou aux places scindées dans $E$ (dû à  Zhang, cf. \cite{Z1}) ou aux autres places archimédiennes (qui n'est pas connu en général mais voir \cite{xue}) fournit \og beaucoup\fg{} de fonctions  qui sont transfert l'une de l'autre. Pour le lecteur intéressé par des progrès sur la conjecture de Gan-Gross-Prasad, nous ne proposons dans cet article que le modeste résultat suivant. C'est une légère 
variation autour des résultats et méthodes de Zhang (cf. \cite{Z1}), Xue (cf. \cite{xue}) et Beuzart-Plessis \cite{RBP} qui repose également sur l'extension aux groupes unitaires de la classification d'Arthur (cf. \cite{Mok} et \cite{KMSW}).  Les articles de Zhang et Xue imposent essentiellement deux conditions à la représentation $\pi\otimes\tilde{\pi}$ en deux places finies distinctes. En l'une de ces places, la condition est d'être tempérée : il s'agit de pouvoir prendre des fonctions test en cette place à support régulier. Par notre théorème de transfert général, il est inutile d'imposer cette condition de régularité et donc cette hypothèse n'a plus lieu d'être.  En l'autre place, notons-la $v$, qui est supposée être scindée dans $E$,  ces auteurs imposent à la représentation d'être supercuspidale de façon à écarter le spectre non cuspidal. En fait, si l'on tient compte de résultats de \cite{RBP} (lemme 2.3.1, corollaire 4.5.1 et preuve du théorème 3.5.7), on voit que la place $v$ n'a pas à être scindée dans $E$ et qu'une hypothèse suffisante pour éliminer le spectre non cuspidal est que le changement de base à $G'_E(F_v)$ de la composante locale $\pi_v\otimes\tilde{\pi}_v$ soit supercuspidale. On laisse les détails au lecteur. Pour une application qui utilise toute la puissance du théorème  \ref{thm:int-main}, il faudra probablement attendre une description fine de la partie discrète du côté spectral de la formule des traces.

\begin{theoreme}
  Soit $\pi\otimes\tilde{\pi}$ une représentation automorphe de $U'_\Phi(\AAA)$. Supposons qu'il existe une place finie $v$ de $F$  telle que le changement de base de $\pi_v\otimes\tilde{\pi}_v$ à $G_E'(F_v)$ soit supercuspidal. Alors les assertions 1 et 2 du §\ref{S:GGP} sont équivalentes.
\end{theoreme}
\end{paragr}

\begin{paragr}[Résultats annexes.] --- Le théorème \ref{thm:int-main} implique le résultat de densité suivant qui, logiquement dans l'article, le précède.

\begin{theoreme}\label{thm:int-dens} (cf. théorème \ref{thm:densiteH})
   Soit $v$ une place de $F$ et $f=f_v\otimes f^v  \in \Cc(G_E'(\AAA))$ avec $f_v\in \Cc(G_E'(F_v))$. Si les intégrales orbitales locales semi-simples régulières \eqref{eq:IOP-int-G} de $f_v$ s'annulent alors on a 
 $$I_a^\eta(f)=0$$
pour tout $a\in \tilde{\Ac}(F)$.
\end{theoreme}

L'analogue pour le groupe $U'$ de ce théorème vaut aussi.

\begin{theoreme}\label{thm:int-densU} (cf. théorème \ref{thm:densiteUU})
   Soit $\Phi\in |\hc|$, soit $v$ une place de $F$ et $f=f_v\otimes f^v  \in \Cc(U'_\Phi(\AAA))$ avec $f_v\in \Cc(U'_\Phi(F_v))$. Si les intégrales orbitales locales semi-simples régulières \eqref{eq:IOP-int-U} de $f_v$ s'annulent alors on a 
 $$I_a^\Phi(f)=0$$
pour tout $a\in \tilde{\Ac}(F)$.
\end{theoreme}

Notons que ce dernier théorème couplé avec un cas facile d'annulation dans le lemme fondamental montre que, dans le théorème \ref{thm:int-main}, le membre de droite de \eqref{eq:identite-tr} ne comporte qu'un nombre fini de termes non nuls.
\end{paragr}

\subsection{Stratégie et plan de l'article}

\begin{paragr}
Donnons tout d'abord la stratégie générale.   Les théorèmes \ref{thm:int-main}, \ref{thm:int-dens} et \ref{thm:int-densU} se réduisent tous à un énoncé analogue dans une situation infinitésimale ; dans ce cas, on a un $F$-groupe réductif $G$ qui agit linéairement sur un certain $F$-espace vectoriel $\tggo$, le groupe et l'espace pouvant varier. On dispose de distributions globales $I_a$ attachées à des éléments rationnels du quotient catégorique $\tggo//G$. On montre ensuite que pour  plupart des $a$, les distributions globales $I_a$ se descendent à une situation où la dimension du groupe qui agit diminue. Pour ces $a$-là, on  obtient les théorèmes voulus par une hypothèse de récurrence. 

Le cas essentiel qui reste à traiter est la classe $a$ nilpotente correspondant à $0\in \tggo$. On utilise alors un analogue dans la situation infinitésimale de la formule des traces de Jacquet-Rallis  tel qu'établi par Zydor dans \cite{leMoi} et \cite{leMoi2}. Le côté géométrique comprend toutes les distributions $I_a$ alors que le côté spectral fait apparaître leurs transformées de Fourier. Il résulte de théorèmes de Zhang et Xue que tous les énoncés envisagés sont stables par transformation de Fourier. En utilisant la formule des traces infinitésimale et l'hypothèse de récurrence, on montre qu'une certaine  distribution est  invariante, à support dans le cône nilpotent (la fibre de $a=0$) et sa transformation de Fourier a la même propriété de support. Notons que dans notre situation la représentation $\tggo$ de $G$ n'est pas irréductible et il nous faut considérer plusieurs transformées de Fourier. Un théorème d'Aizenbud assure qu'une telle distribution est nécessairement nulle. Cela donne alors le résultat cherché dans le cas des théorèmes \ref{thm:int-main}, \ref{thm:int-dens} et \ref{thm:int-densU}.
\end{paragr}

\begin{paragr}
  Illustrons cette stratégie dans le cas très simple suivant. Les notations $E/F$, $\eta$ etc. sont celles de la section \ref{ssec:ppauxres}. On considère d'une part l'action du groupe multiplicatif $\mathbb{G}_{m,F}$ sur le plan affine par $t\cdot(x,y)=(tx,ty^{-1})$. Le quotient catégorique s'identifie par $(x,y)\mapsto xy$ à la droite affine. Soit $a\in F$ et $f\in \Sc(\AAA^2)$ une fonction de Schwartz. Si $a\not=0$, on définit alors l'intégrale orbitale
$$I_a^\eta(f) =\int_{\AAA^\times  } f(tx,t^{-1}y) \eta(t) \,dt$$
où $(x,y)\in F^2$ est tel que  $a=xy$. Pour $a=0$, on pose
$$I_0^\eta(f)=\int_{\AAA^\times }    f(t,0) \eta(t) \,dt + \int_{\AAA^\times }    f(0,t) \eta(t) \,dt$$
où chacune des  intégrales est prise comme valeur en $s=0$ du prolongement analytique à la Tate de l'intégrale où l'on remplace $dt$ par $|t|^s dt$. On obtient en tout cas une distribution $I_a^\eta$ qui est $\eta$-équivariante. Le choix d'un caractère additif non trivial sur $F\back \AAA$ et de la forme quadratique invariante $(x,y)\mapsto xy$ permet de définir une transformation de Fourier $f\mapsto \hat{f}$ sur  $S(\AAA^2)$  et donc d'une transformation de Fourier duale sur l'espace des distributions $\eta$-équivariantes. La formule sommatoire de Poisson implique l'identité suivante entre distributions invariantes
\begin{equation}
  \label{eq:RTFelem}
  \sum_{a\in F} I_a^\eta = \sum_{a\in F} \hat{I}_a^\eta.
\end{equation}

L'autre cas est celui de l'action du groupe unitaire $U(1)=\Res_{E/F}(\mathbb{G}_{m,E})^1$ sur la droite affine sur $E$ vue comme $F$-espace. L'application norme $x\mapsto N_{E/F}(x)$ identifie le quotient catégorique à la droite affine sur $F$. Soit $\nu \in F^\times/N_{E/F}(E^\times)$ (on identifie ce quotient à un système de représentants) et $f_\nu\in  \Sc(\AAA_E)$ une fonction de Schwartz. Pour $a\in F$, on pose 
$$I_a^\nu(f_\nu)=\int_{U_1(\AAA)} f_\nu(t x)\,dt$$
s'il existe $x\in E$ tel que  $N_{E/F}(x)=\nu a$ et $I_a^\nu(f_\nu)=0$ sinon.  Cela  définit une distribution $U_1(\AAA)$-invariante. Comme précédemment en utilisant la forme hermitienne $x\mapsto \nu N_{E/F}(x)$, on définit une transformation de Fourier et on a la formule
\begin{equation}
  \label{eq:RTFelemU} \sum_{a\in F} I_a^{\nu} = \sum_{a\in F} \hat{I}_a^\nu.
\end{equation}
Soit $\vc$ l'ensemble des places de $F$ et pour $v\in\vc$ soit $F_v$ le complété de $F$ et $E_v=E\otimes_F F_v$. Supposons que $f$ et $(f_\nu)_{ \nu \in F^\times/N_{E/F}(E^\times )}$ se correspondent au sens où ces fonctions sont des tenseurs purs $f=\otimes_{v} f_v$ et $f_\nu=\otimes_{v} f_{\nu,v}$ tels que pour tout $v\in \vc$ on ait pour tout $x\in E_v^\times$
\begin{equation}
  \label{eq:transfert}
  \int_{F_v^\times  } f_v(t,t^{-1} \nu N_{E/F}(x)) \eta(t) \,dt=\int_{U_1(F_v)} f_{\nu,v}(t x)\,dt.
\end{equation}
En dehors d'un ensemble fini $S$ de places, $v$ est non-archimédienne et non ramifiée dans $E$ et $f_v$ est la fonction caractéristique de $\oc_v^2$ (où $\oc_v$ est l'anneau des entiers de $F_v$) ; dans ce cas, un calcul facile montre que si $\nu$ n'est pas une norme pour $E_v/F_v$ alors le membre de gauche de \eqref{eq:transfert} est nul  et donc
\begin{equation}
  \label{eq:annul-nu}
  \int_{U_1(F_v)} f_{\nu,v}(t x)\,dt=0
\end{equation}
pour tout $x\in E_v$.
L'analogue du théorème \ref{thm:int-main} est l'énoncé suivant : pour tout $a\in F$, on a
\begin{equation}
  \label{eq:trans-elem}
       I_a^\eta(f)=          \sum_{  \nu \in F^\times/N_{E/F}(E^\times ) } I_a^\nu(f_\nu).
     \end{equation}
     Cet énoncé est dû à Jacquet (\cite{JAENS}, cf. aussi \cite{jacqrall}). L'argument de Jacquet repose sur un développement en germes des intégrales locales et le fait que les distributions globales s'expriment assez aisément en termes des intégrales locales. La généralisation de ces faits en dimension supérieure ne semble pas facile. On donne ici un argument, autre que celui de Jacquet,  qu'on a su généraliser.

   Notons que la condition \eqref{eq:annul-nu} implique que $I_a^\nu(f_\nu)=0$ sauf si $\nu$ est une norme pour $E_v/F_v$ pour $v\notin S$. Il n'y a donc qu'un nombre fini de $\nu$ qui interviennent effectivement dans \eqref{eq:trans-elem}. Si $a \not=0$, tous les termes $ I_a^\nu$  sont nuls sauf éventuellement pour $\nu =aN_{E/F}(E^\times)$. Dans ce cas, l'égalité \eqref{eq:trans-elem} résulte immédiatement de \eqref{eq:transfert}.
Les formules \eqref{eq:RTFelem} et  \eqref{eq:RTFelemU} impliquent qu'on a 
 
 $$\sum_{a\in F}\left( I_a^\eta(f) -     \sum_{  \nu \in F^\times/N_{E/F}(E^\times ) } I_a^\nu(f_\nu) \right)= \sum_{a\in F} \left(\hat{I}_a^\eta(f) -     \sum_{  \nu \in F^\times/N_{E/F}(E^\times ) } \hat{I}_a^\nu(f_\nu)\right)$$
Un théorème de Zhang et Xue assurent que $\hat{f}$ et $(\hat{f}_\nu)_{ \nu \in F^\times/N_{E/F}(E^\times )}$ se correspondent aussi. Dans la formule ci-dessus tous les termes associés à $a\not=0$ disparaissent. Il reste donc 
\begin{equation}
  \label{eq:ident-finale}
  I_0^\eta(f) -     \sum_{  \nu \in F^\times/N_{E/F}(E^\times ) } I_0^\nu(f_\nu) =\hat{I}_0^\eta(f) -     \sum_{  \nu \in F^\times/N_{E/F}(E^\times ) } \hat{I}_0^\nu(f_\nu).
\end{equation}
Soit $v$ une place finie scindée dans $E$. Après extension des scalaires à $F_v$, les situations deviennent isomorphes. On peut prendre $f_{\nu,v}=f_v$ et regarder 
$$f_v \mapsto I_0^\eta(f) -     \sum_{  \nu \in F^\times/N_{E/F}(E^\times ) } I_0^\nu(f_\nu) $$
comme une forme linéaire sur $\Cc(F_v^2)$. Cette distribution est invariante, à support dans le cône nilpotent $xy=0$ et l'identité \ref{eq:ident-finale} implique qu'il est de même de sa tranformée de Fourier (locale). Une telle distribution est nécessairement nulle (cf. \cite{rallSchiff} lemme 8.1) donc les deux membres de \eqref{eq:ident-finale} sont nuls ce qui conclut. 
\end{paragr}

\begin{paragr}[Plan de l'article.] --- Les deux premières parties sont consacrées aux énoncés infinitésimaux analogues des théorèmes \ref{thm:int-dens} et \ref{thm:int-densU} ci-dessus  et ont des structures similaires.  La partie \ref{partie:1} est consacrée aux groupes linéaires alors que la partie \ref{partie:2} se concentre sur les groupes unitaires. Les sections \ref{sec:prelim} et \ref{sec:PrelimU} préparent la descente d'un point de vue algébrique alors que les sections \ref{sec:combi} et \ref{sec:combiU}, qui s'inspirent grandement des travaux d'Arthur, la préparent du point de la troncature qui intervient dans la construction des distributions globales. Les sections \ref{sec:RTFinf} et \ref{sec:RTFinfU} rappellent les propriétés des distributions globales ainsi que la formule des traces infinitésimale. Les sections \ref{sec:densite} et \ref{sec:densiteU} énoncent les théorèmes de densité (analogues  aux  théorèmes \ref{thm:int-dens} et \ref{thm:int-densU} ci-dessus) et les démontrent sous l'hypothèse de la descente des distributions globales. Celle-ci est démontrée dans les sections \ref{sec:desc} et \ref{sec:descU}.
  
Dans la section \ref{sec:lecasinf}, on démontre l'énoncé de transfert infinitésimal (analogue du théorème \ref{thm:int-main}). La section \ref{sec:densiteS} est consacrée aux groupes linéaires. On y rappelle la construction des termes géométriques dans la formule des traces. On prépare la comparaison avec l'espace tangent au moyen d'une application de Cayley. On y démontre au passage le théorème \ref{thm:int-dens}. 
La section \ref{sec:densiteGpU} remplit les mêmes objectifs pour les groupes unitaires. Enfin dans la section \ref{sec:factorisationGpU} on démontre le théorème  \ref{thm:int-main}.

Un certain nombre de notations générales se trouvent dans la section \ref{sec:notations-generales}. On a laissé en appendice un énoncé géométrique sur les \og familles convexes\fg{} d'intérêt indépendant.
\end{paragr}

\begin{paragr}[Remerciements.] ---  Les deux auteurs remercient le projet Ferplay ANR-13-BS01-0012 de l'ANR pour son soutien. Le premier auteur nommé souhaite aussi remercier le  projet  Vargen ANR-13-BS01-0001-01 de l'ANR  dont il fait partie  et plus particulièrement l'Institut Universitaire de France qui lui a fourni d'excellentes conditions de travail. Le deuxième auteur nommé a été partiellement soutenu par le projet \#711733 de la fondation Minerva.
  \end{paragr}

\section{Notations générales}\label{sec:notations-generales}

\subsection{Groupes et quotients}
\label{ssec:notations-combi}

\begin{paragr}
  Soit $F$ un corps de caractéristique $0$ et $G$ un groupe algébrique défini  sur $F$. Suivant les notations d'Arthur (cf. \cite{ar1}), soit $X^*(G)$ le groupe des caractères rationnels de $G$ définis sur $F$ et  $\ago_G^*=X^*(G)\otimes_\ZZ \RR$. Soit $\ago_G$ son espace dual. Pour tout groupe $H$ (algébrique, défini  sur $F$, bien souvent un sous-groupe de $G$), on introduit le signe
$$\eps_H^G=(-1)^{\dim_\RR(\ago_H)-\dim_\RR(\ago_G)}.
$$
\end{paragr}

\begin{paragr}
  Supposons de plus $G$ est réductif. Soit $A_0$ un sous-tore déployé sur $F$ et maximal pour cette propriété. Dans la suite, sauf mention contraire, un sous-groupe parabolique de $G$ est supposé défini sur $F$ et contenant $A_0$ (on dira aussi \og semi-standard\fg{}). Soit $P$ un tel sous-groupe parabolique. Il admet alors une décomposition de Levi $P=M  N_P$ où $N_P$ est le radical unipotent de $P$ et $M=M_P$ est l'unique facteur de Levi contenant $A_0$. Un tel groupe $M_P$ est appelé simplement sous-groupe de Levi de $G$ dans la suite. Soit $A_P=A_M$ le tore central $F$-déployé maximal de $M$. Soit $\Sigma_P$ l'ensemble des racines du tore  $A_P$ dans $P$. Soit $\Delta_P\subset \Sigma_P$ le sous-ensemble des racines simples et $\hat{\Delta}_P$ l'ensemble des poids. Ces ensembles s'identifient à des parties de $\ago_P^*= \ago_M^*$ qui engendrent un sous-espace noté $(\ago_P^G)^*=(\ago_M^G)^*$.

Soit $M_0$  le centralisateur de $A_0$ dans $G$ (c'est un sous-groupe de Levi minimal de $G$) et $\ago_0^*=\ago_{M_0}^*$. On a une inclusion naturelle  $\ago_P^*\hookrightarrow \ago_0^*$ et une somme directe
$$\ago_0^*=(\ago_0^P )^*\oplus (\ago_P)^*.
$$
où l'on pose $(\ago_0^P )^*=(\ago_0^M )^*$.
\end{paragr}

\begin{paragr}
Pour tout sous-groupe $H\subset G$, soit $\fc^G(H)$ l'ensemble des sous-groupes paraboliques de $G$ qui contiennent $H$. Si, de plus, $H$ est un sous-groupe de Levi, soit $\pc^G(H)\subset \fc^G(H)$ le sous-ensemble des $P$ tels que $M_P=H$. Si le contexte est clair, on pourra omettre l'exposant $G$. 
\end{paragr}

\begin{paragr}[Fonctions $\tau$ et $\hat{\tau}$.] \label{S:tau} --- Plus généralement pour des sous-groupes paraboliques $P\subset Q$, on dispose d'ensembles $\Delta_P^Q\subset \Delta_P$ et $\hat{\Delta}_P^Q$, d'une décomposition
$$\ago_P^*=(\ago_P^Q )^*\oplus (\ago_P)^*
$$
et une décomposition duale (notée sans exposant $*$). Soit $\tau_P^Q$ et $\hat{\tau}_P^Q$ les fonctions caractéristiques respectives  des cônes
$$\{H\in \ago_{0} \mid \al(H) > 0, \  \forall \ \al \in  \Delta_{P}^{Q} \}$$
et
$$\{H\in \ago_{0} \mid \varpi(H) > 0, \  \forall \ \varpi \in  \hat{\Delta}_P^Q \}.$$
\end{paragr}

\begin{paragr}[Groupe de Weyl.] ---
 Le groupe de Weyl $W$ de $(G,M_0)$ agit sur $\ago_0$. On fixe alors un produit scalaire $W$-invariant sur $\ago_0$ et ce dernier est muni de la mesure euclidienne ainsi que tous ses sous-espaces.
\end{paragr}

\begin{paragr}[Algèbres de Lie.] ---  Sauf mention expresse du contraire,  l'algèbre de Lie d'un sous-groupe de $G$ est noté par la même lettre en gothique, ici $\ggo$.
  \end{paragr}

  \begin{paragr}[Quotient catégorique.] --- \label{S:quotient-categor} Soit $X$ une variété affine sur $F$ sur laquelle le  groupe $G$ agit. Soit $X//G$ le quotient catégorique ; celui-ci existe et c'est le spectre de l'algèbre $F[X]^G$ des fonctions régulières sur $X$ et $G$-invariantes. Soit
$$a: X \to X//G$$
le morphisme canonique. Pour tout $a\in X//G$, soit $X_a$ la fibre du morphisme canonique au-dessus de $a$.  
    Soit $x\in X$ un point géométrique. On dit que $x$ est ($G$-)\emph{semi-simple} si l'orbite de $x$ sous $G$ est fermée dans $X$. On dit que $x$ est ($G$)-\emph{régulier} si son orbite est de dimension maximale. Les points réguliers forment un ouvert de $X$.
  \end{paragr}

\subsection{Sur les corps de nombres}\label{ssec:cdn}

\begin{paragr}
  Soit $F$ est un corps de nombres et $\AAA=\AAA_F$ l'anneau des adèles de $\AAA$. Soit $|\cdot|_\AAA$ la valeur absolue adélique normalisée. Soit $\vc$ l'ensemble des places de $F$ et $\vc_\infty\subset \vc$ l'ensemble des places archimédiennes. Pour $v\in \vc$, soit $F_v$ le complété de $F$ en $v$. Si $v$ est non-archimédienne soit $\oc_v\subset F_v$ l'anneau des entiers. Pour tout ensemble $S\subset \vc$, on note $\AAA_S$ et $\AAA^S$ respectivement l'anneau des adèles dans $S$ et hors $S$. Si $S=\{v\}$, on a $\AAA_S=F_v$. Soit $\oc^S=\prod_{v\notin S} \oc_v$.
\end{paragr}

\begin{paragr}\label{S:Gcrochet}
Soit $G$ un groupe réductif connexe défini sur $F$. Soit 
$$G(\AAA)^1=\bigcap_{\chi \in X^*(G)} \ker |\chi|_\AAA \ ;$$
c'est un  sous-groupe de $G(\AAA)$. Soit $A_G^\infty$ le sous-groupe central obtenu comme la composante neutre du groupe des $\RR$-points du sous-tore $\QQ$-déployé maximal de la restriction des scalaires de $F$ à $\QQ$ de $A_G$.  Le groupe $G(\AAA)$ est alors isomorphe au produit $A_G^\infty\times G(\AAA)^1$. On définit
$$[G]=G(F)\back G(\AAA).$$
  \end{paragr}

\begin{paragr}[Application $H_P$.] --- \label{S:HP} Soit $K=\prod_{v\in \vc}K_v$ un sous-groupe compact maximal de $G(\AAA)$ admissible vis-à-vis du sous-groupe de Levi minimal $M_0$ (au sens de \cite{arthur2} section 1). Pour tout sous-groupe parabolique $P$ de $G$, on dispose de la décomposition d'Iwasawa et de l'application 
$$H_P : G(\AAA)\to \ago_P
$$
qui vérifie
$$\chi(H_{P}(g))=\log |\chi(p)|_{\AAA}$$
pour tout $g\in pK$, $p\in P(\AAA)$ et $\chi\in X^*(P)$.
\end{paragr}

\begin{paragr}[Espaces de Schwartz-Bruhat.] ---\label{S:Schw} 
Soit $V$ un $F$-espace vectoriel. 
  Soit $\Sc(V(\AAA))$ l'espace de Bruhat-Schwartz des fonctions sur $V(\AAA)$ à valeurs complexes. On définit également pour $S$ fini l'espace $\Sc(V(\AAA_S))$. Soit $S_\infty =S\cap \vc_\infty$ et $ S^\infty=S\setminus S_\infty$. On a alors $\Sc(V(\AAA_S)=\Sc(V(\AAA_{S_\infty}))\otimes_\CC \Cc(V(\AAA_{S^\infty}))$. L'espace $\Sc(V(\AAA_{S_\infty}))$ est l'espace de Schwartz usuel de sa topologie habituelle.

Donnons-nous une action linéaire algébrique de $G$ sur $V$. Pour tout $g\in G(\AAA)$ et $f\in \Sc(V(\AAA))$ soit $f^g\in \Sc(V(\AAA))$ définie par
$$f^g(X)=f(g\cdot X)
$$
pour tout $X\in V(\AAA)$. Cela définit une action à droite de $G(\AAA)$ sur  $\Sc(V(\AAA))$. De même,  $G(\AAA_S)$ agit à droite sur   $\Sc(V(\AAA_S))$.
\end{paragr}

\begin{paragr}[Mesures de Haar.] --- \label{S:Haar}Pour tout $F$-espace vectoriel $V$, le groupe $V(\AAA) $ est muni de la mesure de Haar qui donne le volume $1$ au quotient $ V(F)\back V(\AAA)$.
\end{paragr}

\begin{paragr}[Transformation de Fourier partielle.] --- \label{S:TFPV} Soit $\psi$ un caractère additif, continu et non-trivial
$$\psi: F\back \AAA \to \CC^\times.
$$
Soit 
$$\bg \cdot,  \cdot \bd : V \times V \to F
$$
une forme bilinéaire symétrique, non dégénérée et $G$-invariante. Pour tout sous-espace $V_1\subset V$ qui est $G$-invariant et non dégénéré pour $\bg \cdot,\cdot \bd$, soit $V_2$  son orthogonal. Tout $X\in V$ s'écrit $X_1+X_2$ selon la décomposition $V=V_1\oplus V_2$. On définit alors la transformée de Fourier partielle par
\begin{equation}
  \label{eq:TFPV1}
  \hat{f}_{V_1}(X_1+X_2)=\int_{V_1(\AAA) } f(Y+X_2) \psi(\bg X_1,Y\bd) \, dY.
\end{equation}
On obtient ainsi  un automorphisme de $\Sc(V(\AAA))$. 

On définit de même un analogue local (sur $\AAA_S$ pour un ensemble fini $S$ de places) de cette transformation en prenant comme mesure la mesure auto-duale.
\end{paragr}

\part{Le cas infinitésimal linéaire}\label{partie:1}
 
\section{Préliminaires algébriques}\label{sec:prelim}

\subsection{Stratification} 

\begin{paragr}\label{S:tggo}
Soit $F$ un corps de caractéristique $0$ et $V$ un $F$-espace vectoriel de dimension $n$. Soit $\tgl_F(V)$ l'espace vectoriel défini comme la somme directe
$$\tgl_F(V)=\gl_F(V)\oplus V\oplus V^*
$$
où $\gl_F(V)$  est l'espace vectoriel des endomorphismes de $V$ et  $V^*$ est l'espace vectoriel dual de $V$.  Lorsque le contexte est clair, on omet l'indice $F$. Le groupe $GL(V)$ des automorphismes de $V$ agit à gauche sur  $\tgl(V)$ de la manière suivante : soit $X=(A,b,c)\in \tgl(V)$ alors
$$g\cdot X= (gAg^{-1},gb,cg^{-1}). 
$$
\end{paragr}

\begin{paragr}\label{S:AV}
Suivant §\ref{S:quotient-categor}, on introduit le quotient catégorique
$$\Ac_V=\tgl(V)//GL(V),$$
le morphisme canonique
\begin{equation}
  \label{eq:a}
  a:\tgl(V)\to \Ac_V.
\end{equation}
Par abus, on utilisera souvent la lettre $a$ pour désigner un point de $\Ac$. Dans  ce cas,  $\tgl(V)_a$ sera la fibre en $a$ du morphisme \eqref{eq:a}. Lorsque le contexte est clair, on note simplement $\Ac=\Ac_V$.
\end{paragr}

\begin{paragr}\label{S:A}
Soit $X=(A,b,c)\in \tgl(V)$. Soit $t^n+a_1t^{n-1}+\ldots+a_n\in F[t]$ le polynôme caractéristique de $A$. Pour $1\leq i\leq n$, les fonctions $a_i$ et $b_i=cA^{i-1}b$ sont algébriquement indépendantes, invariantes sous $GL(V)$ et engendrent l'algèbre  $F[\tgl(V)]^{GL(V)}$ (cf. \cite{Z1} lemme 3.1). On utilise ces fonctions pour l'identification
$$\Ac \simeq \AAA_{2n}$$
de $\Ac$ avec l'espace affine de dimension $2n$.
\end{paragr}

\begin{paragr} \label{S:Areg}Pour tout entier $r\geq 1$ et tout $X=(A,b,c)\in \tgl(V)$  soit  la matrice carrée de taille $r$
$$\Delta_r(X)= (cA^{i+j}b)_{0\leq i,j\leq r-1}
$$
et $d_r(X)=\det(\Delta_r(X))$. Par commodité, on pose $d_0=1$. Ce sont des éléments de $F[\tgl(V)]^{GL(V)}$, qui sont nuls pour $r>n$.  Pour tout $r\geq 0$, on définit une partie localement fermée   $\Ac^{(r)}$ par la condition $d_r\not=0$ et $d_i=0$ pour tout $i>r$. Deux telles parties distinctes sont d'intersection vide. Bien sûr, $\Ac^{(r)}\not=\emptyset$ si et seulement si $r\leq n$. La partie   $\Ac^{(0)}$ est un fermé et $\Ac^{(n)}$ est un ouvert dense de $\Ac$ : c'est l'ouvert \og régulier semi-simple\fg{} noté
$$\Ac^{\rs}=\Ac^{(n)}.$$
 On définit des fermés 
$$\Ac^{(\leq r)}=\bigcup_{0\leq i\leq r} \Ac^{(i)}.
$$
qui sont les adhérences des $\Ac^{(r)}$  et des ouverts
$$\Ac^{(\geq r)}=\bigcup_{r \leq i} \Ac^{(i)}.
$$
On définit de même (ou si l'on préfère par image inverse par $a$) une partie localement fermée  $\tgl(V)^{(r)}$ et un ouvert $\tgl(V)^{(\geq r)}$ de $\tgl(V)$.  L'ouvert dense $\tgl(V)^{\rs}=\tgl(V)^{(n)}$ est formé des éléments réguliers et semi-simples.

\begin{lemme}
  \label{lem:A0}
Le fermé $\tgl(V)^{(0)}$ est formé  des triplets $(A,b,c)$ tels que $cA^ib=0$ pour tout $i$.
\end{lemme}

\begin{preuve}
  La condition $d_1=0$ donne $cb=0$. On vérifie immédiatement que la condition  $cA^ib=0$ pour $0\leq i < r-1$ implique la relation $d_r=(cA^{r-1}b)^r$. La conclusion est alors claire.
\end{preuve}
\end{paragr}

\subsection{Sommes directes}
\begin{paragr}[Somme directe canonique $V^+\oplus V^-$ associé à $X$.]\label{S:Xr} --- 
Soit $0\leq r\leq n$ et $X=(A,b,c)\in \tgl(V)^{(r)}$. Soit  $V^+$  le sous-espace de $V$ engendré par $A^ib$ pour $0\leq i\leq r-1$ et $V^-$ l'orthogonal dans $V$ de la famille des  $cA^i$ pour $0\leq i\leq r-1$.  Il résulte de la définition de $r$, qu'on a $\dim(V^+)=\mathrm{codim}(V^-)=r$ et que $V^-\cap V^+=(0)$. On a donc 
 \begin{equation}
    \label{eq:direct2}
    V=V^+\oplus V^-.
  \end{equation}
C'est la \emph{somme directe canonique} associée à $X$. 
\end{paragr}

\begin{paragr}[Morphisme associé à une somme directe.] --- \label{S:somme-directe}Soit
  \begin{equation}
    \label{eq:direct1}
    V=V^+\oplus V^-
  \end{equation}
un décomposition en somme directe.  Soit
$$V^*=(V^+)^*\oplus (V^-)^*$$
la décomposition duale où, par exemple, on identifie le dual $(V^+)^*$ de $V^+$ à l'orthogonal $(V^-)^\perp$ de $V^-$. Soit $r=\dim(V^+)$ et $\tsgo(V^+,V^-)\subset \tgl(V)$ la partie localement fermée formée des triplets  $X=(A,b,c)\in \tgl(V)$ tels que $V^+$ soit le sous-espace de $V$ engendré par $A^ib$ pour $0\leq i\leq r-1$ et $V^-$ soit l'orthogonal dans $V$ de la famille des  $cA^i$ pour $0\leq i\leq r-1$.
  
Lorsque $r\geq 1$ et $(A,b,c)\in\tsgo(V^+,V^-)$, on définit des vecteurs $c'\in (V^+)^*$ et  $b'\in V^+$ par les conditions 
\begin{equation}
  \label{eq:cprime}
  c'A^ib=\left\lbrace
  \begin{array}{l}
    0 \text{  si  } 0\leq i<r-1 \\
1 \text{  si  } i=r-1.
  \end{array}\right.
\end{equation}
et
\begin{equation}
  \label{eq:bprime}
cA^ib'=\left\lbrace
  \begin{array}{l}
    0 \text{  si  } 0\leq i<r-1 \\
1 \text{  si  } i=r-1.
  \end{array}\right.
\end{equation}
 Suivant la décomposition \eqref{eq:direct2}, on écrit matriciellement
$$
A=
\begin{pmatrix}
  A^+ & L^+ \\ L^- & A^{-}
\end{pmatrix}.
$$

Observons que $(A^+,b,c) \in \tgl(V^+)^{(r)}$
Regardons $L^-$ comme un élément de $\Hom(V^+,V^-)$. Pour $0\leq i \leq r-2$, on a $L^-A^ib=0$. Par conséquent, on a  $L^-=vc'$ pour $c'\in (V^+)^*$ défini par la condition \eqref{eq:cprime} et un unique $v\in V^-$. De même, $L^+=b'w$ pour $b'\in V^+$  défini par la condition \eqref{eq:bprime} et un unique $w\in (V^-)^*$.
  On définit alors un isomorphisme 
  \begin{equation}
    \label{eq:iota-s}
    \iota=\iota_{V^+\oplus V^-}: \tgl(V^+)^{(r)} \times \tgl(V^-) \to \tsgo(V^+,V^-)
  \end{equation}
  
par
$$
\iota((A,b,c),(A',v,w))= (   \begin{pmatrix}
  A & b'w \\ vc' & A'
\end{pmatrix}, b,c)
$$
où $b'$ et $c'$ sont les vecteurs définis par les conditions \eqref{eq:bprime} et  \eqref{eq:cprime}. Si $r=n$ alors $V^+=V$ et $\iota$ est l'immersion ouverte de  $\tgl(V)^{(n)}$ dans $\tgl(V)$. Si $r=0$ on a $V^-=V$ et le morphisme est alors par définition l'identité.

Pour tout $X\in \tgl(V)^{(r)}$, on a une  décomposition associée $V=V^+\oplus V^-$, cf. \eqref{eq:direct2}. On pose alors
$$\iota_X=\iota_{V^+\oplus V^-}.$$

 L'inverse de $\iota$ est donné par 

 $$   (\begin{pmatrix}
  A^+ & L^+\\ L^- & A^-
\end{pmatrix}, b,c) \mapsto ((A^+,b,c),(A^-, L^-(A^+)^{r-1}b, c(A^+)^{r-1}L^+)
$$
Le morphisme $\iota$ est $GL(V^+)\times GL(V^-)$-équivariant si l'on identifie $GL(V^+)\times GL(V^-)$ au sous-groupe de de $GL(V)$ qui stabilise à la fois $V^+$ et $V^-$. Il passe donc au quotient catégorique : on obtient un morphisme encore noté $\iota$
\begin{equation}
  \label{eq:iota}
  \iota: \Ac_{V^+}^{(r)}\times \Ac_{V^-} \to \Ac_V.
\end{equation}

\begin{lemme}
 \label{lem:iso-r}
Le morphisme $\iota$ induit un isomorphisme de   $\Ac_{V^+}^{(r)}\times \Ac_{V^-}$ sur l'ouvert $\Ac_V^{(\geq r)}$.
\end{lemme}

\begin{preuve}
Tout d'abord   $\tsgo(V^+,V^-) $ est un fermé de $\tgl(V)^{(\geq r)}$. Compte tenu de l'isomorphisme  \eqref{eq:iota-s}, il suffit de voir que l'inclusion de  $\tsgo(V^+,V^-) $ dans $\tgl(V)^{(\geq r)}$ induit un isomorphisme 
$$\tsgo(V^+,V^-)// GL(V^+)\times GL(V^-) \to \tgl(V)^{(\geq r)}// GL(V)=\Ac_V^{(\geq r)}.$$
Il est facile de voir que chaque $GL(V)$-orbite dans  $\tgl(V)^{(\geq r)}$ rencontre  $\tsgo(V^+,V^-) $ en une unique $GL(V^+)\times GL(V^-)$-orbite. La conclusion s'ensuit aisément (cf. \cite{LuRi} théorème 2.2).

\end{preuve}

\begin{lemme}\label{lem:entierk} Soit $(X,Y)\in \tgl(V^+)^{(r)} \times \tgl(V^-)$ et $k\in \NN$.  On a
  \begin{equation}
    \label{eq:mult}
    d_{r+k}( \iota(X,Y))= d_r(X) d_k(Y)
  \end{equation}
  En particulier,  $Y\in  \tgl(V^-)^{(k)}$ si et seulement si $\iota(X,Y)\in  \tgl(V)^{(r+k)}$.
\end{lemme}

\begin{preuve} Soit $X=(A,b,c)\in \tgl(V^+)^{(r)}$ et $Y=(B,v,w)\in \tgl(V^-)$. Alors 
$$\iota(X,Y)=( Z   , b,c)
$$
où l'endomorphisme $Z$ s'écrit matriciellement relativement à la décomposition \eqref{eq:direct2}
$$
Z=
\begin{pmatrix}
  A  & b'w   \\ vc'& B
\end{pmatrix}.
$$

Pour $0\leq i \leq r-1$, on a $cZ^ib=cA^ib$. Par ailleurs, on a 
\begin{eqnarray*}
  Z^{r}b&=& Z Z_{}^{r-1}b\\
&=& Z A^{r-1}b\\
&\in & v c'A_{}^{r-1}b + V^+\\
&\in & v  + V^+\\
\end{eqnarray*}
On en déduit que $Z(V^+) \subset V^+\oplus \vect(v)$.
En réitérant, on obtient
\begin{eqnarray*}
  Z^{r+1}b&\in & B v  + V^+\oplus \vect(v)
\end{eqnarray*}
puis par récurrence pour $i\geq 0$
\begin{eqnarray*}
  Z^{r+i}b&\in &   B^i v   + V^+\oplus \vect(v, B v,\ldots, B^{i-1}v). 
\end{eqnarray*}

Il s'ensuit que par des manipulations sur les colonnes on a

\begin{eqnarray*}
  d_{r+k}(Z,b,c)&=&\det((cZ^{i+j}b)_{0\leq i,j \leq r+k-1} ) \\
 &=& \det\begin{pmatrix}  (c A^{i+j}b)_{0\leq i,j\leq r-1} & 0 \\
(c Z^{i+j}b)_{\scriptsize{
  \begin{array}{l} r\leq i\leq r+k-1\\
    0\leq j \leq r-1
  \end{array}}} & (c Z^{i} B^{j}v)_{\scriptsize{
  \begin{array}{l} r\leq i\leq r+k-1\\
    0\leq j \leq k-1
  \end{array}}}
\end{pmatrix}\\
&=& \det((c A^{i+j}b)_{0\leq i,j\leq r-1})\cdot \det(   (c Z^{i} B^{j}v)_{\scriptsize{
  \begin{array}{l} r\leq i\leq r+k-1\\
    0\leq j \leq k-1
  \end{array}}})
\end{eqnarray*}

Dualement, par un raisonnement analogue, on obtient pour $i\geq 0$
$$ cZ^{r+i}\in  w B^{i} + (V^+)^* \oplus\vect(w,w B,\ldots wB^{i-1}). 
$$
Par des manipulations sur les lignes, on voit que  
\begin{eqnarray*}
  \det(   (c Z^{i} B^{j}v)_{\scriptsize{
  \begin{array}{l} r\leq i\leq r+k-1\\
    0\leq j \leq k-1
  \end{array}}})&=& \det(   (w B^{i} B^{j}v)_{\scriptsize{
  \begin{array}{l} 0\leq i\leq k-1\\
    0\leq j \leq k-1
  \end{array}}})\\
&=& d_k(Y).
\end{eqnarray*}
Cela conclut.
\end{preuve}

\end{paragr}

\subsection{Décomposition de Jordan}\label{ssec:Jordan}

\begin{paragr}
  Le but de cette section est de donner une construction canonique d'une décomposition de Jordan pour tout $X\in \tgl(V)$. On a la notion d'élément semi-simple de $\tgl(V)$ (cf. §\ref{S:quotient-categor}. On dit que $X$ est  \emph{nilpotent} si son invariant $a$ est nul. On a le lemme suivant.

  \begin{lemme}(Rallis-Schiffmann, cf. \cite{rallSchiff} théorème 6.2).\label{lem:ss}
    Un élément $X=(A,b,c)\in \tgl(V)$ est semi-simple si et seulement si les deux conditions sont satisfaites :
  \begin{enumerate}
  \item Dans la somme directe canonique de $X$, cf. \eqref{eq:direct2},
    \begin{equation*}
      \label{eq:directe}
      V=V^{+}\oplus V^{-},
    \end{equation*}
les espaces $V^{\pm}$ sont stables par $A$.
\item L'endomorphisme de $V^-$ induit par $A$ est semi-simple au sens usuel.
\end{enumerate}
\end{lemme}

\begin{remarque}
  Pour un élément $X=(A,b,c)\in \tgl(V)$ semi-simple, l'espace $V^+$ est engendré par la famille des $A^ib$ pour $i\geq 0$ et l'espace $V^-$ est l'orthogonal de la famille des $cA^i$ pour $i\geq 0$.

L'ouvert $\tgl(V)^{\rs}$ est l'ensemble des éléments à la fois semi-simples et réguliers (cf. \cite{rallSchiff} théorème 6.1) .
\end{remarque}

\end{paragr}

\begin{paragr}
  
  \begin{lemme}\label{lem:fibress}
    Soit $a\in \Ac(F)$. L'ensemble des éléments semi-simples dans la fibre $\tgl(V)_a(F)$ est non-vide et est formé d'une unique orbite sous $G(F)$.
  \end{lemme}

  \begin{preuve}
Deux éléments semi-simples dans $\tgl(V)_a(F)$ sont géométriquement conjugués. Deux éléments semi-simples de $\tgl(V)(F)$ géométriquement conjugués le sont sur $G(F)$. Pour cela il suffit de vérifier que pour tout $X$ semi-simple, l'ensemble de cohomologie galoisienne $H^1(F,G_X)$ associé au centralisateur $G_X$ dans $G$ trivial. Or, d'après le lemme \ref{lem:ss}, un tel élément semi-simple s'écrit $X=(A^++A^-,b,c)$ avec $V=V^+\oplus V^-$, $(A^+,b,c)\in \tgl(V^+)^{\rs}$ et $A^-\in \End(V^-)$ est semi-simple. Le groupe $G_X$ s'identifie alors au centralisateur de $A^-$ dans $GL(V^-)$ et ce dernier groupe est à restriction des scalaires près un produit de groupes linéaires d'où la trivialité du $H^1$ associé.

Il nous reste à voir que $\tgl(V)_a(F)$ possède des éléments semi-simples.     Soit $r$ tel que  $a\in \Ac^{(r)}(F)$. D'après le §\ref{S:somme-directe},  pour toute décomposition $V^+\oplus V^-$ en sous-$F$-espaces, on a un diagramme commutatif

$$
\xymatrix{\tgl(V^+)^{(r)} \times \tgl(V^-)^{(0)}   \ar[rr]_{}   \ar[d]   & & \tgl(V)^{(r)}    \ar[d]_{} \\ \Ac_{V^+}^{(r)}\times \Ac_{V^-}^{(0)} \ar[rr]   & & \Ac_V^{(r)}}
$$
de $F$-morphisme, la flèche horizontale du bas étant un isomorphisme (cf. les lemmes \ref{lem:iso-r} et \ref{lem:entierk}). On est donc ramené aux deux cas suivants   $r=0$ et $r=\dim(V)$. Le premier cas est évident puisque tout polynôme est polynôme caractéristique d'un endomorphisme semi-simple de $V$. Le second résulte de ce que pour tout $a\in \Ac(F)$, la fibre $\tgl(V)_a(F)$ est non vide (il y a en fait une section de $a$ définie sur $F$, cf. \cite{Z1} preuve du lemme 3.1) et que pour $a\in\Ac^{\rs}(F)$ cette fibre est composée d'éléments semi-simples réguliers.
  \end{preuve}

\end{paragr}

\begin{paragr}[Décomposition de Jordan : cas extrêmes.] --- \label{S:Jor1} Soit $X\in \tgl(V)$. Si $X$ appartient à la strate ouverte $\tgl(V)^{(\rs)}$ alors $X$ est semi-simple régulier. On pose dans cas $X_s=X$ et $X_n=0$.  Supposons à l'opposé que $X=(A,b,c)\in \tgl(V)^{(0)}$.  Soit $A=A_s+A_n$ la décomposition de Jordan usuelle de $A$. On  pose alors 
$$X_s=(A_s,0,0)$$
et
$$X_n=(A_n,b,c)
$$
\begin{lemme}\label{lem:J1}
  On a $a(X)=a(X_s)$ et les éléments $X_s$ et $X_n$ sont respectivement semi-simples et nilpotents.
\end{lemme}

\begin{preuve}
  Le fait que $X_s$ est semi-simple résulte du lemme \ref{lem:ss}, l'égalité $a(X)=a(X_s)$ du lemme \ref{lem:A0} et du fait que $A$ et $A_s$ ont même polynôme caractéristique. Les invariants $cA_n^ib$ sont des combinaisons linéaires de $cA^ib$, donc nuls d'après le  lemme \ref{lem:A0}. L'égalité $a(X_n)=0$ en résulte.
\end{preuve}
\end{paragr}

\begin{paragr}[Décomposition de Jordan : cas intermédiaires.] --- \label{S:Jor2} Soit $1\leq r <n$ et $X=(A,b,c)\in \tgl(V)^{(r)}$. Reprenons les constructions du §\ref{S:Xr}. Il existe alors un unique $X^+\in  \tgl(V^+)^{(r)}$ et $X^-\in  \tgl(V^-)$ tels que $X=\iota_X(X^+,X^-)$. D'après le lemme \ref{lem:entierk}, on a nécessairement $X^-\in  \tgl(V^-)^{(0)}$. D'après §\ref{S:Jor1}, on sait définir $X^+_s=X^+$ et $X^-_s$. On pose alors
  \begin{equation}
    \label{eq:ss}
  X_s=\iota_X(X^+_s,X^-_s)
\end{equation}
et
 \begin{equation}
    \label{eq:nilp}
X_n=X-X_s.
\end{equation}

\begin{lemme}\label{lem:J2}
  On a $a(X)=a(X_s)$ et les éléments $X_s$ et $X_n$ sont respectivement semi-simples et nilpotents.
\end{lemme}

\begin{preuve}
Écrivons $X^-=(B,v,w)$. On a alors $X^-_s=(B_s,0,0)$ où $B=B_s+B_n$ est la décomposition de Jordan usuelle dans $\gl(V^-)$. Il résulte du lemme \ref{lem:ss} que $X_s=\iota(X^+,(B_s,0,0))$ est semi-simple. L'égalité $a(X_s)=a(X)$ résulte immédiatement des égalités $a(X^+)=a(X^+_s)$ et $a(X^-)=a(X^-_s)$ (cf. lemme \ref{lem:J1}). 
On a $X_n=( Z, 0,0)$ où  $Z=\begin{pmatrix}
  0 & b'w \\ vc' & B_n
\end{pmatrix}$. Il s'agit de voir que $Z$ est nilpotent au sens usuel. Mais cela résulte  immédiatement  de la formule
$$
Z^k=
\begin{pmatrix}
  0 & b'w B_n^{k-1} \\ B_n^{k-1} vc' & B_n^k + (c'b') \sum_{i=0}^{k-2} B_n^i vw B_n^{k-2-i}
\end{pmatrix}.
$$

\end{preuve}
\end{paragr}

\begin{paragr} On a le lemme suivant.

  \begin{lemme}\label{lem:action}
    Pour tout $\delta\in GL(V)$ et tout $X\in \tgl(V)$ on a $(\delta\cdot X)_s=\delta\cdot X_s  $ et $(\delta\cdot X)_n=\delta\cdot X_n$.
  \end{lemme}

  \begin{preuve}
    Soit $V=V^+\oplus V^-$ la somme canonique associée à $X$. Alors $\delta V^+\oplus \delta V^-$ est la décomposition canonique de $\delta\cdot X$ et
    \begin{equation}
      \label{eq:iota-delta}
      \iota_{\delta\cdot X}=\delta \iota_X(\delta^{-1}\cdot,\delta^{-1}\cdot). 
    \end{equation}
On est alors ramené à prouver la propriété dans les cas extrêmes $X\in  \tgl(V)^{(n)}$ ou $X\in  \tgl(V)^{(0)}$ qui sont évidents.
  \end{preuve}
\end{paragr}

\subsection{Sections transverses}\label{ssec:transv}

\begin{paragr}\label{S:transv}
Soit $F$ un corps de caractéristique $0$. Soit $I$ un ensemble fini et pour tout $ i\in I$ soit  $F_i$ une extension finie de $F$ de degré notée $d_i\geq 1$. Soit $V_i$ un $F_i$-espace vectoriel de dimension $n_i\geq 1$. Soit $\thgo_-= \prod_{i\in I} \tgl_{F_i}(V_i)$ qu'on verra comme un $F$-espace vectoriel. Soit  $V^-=\oplus_{i\in I}  V_i$ vu comme $F$-espace vectoriel de dimension 
$$n_-=\sum_{i\in I} n_id_i.
$$
On a une application $F$-linéaire
\begin{equation}
  \label{eq:appli0}
\iota_{H_-}:  \thgo_-\to \tgl_F(V^-)
\end{equation}
donnée par 
$$((A_i,b_i,c_i)_{ i\in I}\mapsto (\oplus_{i\in I} A_i,\oplus_{i\in I} b_i,\oplus_{i\in I} c_i).$$
Par abus, pour la composante $c_i$ on a utilisé implicitement l'identification de $F$-espaces  vectoriels 
$$V_i^*=\Hom_{F_i}(V_i,F_i) \simeq \Hom_{F}(V_i,F)$$
donnée par 
$$c_i\mapsto (v_i\in V_i \mapsto \trace_{F_i/F}(c_i(v_i)).
$$
L'application \eqref{eq:appli0} est équivariante sous l'action du groupe $H_-=\prod_{i\in I}GL_{F_i}(V_i)$ qu'on identifie naturellement à un sous-$F$-groupe de $GL_F(V)$. On en déduit un morphisme encore noté $\iota_{H_-}$
\begin{equation}
  \label{eq:appli}
  \iota_{H_-}:  \Ac_{H_-}=\prod_{i\in I} \Ac_{V_i} \to \Ac_{V^-}.
\end{equation}

Soit $X_{H^-}=(A_i,b_i,c_i)_{ i\in I}\in \thgo_-$ et $X=(A,b,c)=\iota_{H_-}((A_i,b_i,c_i)_{ i\in I})$.  Pour tout $i\in I$, soit $\chi_{A_i,F}$, resp. $\chi_{A_i,F_i}$, le polynôme caractéristique de $A_i$ vu comme $F$-endomorphisme (resp. $F_i$-endomorphisme) de $V_i$. On dispose aussi du polynôme caractéristique $\chi_{A,F}$ de $A$. Pour tous polynômes $P$ et $Q$,  soit $\disc(P)$ le discriminant de $P$ et $\mathrm{Res}(P,Q)$ le résultant de $P$ et $Q$.
Soit $D_{H^-}^{V^-}$ la fonction régulière sur  $\Ac_{H_-}$ définie par
\begin{eqnarray*}
  D_{H^-}^{V^-}(X_{H^-})&=& \frac{ \disc(\chi_{A,F})  }{\prod_{i\in I} \disc(\chi_{A_i,F})}\\
&=& (\pm) \prod_{ (i,j) \in I^2\!,  \, i\not=j} \mathrm{Res}(\chi_{A_i,F},\chi_{A_j ,F})
\end{eqnarray*}
où $(\pm)$ désigne un signe non explicité.

\begin{lemme}
\label{lem:compatibilite}
Avec les notations ci-dessus, on a 
\begin{eqnarray*}
  d_{n^-}(X)&=& (\pm) \disc(\chi_{A,F}) \cdot \prod_{i\in I} N_{F_i/F}(\frac{d_{n_i}(A_i,b_i,c_i)}{\disc(\chi_{A_i,F_i})})\\
&=& (\pm)  D_{H^-}^{V^-}(X_{H^-}) \cdot \prod_{i\in I} \frac{\chi_{A_i,F} }{N_{F_i/F}(\disc(\chi_{A_i,F_i}))} \cdot \prod_{i\in I} N_{F_i/F}(d_{n_i}(A_i,b_i,c_i))
\end{eqnarray*}
où $(\pm)$ est un signe non précisé et dans la second ligne chacun des trois facteurs est régulier.
 \end{lemme}

 \begin{preuve}
Le passage de la première ligne à la seconde est immédiat. Montrons donc la première égalité.

   Traitons d'abord le cas où $I=\{1,2\}$ et $F_1=F_2=F$. Soit $e_1$ et $e_2$ des bases respectives de $V_1$ et $V_2$. Soit $e=e_1\cup e_2$ la base de $V^-$ qui s'en déduit. Pour tout $(A,b,c)\in \tgl(V^-)$, on a alors
$$d_{n_-}(A,b,c)= {\det}_e(b,Ab,\ldots,A^{n_--1}b) {\det}_{e^*}(c,cA,\ldots,cA^{n_--1})
$$ 
où ${\det}_e$ et ${\det}_{e^*}$ désignent  le déterminant pris respectivement dans la base $e$ de $V^-$ et dans la base duale $e^*$ de $(V^-)^*$. Pour $i\in \{1,2\}$, soit $X_i=(A_i,b_i,c_i)\in \tgl(V^-)$. Soit $(A,b,c)=\iota_{H^-}(X_1,X_2)$. Pour alléger les notations, on pose $P_1=\chi_{A_1,F}$. Par des manipulations sur les colonnes,  on voit qu'on a 
\begin{eqnarray*}
  {\det}_e(b,Ab,\ldots,A^{n_--1}b)&=& {\det}_e(b_1+b_2,A_1b_1+A_2b_2,\ldots,A_1^{n_--1}b_1+A_2^{n_--1}b_2)\\
&=& {\det}_e(b_1+b_2,\ldots,A_1^{n_1-1}b_1+A_2^{n_1-1}b_2, P_1(A_2)b_2, \ldots, P_1(A_2)A_2^{n_2-1}b_2)\\
&=& {\det}_{e_1}(b_1,\ldots,A_1^{n_1-1}b_1) {\det}_{e_2}(P_1(A_2)b_2, \ldots, P_1(A_2)A_2^{n_2-1}b_2)\\
&=& {\det}_{e_1}(b_1,\ldots,A_1^{n_1-1}b_1) {\det}_{e_2}(b_2, \ldots, A_2^{n_2-1}b_2) \det(P_1(A_2))
\end{eqnarray*}
Par un calcul analogue sur le déterminant ${\det}_{e^*}(c,cA,\ldots,cA^{n_--1})$, on aboutit immédiatement au résultat vu qu'on a
 $$\det(P_1(A_2))^2=(\pm) \frac{\disc(\chi_{A,F})}{\disc(\chi_{A_1,F}) \disc(\chi_{A_2,F})}.
$$

Par une récurrence immédiate, on obtient le lemme pour un ensemble $I$ quelconque pour lequel $F_i=F$ pour tout $i\in I$.

Traitons ensuite le cas du singleton $I=\{1\}$ et d'une extension $F_1$ de $F$. Fixons une base de $V_1$ de sorte qu'on identifie $V_1$ à $F_1^{n_1}$, de même pour son dual etc. Pour vérifier la formule annoncée, on peut faire un changement de base de $F$ à $\bar{F}$ une clôture algébrique de $F$. Soit $\Gamma=\Hom_{F}(F_1,E)$ (il s'agit d'homomorphismes de corps) et pour tout ensemble $X$ soit $X^\Gamma$ l'ensemble des applications de $\Gamma$ dans $X$. On a un isomorphisme de $\bar{F}$-espaces vectoriels
$$ \tgl_{F_1}( F_1^{n_1})\otimes_F \bar{F} \simeq \tgl_{\bar{F}}( E^{n_1})^{\Gamma}
$$
donné par $(A_1,b_1,c_1)\otimes 1 \mapsto (A^\sigma_1,b^\sigma_1,c^\sigma_1)_{\sigma \in \Gamma}$.
L'application $\iota_{H^-}$ devient 
$$ (A^\sigma_1,b^\sigma_1,c^\sigma_1)_{\sigma \in \Gamma}\mapsto (\oplus_{\sigma\in \Gamma} A^\sigma_1,\oplus_{\sigma\in \Gamma}b^\sigma_1,\oplus_{\sigma\in \Gamma}c^\sigma_1).
$$
D'après ce que l'on vient de faire, on a 
\begin{eqnarray*}
  d_{n^-}(\iota_{H^-}(A_1,b_1,c_1))&=&\disc(\chi_{A,F}) \prod_{\sigma\in \Gamma} \frac{d_{n_1}(A^\sigma_1,b^\sigma_1,c^\sigma_1)}{\disc_{A^\sigma,\bar{F}}  }\\
&=& \disc(\chi_{A,F}) N_{F_1/F }(\frac{d_{n_1}(A_1,b_1,c_1)}{ \disc_{A_1,F_1}}).
\end{eqnarray*}

Le cas général se déduit aussitôt des cas que nous venons de traiter.
\end{preuve}

 Soit $\Ac_{H_-}'$ l'ouvert de  $\Ac_{H_-}$ défini par $D_{H^-}^{V^-}\not=0$. Le morphisme \eqref{eq:appli} est étale sur l'ouvert $\Ac_{H_-}'$ (cf. \cite{Z1} appendice B). Soit $\Ac_{H_-}^{(n^-)}$ l'ouvert de  $\Ac_{H_-}$ défini comme l'image inverse par $\iota_{H^-}$ de $\Ac^{\rs}_{V^-}$. Il résulte du lemme \ref{lem:compatibilite} qu'on  a
 \begin{equation}
   \label{eq:inclusion-}
   \Ac_{H_-}^{(n^-)}\subset \Ac_{H_-}'.
 \end{equation}

\end{paragr}

\begin{paragr}\label{S:thgo}
  Soit $V^+$ un $F$-espace vectoriel de dimension $r\geq 0$. Soit $\thgo_+=\tgl_{F}(V^+)$ muni de l'action de $H_+=GL_F(V^+)$ de quotient $\Ac_{H_+}=\Ac_{V^+}$. Soit le $F$-espace
$$\thgo=\thgo_+ \times \thgo_-
$$
muni de l'action du $F$-groupe $H=H_+\times H_-$ et 
$$a: \thgo \to \Ac_H=  \Ac_{H_+} \times  \Ac_{H_-}
$$
le morphisme canonique de $\thgo$ vers son quotient. Soit l'ouvert $\Ac_H'=\Ac_{H_+}^{\rs} \times \Ac_{H_-}'$ et $\thgo'$ l'ouvert de $\thgo$ obtenu comme l'image inverse de $\Ac_H'$. 
\end{paragr}

\begin{paragr} Soit $V=V^+\oplus V^-$ et $\tggo=\tgl_F(V)$ muni de l'action de $G=GL_F(V)$. Soit $n=r+\sum_{i\in I}n_i d_i$.
En composant \eqref{eq:iota-s} avec \eqref{eq:appli0}, on obtient un morphisme
\begin{equation}
  \label{eq:iotaH}
  \iota_H : \thgo' \to \tggo
\end{equation}
qui est $H$-équivariant lorsqu'on identifie naturellement $H$ à un sous-groupe de $G$. Il induit donc un morphisme homonyme sur les quotients
\begin{equation}
  \label{eq:iotaHquotient}
\iota_H : \Ac_H' \to \Ac.
\end{equation}
Ce dernier morphisme est étale (cf. \cite{Z1} appendice B). Soit $\Ac_H^{\Grs}$ l'ouvert de $\Ac_H'$ défini comme l'image inverse de l'ouvert $\Ac^{\rs}$. Il résulte de l'inclusion \eqref{eq:inclusion-} et des lemme \ref{lem:iso-r} et \ref{lem:entierk} qu'on a 
\begin{equation}
  \label{eq:egalite(n)}
  \Ac_H^{\Grs}=\Ac_{H_+}^{\rs} \times \Ac_{H_-}^{(n_-)}.
\end{equation}

On a ensuite un isomorphisme (\emph{loc. cit.})
\begin{equation}
  \label{eq:isocrucial}
   G \times^H \thgo' \to \tggo \times_{\Ac} \Ac_H'
 \end{equation}
qui est donné par $(g,Y) \mapsto (g\cdot \iota_H(Y), a(Y))$ et  où la source désigne le quotient de  $G \times \thgo'$ par l'action à droite (libre) de $H$ donnée par $(g,Y)\cdot h=(gh,h^{-1}\cdot Y)$.
\end{paragr}

\begin{paragr}\label{S:surA}
On suppose de plus que $F$ est un corps de nombres ; on utilise les notations de la section \ref{ssec:cdn}. On aura besoin des contructions et des résultats précédents sur une base un peu plus générale. On procède comme suit : on considère $A$ l'anneau des entiers de $F$ \og hors $S$ \fg{} où $S$ est un ensemble fini de places de $F$,  assez grand, qui  contient les places archimédiennes. Soit $A_i$ la clôture intégrale de $A$ dans $F_i$. Lorsque $S$ est assez grand, $A$ est principal, $A_i$ est un $A$-module libre de rang $d_i$ et le morphisme $A\to A_i$ est étale. En considérant $V$ un $A$-module libre de rang $n$, on définit de manière évidente un $A$-module libre $\tggo=\tgl_A(V)$ muni de l'action du $A$-schéma en groupes réductifs $GL_{A}(V)$ et on obtient un $A$-schéma quotient $\Ac$. De même, on définit $\tgl_A(V^+)$ et $\tgl_{A_i}(V_i)$ puis $\thgo$. Par restriction des scalaires à la Weil, ce dernier est muni d'une action du $A$-schéma en groupes réductifs $H=H_+\times H_-$ avec $H_+=GL_A(V^+)$ et $H_-=\prod_{i\in I} \Res_{A_i/A}GL_{A_i}(V_i)$. Soit $\Ac_H$ le quotient. On a un ouvert $\Ac_H'$ défini comme avant.
\end{paragr}

\begin{paragr} On continue avec la situation du paragraphe précédent.  Soit $v\notin S$ une place finie de $F$ et  $\oc_v$ l'anneau des entiers du complété $F_v$.
 
  \begin{lemme}\label{lem:integrite}
    Soit $a\in \Ac_H'(\oc_v)$. 
    \begin{enumerate}
    \item Pour tous $Y\in \thgo_a(F_v)$ et $g\in G(F_v)$, les assertions suivantes sont équivalentes :
      \begin{enumerate}
      \item $g^{-1}\cdot \iota_H(Y)\in \tggo(\oc_v)$.
      \item Il existe $h\in H(F_v)$ tel que $g\in hG(\oc_v)$ et $h^{-1}\cdot Y\in \thgo(\oc_v)$
      \end{enumerate}
    \item Pour tout $Y\in \thgo_a(F_v)$, les assertions suivantes sont équivalentes
      \begin{enumerate}
      \item $\iota_H(Y)\in \tggo(\oc_v)$.
      \item $Y\in \thgo(\oc_v)$.
      \end{enumerate}
    \item Pour tout $g\in G(F_v)$ et tout $X\in \tggo^{\rs}(\oc_v)$ les assertions suivantes sont équivalentes
      \begin{enumerate}
      \item $g^{-1}\cdot X\in \tggo^{\rs}(\oc_v)$.
      \item $g\in G(\oc_v)$.
      \end{enumerate}
          \end{enumerate}
  \end{lemme}

  \begin{preuve}
Prouvons d'abord l'assertion 1. Il est clair que 1.(b) implique 1.(a). Supposons 1.(a).   Le couple $(g^{-1},Y)$ définit donc un élément de $(G \times^H \thgo')(F_v)$ dont l'image $(g^{-1}\cdot \iota_H(Y),a)$ par  \eqref{eq:isocrucial} est un $\oc_v$-point de  $\tggo \times_{\Ac} \Ac_H'$.  Comme \eqref{eq:isocrucial} est un $\oc_v$-isomorphisme, cet élément appartient en fait à $(G \times^H \thgo')(\oc_v)$. Or l'application naturelle
    \begin{equation}
      \label{eq:opoint}
       (G(\oc_v)\times \thgo'(\oc_v) ) \to (G \times^H \thgo')(\oc_v)
     \end{equation}
est surjective. En effet, la surjectivité est vraie pour les points à valeurs dans le corps résiduel par le lemme de Lang. La surjectivité de \eqref{eq:opoint} est alors une conséquence de la lissité du morphisme canonique $G\times \thgo' \to G \times^H \thgo'$.  Par conséquent, il existe un couple $(g_0,Y_0)\in  G(\oc_v)\times \thgo'(\oc_v)  $  et $h\in H(F_v)$ tels que $g^{-1}h=g_0$ et $h^{-1}\cdot Y= Y_0$. Cela conclut.
 
L'assertion 2 est essentiellement un cas particulier de l'assertion 1. En effet, si l'on suppose 2.(a), alors d'après l'assertion 1, Il existe $h\in H(F_v)$ tel que $h\in G(\oc_v)$ et $h^{-1}\cdot Y\in \thgo(\oc_v)$. Or $H(F_v)\cap G(\oc_v)=H(\oc_v)$ d'où $Y\in \thgo(\oc_v)$.

Enfin, l'assertion 3 résulte immédiatement de \cite{rallSchiff} proposition 6.2.
  \end{preuve}

\end{paragr}

\begin{paragr}
Soit $v$ une place quelconque de   $F$.

\begin{lemme}
  \label{lem:compacite} 
Soit $a\in\Ac_H'(F_v)$ et  $\Omega\subset \tggo(F_v)$ un ensemble compact. Il existe un compact $C$ de $G(F_v)$ qui vérifie la propriété suivante. Pour tous $g \in G(F_v)$ et $Y\in \thgo_a(F_v)$ tels que
$$g^{-1}\iota_H(Y)\in \Omega$$
on a $g\in H(F_v)C$.
\end{lemme}

\begin{preuve}
  Le couple $(g^{-1},Y)$ définit un élément de $(G \times^H \thgo')(F_v)$ dont l'image $(g^{-1}\cdot \iota_H(Y),a)$ par  \eqref{eq:isocrucial} appartient à un compact qui ne dépend que de $a$ et $\Omega$. Il en donc est de même pour $(g^{-1},Y)$ ainsi que de son image dans $(G/H)(F_v)=G(F_v)/H(F_v)$.
\end{preuve}
\end{paragr}

\section{Combinatoire des cônes}\label{sec:combi}

\subsection{Sous-espaces paraboliques}\label{ssec:parab}

\begin{paragr}\label{S:Gtilde}
  Soit $V$ un espace vectoriel de dimension $n\geq 0$ sur $F$ un corps de caractéristique $0$. Soit $G=GL(V)$. On considère l'espace $\tlV=V\oplus F e_0 $ de dimension $n+1$ (pour un certain vecteur $e_0$) et $\tilde{G}=GL(\tlV)$. On identifie $G$ au sous-groupe de $\tilde{G}$ qui stabilise $V$ et fixe $e_0$.  Contrairement aux notations de la section \ref{ssec:notations-combi}, $\tggo$ ne désigne pas l'algèbre de Lie de $\tlG$ (qui n'apparaîtra pas dans cette section). On réserve la notation  $\tggo$ pour l'espace $\tgl(V)$.
\end{paragr}

\begin{paragr}[Les projections $r_i$ et $\hat{r}_i$.] --- \label{S:a0}   Soit $(e_1,\ldots,e_n)$ une base de $V$ d'où une base  $(e_1,\ldots,e_n,e_0)$ de $\tlV$. Soit $T_0 \subset \tlT_0$ les sous-tores maximaux respectifs de $G$ et $\tlG$ qui stabilisent les droites engendrées par les vecteurs de base. Ce sont donc des tores déployés maximaux respectivement de $G$ et $\tlG$. 

On pose alors   
$$
\ago_0=\ago_{T_0} \text{   et   }  \ago_{\tl}=\ago_{\tlT_0}
$$ 
et on a une inclusion naturelle $\ago_0\subset \ago_{\tl}$. Le choix de la base identifie cette inclusion à celle de $\RR^n$ dans $\RR^{n+1}$ qui consiste à rajouter une coordonnée nulle en position $n +1$. Le produit scalaire choisi sur ces espaces s'identifie au produit scalaire canonique. 

Soit $ \ago_{\tl}^{0} $ l'orthogonal de $ \ago_{0} $ dans $ \ago_{\tl} $.
Les deux décompositions suivantes de $ \ago_{\tl} $ en somme directe (non-orthogonale) vont intervenir :
\[ 
\ago_{\tl} = \ago_{0} \oplus \ago_{\tlG} = \ago_{\tl}^{\tlG} \oplus \ago_{\tl}^{0}.
\]
Soit $ r_{1} $, $ r_{2} $, $\hat r_{1}$  et $ \hat r_{2} $ les projections respectivement sur  $\ago_{0}, \ago_{\tlG},\ago_{\tl}^{\tlG}$ et $ \ago_{\tl}^{0}$ relativement à ces décompositions.
\end{paragr}

\begin{paragr}\label{S:parab}  Soit 
$$(0)=W_0\subsetneq W_1 \subsetneq \ldots \subsetneq W_r=V$$
un drapeau. Soit $i$ et $j$ deux éléments de $\{0,1,\ldots,r\}$ qui vérifient $0\leq j-i\leq 1$. On note $\tilde{P}$ une telle donnée. On lui associe les objets suivants

\begin{itemize}
\item le sous-groupe parabolique $P\subset G$ qui stabilise le drapeau et $\pgo$ son algèbre de Lie.
\item le sous-espace $V_{\tilde{P}}=W_i$ de $V$ ;
\item le sous-espace $V_{\tilde{P}}'=(W_j)^\perp$ de $V^*$ ;
\item le sous-espace $\tpgo=\pgo\oplus (V_{\tilde{P}}')^\perp  \oplus V_{\tilde{P}}^\perp$ de $\tggo$ ;
\item le sous-espace $\tngo=\tngo_{\tilde{P}}=\ngo\oplus V_{\tilde{P}}  \oplus V_{\tilde{P}}'$ de $\tggo$ où $\ngo$ est l'algèbre de Lie du radical unipotent $N_P$ de $P$.
\end{itemize}
Observons que $\tpgo$ et $\tngo$ sont tous deux  stables sous l'action du sous-groupe $P$ et qu'on a toujours $  V_{\tilde{P}}'\subset V_{\tilde{P}}^\perp   $ avec égalité si et seulement $i=j$.
\end{paragr}

\begin{paragr} La donnée $\tilde{P}$ est  équivalente à celle d'un sous-groupe parabolique de $\tilde{G}$ dont l'intersection avec $G$ est un sous-groupe parabolique de $G$. Les groupes $\tlP$ et $P$ contiennent respectivement $\tlT_0$ et $T_0$ si et seulement si le drapeau $W_\bullet$ est semi-standard. Plus précisément, à la donnée  de $\tilde{P}$, on associe le drapeau
$$(0)=W_0\subsetneq W_1 \subsetneq W_i \subsetneq W_{i+1}\oplus Fe_0\subsetneq  \ldots \subsetneq W_r\oplus Fe =\tlV$$
si $j=i+1$ et
$$(0)=W_0\subsetneq W_1 \subsetneq W_i \subsetneq  W_{i}\oplus Fe_0  \subsetneq  W_{i+1}\oplus Fe \subsetneq \ldots \subsetneq W_r\oplus Fe =\tlV$$
si $j=i$. Dans la suite, on passera sans plus de commentaire d'un point de vue à l'autre.
\end{paragr}

\begin{paragr}[Facteur de Levi.] ---\label{S:Levi}  La donnée supplémentaire d'un facteur Levi $ \tilde{M}$ du sous-groupe parabolique  $\tilde{P}$ de $\tlG$ comme ci-dessus est la donnée d'une somme directe 
$$V=V_1\oplus\ldots\oplus V_r$$
telle que $W_i=V_1\oplus\ldots\oplus V_i$. On a donc une décomposition duale 
$$V^*=V_1^*\oplus\ldots\oplus V_r^*$$
et une décomposition de Levi $P=M N_P$ où $M=GL(V_1)\times \ldots\times GL(V_r)$.
On définit un sous-espace $\tmgo_{\tlP}$ stable par $M$ de $\tggo$ tel que 
 $$\tpgo=\tmgo_{\tlP}\oplus \tngo_{\tlP}$$
de la façon suivante :
\begin{itemize}
\item soit $ V_{\tilde{P}}' = V_{\tilde{P}}^\perp$ et l'on pose $\tmgo_{\tlP}=\mgo$ (algèbre de Lie de $M$) ;
\item soit on a des décompositions $( V_{\tilde{P}}')^\perp= V_{\tilde{P}}\oplus V_j$ et $ V_{\tilde{P}}^\perp= V_{\tilde{P}}'\oplus V_j^*$ et on pose 
$$\tmgo_{\tlP}=\mgo \oplus V_j  \oplus V_j^*= \gl(V_1)\oplus\ldots \oplus \tgl(V_j)\oplus\ldots \oplus\gl(V_r).$$
\end{itemize}
\end{paragr}

\begin{paragr}\label{S:nPQ} Si $\tlP\subset \tlQ$, on pose
$$\tngo_{\tlP}^{\tlQ}=\tmgo_{\tlQ}\cap \tngo_{\tlP}.$$
  \end{paragr}

\begin{paragr}[Décomposition d'un facteur de Levi.] --- \label{S:Levi-dec} Le facteur de Levi $\tlM$ de $\tlP$ se décompose  en 
 $$\tlM = \MM_{\tlP} \times \tlG_{\tlP}$$
où l'on pose
\begin{itemize}
\item pour $j=i$,  $\MM_{\tlP} = M=  \prod_{k}GL(V_{k})$   et $\tlG_{\tlP} = GL(Fe_0)$ ;
\item pour  $j = i+1$,  $\MM_{\tlP} = \prod_{k \neq j}GL(V_{k})$ et $\tlG_{\tlP} = GL(V_{j}\oplus Fe_0)$.
\end{itemize}
Cette décomposition ne dépend réellement que de $\tlM$ (et non du choix de $\tlP$ qui admet $\tlM$ comme facteur de Levi). On introduit aussi le sous-groupe $G_{\tlP}\subset \tlG_{\tlP}$ qui fixe $e_0$ et préserve $V_j$ si $j=i+1$. Si $j=i$ le groupe  $G_{\tlP} $ est trivial et on pose $\tggo_{\tlP}=(0)$.  Si $j=i+1$,   le groupe  $G_{\tlP} $ s'identifie à $GL(V_j)$ et l'on pose $\tggo_{\tlP}= \tgl(V_j)$. Dans chaque cas, $G_{\tlP}$ agit sur $\tggo_{\tlP}$.
\end{paragr}

\begin{paragr}[Espace $\zgo_{\tlP}$.] ---\label{S:zgoP}  Désormais, les sous-groupes paraboliques de $G$ (resp. $\tlG$) seront supposés contenir $T_0$ (resp. $\tlT_0$) et les constructions afférentes aux sous-groupes paraboliques seront relatives à ces sous-tores. Soit $\tlP$ un sous-groupe parabolique de $\tlG$ comme au §\ref{S:parab}.  Dans ce cas, $M_{\tlP}$ est l'unique facteur de Levi de $\tlP$ qui contient $\tlT_0$. La décomposition $M_{\tlP} = \MM_{\tlP} \times \tlG_{\tlP}$ induit des décompositions  orthogonales 
\[
\ago_{\tlP} = \zgo_{\tlP} \oplus \ago_{\tlG_{\tlP}}, \quad \ago_{\tlP}^* = \zgo_{\tlP}^* \oplus \ago_{\tlG_{\tlP}}^*,
\]
où, pour alléger les notations, on pose $\zgo_{\tlP}=\ago_{\MM_{\tlP}}$ et  $\zgo_{\tlP}^*=\ago_{\MM_{\tlP}}^*$. Notons qu'on a naturellement  $\zgo_{\tlP}\subset \ago_0$.  Dans les paragraphes suivants, les  sous-groupes paraboliques de $\tlG$ considérés contiendront $\tlT_0$. Pour un sous-groupe de Levi $\tlM$ de $\tlG$ contenant $\tlT_0$, on pose $\zgo_{\tlM}=\zgo_{\tlP}$ où $\tlP$ est un sous-groupe parabolique de $\tlG$ ayant $\tlM$ comme facteur de Levi. Cette définition ne dépend pas du choix de $\tlP$.
\end{paragr}

\subsection{Fonctions caractéristiques de cônes}

\begin{paragr}[Ensembles $\Pi$ et $\widehat{\Pi}$.] ---\label{S:widePi} On continue avec les notations des §§\ref{S:parab} et \ref{S:Levi-dec}. Pour $1\leq k \leq i$ soit $\delta_k\in X^*(\MM_{\tlP})$  le  déterminant de l'action de $\MM_{\tlP}$ sur $W_k$. De même, pour $j\leq k \leq r$ soit  $\delta_k^-\in X^*(\MM_{\tlP})$   le  déterminant de l'action de $\MM_{\tlP}$ sur $W_k^\perp$. Soit  $\widehat{\Pi}_{\tlP}$  l'ensemble des $\delta_k$,   $1\leq k \leq i$,  et $\delta_k^{-}$ pour $j\leq k \leq r$. On identifie naturellement $\widehat{\Pi}_{\tlP}$  à une partie de $\zgo_{\tlP}^*$.

Soit $\tlP\subset \tlQ$ une inclusion de sous-groupes paraboliques. On a $\widehat{\Pi}_{\tlQ} \subset \widehat{\Pi}_{\tlP}$.  
On définit 
\[
\zgo^{\tlQ}_{\tlP} := \{X \in \zgo_{\tlP} | \delta(X) = 0, \ \forall \, \delta \in \widehat{\Pi}_{\tlQ}\}.
\]

Soit  $\Pi_{\tlP}\subset \zgo_{\tlP}^*$ l'ensemble   obtenu par restriction des éléments de $\Delta_{\tlP}$. Soit $\Pi_{\tlP}^{\tlQ}\subset \Pi_{\tlP}$ le sous-ensemble obtenu par restriction des éléments de  $\Delta_{\tlP}^{\tlQ}$.  On a $ \zgo_{\tlQ} \subset \zgo_{\tlP}$ et plus précisément 
\[
\zgo_{\tlQ} = \zgo_{\tlP} \cap \ago_{\tlQ}=\{X \in \zgo_{\tlP} | \, \gamma(X) = 0, \ \forall \, \gamma \in \Pi^{\tlQ}_{\tlP}\}.
\]

On a de plus
\[
\zgo_{\tlP} = \zgo_{\tlQ} \oplus \zgo^{\tlQ}_{\tlP}
\]
et par dualité
\[
\zgo_{\tlP}^* = \zgo_{\tlQ}^* \oplus (\zgo^{\tlQ}_{\tlP})^*.
\]
On définit alors $\widehat{\Pi}^{\tlQ}_{\tlP}$ comme la projection des éléments de $ \widehat{\Pi}_{\tlP} \setminus  \widehat{\Pi}_{\tlQ}$ sur $(\zgo_{\tlP}^{\tlQ})^{*}$.

On a 
\begin{eqnarray*}
  \ago_{\tl}^* &=& \ago_{\tlP}^* \oplus  (\ago_{\tl}^{\tlP})^*\\
&=& (\zgo_{\tlP}^{\tlQ})^* \oplus \zgo_{\tlQ}^* \oplus \ago_{\tlG_{\tlP}}^* \oplus  (\ago_{\tl}^{\tlP})^*\\
\end{eqnarray*}
de sorte que les ensembles $\Pi_{\tlP}^{\tlQ}$ et  $\widehat{\Pi}_{\tlP}^{\tlQ}$ seront  vus comme des parties de  $\ago_{\tl}^*$. 
\end{paragr}

\begin{paragr}[Fonctions $\sigma$ et $\hat{\sigma}$.] --- \label{S:sigma}Soit
$$
\sigma_{\tlP}^{\tlQ} \text{ et } \hat{\sigma}_{\tlP}^{\tlQ}
$$
respectivement les fonctions caractéristiques respectivement des cônes
$$\{H\in \ago_{\tl} \mid \gamma(H) > 0, \  \forall \ \gamma \in  \Pi_{\tlP}^{\tlQ} \}$$
et
$$\{H\in \ago_{\tl} \mid \delta(H) > 0, \  \forall \ \delta \in  \widehat{\Pi}_{\tlP}^{\tlQ} \}.$$
Notons  les relations suivantes pour tout $H\in \ago_{\tl}$ :
\begin{equation}
  \label{eq:sig-tau}
  \tau_{\tlP}^{\tlQ} (H) = \sigma_{\tlP}^{\tlQ} (r_{1}(H)), \quad
\hat \sigma_{\tlP}^{\tlQ} (H) = \hat \tau_{\tlP}^{\tlQ} (\hat r_{1}(H)),
\end{equation}
où $\hat r_{1}$ et $r_{1}$ sont les projections définies au §\ref{S:a0}.
\end{paragr}

\begin{paragr}[Lemmes géométriques.] 

\begin{lemme}\label{lem:newRoots} Soient $\tlQ \subset \tlP$ des sous-espaces paraboliques de $\tlG$. Alors, la base $\Pi_{\tlQ}^{\tlP}$ 
est une base obtuse de $(\zgo_{\tlQ}^{\tlP})^{*}$ et $\widehat{\Pi}_{\tlQ}^{\tlP}$ est une base aiguë de $(\zgo_{\tlQ}^{\tlP})^{*}$. 
De plus, les bases $\widehat{\Pi}_{\tlQ}^{\tlP}$ et $\Pi_{\tlQ}^{\tlP}$ sont duales l'une de l'autre (pour le produit scalaire fixé sur $\ago_{\tl}$).
\end{lemme}

\begin{preuve}
Soit $\tlT_0\subset \tlB\subset \tlQ$ un sous-groupe de Borel.  Les poids de $\tlT_0$ sur les vecteurs de la base $(e_{1}, \ldots, e_{n})$ forment une base de $\ago_0^{*}$ notée $(\eps_i)_{1\leq i \leq n}$. Soit $i$ la dimension du sous-espace $V_{\tlB}\subset V$ (cf. §\ref{S:parab}). On peut  supposer que la base est indexée de sorte qu'on ait
\begin{align*}
\Pi_{\tlB} &=& \{\eps_{1} - \eps_{2}, \ldots, \eps_{i-1} - \eps_{i}, \eps_{i}, -\eps_{i+1}, \eps_{i+1} - \eps_{i+2}, \ldots, \eps_{n-1} - \eps_{n}\}, \\
\widehat{\Pi}_{\tlB} &=& \{\eps_{1}, \eps_{1} + \eps_{2}, \eps_{1} + \cdots + \eps_{i},  -\eps_{n}, -(\eps_{n} +\eps_{n-1}), \ldots, -(\eps_{n} + \cdots + \eps_{i+1})\}
\end{align*}
où $\eps_{k} = 0$ si $k \notin \{1, \ldots, n\}$.  On vérifie immédiatement l'énoncé pour $\tlP=\tlG$ et $\tlQ=\tlB$. L'ensemble $\Pi_{\tlQ}$ est aussi l'ensemble des projections des éléments de $\Pi_{\tlB} \setminus \Pi_{\tlB}^{\tlQ}$ 
sur $\zgo_{\tlQ}^{*}$.  Il en résulte que, en vertu du lemme 1.2.4 de \cite{labWal}, la base $\Pi_{\tlQ}$ de $\zgo_{\tlQ}^{*}$ est aussi obtuse de même que  $\Pi_{\tlQ}^{\tlP}\subset \Pi_{\tlQ}$. Les ensembles $\Pi_{\tlB}^{\tlP}$ et $\widehat{\Pi}_{\tlB} \setminus \widehat{\Pi}_{\tlP}$ sont en dualité pour le produit scalaire.
Il en résulte que $\widehat{\Pi}_{\tlQ}^{\tlP}$ et $\Pi_{\tlQ}^{\tlP}$ sont aussi en dualité. 
Finalement, le fait que $\widehat{\Pi}_{\tlQ}^{\tlP}$ soit aiguë résulte maintenant de cette dualité 
et du fait que $\Pi_{\tlQ}^{\tlP}$ soit obtuse, comme le démontre le lemme 1.2.6. de \emph{loc. cit.}. 
\end{preuve}

\begin{lemme}\label{lem:langlands}
 Pour tous sous-espaces paraboliques $\tlQ \subset \tlP$ et tout $H \in \ago_{\tl}$ 
on a 
\[
\sum_{\tlQ \subset \tlS \subset \tlP}\varepsilon_{\tlQ}^{\tlS}
\sigma_{\tlQ}^{\tlS}(H)
\hat{\sigma}_{\tlS}^{\tlP}(H) = 
\begin{cases}
0 \text{ si }\tlQ \neq \tlP, \\
1 \text{ si }\tlQ = \tlP.
\end{cases}
\]
\end{lemme}

\begin{preuve}
Il suffit de considérer le cas $H\in \zgo_{\tlQ}^{\tlP}$, restriction que l'on fait dans la suite. Rappelons qu'on a une décomposition $M_{\tlP}=\MM_{\tlP}\times \tlG_{\tlP}$.  Soit  $Q_1=   \MM_{\tlP}\cap \tlQ$ et $\tlQ_2= \tlG_{\tlP}\cap \tlQ$. On a alors une décomposition 
$$\zgo_{\tlQ}^{\tlP}=\ago_{Q_1}^{\MM_{\tlP}} \oplus \zgo_{\tlQ_2}
$$
où la notation $\zgo_{\tlQ_2}$ est celle du §\ref{S:zgoP} relativement au couple $(G_{\tlP},\tlG_{\tlP})$. Suivant cette décomposition, on écrit $H=H_1+H_2$.
L'application   $\tlR\mapsto  (R_1,\tlR_2)=(\MM_{\tlP}\cap \tlR,\tlG_{\tlP}\cap \tlR)$ induit une bijection de l'ensemble des $\tlR$ tels que $\tlQ \subset \tlR \subset \tlP$ sur  le produit 
$$\{R_1\mid Q_1\subset R_1\subset   \MM_{\tlP}    \} \times \{\tlR_2\mid  \tlQ_2\subset  \tlR_2\subset \tlG_{\tlP} \}.
$$
  Pour cette bijection, pour tout $\tlQ \subset \tlS \subset \tlR \subset \tlP$, on a  
$$\sigma_{\tlS}^{\tlR}(H)=\tau_{S_1}^{R_1}(H_1) \sigma_{\tlS_2}^{\tlR_2}(H_2).
$$
et 
$$\hat{\sigma}_{\tlS}^{\tlR}(H)=\hat{\tau}_{S_1}^{R_1}(H_1) \hat{\sigma}_{\tlS_2}^{\tlR_2}(H_2).$$ 

Il s'ensuit que la somme en question est égale à 
\[
\big(  \sum_{Q_1 \subset R_1' \subset \MM_{\tlP}}\varepsilon_{Q_1}^{R_1}
\tau_{Q_1}^{R_1}(H_1)
\hat{\tau}_{R_1}^{\MM_{\tlP}}(H_1) \big).
\big( \sum_{\tlQ_2\subset \tlR_2  \subset \tlG_{\tlP}}\varepsilon_{\tlQ_2}^{\tlR_2}
\sigma_{\tlQ_2}^{\tlR_2}(H_2)
\hat{\sigma}_{\tlR_2}^{\tlG_{\tlP}}(H_2) \big)
\]

Le premier facteur vaut $1$ si $Q_1=\MM_{\tlP}$ et $0$ sinon (lemme de Langlands  cf. \cite{labWal}, proposition 1.7.2). L'autre facteur  vaut $1$ si $\tlQ_2=\tlG_{\tlP}$ et $0$ sinon ; la preuve est la même que celle du lemme de Langlands, il suffit de remplacer les ensembles de racines simples par les ensembles $\Pi$ correspondant, la preuve de la proposition 1.7.2 de \emph{loc. cit.} reposant seulement sur 
le lemme 1.2.6 de \emph{loc. cit.} qui est le lemme \ref{lem:newRoots} ci-dessus.
\end{preuve}
\end{paragr}

\begin{paragr}[Les fonctions $ \upla $.] --- Soit $ \upla_{\tlP} $ le produit du  déterminant de l'action du tore $ A_{\MM_{\tlP}}$ sur $ V_{\tlP} $ et du déterminant de son action sur $ V_{\tlP}' $ (à ce propos, rappelons que $G$ agit à gauche  sur $V^*$  par $g\cdot c=c g^{-1}$). On obtient ainsi un élément de $\zgo_{\tlP}^{*}$ qui vérifie l'identité (écrite de manière additive)
$$ \upla_{\tlP} = 2\rho_{\tlP} - 2\rho_{P} $$
où, comme d'habitude, $\rho_{\tlP} $ et $\rho_{P} $ désignent la demi-somme des éléments de $\Sigma_{\tlP}$ et $\Sigma_P$ respectivement.
\end{paragr}

\begin{paragr}
Pour tout sous-espace parabolique $ \tlP $ et tous $H,X\in \ago_{\tl}$, soit
\begin{equation}
  \label{eq:B}
  \mathrm{B}_{\tlP}^{\tlG}(H,X) = 
\sum_{\tlR\in \fc^{\tlG}(\tlP) }\varepsilon_{\tlP}^{\tlR}
\hat{\sigma}_{\tlR}(H - X) \sigma_{\tlP}^{\tlR}(H).
\end{equation}
Lorsque le contexte est clair, on omet l'exposant $\tlG$.

\begin{lemme}\label{lem:upGamma}
\begin{enumerate}
\item On a
\[
\hat{\sigma}_{\tlP}(H-X) = \sum_{\tlR\in \fc^{\tlG}(\tlP) } \varepsilon_{\tlR}^{\tlG}
\hat{\sigma}_{\tlP}^{\tlR}(H)
\mathrm{B}_{\tlR}(H,X), \quad X,H \in \ago_{\tl}.
\]
\item Pour tout $X \in\ago_{\tl}$, la fonction $H\in  \zgo_{\tlP} \mapsto \mathrm{B}_{\tlP}(H,X)$ est à support compact. 
\item Pour tous $\tlP\subset \tlR$, il existe un  polynôme $ p_{\tlP,\tlR} $ sur $ \zgo_{\tlR} $ tel que pour tout $X\in \ago_{\tl}$ on ait 
\[ 
\int_{\zgo_{\tlP}}e^{\upla_{\tlP}(H)}\mathrm{B}_{\tlP}(H,X)dH = 
\sum_{\tlP \subset \tlR} e^{\upla_{\tlR}(X)}p_{\tlP,\tlR}(X_{\tlR})
\]
où $ X_{\tlR} $ désigne la projection orthogonale de $ X $ sur $ \zgo_{\tlR} $. 
\end{enumerate}
\end{lemme}

\begin{preuve}
  Le point 1 découle immédiatement du lemme \ref{lem:langlands}. 
  Comme on l'a déjà remarqué, le lemme \ref{lem:newRoots} permet d'appliquer l'analyse de \cite{labWal} 
  sans changement, il suffit de remplacer $\Delta$ par $\Pi$. 
  Le point 2 résulte alors du lemme 1.8.1 de loc. cit.
  Le dernier s'obtient de même façon que le lemme 4.3 de \cite{leMoi}. 
\end{preuve}
\end{paragr}

\begin{paragr} La fonction $\mathrm{B}$ introduite en \eqref{eq:B} est le pendant de la fonction $\Gamma'$ d'Arthur définie dans \cite{arthur2} section 2. Rappelons que cette fonction $\Gamma'$ est définie par la relation \eqref{eq:B} où l'on remplace $\sigma$ et $\hat \sigma$ respectivement par $\tau$ et $\hat{\tau}$. Ces deux fonctions sont étroitement reliées comme le montre le lemme suivant.

\begin{lemme}\label{lem:GammaB}
  Pour tout $ H \in \zgo_{\tlP} $ et $T\in \ago_{\tlP}^{\tlG}$ on a :
\[ 
\Gamma_{\tlP}'(H-T, \hat r_{2}(H) ) = 
\mathrm{B}_{\tlP}(H - r_{1}(T), r_{2}(T)).
\]
\end{lemme}

\begin{preuve}
En utilisant la définition de $\Gamma'$ et les formules \eqref{eq:sig-tau}, on a
\begin{eqnarray*}
  \Gamma_{\tlP}'(H-T, \hat r_{2}(H) ) &=& 
\sum_{ \tlP \subset \tlR }\varepsilon_{\tlP}^{\tlR}\hat \tau_{\tlR}(H - T - \hat r_{2}(H))
\tau_{\tlP}^{\tlR}(H -T) \\ 
&=& \sum_{ \tlP \subset \tlR } \varepsilon_{\tlP}^{\tlR}
\hat \tau_{\tlR}(\hat r_{1}(H- T))
\tau_{\tlP}^{\tlR}(H -T)\\
&=& \sum_{ \tlP \subset \tlR } \varepsilon_{\tlP}^{\tlR}
\hat \sigma_{\tlR}(H- T)
\sigma_{\tlP}^{\tlR}(r_1(H -T))\\
&=& \sum_{ \tlP \subset \tlR } \varepsilon_{\tlP}^{\tlR}
\hat \sigma_{\tlR}(H- r_1(T)-r_2(T))
\sigma_{\tlP}^{\tlR}(H -r_1(T))
\end{eqnarray*}
d'où le lemme.
\end{preuve}

On utilisera la variante suivante du lemme \ref{lem:upGamma}.

  \begin{lemme}\label{lem:upGamma2} Soit $\tlP\subsetneq \tlG$ un sous-groupe parabolique.
Pour tout sous-groupe parabolique $\tlP\subset \tlR\subset \tlG$, il existe un polynôme $ p_{\tlP,\tlR} $ sur $ \ago_{\tlG} $ et un élément $\la_{\tlR}$ \emph{non nul} de $\ago_{\tlP}^{\tlG,*}$ tels que pour tout $T\in \ago_{\tlP}^{\tlG}$ on ait 
\begin{equation}\label{eq:purExpDesc0}
\int_{\zgo_{\tlP}}e^{\upla_{\tlP}(H)}
\Gamma_{\tlP}'(H-T, \hat r_{2}(H) )dH= \sum_{\tlP \subset \tlR} e^{\la_{\tlR}(T)}p_{\tlP,\tlR}(r_2(T))
\end{equation}
où l'intégrande dans l'intégrale de gauche est à support compact.
\end{lemme}

\begin{preuve}
La compacité du support est immédiate d'après les lemmes \ref{lem:GammaB} et \ref{lem:upGamma} assertion 2.
En utilisant le lemme \ref{lem:GammaB} et un changement de variables, l'expression considérée est égale à 
$$  e^{\upla_{\tlP}(r_1(T))} \int_{\zgo_{\tlP}}e^{\upla_{\tlP}(H)}\mathrm{B}_{\tlP}(H,r_2(T))dH.
$$
Le lemme \ref{lem:upGamma} donne alors la forme voulue avec $\la_{\tlR}(T)= \upla_{\tlP}(r_1(T))+\upla_{\tlR}(r_2(T))$. Il reste à voir que $\la_{\tlR}$ induit une forme linéaire non triviale sur  $\ago_{\tlP}^{\tlG}$. Pour cela, on observe qu'on a 

$$\la_{\tlR}(T)= \upla_{\tlP}(T) +(\upla_{\tlR}-\upla_{\tlP})(r_2(T)).
$$
Notons que la projection orthogonale $T\mapsto T_{\tlP}^{\tlG}$ sur $a_{\tlP}^{\tlG}$ induit un isomorphisme de $\zgo_{\tlP}$ sur  $\ago_{\tlP}^{\tlG}$. On va en fait prouver que $T\mapsto \la_{\tlR}( T_{\tlP}^{\tlG})$ est non identiquement nulle sur $\zgo_{\tlP}$. Soit $\delta$ et $\delta'$ les éléments de $\zgo_P^*$ correspondant aux déterminants de $A_{\tlP}$ agissant respectivement sur $V_{\tlP}$ et   $V_{\tlP}'$. Soit $d$, $d'$, $e$ et $e'$ les dimensions respectives de $V_{\tlP}$, $V_{\tlP}'$, $V_{\tlR}$ et $V_{\tlR}'$. Soit $n=\dim(V)$. On a donc pour $T\in \zgo_{\tlP}$
\begin{eqnarray*}
  \upla_{\tlP}(T_{\tlP}^{\tlG})&=& (\delta+\delta'- \frac{d-d'}{n+1}(\delta-\delta'))(T).
\end{eqnarray*}
Par ailleurs, on a 
\begin{eqnarray*}
  (\upla_{\tlR}-\upla_{\tlP})(r_2(T_{\tlP}^{\tlG}))&=& (- \frac{e+d'-(e'+d)}{n+1}(\delta-\delta'))(T)\\
\end{eqnarray*}
et
\begin{eqnarray}\label{eq:forme}
\la_{\tlR}( T_{\tlP}^{\tlG})= (\frac{n+1-e+e'}{n+1}\delta+ \frac{n+1+e-e'}{n+1}\delta')(T).
\end{eqnarray}
On a toujours $n\geq d\geq e$ et $n\geq d'\geq e'$. Si $\delta=0$, on a alors $d=0$ et donc $e=0$. Puisque $\tlP\not=\tlG$, on a alors $\delta'\not=0$ sur $\zgo_{\tlP}$ et le coefficient de $\delta'$, à savoir $\frac{n+1-e'}{n+1}$ est non nul. Donc la forme linéaire \eqref{eq:forme} est non nulle sur $\zgo_{\tlP}$. De même, on traite le cas   $\delta'=0$. Si $\delta\not=0$ et $\delta'\not=0$, ces formes  sont linéairement indépendantes sur $\zgo_{\tlP}$. . Donc la forme linéaire \eqref{eq:forme} ne peut pas être nulle sur $\zgo_{\tlP}$ car, la somme des coefficients de $\delta$ et $\delta'$  dans  \eqref{eq:forme} vaut $2$.
\end{preuve}
\end{paragr}

\subsection{Descente et  combinatoire des cônes}\label{ssec:desc-combi}

\begin{paragr}
 On reprend les notations du §\ref{S:transv}. Pour $i\in I$, on dispose de $F_i$-espaces vectoriels $V_i$ de dimension $n_i$ et d'un $F$-espace vectoriel $V=V^+\oplus V^-$ où $V^-=\oplus_{i\in I} V_i$. Soit $\tilde{G}=GL_F(V\oplus F e_0)$ et son sous-groupe $G=GL_F(V)$.
\end{paragr}

\begin{paragr}[Sous-groupes $\tlM_1$.] ---  \label{S:M1}Pour tout $i\in I$, soit
$$V_i=\oplus_{l=1}^{n_i}  D_{i,l}$$
une décomposition en somme de  sous-$F_i$-espaces $D_{i,l}$ de dimension $1$. Soit 
$$\tlM_1=GL_F( V_+\oplus Fe_0 )\times \prod_{i\in I} \prod_{l=1}^{n_i} GL_F(D_{i,l})$$
qu'on identifie naturellement à un sous-groupe de Levi de $\tlG$. 
\end{paragr}

\begin{paragr}[L'application $\tlQ\mapsto \tlQ^-$.] --- \label{S:QQ-}Soit $\tlQ\in \fc^{\tlG}(\tlM_1)$. Comme on l'a vu au §\ref{S:parab}, à un tel élément est associé un drapeau
\begin{equation}
  \label{eq:drapeau}
  (0)=W_0\subsetneq W_1 \subsetneq \ldots \subsetneq W_s=V
\end{equation}
et deux éléments  $i_0$ et $j_0$ de $\{0,1,\ldots,s\}$ qui vérifient $0\leq j_0-i_0\leq 1$ ; en outre les sous-espaces $W_j$ sont des sommes directes de $D_{i,l}$ et $V^+$.

Soit $\tlH_i=GL_{F_i}(V_i \oplus F_ie_0)$ et $\tlL_{i}$ son groupe de Levi qui s'identifie naturellement à 
$$GL_{F_i}(F_ie_0)\times \prod_{l=1}^{n_i} GL_{F_i}(D_{i,l}).
$$ Soit $\tlQ_i\in \fc^{\tlH_i}(\tlL_i)$ l'élément associé aux données suivantes : 
\begin{itemize}
\item le drapeau  de $F_i$-espaces vectoriels 
\begin{equation}
  \label{eq:drapeau2}
  (0)=W_{i,0}\subsetneq W_{i,1} \subsetneq \ldots \subsetneq W_{i,s_i}=V_i
\end{equation}
où les $W_{i,l}$ sont obtenus par intersection des $W_j$ avec $V_i$ et élimination des doublons ;
\item le couple  d'entiers $(i_0',j_0')$ définis par la condition $W_{i,i_0'}=W_{i_0}\cap V_i$ et  $W_{i,j_0'}=W_{j_0}\cap V_i$. (Notons qu'on a $0 \leq j_0'-i_0'\leq 1$.)
\end{itemize}

Soit 
$$\tlH^-= \prod_{i\in I} \tlH_i \text{  et   } \tlM_0=  \prod_{i\in I} \tlL_i.$$ 
On pose alors 
$$\tlQ^-=  \prod_{i\in I} \tlQ_i
$$
Soit $\fc^{\tlH^-}(\tlM_0)=\prod_{i\in I} \fc^{\tlH_i}(\tlL_i)$ où chaque facteur est défini relativement au corps $F_i$. On peut définir de manière équivalente cet ensemble comme l'ensemble des sous-groupes paraboliques de $\tlH^-$ contenant $\tlM_0$ si l'on voit tous ces groupes comme des groupes sur $F$ par restriction des scalaires. En tout cas, on a $\tlQ^-\in \fc^{\tlH^-}(\tlM_0)$.
\end{paragr}

\begin{paragr}
 On a une décomposition (cf. §\ref{S:Levi-dec})
$$M_{\tlQ}=\MM_{\tlQ}\times \tlG_{\tlQ}
$$
associée à certaine décomposition 
$$V\oplus Fe _0= U_1\oplus \ldots \oplus (V^+\oplus U_{j_0} \oplus Fe_0) \oplus \ldots \oplus U_r
$$
où chaque $U_j$ est une somme directe de $F$-espaces $D_{i,l}$. Par exemple, $\MM_{\tlQ}=\prod_{j\not=j_0} GL_{F}(U_j)$. On  a une décomposition
$$V_i\oplus F_i {  e_0}= (U_1\cap V_i) \oplus \ldots \oplus ((U_{j_0}\cap V_i) \oplus F_ie_0) \oplus \ldots \oplus (U_r\cap V_i)
$$
donc un sous-groupe de Levi $\tlM_i $ de $\tlH_i$ muni de la décomposition
$$\tlM_i =\MM_i \times \tlG_i$$
avec $\MM_i= \prod_{j\not=j_0} GL_{F_i}(U_j\cap V_i)$ et $\tlG_i=  GL_{F_i}((U_{j_0}\cap V_i) \oplus F_i e_0) $. On peut alors identifier $\prod_{i\in I} \MM_i$ à un sous-$F$-groupe de $\MM_{\tlQ}$. On a donc une inclusion naturelle 
$$\ago_{\MM_{\tlQ}} \hookrightarrow \oplus_{i\in I} \ago_{\MM_i}.
$$
L'espace à gauche n'est autre que $\zgo_{\tlQ}$ et celui à droite est $\zgo_{\tlQ^-}$ défini selon  le §\ref{S:zgoP} dans le contexte du groupe produit $\tlH^-$ : en effet, pour tout $\tlR\in  \fc^{\tlH^-}(\tlM_0)$, on définit un espace $\zgo_{\tlR}$ qui ne dépend en fait que du facteur de Levi de $\tlR$ contenant $\tlM_0$ qu'on note $\tlL$. On pose aussi $\zgo_{\tlL}= \zgo_{\tlR}$ (ce sont les notations de  §\ref{S:zgoP} mais dans la situation du produit $\tlH^-$).
\end{paragr}

\begin{paragr}[Ensembles $\fc$.] --- \label{S:combi}  Pour tout $\tlR \in \fc^{\tlH^-}(\tlM_0)$, on définit des ensembles  $\fc_{\tlR}^{0}(\tlM_1)\subset \fc_{\tlR}^{}(\tlM_1)\subset\overline{\fc}_{\tlR}^{}(\tlM_1) \subset \fc^{\tlG}(\tlM_1)$ par les conditions
\begin{eqnarray*}
\overline{\fc}_{\tlR}^{}(\tlM_1) & = & \{\tlP \in \fc^{\tlG}(\tlM_1) | \,    \tlR\subset \tlP^{-} \}\\
\fc_{\tlR}^{}(\tlM_1) & = & \{\tlP \in \fc^{\tlG}(\tlM_1) | \, \tlP^{-} = \tlR\}, \\ 
\fc_{\tlR}^{0}(\tlM_1) & = &\{\tlP \in \fc_{\tlR}^{}(\tlM_1) | \, \zgo_{\tlP} = \zgo_{\tlR}\}.
\end{eqnarray*}
Si $\tlR$ est minimal, c'est-à-dire s'il admet $\tlM_0$ comme facteur de Levi, on a  $\zgo_{\tlM_1} = \zgo_{\tlM_0}=\zgo_{\tlR}$. En général, on a donc une inclusion $\zgo_{\tlR}\subset \zgo_{\tlM_1}$.

\begin{lemme}
  \label{lem:convexite}
Pour tout $\tlR\in  \fc^{\tlH^-}(\tlM_0)$, il existe un unique sous-groupe de Levi $\tlM$ de $\tlG$ contenant $\tlM_1$ tel que
$$ \fc_{\tlR}^{0}(\tlM_1) \subset\pc^{\tlG}(\tlM). 
$$
En outre,  l'ensemble $ \fc_{\tlR}^{0}(\tlM_1) $ est une famille convexe dans $ \pc(\tlM) $   au sens  de l'appendice \ref{App:domConvP}.
\end{lemme}

\begin{preuve}
  Soit $M_{\tlR}$ le facteur de Levi de $\tlR$ qui contient $\tlM_0$. On a alors  $\tlR=\prod_{i\in I} \tlR_i$ et $M_{\tlR}=\prod_{i\in I} \tlM_i$. Chaque $\tlM_i$ est le stabilisateur des sous-espaces qui apparaissent dans la décomposition
$$V_i\oplus F _i e_0=V_{i,1}\oplus \ldots \oplus V_{i,s_i}
$$
et il existe un unique indice $l_i$ tel que $e_0\in V_{i,l_i}$. De plus, $\tlR_i$ est le stabilisateur du drapeau 
$$(0)\subsetneq V_{i,1}\subsetneq V_{i,1}\oplus V_{i,2}\subsetneq \ldots$$
Le sous-groupe de Levi $\tlM$ est nécessairement le stabilisateur des sous-espaces dans la décomposition
\begin{equation}
  \label{eq:decompo}
V \oplus F e_0=\oplus_{(i,l)} V_{i,l} \oplus W
\end{equation}
où l'on somme sur les couples $(i,l)$ tels que $i\in I$, $1\leq l \leq s_i$ et $l\not=l_i$ et $W$ est le sous-espace engendré par $\cup_{i\in I }V_{i,l_i}$ et $V^+$.
Les racines de $A_{\tlM}$ dans $\tlG$ sont alors naturellement indexées par les couples $(V',V'')$ formés de deux éléments distincts parmi les sous-espaces qui apparaissent dans la décomposition \ref{eq:decompo}. Soit $\Sigma(\tlR)$ le sous-ensemble des racines associées à l'un des  couples suivants pour $i\in I$ :
\begin{itemize}
\item   $(V_{i,l},V_{i,l'})$ pour $l<l'$ et $l\not=l_i$, $l'\not=l_i$;
\item  $(V_{i,l},W)$ pour  $l<l_i$ ;
\item $(W,V_{i,l})$ pour $l_i<l$.
\end{itemize}
On a alors 
\[
 \fc_{\tlR}^{0}(\tlM_1) = \bigcap_{\al \in \Sigma(\tlR)}H(\al)^+ 
 \]
avec les notations de l'appendice \ref{App:domConvP}. Le lemme résulte alors du lemme \ref{lem:hfSpCnv}.
\end{preuve}
\end{paragr}

\begin{paragr}[Des lemmes combinatoires.] --- Soit $\tlR\in  \fc^{\tlH^-}(\tlM_0)$.

\begin{lemme}\label{lem:artLem51}
La somme
$$
\sum_{  \tlP\in \overline{\fc}_{\tlR}^{}(\tlM_1)}\eps_{\tlP}^{\tlG} \, \hat{\sigma}_{\tlP},
$$
vue comme fonction sur $\zgo_{\tlM_1}$ est égale à la fonction caractéristique du cône fermé
\[
\{H \in \zgo_{\tlM_1} | \gamma(H) \leq 0, \ \forall \ \gamma \in \widehat{\Pi}_{\tlR}\},
\]
où l'on voit $\widehat{\Pi}_{\tlR}$ comme une partie de $ \zgo^*_{\tlR}\subset  \zgo_{\tlM_0}^*=\zgo_{\tlM_1} ^*$.
\end{lemme}

\begin{preuve}
On reprend les notations du lemme \ref{lem:convexite} et de sa démonstration. On sait alors que $ \fc_{\tlR}^{0}(\tlM_1)$ est une partie convexe de $\pc^{\tlG}(\tlM)$. On a aussi
$$
\overline{\fc}_{\tlR}^{}(\tlM_1)=\{\tlP\in \fc^{\tlG}(\tlM_1) \mid \exists \tlQ \in \fc_{\tlR}^{0}(\tlM_1) \  \tlQ\subset \tlP\}
$$
En utilisant  \eqref{eq:sig-tau} et les notations de l'appendice \ref{App:domConvP} (cf. \eqref{eq:psiAsSumHtau}), on voit qu'on a pour tout $H\in  \zgo_{\tlM_1}$
\[
\sum_{\tlP \in \bar{\fc}_{\tlR}(\tlM_1) }
\varepsilon_{\tlP}^{\tlG}
\hat \sigma_{\tlP}(H) = \psi_{\fc_{\tlR}^{0}(\tlM_{1})}(\hat r_1 (H),  (0)), 
\]
où $(0)$ est la famille \og nulle\fg. Il résulte du lemme \ref{lem:arthGenrlz} et d'une nouvelle application des relations  \eqref{eq:sig-tau} que l'application $H\in \zgo_{\tlM_1} \mapsto  \psi_{\fc_{\tlR}^{0}(\tlM_{1})}(\hat r_1 (H),  (0))$  est égale à la fonction caractéristique de l'ensemble 
\[
\{H \in \zgo_{\tlM_1}| \delta(H) \le 0, \ \forall \ \tlP \in \fc_{\tlR}^{0}(\tlM_1) \ \forall \ \delta \in \widehat \Pi_{\tlP}\}.
\]

Il est facile de voir que 
\[
\widehat \Pi_{\tlR} \subset \bigcup_{\tlP \in \fc_{\tlR}^{0}(\tlM_1)}\widehat \Pi_{\tlP}.
\]
D'autre part, tout $\delta \in \bigcup_{\tlP \in \fc_{\tlR}^{0}(\tlM_1)}\widehat \Pi_{\tlP}$ s'écrit comme une somme d'éléments de $\widehat \Pi_{\tlR}$ à coefficients positifs. Le résultat s'en déduit. 

\end{preuve}

\begin{lemme}
  \label{cor:artLem52}
Pour tout $H\in \zgo_{\tlM_0}=\zgo_{\tlM_1}$, on a  
$$
\sum_{  \tlP\in \fc_{\tlR}^{}(\tlM_1)}\eps_{\tlP}^{\tlG} \, \hat{\sigma}_{\tlP}(H)=\eps_{\tlR}^{\tlH^-} \hat{\sigma}_{\tlR}(H).
$$
\end{lemme}

\begin{preuve}
La somme de gauche égale :
$$
\sum_{\tlP \in \overline{\fc}_{\tlR}(\tlM_1)} \eps_{\tlP}^{\tlG}  \hat{\sigma}_{\tlP} \big( \sum_{\tlR \subset \tlS \subset \tlP^-  }\eps_{\tlR}^{\tlS}\big)= \sum_{ \tlS \in  \fc^{\tlH^-}(\tlR) }\eps_{\tlR}^{\tlS} \sum_{\tlP \in \overline{\fc}_{\tlS}(\tlM_1)} \eps_{\tlP}^{\tlG}  \hat{\sigma}_{\tlP} 
$$
Le résultat se déduit alors immédiatement du lemme \ref{lem:artLem51}.
\end{preuve}

\begin{lemme}
  \label{lem:artLem42}  Soit  $\tlP \in  \overline{\fc}_{\tlR}^{}(\tlM_1)$. Pour tout $X\in \zgo_{\tlR}$ on a:
\[
\sum_{\{\tlQ \in \overline{\fc}_{\tlR}^{}(\tlM_1) \mid \tlQ\subset \tlP\}}
\eps_{\tlQ}^{\tlP}\sigma_{\tlR}^{\tlQ^{-} }(X) \hat{\sigma}_{\tlQ}^{\tlP}(X) = 
\begin{cases} 
1 \text{ si } \tlP \in \fc_{\tlR}^{}(\tlM_1) \ \text{ et }\ X \in \zgo_{\tlP}, \\
0 \text{ sinon}.
\end{cases}
\]

\end{lemme}

\begin{preuve}
Démontrons d'abord le cas $\tlP = \tlG$. Dans le membre de gauche, on groupe les $\tlQ$ suivant leur $\tlQ^-$ et on utilise le  corollaire \ref{cor:artLem52} : 
\[
\sum_{\tlS \in \fc^{\tlH^-}(\tlR)} \sigma_{\tlR}^{\tlS}(X)
\sum_{\tlQ \in \fc_{\tlS}(\tlM_{1})} \varepsilon_{\tlQ}^{\tlG} \hat \sigma_{\tlQ}(X)
=\sum_{\tlS \in \fc^{\tlH^-}(\tlR)}
\varepsilon_{\tlS}^{\tlH^-}
\sigma_{\tlR}^{\tlS}(X)\hat \sigma_{\tlS}(X).
\]
On conclut en invoquant le lemme \ref{lem:langlands}. Dans ce cas, la somme est nulle sauf si $\tlR=\tlH^-$ (dans cette situation $\zgo_{\tlR}=(0)$). 

Passons au cas d'un sous-groupe parabolique $\tlP$ général. On a la décomposition $M_{\tlP} = \MM_{\tlP}\times  \tlG_{\tlP} $ qui induit une décomposition de l'espace ambiant $V=W_1\oplus W_2$. On a $W_1=\oplus_{i\in I} (W_1\cap V_i)$ et $W_2=\oplus_{i\in I} (W_2\cap V_i) \oplus V^+$.  À  $M_{\tlP}$ on associe le sous-groupe de Levi $M_{\tlP^-}$ de $\tlH^-$. Celui-ci se décompose en $H_1\times \tlH_2$ avec 
$H_1\subset \prod_{i\in I}GL_{F_i}(W_1\cap V_i)$ et  $\tlH_2\subset \prod_{i\in I}GL_{F_i}((W_2\cap V_i)\oplus F_i e_0)$.  Plus généralement, pour tout $\tlS\in \fc^{\tlH^-}(\tlR)$, le sous-groupe parabolique $  M_{\tlP^-}\cap \tlS$ est un produit $S_1\times \tlS_2$. Le sous-groupe de Levi $\tlM_1$ se décompose en $M_{11}\times \tlM_{12}$. Ainsi on a une bijection naturelle 
$$\{\tlQ \in \overline{\fc}_{\tlR}^{}(\tlM_1) \mid \tlQ\subset \tlP\} \to  \overline{\fc}^{\MM_{\tlP}}_{R_1}(M_{11})  \times \overline{\fc}_{\tlR_2}^{\tlG_{\tlP}}(\tlM_{12}) $$
où l'ensemble $\overline{\fc}^{\MM_{\tlP}}_{R_1}(M_{11})$ est l'ensemble des $Q\in \fc^{\MM_{\tlP}}(M_{11})$ tel que $R_1\subset Q$ (on peut voir $R_1$ comme un sous-groupe de $\MM_{\tlP}$). On a de plus
$$\zgo_{\tlR}= \ago_{R_1}    \oplus \zgo_{\tlR_2}$$
et on écrit $X=X_1+X_2$ selon cette décomposition.

En raisonnant \emph{mutatis mutandis} comme dans la preuve du lemme \ref{lem:langlands}, on est ramené à la situation produit des deux cas liés à $\MM_{\tlP}$ et $\tlG_{\tlP}$. Le premier cas est traité dans l'article d'Arthur (cf. lemme 5.3 de \cite{arthur7}) et le second cas est le cas traité en début de démonstration. On trouve que la somme en question est toujours nulle sauf  si les trois conditions suivantes sont satisfaites (auquel cas on obtient $1$) :
\begin{enumerate}
\item $X_1\in \zgo_{\tlP}=\ago_{\MM_{\tlP}}$ ;
\item $R_1=H_1$  ;
\item   $\tlR_2=\tlH_{2}$.
\end{enumerate}
 Les deux dernières conditions sont équivalentes au fait $\tlP\in \fc_{\tlR}(\tlM_1)$. La condition 3 implique $ \zgo_{\tlR_2}=0$. La condition 1 implique  donc $X\in \zgo_{\tlP}$.
\end{preuve}

\end{paragr}

\begin{paragr}Soit $\tlR \in \fc^{\tlH^{-}}(\tlM_{0})$ et $\tlM$ le sous-groupe de Levi de $\tlG$ tel que $\fc_{\tlR}^{0}(M_{1}) \subset \pc(\tlM)$ (cf. lemme \ref{lem:convexite}). Soit $\tlA=A_{\tlM}$. Soit $\yc_{\tlR} = (Y_{\tlP})_{\tlP \in  \fc_{\tlR}^{0}(\tlM_{1})}$ une famille des vecteurs dans $\ago_{\tlM}$ qui est  $\tlA$-orthogonale positive au sens de la définition \eqref{eq:orthPosit}.  Pour tout $\tlQ \in \overline{\fc}_{\tlR}(\tlM_{1})$ soit $\tlP \in \fc_{\tlR}^{0}(\tlM_{1})$ contenu dans $\tlQ$ et $Y_{\tlQ} \in \ago_{\tlQ}$ la projection orthogonale de $Y_{\tlP}$ sur $\ago_{\tlQ}$. Cette  définition ne dépend du choix de $\tlP \subset \tlQ$. De cette manière, on obtient pour tout $\tlS \in \fc^{\tlH^{-}}(\tlR)$  une  famille positive $\yc_{\tlS} := (Y_{\tlQ})_{\tlQ \in \fc_{\tlS}^{0}(\tlM_{1})}$.

Pour tout $ H \in \zgo_{\tlR}$, soit
\[
\mathrm{B}_{\tlR}^{\tlG}(H, \yc_{\tlR}) = \sum_{\tlS \in \fc^{\tlH^{-}}(\tlR)}
\sigma_{\tlR}^{\tlS}(H) 
\big(
\sum_{\tlQ \in \fc_{\tlS}(\tlM_{1})}
\varepsilon_{\tlQ}^{\tlG}
\hat \sigma_{\tlQ}(H - Y_{\tlQ})
\big).
\]
En vertu du lemme \ref{lem:langlands} on a alors:
\begin{equation}\label{eq:42star}
\sum_{ \tlP \in \fc_{\tlR}(\tlM_{1})}
\varepsilon_{\tlP}^{\tlG}
\hat \sigma_{\tlP}(H - Y_{\tlP}) = 
\sum_{\tlS \in \fc^{\tlH^{-}}(\tlR)}
\varepsilon_{\tlR}^{\tlS}
\hat \sigma_{\tlR}^{\tlS}(H)
\mathrm{B}_{\tlS}^{\tlG}(H, \yc_{\tlS}).
\end{equation}

\begin{lemme}\label{lem:gamSumGam0} À un ensemble de mesure $0$ près, on a l'égalité des fonctions sur $\zgo_{\tlR}$ suivante:
\[
\mathrm{B}_{\tlR}^{\tlG}(\cdot, \yc_{\tlR}) = 
 \sum_{\tlP \in \fc_{\tlR}^{0}(\tlM_{1})} 
 \mathrm{B}_{\tlP}(\cdot, Y_{\tlP}).
\]
\end{lemme}

\begin{preuve}
La preuve est identique à la preuve du lemme 4.1 de \cite{arthur7}, 
le point clé étant l'application du lemme \ref{lem:artLem42}.
\end{preuve}
\end{paragr}

\section{Une formule des traces infinitésimale}\label{sec:RTFinf}

\subsection{Intégrales tronquées}

\begin{paragr}[Situation.] ---  Soit $F$ est un corps de nombres et  $V$ un  $F$-espace vectoriel de dimension $n\geq 1$.  Soit $G=GL_F(V)$,  $\tlG=GL_F(V\oplus Fe_0)$ et $\tggo=\tgl_F(V)$ (cf. §\ref{S:tggo}). On identifie $G$ à un sous-groupe de $\tlG$ (cf. §\ref{S:Gtilde}). On reprend les notations de la section \ref{sec:prelim}.   
\end{paragr}

\begin{paragr}[Choix auxiliaires.] --- \label{S:auxil}On fixe une base de $V$. Ce choix détermine des identifications $G\simeq GL(n)$ et $\tlG\simeq GL(n+1)$. Les groupes $G(\AAA)$ et $\tlG(\AAA)$ sont alors munis des sous-groupes compacts maximaux usuels  notés  respectivement $K$ et $\tlK$ : par exemple, on a $K=\prod_{v\in \vc} K_v$  et, via l'identification de $G$ à $GL(n)$, le groupe  $K_v$ est le groupe orthonal $O(n,\RR)$ en une place réelle,  le groupe unitaire $U(n,\RR)$ en une place complexe et $K_v=GL(n,\oc_v)$ si la place $v$ est non-archimédienne.

On dispose pour tout sous-groupe parabolique $\tlQ$ de $\tlG$ de l'application
$H_{\tlQ} : \tlG(\AAA)\to \ago_{\tlQ} $
définie relativement au sous-groupe $\tlK$ (cf. §\ref{S:HP}).

Le choix de la base de $V$  détermine aussi une inclusion  $T_0\subset \tlT_0$ de sous-tores déployés maximaux respectivement de $G$ et $\tlG$ (cf. §\ref{S:a0}). Soit $B$ un sous-groupe de Borel de $G$ contenant $T_0$. Soit $\tlB$ un sous-groupe de Borel de $\tlG$ tels que $\tlB\cap G=B$. Notons que $\tlT_0\subset \tlB$. L'épithète \emph{relativement standard}, resp. \emph{relativement semi-standard}, distinguera parmi les sous-groupes paraboliques de $\tlG$ ceux contenant $B$, resp. contenant $T_0$. De tels sous-groupes paraboliques $\tlP$ possèdent un unique facteur de Levi noté $M_{\tlP}$ contenant $T_0$. 

On fixe une mesure de Haar sur $G(\AAA)$.
\end{paragr}

\begin{paragr}[Noyaux paraboliques.] ---\label{S:noyau}
  Soit  $\tilde{P}$ un sous-groupe parabolique relativement standard de $\tlG$ et $M_{\tlP}N_{\tlP}$   sa décomposition de Levi standard. À cette donnée, on associe un sous-espace $\tpgo=\tmgo_{\tlP}\oplus \tngo_{\tlP}$ de $\tggo$  (cf. section \ref{ssec:parab}). Pour alléger, on omet en  indice $\tlP$. Pour tout $a\in \Ac$, on pose
 $$\tmgo_{a}= \tmgo \cap \tggo_a.$$
 
Soit $f\in \Sc(\tggo(\AAA))$ et  $a\in \Ac(F)$. Pour tout $g\in G(\AAA)$  on définit ce qu'on appellera improprement un \og noyau\fg{} 
\begin{equation}
  \label{eq:k}
  k_{\tilde{P},a}(f,g)=\sum_{X\in \tmgo_{a}(F)}  \int_{\tngo(\AAA)} f(g^{-1}\cdot (X+U))dU,
\end{equation}
(pour le choix de la mesure sur $\tngo(\AAA)$, cf. §\ref{S:Haar}).
\end{paragr}

\begin{paragr}[Noyau tronqué.] --- \label{S:noy-tronque}Pour tout sous-groupe de Borel $\tlB'$ contenant $\tlT_0$, soit $\ago_{\tlB'}^+\subset \ago_{\tl}^{\tlG}$ la chambre de Weyl aiguë et ouverte associée à $\tlB'$. On note simplement  $\ago_{\tl}^+=   \ago_{\tlB}^+$. Soit $T\in \ago_{\tl}^+$.   Par l'action du groupe de Weyl de $(\tlG,\tlT_0)$, on obtient des points $T_{\tlB'}\in \ago_{\tlB'}^+$ indexés par les sous-groupes de Borel semi-standard.  Soit $\tlB'$ un tel sous-groupe de Borel contenu dans $ \tlP$.  Par projection orthogonale de $T_{\tlB'}$  sur $\ago_{\tlP}$, on obtient un point $T_{\tlP}$ indépendant d'ailleurs du choix de  $\tlB'\subset \tlP$.

Avec les notations du § \ref{S:noyau}, on introduit le  \og noyau  tronqué\fg{} 
\begin{equation}
  \label{eq:kT}k^T_a(f,g)=\sum_{\tilde{P}\in \fc^{\tlG}(B)} \eps_{\tilde{P}}^{\tilde{G}} \sum_{\delta\in P(F)\back G(F)} \hat{\tau}_{\tilde{P}}(H_{\tilde{P}}(\delta g)-T_{\tlP}) \,  k_{\tilde{P},a}(f,\delta g)
\end{equation}
où $P=\tlP\cap G$ est un sous-groupe parabolique (standard) de $G$.
\end{paragr}

\begin{paragr}[Convergence d'intégrales.] ---\label{S:cv}

Soit
$$\eta : F^\times \back \AAA \to \CC^\times$$
un caractère \emph{quadratique} continu.
Par abus, on note encore $\eta$ le caractère quadratique 
$$\eta: G(\AAA) \to \CC^\times$$
obtenu par composition avec le déterminant. Le théorème suivant est une légère variante des théorèmes 3.6 et 4.8 de \cite{leMoi2}. 

\begin{theoreme}
  \label{thm:cv}
  \begin{enumerate}
  \item Il existe un point $T_+\in  \ago_{\tl}^+$ tel que pour tout $T\in T_++ \ago_{\tl}^+$, l'intégrale 
$$I^{\eta,T}_a(f)=\int_{[G]} k^T_a(f,g) \, \eta(g)\, dg
$$
converge absolument.
\item L'application $T\mapsto I^{\eta,T}_a(f)$ est la restriction d'une fonction exponentielle-polynôme en $T$.
\item Le terme purement polynomial de cette exponentielle-polynôme est constant.
  \end{enumerate}
\end{theoreme}
\end{paragr}

\subsection{Distributions globales}

\begin{paragr}[Distributions $I^{\eta}_a$.] --- \label{S:Ia}On continue avec les notations des sections précédentes. 
Soit  $I^\eta_a(f)$ le terme constant de l'intégrale $I^{\eta,T}_a(f)$ dans le théorème \ref{thm:cv}. On obtient ainsi une distribution $I^{\eta}_a$ (simplement une forme linéaire sur $\Sc(\tggo(\AAA))$) qui vérifie les propriétés suivantes.

  \begin{theoreme}
\label{thm:I}
\begin{enumerate}
\item La distribution $I^\eta_a$ ne dépend que du choix de la mesure de Haar sur $G(\AAA)$ mais pas du choix de $B$, $\tlB$ et $\tlK$.
\item La distribution $I^\eta_a$ est $\eta$-équivariante au sens où pour tous $f\in\Sc(\tggo(\AAA))$ et $g\in G(\AAA)$, on a 
$$I^\eta_a(f^g)=\eta(g) \, I^\eta_a(f).
$$
\item Le support de $I^{\eta}_a$ est inclus dans $\tggo_a(\AAA)$.
\item Soit $S\subset \vc_\infty$ un ensemble de places archimédiennes. Soit $f^S\in \Sc(\tggo(\AAA^S))$. La forme linéaire sur $\Sc(\tggo(\AAA_S))$ donnée par 
$$f_S \mapsto I^\eta_a(f_S\otimes f^S)
$$
est continue pour la topologie usuelle sur $\Sc(\tggo(\AAA_S))$.
\end{enumerate}
      \end{theoreme}

Ce théorème est une légère variante de résultats de  \cite{leMoi2} (les références renvoient à cet article) : l'assertion 1 vient de la section 4.5, l'assertion 2 du  théorème  4.11, l'assertion 3 est évidente par construction et enfin l'assertion 4 résulte de la démonstration   des théorèmes 3.6 et 3.7.
\end{paragr}

\begin{paragr}[Transformation de Fourier partielle.] --- \label{S:TFP} On suit les choix et les notations de §\ref{S:TFPV}$$\psi: F\back \AAA \to \CC^\times.
$$
Soit 
$$\bg \cdot,  \cdot \bd : \tggo \times \tggo \to F
$$
la forme bilinéaire symétrique, non dégénérée et $G$-invariante définie pour $X=(A,b,c)\in \tggo$ par
\begin{equation}
  \label{eq:fbs}
  \bg X,X\bd= \trace(A^2)+2cb.
\end{equation}
Pour tout sous-espace $\tggo_1\subset \tggo$ qui est $G$-invariant et non dégénéré, on dispose de la transformée de Fourier partielle 
\begin{equation}
  \label{eq:TFPg1}
 f \mapsto  \hat{f}_{\ggo_1}.
\end{equation}
Les seuls espaces $\tggo_1$ possibles sont les suivants et leurs sommes directes 
\begin{itemize}
\item $F$ vu comme le sous-espace des homothéties de $V$ ;
\item $\mathfrak{sl}_F(V)$ le sous-espace des endomorphismes de $V$ de trace nulle ;
\item $V\oplus V^*$.
\end{itemize}

\begin{remarque}
  \label{rq:TFP}
De notre point de vue, les sous-espaces $\tggo_1$ les plus intéressants seront $\mathfrak{sl}_F(V)$, $V\oplus V^*$ et $\mathfrak{sl}_F(V) \oplus V\oplus V^*$.
\end{remarque}
\end{paragr}

\begin{paragr}[Formule des traces infinitésimale.] --- \label{S:RTFinf}Pour toute  transformation de Fourier partielle $f\mapsto \hat{f}$ (cf. §\ref{S:TFP}), on note $D\mapsto \hat{D}$ la transformation duale au niveau des distributions. Notons que celle-ci préserve le fait d'être $\eta$-équivariante. L'application qui, à $(A,b,c)\in \tggo$ associe la trace de l'endomorphisme $A$, est bien sûr $G$-invariante et donc définit un morphisme $\trace$ de $\Ac$ dans la droite affine. Pour tout $\al\in F$, soit  $\Ac_{\al}$ la fibre de ce morphisme.

  \begin{theoreme}
    \label{thm:RTFinf}
    \begin{enumerate}
    \item Pour tout $f\in \Sc(\tggo(\AAA))$, la somme $\sum_{a\in \Ac(F)} I_a^{\eta}(f)$ converge absolument ce qui définit une distribution $\sum_{a\in \Ac(F)} I_a^{\eta}$.
    \item On a 
$$\sum_{a\in \Ac(F)} I_a^{\eta}=\sum_{a\in \Ac(F)} \hat{I}_a^{\eta}
$$
pour toute transformation de Fourier partielle au sens du §\ref{S:TFP}.
\item Si la transformation de Fourier partielle est associée à l'un des trois espaces décrits dans la remarque \ref{rq:TFP}, alors pour tout $\al\in F$, on a 
$$\sum_{a\in \Ac_\al(F)} I_a^{\eta}=\sum_{a\in \Ac_\al(F)} \hat{I}_a^{\eta}.
$$
    \end{enumerate}
  \end{theoreme}
  
Ce théorème est une légère variante de théorèmes de \cite{leMoi2} (les références qui suivent renvoient à cet article) : l'assertion 1 vient du théorème 3.6, l'assertion 2 du théorème 5.1 et l'assertion 3 en est juste une variante (il suffit de considérer une formule sommatoire de Poisson limitée à un sous-espace invariant de sous-espace $\mathfrak{sl}_F(V)\oplus V\oplus V^*$).
\end{paragr}

\subsection{Le cas d'un produit}\label{ssec:produit}

\begin{paragr}
C'est le cas considéré à la section \ref{ssec:transv} de $\thgo$ muni de l'action du groupe  $H$.
\end{paragr}

\begin{paragr}
  Les constructions et les théorèmes qu'on vient d'évoquer aux sections ci-dessus se généralisent sans mal au cas de $H$ agissant sur $\thgo$. La restriction du caractère $\eta$ au sous-groupe $H(\AAA)$ de $G(\AAA)$ est encore noté $\eta$. Pour tout $a\in \Ac_H(F)$, on dispose donc d'une distribution $I^{H,\eta}_a$ sur $\Sc(\thgo(\AAA))$ qui est $\eta$-équivariante et à support dans la fibre $\thgo_a(\AAA)$ du morphisme canonique $\thgo\to \Ac_H$. On laisse au lecteur le soin de formuler  les analogues des théorèmes \ref{thm:I} et \ref{thm:RTFinf} dans ce cadre. Notons que si $f$ est un produit de fonctions $f_+\in \Sc(\tgl(V^+,\AAA))$ et $f_i\in \Sc(\tgl(V_i, \AAA_{F_i}))$ et si $a$ correspond à un couple $(a_+,(a_i)_{i\in I})$ alors $I^{H,\eta}(f)$ est simplement le produit des distributions définies plus haut associés à $a_i,f_i$ et $\tgl_{F_i} (V_i)$ (ou   $a_+,f_+$ et $\tgl_{F}(V^+) $) et le caractère obtenu par restriction de $\eta$ sur le facteur correspondant.
\end{paragr}

\section{Le théorème de densité}\label{sec:densite} 

\subsection{Intégrales orbitales locales}\label{ssec:IOloc}

\begin{paragr}[Facteur $\tilde{\eta}$.] ---   \label{S:eta} Dans toute cette section, on reprend les notations de la section \ref{sec:RTFinf}. Soit $X=(A,b,c) \in \tggo^{\rs}$. Rappelons qu'on a fixé une base de l'espace $V$ ce qui identifie $G$ à $GL(n)$.  La matrice de la famille $(b,Ab,\ldots,A^{n-1}b)$ dans cette base définit donc un élément noté $\delta_X \in G$. On a, de plus, $\delta_{g\cdot X}=g\delta_X$. On a ainsi obtenu un morphisme équivariant $  \tggo^{\rs} \to G$.

Pour tout ensemble $S$ de places de $F$ et tout $X \in \tggo^{\rs}(\AAA_S)$, on obtient un élément $\delta_X \in G(\AAA_S)$. En regardant ce dernier groupe comme un sous-groupe de $G(\AAA)$,  on définit
$$\tilde{\eta}(X)=\eta(\delta_X)^{-1}$$
où $\eta$ est le caractère fixé au §\ref{S:cv}.
On a donc 
$$\tilde{\eta}(g^{-1} \cdot X)=\eta(g) \tilde{\eta}(X).$$
\end{paragr}

\begin{paragr}[Mesure de Haar.] --- \label{S:ensS}  Soit $S\subset \vc$ un ensemble fini  de places de $F$. Soit $dg_S$ et $dg^S$ des mesures de Haar sur $G(\AAA_S)$ et $G(\AAA^S)$ de sorte que $dg=dg_S\otimes dg^S$ soit la mesure de Haar choisie au §\ref{S:auxil}.
\end{paragr}

\begin{paragr}[Intégrales orbitales locales.] --- \label{S:IOPloc} Soit $f \in \Sc(\tggo(\AAA_S))$ et $a\in \Ac^{\rs}(\AAA_S)$. D'après le lemme \ref{lem:fibress}, il existe $X\in \tggo^{\rs}(\AAA_S)$ tel que $\ggo_a(\AAA_S)$ soit l'orbite sous $G(\AAA_S)$ de $X$. On introduit alors \emph{l'intégrale orbitale semi-simple régulière locale}
  \begin{eqnarray}
    I_a^{\tilde{\eta}}(f)&=&\int_{G(\AAA_S)} f(g^{-1}\cdot X) \, \tilde{\eta}(g^{-1}\cdot X)\, dg_S \\
&=&  \tilde{\eta}(X) \int_{G(\AAA_S)} f(g^{-1}\cdot X) \, \eta(g)\, dg_S.\label{eq:IOPloc2}
  \end{eqnarray}
  
\begin{remarque}
  Cette construction dépend implicitement du choix de la mesure de Haar $dg_S$. Cependant, elle ne dépend pas du choix de $X$.
\end{remarque}

\begin{remarque}
  \label{rq:decomposition}
  Soit $S'\subset \vc$ un ensemble fini de places disjoint de $S$ et $S''=S\cup S'$. On suppose les mesures de Haar vérifient l'égalité $dg_S\otimes dg_{S'}=dg_{S''}$. Soit $a=(a_S,a_{S'})\in \Ac^{\rs}(\AAA_{S''})$. Pour toutes fonctions  $f _{S}\in \Sc(\tggo(\AAA_S))$ et $f _{S'}\in \Sc(\tggo(\AAA_{S'}))$, on a 
 $$I_a^{\tilde{\eta}}(f_S\otimes f_{S'})=I_{a_{S}}^{\tilde{\eta}}(f_{S})\cdot I_{a_{S'}}^{\tilde{\eta}}(f_{S'}).
$$
\end{remarque}
\end{paragr}

\subsection{Distributions stables}

\begin{paragr}[Fonctions instables.] --- \label{S:inst} Soit 
$$\Sc(\tggo(\AAA_S))_\eta\subset \Sc(\tggo(\AAA_S))$$
 le sous-espace des fonctions $\eta$-instables c'est-à-dire des fonctions $f\in  \Sc(\tggo(\AAA_S))$ telles $I^{\tilde{\eta}}_a(f)=0$ pour tout $a\in \Ac^{\rs}(\AAA_S)$. Comme le montre l'expression \eqref{eq:IOPloc2}, cette définition ne dépend que du caractère $\eta$ ; elle ne dépend ni du choix de la mesure de Haar ni du facteur $\tilde{\eta}$.
\end{paragr}

\begin{paragr}
  On dit qu'une distribution, c'est-à-dire une forme linéaire sur $\Sc(\tggo(\AAA_S))$, est $\eta$-stable  si elle s'annule sur le sous-espace 
des fonctions $\eta$-instables. Toute distribution $\eta$-stable est évidemment $\eta$-équivariante. Il est fort possible que la réciproque soit vraie (avec vraisemblablement des conditions de continuité pour la distribution aux places archimédiennes) mais, à notre connaissance, cela n'est pas connu (pour un résultat lorsque $\dim(V)\leq 2$, on pourra consulter \cite{Zhrang3}). L'objet principal de cette section et de cette partie est de démontrer le théorème suivant.

  \begin{theoreme}
    \label{thm:densite}
Soit $a\in \Ac(F)$ et $f^S\in \Sc(\tggo(\AAA^S))$. La distribution
$$f \in \Sc(\tggo(\AAA_S)) \mapsto I_a^{\eta}(f\otimes f^S)$$
est $\eta$-stable.
  \end{theoreme}

Avant de donner la démonstration de ce théorème à la section \ref{ssec:demo-densite}, nous aurons besoin de deux résultats qui sont énoncés dans les deux section qui suivent.
\end{paragr}

\subsection{Stabilité et transformée de Fourier}

\begin{paragr}[Stabilité et transformation de Fourier.] --- \label{S:stab} Comme dans le cas global (cf. §\ref{S:TFP}), le choix de la forme $\bg\cdot,\cdot\bd$ et du caractère $\psi$ permet de définir des transformées de Fourier partielles relatives à des sous-espaces non dégénérées de $\tggo$ (cf. §\ref{S:TFPV}). On en fixe une notée $f\mapsto \hat{f}$. 

  \begin{theoreme}\label{thm:TFstable}
    L'espace $\Sc(\tggo(\AAA_S))_\eta$ est invariant par la transformée de Fourier  $f\mapsto \hat{f}$.
  \end{theoreme}
et son corollaire évident :

\begin{corollaire}
  \label{cor:TFstable}
Si $D$ est $\eta$-stable, il en est de même de $\hat{D}$.
\end{corollaire}
\end{paragr}

\begin{paragr}[Démonstration du théorème \ref{thm:TFstable}.] --- Lorsque $\eta$ est non-trivial le théorème \ref{thm:TFstable} résulte du théorème 4.17 (plus général) de \cite{Z1} et de sa variante archimédienne due à H.Xue (cf. \cite{xue}, combinaison des lemmes 9.2, 9.3 et 9.4).  Lorsque $\eta$ est quadratique (ou même unitaire) quelconque, esquissons pour la commodité du lecteur une preuve simple directement inspirée de   \cite{Z1}. Soit $\tggo_1\subset \tggo$ un sous-espace $G$-invariant et non dégénéré pour la forme $\bg \cdot,\cdot \bd$ (cf. §\ref{S:TFP}) et $\tggo_2$ son orthogonal. La transformée de Fourier partielle est relative à $\tggo_1$. Soit $X_2 \in \tggo_2(\AAA_S)$ un élément régulier semi-simple (au sens de l'action de $G$ sur $\tggo_2$, c'est-à-dire sa $G$-orbite géométrique est fermée et de dimension minimale). Pour tous $f_1,f_2\in \Sc(\tggo(\AAA_S))$ soit
$$T(f_1,f_2)=\int_{ \tggo_1(\AAA_S)  }  \left(\int_{G(\AAA_S)} f(g^{-1}\cdot (X_1+X_2)) \, \eta(g)\, dg_S \right)f_2(X_1+X_2) dX_1.
$$
L'intégrale converge dans l'ordre indiqué. De plus, on a \og la formule des traces locale\fg{} suivante (cf. \cite{Z1} théorème 4.6 et \cite{xue} théorème 6.1)
$$T(f_1,f_2)=T(\hat{f}_1,\check{f}_2)
$$
où $\check{f}_2(X)=\hat{f_2}(-X)$. Pour $f_1\in \Sc(\tggo(\AAA_S))_{\eta}$, on a clairement $T(f_1,f_2)=0$ pour tout $f_2\in  \Sc(\tggo(\AAA_S))$. On déduit donc de la formule des traces qu'on a aussi $T(\hat{f}_1,f_2)=0$ pour tout $f_2\in  \Sc(\tggo(\AAA_S))$. En utilisant la lissité de l'application

$$X\mapsto   \int_{G(\AAA_S)} f(g^{-1}\cdot (X_1+X_2)) \, \eta(g)\, dg_S$$
au voisinage d'éléments réguliers, on en déduit que les intégrales orbitales de $f_1$ sont nulles pour tous les éléments réguliers semi-simples de la forme $X_1+X_2$ avec $X_2$ lui-même régulier semi-simple comme élément de $\tggo_2$. L'ensemble de tels éléments étant dense dans celui des réguliers semi-simples, la même lissité implique que \emph{toutes} les intégrales orbitales de $f_1$ sont nulles. 
\end{paragr}
\subsection{Descente}\label{ssec:descente}

\begin{paragr}[Situation.] --- \label{S:situation-desc}Soit $0\leq r\leq n$ et $a\in \Ac^{(r)}(F)$.  Soit $V=V^+\oplus V^-$ une décomposition de l'espace vectoriel sous-jacent aux données telle que $\dim(V^+)=r$. D'après les lemmes \ref{lem:iso-r} et \ref{lem:entierk}, il existe $a_+\in \Ac_{V^+}^{(r)}(F)$ et $a_- \in \Ac_{V^-}^{(0)}(F)$ tel que $a$ soit l'image du couple $(a_+,a_-)$ par l'isomorphisme \eqref{eq:iota}. La donnée $a_-$ est équivalente à la donnée d'un polynôme $P$, unitaire, de degré $n-r$, à coefficient dans $F$. Soit
$$P=\prod_{i\in I} P_i^{n_i}$$
la décomposition de $P$ en facteurs irréductibles $P_i$ deux à deux distincts. Pour tout $i\in I$, on introduit l'extension $F_i=F[t]/P_i$ de $F$. Soit $\al_i\in F_i$ la classe de l'indéterminée $t$. Pour $i\in I$ soit $V_i$ un $F_i$-espace vectoriel de dimension $n_i$. Fixons un $F$-isomorphisme $V^-=\oplus_{i\in I}V_i$. Soit $\Ac_{V_i}$ l'espace affine sur $F_i$ attaché au $F_i$-espace vectoriel $V_i$ (cf. section \ref{sec:prelim}). Soit $a_i\in \Ac_{V_i}^{(0)}$ l'image du triplet $(\al_i,0,0)\in \tgl_{F_i}(V_i)$ par le morphisme canonique (on voit $\al_i$ comme une homothétie de $V_i$). Avec nos données $V^+, V_i, $ etc., on reprend les notations de la section \ref{ssec:transv}. On dispose donc d'un espace $\thgo$ muni d'une action du groupe $H$ et du quotient  $\Ac_H$. La famille $(a_+,(a_i)_{i\in I})$ définit un élément $a_H$ de l'ouvert $\Ac_{H}'(F)$ dont l'image par le morphisme $\Ac_{H}'\to \Ac$ (cf. \eqref{eq:iotaHquotient}) est l'élément $a$ de départ.
\end{paragr}

\begin{paragr}
Soit $dh=dh_S\otimes dh^S$ une mesure de Haar sur $H(\AAA)=H(\AAA_S)\times H(\AAA^S)$ qui se décompose en mesures de Haar sur les facteurs.
\end{paragr}

\begin{paragr}[Espace de Schwartz.]  --- \label{S:Sch-salg}Pour toute variété réelle semi-algébrique $X$, on sait définir un espace de Schwartz, $\Sc(X)$ muni d'une topologie naturelle (cf. par exemple \cite{ducloux} définition 1.2.1). Lorsque $X$ est un $\RR$-espace vectoriel, on retrouve la notion usuelle.
  \end{paragr}

\begin{paragr}\label{S:hypsurS}
  On a fixé au §\ref{S:ensS} un ensemble $S$ fini de places. Quitte à l'agrandir $S$, on suppose vérifiées les propriétés supplémentaires du §\ref{S:surA} et les propriétés suivantes. On dispose donc d'un anneau $A\subset F$ et tous les objets viennent avec une structure sur $A$. On suppose également que l'élément $a_H$ construit ci-dessus appartient en fait à $\Ac_H'(A)$. Le caractère $\eta$ est trivial sur $G(\oc^S)$ et $H(\oc^S)$.
\end{paragr}

\begin{paragr}\label{S:Omega}  Le morphisme $\iota_H$ (cf. \eqref{eq:iotaHquotient})  est étale.  Par conséquent (cf. \cite{BCR} proposition 8.1.2) pour le cas archimédien) pour tout $v\in S$, on peut trouver des ouverts $\om_v\subset \Ac(F_v)$ et $\om_{H,v}\subset \Ac_H'(F_v)$ voisinages respectifs  de $a$ et $a_H$ de sorte qu'on a 
\begin{itemize}
\item si $v$ est archimédienne, les ouverts  $\om_v$ et $\om_{H,v}$ sont semi-algébriques ;
\item $\iota_H$ induit une bijection $\om_{H,v}\to \om_v$ qui est un difféomorphisme de Nash si $v$ est archimédienne et un isomorphisme de variété analytique  (au sens de  \cite{LaLg}) si $v$ est non-archimédienne.
\end{itemize}
  Soit  $\om=\prod_{v\in S}\om_v \subset \Ac(\AAA_S)$ et $\om_H=\prod_{v\in S } \om_{H,v}\subset \Ac_H'(\AAA_S)$. Le morphisme $\iota_H$ induit donc une bijection
  \begin{equation}
    \label{eq:iso-om}
     \om_H\to\om. 
   \end{equation}
   Soit $\Om=\prod_{v\in S} \Om_v \subset \tggo(\AAA_S)$ et $\Om_H=\prod_{v\in S} \Om_{H,v} \subset \thgo'(\AAA_S)$ les ouverts obtenus comme images réciproques respectives par le morphisme $a$  de $\om$ et $\om_H$.

Pour tout $v\in S$, on déduit de l'isomorphisme \eqref{eq:isocrucial} une bijection  qui est un difféomorphisme de Nash ou un isomorphisme de variété analytique selon que la place $v$ est archimédienne ou pas  (la surjectivité résulte du fait  que la cohomologie galoisienne en degré $1$ du groupe $H\times_F F_v$ est triviale)
$$G(F_v) \times ^{H(H_v)}\Om_{H,v} \to \Om_v
$$
donnée par $(g,Y)\mapsto g\cdot \iota_H(Y)$.  Soit $S_\infty=S\cap \vc_\infty$. On dispose alors de l'espace de Schwartz $\Sc(\Om)$  qui est par définition l'espace engendré par les fonctions du type $f_\infty \otimes f^\infty$ où $f_\infty\in \Sc(\prod_{v\in S_\infty} \Om_v)$ (cf. \ref{S:Sch-salg}, le produit $\prod_{v\in S_\infty} \Om_v$ étant une variété semi-algébrique) et $f^\infty$ est une fonction localement constante à support compact sur $\prod_{v\in S\setminus S_\infty} \Om_v$.

On en déduit alors une application linéaire surjective qui est continue sur sa composante archimédienne
$$\Sc(G(\AAA_S) \times \Om_H)\to \Sc(\Om)
$$
donnée par $\al \mapsto f_\al$ où $f_\al$ est déterminé par l'égalité
\begin{equation}
  \label{eq:falpha}
  f_\al(g\cdot \iota_H(Y))=\int_{H(\AAA_S)} \al(gh,h^{-1}\cdot Y) \, dh
\end{equation}
pour tout $g\in G(\AAA_S)$ et $Y\in \Om_H$. La surjectivité aux places archimédiennes résulte du fait que  l'application composée
$$\prod_{s\in S_\infty }G(F_v) \times \Om_{H,v} \to \prod_{s\in S_\infty } G(F_v) \times ^{H(H_v)}\Om_{H,v} \to \prod_{s\in S_\infty } \Om_v
$$
est une $\prod_{v\in S}G(F_v)$-fibration localement triviale et triviale  sur un recouvrement fini de   $\prod_{s\in S_\infty } \Om_v$ par des ouverts semi-algébriques  (cf. \cite{AGdeRham} proposition 4.0.6) et de l'existence de partitions de l'unité subordonnées à un tel recouvrement (cf. \cite{schwFons} théorème 4.4.1).

On utilisera également l'application continue et surjective
$$\Sc(G(\AAA_S) \times \Om_H)\to \Sc(\Om_H)
$$
donnée par $\al\mapsto f^{H,\eta}_\al$ où $f^{H,\eta}_\al$ est définie par
\begin{equation}
  \label{eq:fHalpha}
f^{H,\eta}_\al(Y)=\int_{G(\AAA_S)} \al(g, Y) \, \eta(g)^{-1}\, dg
\end{equation}
pour tout $Y\in \Om_H$.
\end{paragr}

\begin{paragr}[Descente des distributions globales.] --- Il s'agit du théorème suivant qui sera démontré à la toute fin de la section \ref{sec:desc}. Les notations et les hypothèses sont celles du paragraphe précédent. On considère les distributions globales $I^{G,\eta}_a$ et $I^{H,\eta}_{a_H}$ définies dans la section \ref{sec:RTFinf} relatives respectivement à l'action de $G$ sur $\tggo$ et $H$ sur $\thgo$ (pour ce dernier cas, cf. §\ref{ssec:produit}). Précisons qu'on voit $H(\AAA)$ comme un sous-groupe de $G(\AAA)$. Le caractère $\eta$ sur $H(\AAA)$ est simplement la restriction du caractère $\eta$. 

  \begin{theoreme}
    \label{thm:desc}
Pour tout $\al \in \Sc(G(\AAA_S) \times \Om_H)$, on a
$$I^{G,\eta}_a(f_\al\otimes \mathbf{1}_{\tggo(\oc^S)})= \frac{\vol(G(\oc^S))}{\vol(H(\oc^S))}\cdot I^{H,\eta}_{a_H}(f^{H,\eta}_\al\otimes \mathbf{1}_{\thgo(\oc^S)}).$$
  \end{theoreme}
\end{paragr}

\begin{paragr}\label{S:desc-locale} Soit $\Ac_H^{\rs}=\Ac_{V^+}^{\rs}\times \prod_{i\in I} \Ac_{V_i}^{\rs}$. On ne confondra pas cet ouvert avec l'ouvert $\Ac_H^{\Grs}$ de \eqref{eq:egalite(n)}. D'après le lemme \ref{lem:compatibilite}, on a  $\Ac_H^{\Grs}\subset \Ac_H^{\rs}$.

Soit $\thgo^{\rs}$ l'image inverse de $\Ac_H^{\rs}$ par le morphisme canonique. En travaillant composante par composante, on définit comme au §\ref{S:eta} un morphisme $H$-équivariant  $\delta^H:\thgo^{\rs}\to H$. Rappelons qu'on identifie $H$ à un sous-groupe de $G$. On note encore $\det$ le morphisme induit par la restriction du morphisme déterminant. On en déduit alors un morphisme
$$\thgo^{\Grs} \to \mathbb{G}_{m,F}$$
donné par $\det(\delta^H_Y)  \det(\delta_{\iota_H(Y)})^{-1}$. On vérifie par un calcul explicite que ce morphisme se prolonge en un morphisme $\thgo^{'} \to \mathbb{G}_{m,F}$ où $\thgo'$ est l'ouvert défini au §\ref{S:thgo} : il est d'abord facile de se ramener au cas $r=0$ (l'entier $r$ est introduit au §\ref{S:situation-desc}) et dans ce cas on utilise les calculs donnés dans la preuve du lemme \ref{lem:compatibilite}. Ce morphisme est évidemment $H$-invariant. On en déduit un morphisme 
$$\delta_G^H : \Ac_H' \to  \mathbb{G}_{m,F}.$$
En particulier, l'application
$$a_H'\in \Ac_H(\AAA_S) \mapsto \eta(\delta_G^H(a_H'))$$
est localement constante. Comme $a_H\in \Ac_H'(F)$ et que $\eta$ est trivial sur $F^\times$, quitte à agrandir $S$, on peut et on va supposer que 
$$  \eta(\delta_G^H(a_H))=1.$$
Quitte à restreindre les ouverts $\om$ et $\om_H$, on va supposer que $\eta(\delta_G^H(a_H'))=1$ pour tout $a_H'\in \om_H$.

On définit alors un facteur $\tilde{\eta}_H$ par 
$$\tilde{\eta}_H(Y)=\eta(\delta^H_Y)$$
pour tout $Y \in \Ac_H^{\Grs}(\AAA_S)$.

Il résulte de ce qui précède que pour $Y\in \Ac_H^{\Grs}(\AAA_S)\cap  \Om_H $, on a 
$$\tilde{\eta}(\iota_H(Y) )=\tilde{\eta}_H(Y ).$$

On définit sans peine les intégrales orbitales locales $I_{a_H'}^{H,\tilde{\eta}_H}$ pour $a_H'\in \Ac_H^{\rs}(\AAA_S)$.  L'intérêt de considérer les fonctions $f_\al$ et $f^{H,\eta}_\al$ vient en partie du lemme suivant. 

\begin{lemme}
   \label{lem:desc-locale}
Soit  $\al \in \Sc(G(\AAA_S) \times \Om_H)$ , $a'\in\Ac^{\rs}(\AAA_S)$ et $a'_H\in  \Ac_H^{\rs}(\AAA_S)$. 
\begin{enumerate}
\item Si $a'\notin \om$, l'intégrale orbitale locale de $f_\al$ est nulle
$$I_{a'}^{G,\tilde{\eta}}(f_\al)=0.
$$
\item Si $a'_H\notin \om_H$, l'intégrale orbitale locale de $f_\al^{H,\eta}$ est nulle
$$I_{a'_H}^{H,\tilde{\eta}_H}(f_\al^{H,\eta})=0.
$$
\item Si $a'\in \om$ et $a'_H\in \om_H$ se correspondent par l'isomorphisme \eqref{eq:iso-om}, on a
  \begin{equation}
    \label{eq:desc-loc}
    I_{a'}^{G,\tilde{\eta}}(f_\al)= I_{a'_H}^{H,\tilde{\eta}_H}(f_\al^{H,\eta}).
  \end{equation}
\end{enumerate}
\end{lemme}

\begin{preuve}
  Les assertions $1$ et $2$ résultent immédiatement des propriétés de support des fonctions. L'assertion 3 est une manipulation élémentaire.
\end{preuve}
  \end{paragr}

\subsection{Démonstration du théorème \ref{thm:densite}}\label{ssec:demo-densite}

\begin{paragr}[Hypothèse de récurrence.] --- On rappelle que $n$ désigne la dimension de l'espace $V$. Pour démontrer le  théorème \ref{thm:densite}, on raisonne par récurrence sur $n$ en supposant le théorème vrai pour l'espace $\tgl_F(V)$ pour tout $F$-espace vectoriel de dimension $<n$. On a même besoin d'une hypothèse de récurrence un peu plus forte : le théorème  \ref{thm:densite} est vrai pour tout produit fini $\prod_{i\in I} \tgl_{F_i}(V_i)$ où $F_i$ est une extension finie de $F$ et $V_i$ est un $F_i$-espace vectoriel de dimension $n_i<n$. En principe, cela nous oblige à traiter le cas de tels produits aussi dans la récurrence. Il n'y a en fait aucune difficulté supplémentaire autre que notationnelle. Par souci de simplicité, on se borne à rédiger le cas de $\tgl_F(V)$.
  \end{paragr}

\begin{paragr}
  Soit $a\in \Ac(F)$. Soit $f\in \Sc(\tggo(\AAA_S))_{\eta}$ et $f^S\in \Sc(\tggo(\AAA^S))$. Comme les intégrales locales d'un produit sont elles-mêmes des produits (cf. remarque \ref{rq:decomposition}), il suffit de prouver le théorème \ref{thm:densite} pour un ensemble fini $S$ de places aussi grand que l'on veut. Quitte à agrandir $S$, on peut et on va supposer que $S$ satisfait toutes les propriétés requises par la section \ref{ssec:descente}. On peut également supposer qu'on a $f^S=\mathbf{1}_{\tggo(\oc^S)}$ est la fonction caractéristique de $\tggo(\oc^S)$, que $\eta$ est trivial sur $G(\oc^S)$ et que le morphisme $X\mapsto \delta_X$ est défini sur  $\oc^S$.
\end{paragr}

\begin{paragr}[Cas régulier semi-simple.] --- \label{S:casrss} C'est le cas le plus immédiat, pour lequel $a\in \Ac^{\rs}(F)$. Soit $X\in \tggo_a(F)$. On voit cet élément comme un adèle et on écrit avec des notations évidentes $X=(X_S,X^S)$. De même, on écrit $a=(a_S,a^S)$. On peut supposer et on suppose que $X^S\in  \tggo_a^{\rs}(\oc^S)$. On a alors 
$$\tilde{\eta}(X_S)=\tilde{\eta}(X)=1.
$$
La distribution globale correspondant à $a$ est donnée par une intégrale orbitale qu'on relie aisément à une intégrale locale à l'aide du lemme \ref{lem:integrite}. On a alors
\begin{eqnarray*}
  I_a^{\eta}(f\otimes  \mathbf{1}_{\tggo(\oc^S)}  )&=&\int_{G(\AAA)}   (f\otimes  \mathbf{1}_{\tggo(\oc^S)})(g^{-1}X)\eta(g)\, dg\\
&=& \vol(G(\oc^S),dg^S)\cdot I_{a_S}^{G,\tilde{\eta}}(f)
\end{eqnarray*}
qui est nulle puisque $f\in \Sc(\tggo(\AAA_S))_{\eta}$ par hypothèse.  
\end{paragr}

\begin{paragr}[Cas de descentes.] ---\label{S:cas-desc} On suppose désormais $a\notin  \Ac^{\rs}(F)$. On reprend les notations de la section \ref{ssec:descente}. On dispose donc de $H$, $\thgo$ etc. L'entier $r$ qui y est défini vérifie $0\leq r <n$. On suppose qu'on est dans le cas  $n_i<n$ pour tout $i\in I$. L'hypothèse de récurrence s'applique donc à $\thgo$. 

Soit $\theta \in \Cc(\om)$ qui vaut $1$ dans un voisinage de $a_S$ et $f'=(\theta\circ a) f$. Observons que la fonction $f'$ d'une part est à support inclus dans $\Om$ et d'autre part qu'elle appartient aussi à  $\Sc(\tggo(\AAA_S))_{\eta}$. Comme la fonction $f'-f$ est nulle au voisinage de $\tggo_a(\AAA_S)$, on a, par la propriété de support (cf. assertion 3 du théorème \ref{thm:I}), l'égalité
$$I_a^\eta(f \otimes  \mathbf{1}_{\tggo(\oc^S)}  )= I_a^\eta(f' \otimes  \mathbf{1}_{\tggo(\oc^S)}  ).$$
Il suffit donc de prouver la nullité du membre de droite. Quitte à remplacer $f$ par $f'$, on peut et on va supposer que la fonction $f$ est de la forme $f_\al$ pour $\al\in \Sc(G(\AAA_S)\times \Omega_H)$. En tenant compte du théorème \ref{thm:desc} et de l'hypothèse de récurrence, on voit qu'il suffit de montrer que la fonction $f_\al^{H,\eta}$ appartient à  $\Sc(\thgo(\AAA_S))_{\eta}$. Le lemme \ref{lem:desc-locale} donne la nullité des intégrales orbitales de cette fonction pour tout $a_S\in\Ac^{\Grs}_H(\AAA_S)$. Il en est alors de même pour tout $a_S\in\Ac^{\rs}_H(\AAA_S)$. En effet, d'une part, pour tout $g\in \Sc(\thgo(\AAA_S))$, la fonction
$$
a_S\in\Ac^{\rs}_H(\AAA_S)\mapsto I^{H,\tilde{\eta}_H}_{a_S}(g)
$$
est lisse  et d'autre part $\Ac^{\Grs}_H(\AAA_S)$ est dense dans $\Ac^{\rs}_H(\AAA_S)$.  
\end{paragr}

\begin{paragr}[Fin de la démonstration.] --- \label{S:findemodensite}Soit $a_0 \in \Ac(F)$ qui échappe aux deux cas précédents. Alors $a_0$ est l'image d'un élément de la forme $(A,0,0)$ où $A$ est une homothétie de $V$. Soit $\al\in F$ tel que $a_0\in \Ac_\al(F)$ (cf.  §\ref{S:RTFinf}, $\al$ est la trace de $A$). Pour tout $g\in \Sc(\tggo(\AAA))$ on a 
$$I_{a_0}^{G,\eta}(g)=I_{0}^{G,\eta}(g_0)
$$
où $g_0(X)=g(X+(A,0,0))$ et $0$ est la classe du triplet $(0,0,0)$. Localement, la translation par un tel $(A,0,0)$ préserve $\Sc(\tggo(\AAA_S))_\eta$. Il s'ensuit qu'il suffit de considérer le cas $a_0=0$ c'est-à-dire $\al=0$.

Tout élément $a\in \Ac^0(F)$ avec $a\not=0 $ tombe dans l'un des deux cas traités aux §§\ref{S:casrss} et \ref{S:cas-desc}. Soit $f^S\in \Sc(\tggo(\AAA_S)$. La distribution 
$$f\in \Sc(\tggo(\AAA_S))\mapsto I_a^{G,\eta}(f\otimes f^S)
$$ 
est donc $\eta$-stable pour $a\not= 0$. D'après le corollaire \ref{cor:TFstable}, il en est de même pour la distribution 
$$f\in \Sc(\tggo(\AAA_S))\mapsto \hat{I}_a^{G,\eta}(f\otimes f^S)
$$ 
pour n'importe quelle transformée de Fourier partielle de $\hat{I}$ de $I$. Supposons de plus qu'on a $f\in \Sc(\tggo(\AAA_S))_\eta$. On suppose que la transformée de Fourier partielle $I\mapsto \hat{I}$ est l'une des trois considérées dans la remarque \ref{rq:TFP}. La formule des traces infinitésimales (théorème \ref{thm:RTFinf} assertion 3 pour $\al=0$) se simplifie alors en l'égalité
\begin{equation}
  \label{eq:RTF0}
I_{0}^{G,\eta}(f\otimes f^S)=\hat{I}_{0}^{G,\eta}(f\otimes f^S).
\end{equation} 
Soit $u\notin S$  place non-archimédienne auxiliaire et  $S'=S\cup \{u\}$. Soit $f^{S'}\in \Sc(\tggo(\AAA^{S'}))$. Soit $\nc=\tggo_{0}(F_u)$ le cône nilpotent. Soit $D$ la distribution sur $\Cc(\tggo(F_u))$ définie par 
$$D(f_u)=I_{0}^{G,\eta}(f\otimes f_u \otimes f^{S'}).
$$
pour $f_u\in \Cc(\tggo(F_u))$.
Alors $D$ vérifie les trois hypothèses du théorème \ref{thm:incertitude} ci-dessous pour le caractère de $G(F_u)$ obtenu par restriction de $\eta$. Les deux premières propriétés résultent immédiatement des propriétés de la distribution globale   $I_{0}^{G,\eta}$ (cf. théorème \ref{thm:I}). Vérifions la propriété 3. Pour tout $f_u\in \Cc(\tggo(F_u))$, on a, en utilisant \eqref{eq:RTF0},
\begin{eqnarray*}
  \hat{D}(f_u)&=& I_{0}^{G,\eta}(f\otimes \hat{f}_u \otimes f^{S'})\\
&=&  \hat{I}_{0}^{G,\eta}(f\otimes \hat{f}_u \otimes f^{S'})\\
&=& I_{0}^{G,\eta}(\hat{f}\otimes \hat{\hat{f}}_u \otimes \hat{f}^{S'}).
\end{eqnarray*}
En utilisant le fait que $\hat{\hat{f}}_u=f_u(-X)$, on voit que la propriété de support de $\hat{D}$ résulte elle-aussi du contrôle du support de la distribution globale   $I_{0}^{G,\eta}$. Le   théorème \ref{thm:incertitude} implique donc qu'on a $D=0$ ce qu'on voulait démontrer.
\end{paragr}

\begin{paragr}[Le principe d'incertitude d'Aizenbud.] --- Dans ce paragraphe, $F_u$ est  un corps local non-archimédien.   Soit  $V$ est un $F_u$-espace vectoriel de dimension $n$  et $\tgl(V)=\tgl_{F_u}(V)$ muni de l'action de $GL(V)=GL_{F_u}(V)$.  Soit $\nc\subset \tgl(V)$ le cône nilpotent, c'est-à-dire le fermé de  $\tgl(V)$ formé des éléments dont l'invariant $a$ est nul. On note $\supp$ le support d'une distribution. Le théorème suivant est dû à Aizenbud  (cf. \cite{Aiz} section 6.3 où il est énoncé pour $\eta$ trivial ; un point crucial étant le lemme 8.1 de \cite{rallSchiff}). Le cas $\eta$ trivial nous suffirait d'ailleurs (il suffirait de prendre ci-dessus la place $u$ déployée dans $F$). 

  \begin{theoreme} (Aizenbud) \label{thm:incertitude}
Soit $\eta$ un caractère quadratique  de $GL(V)$.  Toute distribution $D$ sur $\tgl(V)$ qui vérifie les trois propriétés suivantes :
\begin{enumerate}
\item $D$ est $\eta$-équivariante ;
\item $\supp(D)\subset \nc$ ;
\item  $\supp(\hat{D})\subset \nc$ pour l'une des trois transformées de Fourier considérées dans la remarque \eqref{rq:TFP}
\end{enumerate}
est nulle.
  \end{theoreme}
\end{paragr}

\section{La descente des distributions globales} \label{sec:desc}


\label{ssec:intro-desc}\subsection{Introduction}

\begin{paragr}
  Cette section a pour finalité de démontrer le théorème \ref{thm:desc}. On reprend les notations des sections précédentes.
\end{paragr}

\begin{paragr}Soit $\tlP$ un sous-groupe parabolique de $\tlG$ contenant $\tlT_0$. Soit $\tlP=\tlM \tlN$ sa décomposition de Levi associée. Soit $\tmgo$ l'espace associé (cf. §\ref{S:Levi}).  Soit $X\in \tmgo$ et $V=V^+\oplus V^-$ la somme directe qui lui est associée (cf. §\ref{S:Xr}). Soit $X_s$ sa partie semi-simple (cf. §§\ref{S:Jor1} et \ref{S:Jor2}). Si  $V_{\tilde{P}}'= V_{\tilde{P}}^\perp$ (cf.§\ref{S:parab}), on obtient que $X\in \mgo$, que $X_s$ est la partie semi-simple usuelle de $X$ et donc qu'on a $X_s\in \mgo$. Dans ce cas, on a, bien sûr, $V^+=0$. Si  $V_{\tilde{P}}'\subsetneq V_{\tilde{P}}^\perp$, on écrit $X=(X_1,\ldots,X_r)$ selon la décomposition $\tmgo=\gl(V_1)\oplus\ldots \oplus \tgl(V_j)\oplus\ldots \oplus\gl(V_r)$. On a $X_s=(X_{1,s},\ldots, X_{r,s})$ où, sur chaque facteur d'indice distinct de $j$, on prend la partie semi-simple usuelle et sur le facteur $\tgl(V_j)$ on prend la partie semi-simple relative de $X_j$ telle que définie en   §§\ref{S:Jor1} et \ref{S:Jor2}. En particulier, $X_s\in \tmgo$. Soit $V_j=V_j^+\oplus V_j^-$ la somme directe associée à $X_j$. Elle est  reliée à celle associée à  $X$  de la manière suivante : on a $V^+=V_j^+$ et $V^-=V_1\oplus \ldots \oplus V_j^- \oplus \ldots \oplus V_r$. On a $V_{\tilde{P}} \subset V^-$ et $V_{\tilde{P}}'\subset (V^-)^*$.
\end{paragr}

\begin{paragr}[Centralisateurs.] ---\label{S:centra}  On écrit ensuite $X=\iota_X(X^+,X^-)$ avec $X^+\in \tgl(V^+)^{(r)}$ et $X^-\in \tgl(V^-)^{(0)}$. La partie semi-simple $X_s^-$ est donc de la forme $(Y,0,0)\in  \tgl(V^-)$. Soit $\gl(V^-,X_s^-)$ la sous-algèbre de Lie des $Z\in \gl(V^-)$ qui vérifient $[Z,Y]=0$. On identifie d'ailleurs $\gl(V^-)$ à la sous-algèbre de $\gl(V)$ des endomorphismes qui laissent stable $V^-$ et sont nuls sur $V^+$. 

Soit $N$ le radical unipotent de $P=\tlP\cap G$ et $\ngo$ son algèbre de Lie. Soit $N_{X_s}$ le centralisateur de $X_s$ dans $N$ et $\tngo(X_s)$ le sous-espace de $\tgl(V^-)=\mathfrak{gl}(V^-)\oplus  V^-\oplus (V^-)^*$ défini par
$$\tngo(X_s)=[\ngo\cap \gl(V^-,X_s^-)]\oplus V_{\tilde{P}}\oplus V_{\tilde{P}}'.$$

\begin{lemme}\label{lem:ss-U}
  Pour tout $U\in \tngo(X_s)$, on a 
$$\iota_X(X^+,X^-+U)_s=X_s.$$
\end{lemme}

\begin{preuve}
On a par définition   $\iota_X(X^+,X^-+U)_s=\iota_X(X^+_s,(X^-+U)_s)$. Il s'agit donc de vérifier l'égalité $X_s^-=(X^-+U)_s$. Or $X^-+U$ et $X^-$ ont même image dans $\Ac_{V^-}$ et cette image est d'ailleurs dans  $\Ac_{V^-}^{(0)}$. La construction du §\ref{S:Jor1} donne alors  $(X^-+U)_s=X_s^-$.
\end{preuve}
\end{paragr}

\subsection{Un premier noyau auxiliaire}

\begin{paragr}
  Le but de cette section est de prouver qu'on peut remplacer dans la contruction des distributions globales (cf. section \ref{sec:RTFinf}) le noyau $k^T_a$ (cf. \eqref{eq:k}) par un noyau plus adapté à la descente.
\end{paragr}

\begin{paragr}[Un nouveau noyau tronqué.] --- Soit $f\in \Sc(\tggo(\AAA))$,  $a\in \Ac(F)$ et $\tlP$ un sous-groupe parabolique relativement standard de $\tlG$. On reprend les notations des §§\ref{S:noyau} et \ref{S:noy-tronque}. On introduit la variante suivante de \eqref{eq:k}
  \begin{equation}
    \label{eq:j}
    j_{\tilde{P},a}(f,g)=\sum_{X\in \tmgo_{\tilde{P},a}(F)} \sum_{\nu \in N_{X_s}(F)\back N(F)} \int_{\tilde{\ngo}(X_s,\AAA)} f((\nu g)^{-1}\cdot \iota_X(X^+,X^-+U))dU. 
  \end{equation}

On a également la variante suivante de \eqref{eq:kT} pour $T\in \ago_{\tl}^+$
\begin{equation}
  \label{eq:jT}
  j^T_a(f,g)=\sum_{\tilde{P}\in \fc^{\tlG}(B)} \eps_{\tilde{P}}^{\tilde{G}} \sum_{\delta\in P(F)\back G(F)} \hat{\tau}_{\tilde{P}}(H_{\tilde{P}}(\delta g)-T_{\tlP})j_{\tilde{P},a}(f,\delta g).
\end{equation}

\end{paragr}

\begin{paragr}[Théorème de comparaison.] --- Il s'agit du théorème suivant. 

  \begin{theoreme}\label{thm:cvj}
Soit $a\in \Ac(F)$ et $f\in \Cc(\tggo(\AAA))$. Il existe $T_f$ tel que pour tout $T\in T_f+\ago_{\tl}^+$   on ait 
\begin{equation}
  \label{eq:cv}
  \int_{[G]} |j^T_a(f,g)  \eta(g)|\, dg <\infty
\end{equation}
et
\begin{equation}
  \label{eq:j=k}
   \int_{[G]} j^T_a(f,g)  \eta(g) \, dg = \int_{[G]} k^T_a(f,g)   \eta(g) \, dg.
\end{equation}

\end{theoreme}

La démonstration du théorème \ref{thm:cvj} va nous occuper jusqu'à la fin de cette section. On commence par un résultat annexe.

\end{paragr}

\begin{paragr}[Un lemme sommatoire.] --- On reprend les notations de la section \ref{ssec:intro-desc}.
  
\begin{lemme}\label{lem:resommation}Soit $\tlP$ un sous-groupe parabolique contenant $\tlT_0$. Soit $\tmgo=\tmgo_{\tlP}$ et $\tngo=\tngo_{\tlP}$.

  \begin{enumerate} 
  \item Soit $X\in \tmgo(F)$. Pour toute  fonction  $f$ absolument sommable sur $\tngo(F)$,
$$
\sum_{U\in \tngo(F)} f(U)=\sum_{\delta \in N_{X_s}(F)\back N(F)} \sum_{U\in \tngo(X_s,F)} f(\delta^{-1}\cdot \iota_X(X^+,X^-+U) -X)
$$
\item Soit $X\in \tmgo(\AAA)$. Pour toute fonction $f\in L^1(\tngo(\AAA))$, on a
$$
\int_{\tngo(\AAA)} f(U)\, dU=\int_{N_{X_s}(\AAA)\back N(\AAA)} \int_{\tngo(X_s,\AAA)} f(n^{-1}\cdot \iota_X(X^+,X^-+U) -X)\, dUdn.
$$
 
  \end{enumerate}
\end{lemme}

\begin{preuve}
On se contente de prouver la première assertion, la seconde se démontrant de manière similaire.

  Traitons d'abord le cas où $X\in \tgl(V)^{(0)}(F)$.   On écrit $X=(Y,b,c)$. Dans ce cas $V^-=V$, $X_s=(Y_s,0,0)$, $\tngo(X_s)= \ngo_{Y_s}\oplus V_{\tilde{P}} \oplus V_{\tilde{P}}'$ et le membre de droite de l'égalité à prouver devient
$$\sum_{\delta \in N_{Y_s}(F)\back N(F)} \sum_{U\in \ngo_{Y_s}(F)} \sum_{v\in V_{\tilde{P}}(F)}\sum_{w\in V_{\tilde{P}}'(F)}  f(\delta^{-1}(Y+U)\delta -Y, \delta^{-1}(v+b)-b, (w+c)\delta -c).
$$
On utilise ensuite le fait que $v\mapsto  \delta^{-1}(v+b)-b$ et  $w\mapsto (w+c)\delta -c$ induisent des automorphismes respectifs de $V_{\tilde{P}}$ et $V_{\tilde{P}}'$. On obtient alors
$$\sum_{v\in V_{\tilde{P}}(F)}\sum_{w\in V_{\tilde{P}}'(F)} \sum_{\delta \in N_{Y_s}(F)\back N(F)} \sum_{U\in \ngo_{Y_s}(F)}  f(\delta^{-1}(Y+U)\delta -Y, v, w).
$$
La double somme intérieure se simplifie en une somme sur $\ngo(F)$ (cf. \cite{PH1} corollaire 2.4) si bien qu'on obtient
$$\sum_{v\in V_{\tilde{P}}(F)}\sum_{w\in V_{\tilde{P}}'(F)}  \sum_{U\in \ngo(F)}  f(U, v, w)
$$
qui est bien le membre de gauche  dans l'assertion 1.

Traitons ensuite le cas où $X\in \tgl(V)^{(r)}(F)$ avec $r\geq 1$. Comme d'habitude, on a  $X=\iota(X^+,X^-)$ (pour alléger on note simplement $\iota=\iota_X$). Dans la suite, on identifie  $GL(V^-)$ au sous-groupe de $GL(V)$ qui stabilise $V^-$ et agit trivialement sur $V^+$.

On dispose du drapeau de $V^-$
$$(0)\subsetneq V_1 \subsetneq \ldots\subsetneq  (V_1\oplus\ldots \oplus V_{j-1}) \subset (V_1\oplus\ldots \oplus V_{j-1}\oplus V_j^-) \subsetneq (V_1\oplus\ldots \oplus V_{j-1}\oplus V_j^-\oplus V_{j+1}) \subsetneq \ldots \subsetneq V^-.
$$
On fait le choix du couple $(j-1,j)$ si $V_j^-\not=0$ et $(j-1,j-1)$ sinon.  Suivant les conventions de §\ref{S:parab}, on note $\tlQ$ de telles données. On a donc un  sous-groupe parabolique $Q$ de $GL(V^-)$ qui stabilise le drapeau ci-dessus et un espace $\tngo_{\tilde{Q}}=\ngo_Q\oplus V_{\tilde{P}}\oplus V_{\tilde{P}}'$.  

On a défini au §\ref{S:centra}, un sous-espace  $\tngo_{\tilde{Q}}(X_s^-)$ relatif à $X_s^-$ : on a  $\tngo_{\tilde{Q}}(X_s^-)=\ngo_{Q,X_s^-}\oplus V_{\tilde{P}}\oplus V_{\tilde{P}}'$. Observons  que le centralisateur  $N_{X_s}$ de $X_s$ dans $N$  n'est autre que le centralisateur  de $X^-_s$ dans $N_Q$. Il s'ensuit qu'on a $\ngo_{Q,X_s^-}=\ngo_{X_s}$. Par conséquent, avec les notations de  §\ref{S:centra}, on a aussi $\tngo(X_s)=\tngo_{\tilde{Q}}(X_s^-)$.

 Le membre de droite de l'assertion 1 s'écrit donc
$$
\sum_{\delta \in N_{Q}(F)\back N(F)}\sum_{\nu\in N_{Q,X_s^-}(F)\back N_Q(F)} \sum_{U\in \tngo(X_s,F)} f(\delta^{-1}\cdot \iota(X^+, \nu^{-1 }(X^-+U)\nu)-X).
$$
En utilisant l'assertion 1 pour $X^-\in \tgl(V^-)^{(0)}$ cas qu'on vient de démontrer, on obtient
\begin{equation}
  \label{eq:mb}
  \sum_{\delta \in N_{Q}(F)\back N(F)} \sum_{U\in \tngo_{\tlQ}(F)} f(\delta^{-1}\cdot \iota(X^+, X^-+U) -X).
\end{equation}

Soit $(Y,v,w)\in X+\tngo$. Un tel élément a même invariant $a$ que $X$ donc il appartient à $\tgl(V)^{(r)}$. Donc l'espace engendré par $Y^iv$ pour $0\leq i <r$ est un sous-espace de $V^+\oplus V_{\tilde{P}}$ de dimension $r$. De même, l'espace engendré par $wY^i$ pour $0\leq i <r$ est un sous-espace de $(V^+)^*\oplus V_{\tilde{P}}'$ de dimension $r$. Quitte à faire agir $N$, on peut supposer ces espaces engendrés sont exactement  $V^+$ et $(V^+)^*$. Mais alors on a $(Y,v,w)=\iota(X^+,X^-+U)$ pour un unique élément $U\in \tngo_{\tilde{Q}}$. Maintenant si deux éléments de $\iota(X^+,X^-+ \tngo_{\tilde{Q}})$ se déduisent l'un de l'autre par l'action de $n\in N$ alors $n\in N_Q$. Il s'ensuit d'une part que le centralisateur $N_Z$ de $Z\in  \iota(X^+,X^-+ \tngo_{\tilde{Q}})$ dans $N$ est inclus dans $N_Q$ et d'autre part  que les orbites de $X+\tngo(F)$ sous $N(F)$ sont en bijection avec les orbites de  $\iota(X^+,X^-+ \tngo_{\tilde{Q}}(F))$ sous $N_Q(F)$.

En partant du membre de droite dans l'assertion 1, on peut donc écrire avec l'abus de  confondre un quotient avec un système de représentants
\begin{eqnarray*}
  \sum_{Z\in  X+\tngo(F) } f(Z-X) &=&  \sum_{ Z  \in  X+\tngo(F)/N(F) }  \sum_{\nu \in N_Z(F)\back N(F) } f(\nu^{-1}\cdot Z-X)\\
&=&   \sum_{ Z  \in  \iota(X^+,X^-+ \tngo_{\tilde{Q}}(F))/N_Q(F) }  \sum_{\nu \in N_Z(F)\back N(F) } f(\nu^{-1}\cdot Z-X)\\
&=&  \sum_{\delta \in N_Q(F)\back N(F) }  \sum_{ Z  \in  \iota(X^+,X^-+ \tngo_{\tilde{Q}}(F))/N_Q(F) }  \sum_{\nu \in N_Z(F)\back N_Q(F) } f((\delta \nu^{-1})\cdot Z-X)\\
&=&  \sum_{\delta \in N_Q(F)\back N(F) }  \sum_{ Z  \in  \iota(X^+,X^-+ \tngo_{\tilde{Q}}(F))}  f(\delta \cdot Z-X)\\
\end{eqnarray*}
qui est bien \eqref{eq:mb}.
\end{preuve}
\end{paragr}

\begin{paragr}[Première réduction.] --- Dans toute la suite, on fixe $a\in \Ac(F)$ et $f\in \Cc(\tggo(\AAA))$.

Par les manipulations de \cite{leMoi2}, preuve du théorème 3.6,  il suffit de prouver pour $B\subset \tilde{P}_1 \subsetneq    \tilde{P}_2$ qu'on a

$$ \int_{P_1(F)\back G(\AAA)} \chi_{ \tilde{P}_1, \tilde{P}_2  }^T(g) |j_a^{1,2}(f,g)  |\, dg <\infty
$$
où
$$
j_a^{1,2}(f,g)=\sum_{  \tilde{P}_1 \subset   \tilde{P}  \subset    \tilde{P}_2 }\eps_{ \tilde{P}   }^{ \tilde{G}   }\cdot  j_{\tilde{P},a}(f,g)
$$
et  $\chi_{ \tilde{P}_1, \tilde{P}_2  }^T$ est la fonction définie en \emph{loc. cit.} (2.4) : elle vaut $0$ ou $1$. Désormais, on fixe  $\tilde{P}_1 \subsetneq    \tilde{P}_2$.
\end{paragr}

\begin{paragr}[Deuxième réduction.] --- 

  \begin{lemme}
     Pour tout  $T$ dans un translaté de $\ago_{\tl}^+$,    tout $g\in P_1(F)\back G(\AAA)$ tel que $\chi_{ \tilde{P}_1, \tilde{P}_2  }^T(g)=1$ et tout $\tilde{P}_1 \subset    \tilde{P} \subset   \tilde{P}_2$, on a 
$$
j_{\tilde{P},a}(f,g)=\sum_{X\in (\tmgo_{\tilde{P},a}\cap \tpgo_1)(F)} \sum_{\nu \in N_{X_s}(F)\back N(F)} \int_{\tilde{\ngo}(X_s,\AAA)} f((\nu g)^{-1}\cdot \iota_X(X^+,X^-+U))\,dU.
$$
  \end{lemme}

  \begin{preuve}
    \emph{A priori}, le membre de gauche de l'égalité à prouver est donné par
$$
j_{\tilde{P},a}(f,g)=\sum_{X\in \tmgo_{\tilde{P},a}(F)} \sum_{\nu \in N_{X_s}(F)\back N(F)} \int_{\tilde{\ngo}(X_s,\AAA)} f((\nu g)^{-1}\cdot \iota_X(X^+,X^-+U))dU.
$$ 
Il résulte du corollaire 2.5 de \cite{leMoi2} et de \cite{ar1} pp.943-944 qu'il existe un compact $\om\subset (N_B\cap M_2)(\AAA) T_0(\AAA)^1$ et une constante $ c $ tous deux indépendants de $T$ tels que sous l'hypothèse   $\chi_{ \tilde{P}_1, \tilde{P}_2  }^T(g)=1$,  il existe un sous-groupe de Borel $ \tlB \in \fc^{\tlG}(B) $ contenu dans $ \tlP_{1} $ tel que $g\in B(F) N_2(\AAA) A^{\infty}_{1,2}(\tlB,T)\om  K$
où  l'on définit 
$$A^{\infty}_{1,2}(\tlB,T)=\{ a\in A^\infty_{T_0} \mid \al(H_B(a)-T_{\tlB}) >0 \ \forall \al\in \Delta_{\tilde{B}}^{\tilde{P}_2 }\smallsetminus\Delta^{\tilde{P}_1}_{\tilde{B}}, 
\  \al(H_B(a) > c \ \forall \al \in \Delta_{\tlB}^{\tlP_{1}} 
\}.
$$
Quitte à translater, on peut et on va supposer que $g=na y$ avec $ n\in N_1^2(\AAA)$, $a\in A^{\infty}_{1,2}(\tlB,T)$ et $y\in\om  K$. 
Soit $X\in \tmgo_{\tilde{P},a}(F)$, $\nu \in N(F)$ et $U\in\tilde{\ngo}(X_s,\AAA)$ tel que $f((\nu g)^{-1}\cdot \iota_X(X^+,X^-+U))$ soit non nul. La composante sur $\tmgo_{\tilde{P}}$ de $(\nu na )^{-1}\cdot \iota_X(X^+,X^-+U)$ est donc $a^{-1}\cdot X$. Cette dernière est astreinte à rester dans un compact qui ne dépend que de $\om$, $K$ et du support de $f$. Si $X\notin \tpgo_1(F)$, cet élément  admet au moins une composante non nulle sur un espace propre pour $A_{T_0}^\infty$ et pour un caractère qui est la restriction à $A_{T_0}^\infty$ d'une somme à coefficients négatifs de racines dans  $\Delta_{\tilde{B}}^{\tilde{P}_2 }$ dont une composante au moins sur $\Delta_{\tilde{B}}^{\tilde{P}_2 }\smallsetminus\Delta^{\tilde{P}_1}_{\tilde{B}}$ est non nulle. Par conséquent,  si $T$ est assez positif, on voit que $a^{-1}\cdot X$  sort de ce compact. Cela conclut.
  \end{preuve}

\end{paragr}

\begin{paragr} Pour alléger, on écrit parfois en indice ou exposant $1$ ou $2$ à la place de $\tlP_1$ et $\tlP_2$. Avec la notation   $\tngo_{ \tilde{P}_1}^{ \tilde{P}}$ du §\ref{S:nPQ}, on introduit
$$
j_{\tilde{P},X}(f,g)= \sum_{Z\in  \tngo_{ \tilde{P}_1}^{ \tilde{P}}(F)} \sum_{\nu \in N_{Y_s}(F)\back N(F)} \int_{\tilde{\ngo}(Y_s,\AAA)} f((\nu g)^{-1}\cdot\iota_Y(Y^+,Y^-+U))dU 
$$
où pour un $Z$ comme dans la somme ci-dessus, on pose $Y=X+Z$. On pose alors
   $$
 j'_{\tilde{P},a}(f,g)=\sum_{X\in  \tmgo_{\tilde{P}_1,a}(F)}  j_{\tilde{P},X}(f,g).
$$

On utilise ensuite le lemme \eqref{lem:resommation} (ou plutôt une variante évidente adaptée à $\tmgo$) pour réécrire la somme sur $Z$. On obtient
$$
j_{\tilde{P},X}(f,g)= \sum_{\delta\in N_{1,X_s}^P(F)\back N_1^P(F)}  \sum_{Z\in  \tngo_{ \tilde{P}_1}^{ \tilde{P}}(X_s,F)} \sum_{\nu \in N_{Y_s}(F)\back N(F)} \int_{\tilde{\ngo}(Y_s,\AAA)} f((\nu g)^{-1}\cdot \iota_Y(Y^+,Y^-+U))dU 
$$
où, cette fois-ci, $Y=\delta^{-1}\cdot \iota_X(X^+,X^-+Z)$ et où l'on note  $\tngo_{ \tilde{P}_1}^{ \tilde{P}}(X_s)= \tngo_{ \tilde{P}_1}(X_s)\cap \tmgo_{\tlP}$.
Observons qu'on a
\begin{enumerate}
\item $Y_s=\delta^{-1}\cdot X_s$ (cf. lemmes \ref{lem:action} et \ref{lem:ss-U}) ;
\item $Y^+=\delta^{-1}\cdot X^+$ ;
\item $Y^-=\delta^{-1}\cdot X^-$ ;
\item $\iota_Y(Y^+,Y^-+U)=\delta^{-1}\cdot\iota_X(X^+,X^-+Z+\delta\cdot U)$.
\end{enumerate}
De là, on déduit via un changement de variables qu'on a 
$$
j_{\tilde{P},X}(f,g)=  \sum_{\nu \in N_{1,X_s}(F)\back N_1(F)}  \sum_{Z\in  \ngo_{ \tilde{P}_1}^{ \tilde{P}}(X_s,F)} \int_{\tilde{\ngo}(X_s,\AAA)} f((\nu g)^{-1}\cdot \iota_X(X^+,X^-+Z+U))dU.
$$
On applique ensuite la formule sommatoire de Poisson à la somme sur $Z$. Il vient
$$
j_{\tilde{P},X}(f,g)= \sum_{\nu \in N_{1,X_s}(F)\back N_1(F)}  \sum_{Z\in  \bar{\ngo}_{ \tilde{1}}^{ \tilde{P}}(X_s,F)}\Phi(\nu g,X,Z)
$$
où
\begin{itemize}
\item on note $\bar{\ngo}_{ \tilde{P}_1}^{ \tilde{P}}(X_s)= \tngo_{ \overline{\tilde{P}_1}}(X_s)\cap \tmgo_{\tlP}$ avec $\overline{\tilde{P}_1}$ le sous-groupe parabolique opposé à $\tlP_1$ ;
\item on introduit l'expression auxiliaire
$$\Phi(g,X,Z)= \int_{\tngo_{\tilde{P}_1}(X_s,\AAA)} f(g^{-1}\cdot \iota_X(X^+,X^-+U)) \psi(\bg Z,U\bd)\,  dU.
$$
\end{itemize}
Comme dans \cite{leMoi2}, preuve du théorème 3.7, on aboutit à 
$$
j_{a}^{1,2}(f,g)= \eps_{\tilde{P}_2}^{\tilde{G}} \sum_{X\in \tmgo_{\tilde{P}_1,a}(F)} \sum_{\nu \in N_{1,X_s}(F)\back N_1(F)}\sum_{Z\in  \bar{\ngo}_{ \tilde{P}_1}^{ \tilde{P}_2}(X_s,F)'}\Phi(\nu g,X,Z)
$$

où $\bar{\ngo}_{\tilde{P}_1}^{ \tilde{P}_2}(X_s,F)'=\bar{\ngo}_{\tilde{P}_1}^{ \tilde{P}_2}(X_s,F) \setminus \bigcup_{\tilde{P_1}\subset \tilde{R} \subsetneq \tilde{P}_2 }\bar{\ngo}_{\tilde{P}_1}^{ \tilde{R}}(X_s,F)$.

Il suffit donc de prouver qu'on a
\begin{equation}
  \label{eq:maj12}
   \int_{P_1(F)\back G(\AAA)} \chi_{ \tilde{P}_1, \tilde{P}_2  }^T(g)   \sum_{X\in \tmgo_{  \tilde{P}_1,a}(F)} \sum_{\nu \in N_{1,X_s}(F)\back N_1(F)}\sum_{Z\in  \bar{\ngo}_{ \tilde{P}_1}^{ \tilde{P}_2}(X_s,F)'} |\Phi (\nu g, X,Z)  |\, dg <\infty
 \end{equation}
 
\end{paragr}

\begin{paragr} Pour tous $g\in G(\AAA)$ et $X\in  \tmgo_{  \tilde{P}_1}(F)$, on a  en utilisant le fait que $\vol([N_{1,X_s}])=1$

$$
 \int_{N_{1,X_s}(F)\back N_1(\AAA)}\sum_{Z\in  \bar{\ngo}_{ \tilde{P}_1}^{ \tilde{P}_2}(X_s,F)'} |\Phi (ng, X,Z)  |\, dn
$$
$$
= \int_{N_{1,X_s}(\AAA)\back N_1(\AAA)}\sum_{Z\in  \bar{\ngo}_{\tilde{P}_1 }^{ \tilde{P}_2}(X_s,F)'} |\Phi (ng, X,Z)  |\, dn
$$
$$
=\int_{N_{1,X_s}(\AAA)N_2(\AAA)\back N_1(\AAA)} \int_{ N_{2,X_s}(\AAA)\back N_2(\AAA) }\sum_{Z\in  \bar{\ngo}_{ \tilde{P}_1}^{ \tilde{P}_2}(X_s,F)'} |\Phi (n_2 ng, X,Z)  |\, dn_2 dn
$$
Pour tous $X\in  \tmgo_{  \tilde{P}_1}(F)$ et $Z\in  \bar{\ngo}_{ \tilde{P}_1}^{ \tilde{P}_2}(X_s,F)'$, l'application  $g\in G(\AAA)\mapsto \Phi (g, X,Z) $ est invariante par translatation à gauche de $N_{2,X_s}(\AAA)$. On a donc
\begin{eqnarray*}
    \int_{ N_{2,X_s}(\AAA)\back N_2(\AAA) } |\Phi (n_2 g, X,Z)  |\, dn_2\\
\leq  \int_{ N_{2,X_s}(\AAA)\back N_2(\AAA) }  \int_{\tngo_{2}(X_s,\AAA)} | \int_{\tngo^2_{\tilde{P}_1}(X_s,\AAA)}f((n_2g)^{-1}\cdot \iota_X(X^+,X^-+U+U_2)) \psi(\bg Z,U\bd)\,  dU| \, dU_2 dn_2\\
= \int_{\tngo_{2}(\AAA)} | \int_{\tngo^2_{\tilde{P}_1}(X_s,\AAA)}f(g^{-1}\cdot( \iota_X(X^+,X^-+U) +U_2)) \psi(\bg Z,U\bd)\,  dU| \, dU_2 
\end{eqnarray*}
la dernière égalité résultant du  lemme \ref{lem:resommation}. 

Après ces observations, on revient à l'étude de l'intégrale \eqref{eq:maj12}. D'après le § \ref{S:Levi-dec}, on a une décomposition 
$$M_1= \MM_{\tlP_1} \times G_{\tlP_1}.
$$
Soit $M_1(\AAA)'=\MM_{\tlP_1}(\AAA)^1 \times G_{\tlP_1}(\AAA)$. On a donc $M_1(\AAA)= \mathbf{A}_{1} M_1(\AAA)'$ où l'on pose $\mathbf{A}_{1}=A_{\MM_{\tlP_1}}^\infty$.
Par la décomposition d'Iwasawa,  on remplace l'intégrale sur le quotient $P_1(F)\back G(\AAA)$ par une intégrale multiple  sur $[N_1]$,  $M_1(F)\back M_1(\AAA)'$, $\mathbf{A}_1$ et $K$. L'intégration sur $[N_1]$ se combine avec la somme sur $\nu \in N_{1,X_s}(F)\back N_1(F)$.  D'après ce qui précède, on est ramené à majorer l'intégrale suivante

$$\int_K \int_{ M_1(F)\back M_1(\AAA)'}  \int_{\mathbf{A}_1}  e^{-2\rho_1(H_B(t))} \chi_{ \tilde{P}_1, \tilde{P}_2  }^T(tm)  \sum_{X\in  \tmgo_{  \tilde{P}_1}(F)} \int_{N_{1,X_s}(\AAA)N_2(\AAA)\back N_1(\AAA)} $$

$$
\sum_{Z\in  \bar{\ngo}_{ \tilde{P}_1}^{ \tilde{P}_2}(X_s,F)'} \int_{\tngo_{2}(\AAA)} | \int_{\tngo^2_{\tilde{P}_1}(X_s,\AAA)}f((ntmk)^{-1}\cdot( \iota_X(X^+,X^-+U) +U_2)) \psi(\bg Z,U\bd)\,  dU| \, dU_2  dn dt dm dk
$$
Sous la condition 
$$\chi_{ \tilde{P}_1, \tilde{P}_2  }^T(tm)=F^{\tlP_1}(m,T) \sigma_{\tlP_1}^{\tlP_2}(H_{\tlP_1}(tm)-T)=1,$$
on peut remplacer les intégrales sur $m$ et $k$ par la borne supérieure lorsque $m$ parcourt un compact fixé de $M_1(\AAA)'$ et $k$ parcourt $K$ (cf. corollaire 2.5 de \cite{leMoi2}). Observons que $\mathbf{A}_1$ centralise $X$ et $X_s$. La composante sur $\tmgo_1$ de $(nt)^{-1}\cdot( \iota_X(X^+,X^-+U) +U_2)$ est $X$. Par conséquent, il y a un nombre fini de $X$ qui intervienne. On peut et on va donc supposer que $X$ est fixé. Comme $m$ est astreint à rester dans un compact, on peut le supposer fixé et quitte à changer $T$ et $f$ on va supposer $m=1$.

Pour conclure, il suffit de montrer que l'expression ci-dessous vue comme fonction de $t\in \mathbf{A}_1$ est à décroissance rapide sur le cône $\sigma_{\tlP_1}^{\tlP_2}(H_{\tlP_1}(t)-T)=1$

$$
\int_{N_{1,X_s}(\AAA)N_2(\AAA)\back N_1(\AAA)}\sum_{Z\in  \bar{\ngo}_{ \tilde{P}_1}^{ \tilde{P}_2}(X_s,F)'} \int_{\tngo_{2}(\AAA)} | \int_{\tngo^2_{\tilde{P}_1}(X_s,\AAA)}f((nt)^{-1}\cdot( \iota_X(X^+,X^-+U) +U_2)) \psi(\bg Z,U\bd)\,  dU| \, dU_2  dn
$$
Observons qu'on a pour $t\in \mathbf{A}_1$ 

$$t^{-1}\cdot\iota_X(X^+,X^-+U)=t\iota_X(X^+,X^-+ t^{-1}\cdot U).$$

À un jacobien près (qui ne pertube le fait d'être à décroissance rapide), on est ramené, par des changements de varaiables,  à traiter l'expression
$$\int_{N_{1,X_s}(\AAA)N_2(\AAA)\back N_1(\AAA)}\sum_{Z\in  \bar{\ngo}_{ \tilde{P}_1}^{ \tilde{P}_2}(X_s,F)'} \int_{\tngo_{2}(\AAA)} | \int_{\tngo^2_{\tilde{P}_1}(X_s,\AAA)}f(n^{-1}\cdot \iota_X(X^+,X^-+U) +U_2)) \psi(\bg Z,t\cdot U\bd)\,  dU| \, dU_2  dn.
$$
L'intégrande dans l'intégrale sur $n$ est nulle hors d'un compact indépendant de $t$. Il est alors facile de conclure.
\end{paragr}

\begin{paragr}[Preuve de l'égalité \eqref{eq:j=k}.] --- D'après les manipulations ci-dessus, on a, pour $T$ dans un certain translaté de $\ago_{\tl}^+$  (les manipulations sont justifiées par la convergence absolue)
    \begin{eqnarray*}
        \int_{[G]} j^T_a(f,g)  \, \eta(g)dg&=&  \sum_{B\subset \tilde{P}_1 \subsetneq    \tilde{P}_2} \int_{P_1(F)\back G(\AAA)} \chi_{ \tilde{P}_1, \tilde{P}_2  }^T(g) j_a^{1,2}(f,g)  \, \eta(g)dg \\
&=& \int_{K} \int_{[M_1]}  e^{- 2\rho_{P_1}(H_{P_1}(m))} \chi_{ \tilde{P}_1, \tilde{P}_2  }^T(m)   \big( \int_{[N_1]}  j_a^{1,2}(f,nmk) \, dn\big)\,  \eta(mk)dm dk.
    \end{eqnarray*}
Fixons $g\in G(\AAA)$. On a 
\begin{eqnarray*}
  \int_{[N_1]}  j_a^{1,2}(f,ng)\,dn= \sum_{  \tilde{P}_1 \subset   \tilde{P}  \subset    \tilde{P}_2 }\eps_{ \tilde{P}   }^{ \tilde{G}   }\cdot   \int_{[N_1]} j_{\tilde{P},a}(f,ng)\, dn
\end{eqnarray*}
Soit $\tilde{P}$ tel que $\tilde{P}_1 \subset   \tilde{P}  \subset    \tilde{P}_2$. On pose $\tmgo=\tmgo_{\tilde{P}}$ et $N=N_P$ pour soulager les notations. En remarquant que $N_P$ est un sous-groupe distingué de $N_1$, on a 
\begin{eqnarray*}
  \int_{[N_1]} j_{\tilde{P},a}(f,n_1g)\, dn_1 &=&  \int_{[N_1]} j_{\tilde{P},a}(f,n_1g)\, dn_1\\
&=&  \int_{[N] }\int_{[N_1]} j_{\tilde{P},a}(f,n_1n g) \, dn\, dn_1\\
&=& \int_{[N_1] }\int_{[N]} j_{\tilde{P},a}(f,n n_1 g) \, dn\, dn_1\\
&=&  \int_{[N_1] } \sum_{X\in \tmgo_{a}(F)} \int_{N_{X_s}(F)\back N(\AAA)} \int_{\tilde{\ngo}(X_s,\AAA)} f((n n_1 g)^{-1}\cdot \iota_X(X^+,X^-+U))dU \, dn \, dn_1.
\end{eqnarray*}
On applique ensuite le lemme \ref{lem:resommation} ce qui donne
\begin{eqnarray*}
  \int_{[N_1]} j_{\tilde{P},a}(f,n_1g)\, dn_1 &=& \int_{[N_1] } \sum_{X\in \tmgo_{a}(F)} \int_{[N_{X_s}]} \int_{\tilde{\ngo}(\AAA)} f((n n_1 g)^{-1}(X+U))dU \, dn \, dn_1\\
&=& \int_{[N_1] } \sum_{X\in \tmgo_{a}(F)} \int_{\tilde{\ngo}(\AAA)} f((n n_1 g)^{-1}(X+U))dU \, dn_1\\\
&=& \int_{[N_1] } k_{\tilde{P},a}(f,n_1g)\, dn_1.
\end{eqnarray*}
On peut alors rebrousser chemin et retomber sur le membre de droite dans \eqref{eq:j=k}.
\end{paragr}

\subsection{Un deuxième noyau auxiliaire}\label{par:noyTronq}

\begin{paragr}
Soit $ f \in \Sc(\tggo(\AAA) )$ et   $ a \in \Ac(F)$. On va introduire une variante du noyau tronqué \eqref{eq:kT}. On utilise le même noyau $k_{\tlP,a}$ introduit en \eqref{eq:k} mais on remplace les fonctions caractéristiques de cônes $\hat{\tau}_{\tlP}$ par les fonctions $\hat{\sigma}_{\tlP}$ (cf. § \ref{S:sigma}). Pour tout $T \in \ago_{\tl}^{+} $ et tout $x\in G(\AAA)$, on pose
\begin{equation}
  \label{eq:kappaT}
  \kappa_{a}^{T}(x) = \sum_{\tilde{P}\in \fc^{\tlG}(B)}\varepsilon_{\tlP}^{\tlG}
\sum_{\delta \in P(\rmF) \backslash G(\rmF)}
\hat \sigma_{\tlP}(H_{\tlP}(\delta x) - T_{\tlP})
k_{\tlP, a}(\delta x).
\end{equation}
\end{paragr}

\begin{paragr} Rappelons qu'à l'aide du noyau $k_a^T$ et du caractère $\eta$, on définit une distribution $I_a^{\eta}$ (cf. théorème \ref{thm:I}). La proposition suivante indique qu'à cette fin,  le noyau $\kappa_a^T$ joue le même rôle.
  
\begin{proposition}\label{prop:newTroncKer1} Il existe $T_+$ tel que 
\begin{enumerate}
\item Pour tout $T \in T_++ \ago_{\tl}^{+}$, on a
\[
\int_{[G]}
|\kappa_{a}^{T}(x) |dx < \infty.
\]
\item Il existe un unique polynôme-exponentielle en $T$ qui coïncide sur $T_++  \ago_{\tl}^{+}$ avec l'intégrale $\int_{[G]}\kappa_{a}^{T}(x) \eta(x)dx$.  En outre, le terme  purement polynomial  de ce polynôme-exponentielle est constant et égal à  $I_{a}^{\eta}(f)$.
\end{enumerate}
\end{proposition}

\begin{preuve}  On rappelle qu'on a introduit au §\ref{S:a0} des projections $r_i$ et $\hat{r}_i$. En vertu de \eqref{eq:sig-tau}, on a 
$$\hat \sigma_{\tlP}(H_{\tlP}(\delta x) - T_{\tlP}) = \hat \tau_{\tlP} (\hat r_{1}(H_{\tlP}(x)) - T_{\tlP}).$$
En utilisant la formule d'Arthur (cf. \cite{arthur2} section 2)
\begin{equation}\label{eq:GammaPrimeRecurrenceG}
\hat \tau_{\tlP}(H - X) = 
\sum_{\tlQ \supseteq \tlP}
\varepsilon_{\tlQ}^{\tlG}
\hat \tau_{\tlP}^{\tlQ}(H)
\Gamma_{\tlQ}'(H,X), 
\end{equation}
pour  $H = H_{\tlP}(\delta x) - T_{\tlP}$ et $X = H_{\tlP}(\delta x) - \hat r_{1}(H_{\tlP}(\delta x)) = \hat r_{2}(H_{\tlP}(\delta x))$, on voit que l'intégrale en question est égale à la somme sur $\tlQ \in \fc^{\tlG}(B)$ de 
\begin{equation}\label{eq:QsousG}
\int\limits_{Q(\rmF)\backslash  G(\AAA)}
\sum_{\tlP\in \fc^{\tlQ}(B) }
\varepsilon_{\tlP}^{\tlQ}
\sum_{\delta \in (P \cap M_{Q})(\rmF)\backslash M_{Q}(\rmF)}
\Psi^{T}_{\tlP,\tlQ, a}(\delta x)dx
\end{equation}
où l'on introduit
\begin{displaymath}
\Psi^{T}_{\tlP,\tlQ, a}(x) = 
k_{\tlP,a}(x)
\hat \tau_{\tlP}^{\tlQ}(H_{\tlP}(x)- T_{\tlP})
\Gamma_{\tlQ}'(H_{\tlQ}(x)-T_{\tlQ},\hat r_{2}(H_{\tlQ}(x)))
\eta(x).
\end{displaymath}

Soit $\tlQ \in \fc^{\tlG}(B)$. En faisant une analyse analogue à celle de la preuve de la proposition 3.7 de \cite{leMoi2}, on voit que \eqref{eq:QsousG} est le produit de 
\begin{equation}\label{eq:purExpDesc}
\int_{\zgo_{\tlQ}}e^{\upla_{\tlQ}(H)}
\Gamma_{\tlQ}'(H-T_{\tlQ}, \hat r_{2}(H) )dH
\end{equation}
et
\begin{equation}\label{eq:IMQ}
I_{a}^{M_{\tlQ}, \eta',T}(f_{\tlQ}),
\end{equation}
où
\begin{itemize}
\item d'après le lemme \ref{lem:upGamma2}, l'intégrale \eqref{eq:purExpDesc} est absolument convergente et c'est un polynôme-exponentielle en $T_{\tlQ}$ sans terme purement polynomial pour $\tlQ\subsetneq \tlG$ ;
\item pour $\tmgo=\tmgo_{\tlQ}$ et  $\tngo=\tngo_{\tlQ}$, la fonction $f_{\tlQ} \in \Sc(\tmgo(\AAA))$ est définie pour tout $X\in \tmgo(\AAA)$ par la formule
  \begin{equation}
    \label{eq:fQ}
    f_{\tlQ}(X) = \int_{K}\int_{\tngo(\AAA)}f(k^{-1}\cdot(X + U))\eta(k)dUdk \, ;
  \end{equation}
\item la  distribution $I_{a}^{M_{\tlQ}, \eta',T}$ est définie  comme au paragraphe 4.2 de \cite{leMoi2} relativement à l'action de $M_{Q}$ sur $\tmgo$ et d'un certain caractère $\eta'$ de $\MM_{\tlQ}(\AAA)^1 \times G_{\tlQ}(\AAA)$ ; pour $\tlQ\subsetneq \tlG$, d'après \emph{loc. cit.} remarque 4.9,  le facteur \eqref{eq:IMQ} est un polynôme-exponentielle en la projection orthogonale de $ T $ sur $ \ago_{\tl}^{\tlQ} $. 
\end{itemize}

Il s'ensuit que,  pour $\tlQ\subsetneq \tlG$, le terme \eqref{eq:QsousG} est un polynôme-exponentielle sans terme purement polynomial. Comme pour $ \tlQ = \tlG $, le facteur  \eqref{eq:purExpDesc} vaut  $1$ et \eqref{eq:IMQ} est égal à  $ I_{a}^{\eta,T}(f) $, la proposition résulte du théorème \ref{thm:cv} et de la définition de $I_a^\eta$ (cf. §\ref{S:Ia}).
\end{preuve}
\end{paragr}

\subsection{Un troisième noyau auxiliaire}\label{par:noySemDes}

\begin{paragr}
Soit $f\in \Sc(\tggo(\AAA))$. On introduit l'analogue suivant de \eqref{eq:kT} où l'on remplace à la fois les fonctions $\hat\tau$ par leurs analogues $\hat\sigma$ et les noyaux $k_{\tlP,a}$ par leurs analogues $j_{\tlP,a}$ :
\begin{equation}
  \label{eq:tljT}
  \tlj_{a}^{T}(f,x) = 
\sum_{\tlP \in \fc^{\tlG}(B)}\varepsilon_{\tlP}^{\tlG}
\sum_{\delta \in P(\rmF) \backslash G(\rmF)}
\hat\sigma_{\tlP}(H_{\tlP}(\delta x) - T_{\tlP})j_{\tlP,a}(\delta x).
\end{equation}
\end{paragr}

\begin{paragr} Il nous faut aussi l'analogue de la proposition \ref{prop:newTroncKer1} suivant:

\begin{proposition}\label{prop:newTroncKer2}
Pour tout  $f\in \Cc(\tggo(\AAA))$, la proposition \ref{prop:newTroncKer1} vaut encore lorsqu'on remplace dans l'énoncé le noyau $\kappa_{a}^{T}(f,\cdot)$ par $\tlj_{a}^{T}(f,\cdot)$.
\end{proposition}

\begin{preuve}
Elle est semblable à la preuve de la proposition \ref{prop:newTroncKer1} sauf qu'ici on fait appel à l'égalité \eqref{eq:j=k} du théorème \ref{thm:cvj} et qu'on  compare les intégrales   $\int_{[G]} \tlj_{a}^{T}(x) \eta(x)dx$ et $\int_{[G]}j_{a}^{T}(x) \eta(x)dx$. La différence entre ces deux dernières intégrales fait apparaître une somme sur $\tlQ\in \fc^{\tlG}(B)$ de termes 
\begin{equation}\label{eq:noySemDes}
\int\limits_{Q(\rmF)\backslash  G(\AAA)}
\sum_{\fc^{\tlG}(B) \ni \tlP \subseteq \tlQ }
\varepsilon_{\tlP}^{\tlQ}
\sum_{\delta \in (P \cap M_{Q})(\rmF)\backslash M_{Q}(\rmF)}
\Psi^{T}_{\tlP,\tlQ, a}(\delta x)dx
\end{equation}
où l'on définit
\begin{displaymath}
\Psi^{T}_{\tlP,\tlQ, a}(x) = 
j_{\tlP,a}(x)
\hat \tau_{\tlP}^{\tlQ}(H_{\tlP}(x)- T_{\tlP})
\Gamma_{\tlQ}'(H_{\tlQ}(x)-T_{\tlQ}, \hat r_{2}(H_{\tlQ}(x)))
\eta(x).
\end{displaymath}

On remplace l'intégrale sur $Q(\rmF)\backslash G(\AAA)$ 
par l'intégrale sur 
\[
N_{Q}(\rmF) \backslash N_{Q}(\AAA) \times A_{\MM_{\tlQ}}^{ \infty} \times  
(M_{Q}(\rmF) \backslash \MM_{\tlQ}(\AAA)^{1}\times G_{\tlQ}(\AAA)) \times K
\]
ce qui donne $dx = e^{-2\rho_{Q}(H_{Q}(am))}dndadmdk$.

Pour $m \in (\MM_{\tlQ}(\AAA)^{1}\times G_{\tlQ}(\AAA))$, 
$a \in A_{\MM_{\tlQ}}^{ \infty}$, 
$k \in K$, $\delta \in M_{Q}(\rmF)$ et $\tlP \subseteq \tlQ$ on regarde:
\[
\int\limits_{[N_{Q}]} j_{\tlP, a}(na \delta mk)dn = 
\int\limits_{[N_{Q}]}
\sum_{X \in \tmgo_{\tlP,a}(\rmF)}
\sum_{\nu \in N_{P}(X_{s},\rmF) \backslash N_{P}(\rmF)} 
\int\limits_{\tngo_{\tlP}(X_{s},\AAA)}
f((\nu na \delta mk)^{-1} \cdot \iota_{X}(X^{+}, X^{-} + U))dUdn.
\]

On a $\tngo_{\tlP}(X_{s},\cdot) = \tngo_{\tlP}^{\tlQ}(X_{s}, \cdot) \oplus 
\tngo_{\tlQ}(X_{s}, \cdot)$ 
et $N_{P}(X_{s}, \cdot) = N_{P}^{Q}(X_{s}, \cdot) N_{Q}(X_{s}, \cdot)$, 
où $ \tngo_{\tlP}^{\tlQ}(X_{s}, \cdot) $ 
et $ N_{P}^{Q}(X_{s}, \cdot) $ sont les analogues des définitions 
du paragraphe \ref{S:centra} associés à $ \tmgo_{\tlQ} $.
On fixe alors $X \in \tmgo_{\tlP,a}(\rmF)$ et $\nu'' \in N_{P}^{Q}(X_{s},\rmF) \backslash N_{P}^{Q}(\rmF)$, 
on pose $$z = a \nu'' \delta mk$$ et on regarde:
\begin{multline*}
\int_{[N_{Q}]}
\sum_{\nu \in N_{Q}(X_{s},\rmF) \backslash N_{Q}(\rmF)} 
\int_{\tngo_{\tlQ}(X_{s},\AAA)}
\int_{\tngo_{\tlP}^{\tlQ}(X_{s},\AAA)}
f((\nu n z)^{-1} \cdot \iota_{X}(X^{+}, X^{-} + U_{\tlP}^{\tlQ} + U_{\tlQ}))dU_{\tlP}^{\tlQ} dU_{\tlQ} dn = \\
\int_{\tngo_{\tlP}^{\tlQ}(X_{s},\AAA)}
\int_{N_{Q}(X_{s}, \AAA) \backslash N_{Q}(\AAA)}
\int_{\tngo_{\tlQ}(X_{s},\AAA)}
f((n z)^{-1} \cdot \iota_{X}(X^{+}, X^{-} + U_{\tlP}^{\tlQ} + U_{\tlQ}))dU_{\tlP}^{\tlQ} dU_{\tlQ} dn. 
\end{multline*}
On fixe $U_{\tlP}^{\tlQ} \in \tngo_{\tlP}^{\tlQ}(X_{s},\AAA)$ et on applique 
le lemme \ref{lem:resommation} point \textit{2} pour la fonction 
\[
\tngo_{\tlQ}(\AAA) \ni U_{\tlQ} \mapsto 
f( z^{-1} \cdot ( \iota_{X}(X^{+}, X^{-} + U_{\tlP}^{\tlQ} + U_{\tlQ})).
\]
On trouve que l'équation ci-dessus devient:
\[
\int_{\tngo_{\tlP}^{\tlQ}(X_{s},\AAA)}
\int_{\tngo_{\tlQ}(\AAA)}
f(z^{-1} \cdot \iota_{X}(X^{+}, X^{-} + U_{\tlP}^{\tlQ}) + z^{-1}\cdot U_{\tlQ})dU_{\tlP}^{\tlQ} dU_{\tlQ} dn. 
\]

En faisant le changement de variable $(a \nu'' \delta m)^{-1} \cdot U_{\tlQ} \mapsto U_{\tlQ}$ 
on obtient que l'équation \eqref{eq:noySemDes} est le produit de deux facteurs : 
\begin{itemize}
\item le premier facteur est exactement \eqref{eq:purExpDesc} ;
\item le second est :
\begin{equation}\label{eq:jpm}
\int\limits_{
M_{Q}(\rmF) \backslash \MM_{\tlQ}(\AAA)^{1}\times G_{\tlQ}(\AAA)
}
\! \!
\sum_{\fc^{\tlG}(B) \ni \tlP \subseteq \tlQ }
\! \! \! \!
\varepsilon_{\tlP}^{\tlQ}  \! \! \! \!
\sum_{\delta \in (P \cap M_{Q})(\rmF)\backslash M_{Q}(\rmF)} \! \! \! \!
\hat \tau_{\tlP}^{\tlQ}(H_{\tlP}(\delta x) - T_{\tlP})
j_{\tlP \cap M_{\tlQ}, a}(f_{\tlQ}, \delta x)  \eta'(m)dm
\end{equation}
où $f_{\tlQ} \in C_{c}^{\infty}(\tmgo_{\tlQ}(\AAA))$ est définie par \eqref{eq:fQ}, le caractère $\eta'$ est comme dans la preuve de la proposition \ref{prop:newTroncKer1}. 
\end{itemize}
Le noyau $j_{ \tlP \cap M_{\tlQ}, a}(f_{\tlQ}, \cdot)$ est l'analogue de $j_{\tlP, a}(f,\cdot)$ pour le groupe $M_{Q}$ agissant sur l'espace $ \tmgo_{\tlQ} $. Par un théorème analogue au théorème \ref{thm:cvj}, l'intégrale \eqref{eq:jpm} est convergente et égale à l'intégrale $I_a^{M_{\tlQ},\eta',T}$. La conclusion est alors claire.
\end{preuve}
\end{paragr}

\subsection{Combinatoire de la descente} 

\begin{paragr}[Sous-groupes $\tlM_1$ et $\tlP_1$.] --- On reprend les notations et les hypothèses de la section \ref{ssec:descente}, en particulier celles du §\ref{S:situation-desc}. On dispose donc d'une décomposition $V=V^+\oplus (\oplus_{i\in I}V_i)$ Selon les notations de la section \ref{ssec:desc-combi}, en particulier cf. §\ref{S:M1}, par le choix d'une décomposition du $F_I$-espace vectoriel $V_i$, on fixe un sous-groupe de Levi $\tlM_1$ de $\tlG$. Quitte à faire agir $G(F)$, on peut et on va supposer que $\tlM_1$ contient $\tlT_{0}$ et qu'il est le facteur de Levi semi-standard d'un sous-groupe parabolique $\tlP_1$ contenant $B$.
\end{paragr}

\begin{paragr}[Sous-groupe $M_0$ et $P_0$.] --- D'après les notations et la construction du §\ref{S:QQ-}, on dispose d'une application $\tlQ\mapsto \tlQ^-$ de $\fc^{\tlG}(\tlM_1)$ dans $\fc^{\tlH^-}(\tlM_0)$. Le groupe $H^-$ est naturellement un sous-groupe de $\tlH^-$ et l'on pose $\tlP_0=\tlP_1^-$ et $P_0=\tlP_0\cap H^-$. Alors $P_0$ est un sous-groupe parabolique défini sur $F$  minimal de  $H^-$. Le groupe $M_0=\prod_{i \in I} M_i$ avec $M_i=\prod_l GL_{F_i}(D_{i,l})$ est un facteur de Levi de $P_0$.  Soit
$$\fc^{\tlG}_{P_0}(\tlM_1)$$
 l'ensemble des $\tlQ\in \fc^{\tlG}(\tlM_1)$ tel que $P_0\subset \tlQ^-$.
\end{paragr}

\begin{paragr}[Un élément semi-simple $X$.] ---  On dispose également d'un élément $a\in \Ac^{(r)}(F)$ image d'un élément $(a_+,(a_i)_{i\in I}\in \Ac_H'(F)$ par le morphisme \eqref{eq:iotaHquotient}. On fixe alors un élément $X+\in \thgo^{\rs}(F)$ d'image $a_+$ dans le quotient ainsi que $X_i=(\al_i,0,0)\in \thgo_i(F)$ avec $\al_i\in F_i$ d'image $a_i$. On obtient alors un élément $(X_+,(X_i)_{i\in I})\in \thgo'(F)$. Soit  $X\in \tggo^{(r)}(F)$ l'image de  cet élément par le morphisme \eqref{eq:iotaH}. Notons que $X$ est un élément semi-simple de $\tmgo_1(F)$. Soit 
$$\fc_X^{\tlG}(B)$$
 l'ensemble des couples $(\tlQ,\of)$ formés d'un élément $\tlQ\in \fc^{\tlG}(B)$ et d'une classe de $M_Q(F)$-conjugaison d'éléments semi-simples $\of\subset \tmgo_{\tlQ}(F)$ telle que $\of\subset G(F)\cdot X$
\end{paragr}

\begin{paragr} Les ensembles $ \fc_X^{\tlG}(B)$ et $\fc^{\tlG}_{P_0}(\tlM_1)$ sont naturellement en bijection comme le montre le lemme suivant

\begin{lemme}\label{lem:bij}
Pour tout $(\tlQ,\of)\in \fc_X^{\tlG}(B)$, il existe un unique élément $w\in W^Q\back W$ tel que $w^{-1}\tlQ w\in \fc^{\tlG}_{P_0}(\tlM_1)$ et $\of=(M_Q(F)w)\cdot X$.
L'application $(\tlQ,\of)\mapsto w^{-1}\tlQ w$ induit une bijection
$$ \fc_X^{\tlG}(B) \to \fc^{\tlG}_{P_0}(\tlM_1).
$$
\end{lemme}

\begin{preuve}
Soit $(\tlQ,\of)\in \fc_X^{\tlG}(B)$. Selon le §\ref{S:Levi-dec}, on a des décompositions $M_{\tlQ}=\MM_{\tlQ}\times \tlG_{\tlQ}$ et $\tlM_1=\MM_1\times \tlG_{\tlP_1}$.  Soit $Y\in \of$. Le tore déployé $A_{\MM_{\tlQ}}  $, qui est central dans $\MM_{\tlQ}$, est inclus dans le centralisateur $G_Y$ de $Y$.  Le centralisateur de $X$ s'identifie au groupe $\prod_{i\in I} GL_{F_i}(V_i)$ dont $A_{\MM_1}$ est un sous-tore déployé maximal. Soit $g\in G(F)$ tel que $Y=g\cdot X$. Les groupes $G_X$ et $G_Y$ sont donc conjugués par $g$. Il s'ensuit que $g A_{\MM_1} g^{-1}$ est un sous-tore déployé maximal de $G_Y$. Quitte à translater $g$ à gauche par un élément de $G_Y(F)$, on peut et on va supposer qu'on a $A_{\MM_{\tlQ}} \subset g A_{\MM_1} g^{-1}$. En prenant les centralisateurs de ces groupes dans $\tlG$, on a alors $g \tlM_1  g^{-1}\subset M_{\tlQ}$. 

Il est alors facile de voir que, quitte à translater $g$ à gauche par un élément de $M_Q(F)$ et changer $Y$ en conséquence, on a $g\in \Norm_{G(F)}(T_0)$. Soit $w$ la classe de $g$. On a  $w^{-1}\tlQ w\in \fc^{\tlG}(\tlM_1)$ et $\of=(M_Q(F)w)\cdot X$. On peut encore translater $w$ par un élément de  $\Norm_{G_X(F)}(T_0)$ ; on peut donc supposer supposer de plus que $P_0\subset (w^{-1}\tlQ w)^-$. On a ainsi établi l'existence de $w$.

Vérifions l'unicité.  La condition $\of=(M_Q(F)w)\cdot X$ implique que $w\in W^Q\back W$ est bien défini à une translation à droite  par un élément de  $\Norm_{G_X(F)}(T_0)$. Mais la condition   $P_0\subset (w^{-1}\tlQ w)^-$ assure l'unicité de $w\in W^Q\back W$.

 L'application est donc bien définie. Elle est injective car si $(\tlQ,\of)$ et $(\tlQ,\of')$ ont même image, les sous-groupes paraboliques $\tlQ$ et $\tlQ'$ sont standard et conjugués donc égaux et les éléments $w$ et $w'$ associés sont égaux modulo $W^Q$. Elle est enfin surjective car tout élément de  $\fc^{\tlG}(\tlM_1)$ est conjugué sous $W$ à un élément de  $\fc_X^{\tlG}(B)$.
\end{preuve}
\end{paragr}

\begin{paragr} Soit $T\in \ago_{\tl}^+$  Pour tout $g\in G(\AAA)$ et $\tlR\in \fc^{\tlH^-}(P_0)$, soit
  \begin{equation}
    \label{eq:Upsilon}
     \Upsilon^T_{\tlR}(g)= \sum_{\tlP \in \fc^{\tlG}_{\tlR}(\tlM_1)}   \eps_{\tilde{P}}^{\tilde{G}} \,  \hat{\sigma}_{\tilde{P}}(H_{\tilde{P}}(g)-T_{\tlP}),
   \end{equation}
   où  $\fc^{\tlG}_{\tlR}(\tlM_1)$ est défini au §\ref{S:combi}, et soit $\nc_{\tmgo_{\tlR}}$ le cône nilpotent de $\tmgo_{\tlR}$ c'est-à-dire le fermé de $\tmgo_{\tlR}$ obtenu par intersection avec la fibre en $0$ de l'application $\thgo^- \to \Ac_{H^-}$. Soit  $f\in \Sc(\tggo(\AAA))$.  On a introduit en \eqref{eq:tljT} un noyau  $\tlj_{a}^{T}(f,g)$ pour $g\in G(\AAA)$. Le lemme suivant en fournit une expression alternative.

\begin{lemme}\label{lem:formj} 
Pour tout $g\in G(\AAA)$,  on a 
  \begin{eqnarray*}
\tlj^T_a(f,g)=\sum_{\tlR \in \fc^{\tlH^-}(P_0)} \sum_{\delta\in R(F)\back G(F)} \Upsilon^T_{\tlR}(\delta g) \sum_{Y\in \nc_{\tmgo_{\tlR}}(F)} \, \int_{\tilde{\ngo}_{\tlR}(\AAA)} f((\delta g)^{-1}\cdot \iota_X(X^+,X^-+Y+U))dU.
\end{eqnarray*}
où $R=\tlR \cap H^-$. 
\end{lemme}

\begin{preuve} 
Pour tout $(\tlQ,\of)\in \fc_X^{\tlG}(B)$, introduisons
$$  j_{\tilde{Q},\of}(f,g)= \sum_{\{Y\in \tmgo_{\tlQ}(F), Y_s\in \of\}} \sum_{\nu \in N_{Q,Y_s}(F)\back N_Q(F)} \int_{\tilde{\ngo}_{\tlQ}(Y_s,\AAA)} f((\nu g)^{-1}\cdot \iota_Y(Y^+,Y^-+U))dU. 
$$
Par définition, on  a
$$
j_{\tilde{Q},a}(f,g)=\sum_{\{\of \mid(\tlQ,\of)\in \fc_X^{\tlG}(B)\}} j_{\tilde{Q},\of}(f,g)$$
et
$$\tlj^T_a(f,g)=\sum_{\tilde{P}\in \fc(B)} \eps_{\tilde{P}}^{\tilde{G}} \sum_{\delta\in P(F)\back G(F)} \hat{\sigma}_{\tilde{P}}(H_{\tilde{P}}(\delta g)-T)j_{\tilde{P},a}(f,\delta g).
$$

On utilise ensuite la bijection donnée par le lemme \ref{lem:bij}. Soit $\tlP\in \fc^{\tlG}_{P_0}(\tlM_1)$ correspondant à $(\tlQ,\of)$ et $w\in W^Q\back W$ l'élément du groupe de Weyl correspondant. On a alors 
$$  j_{\tilde{Q},\of}(f,g)= \sum_{Y\in \nc_{\tmgo_{\tlP^-}}(F)} \, \sum_{\nu \in P_{X}(F)\back P(F)} \int_{\tngo_{\tlP^-}(\AAA)} f((\nu w g)^{-1}\cdot \iota_X(X^+,X^-+Y+U))dU, 
$$
où $P_X$ est le centralisateur de $X$ dans $P$ et où l'on utilise l'égalité $\tngo_{\tlP}(X)= \tngo_{\tlP^-}$. Il vient alors
\begin{eqnarray*}
\tlj^T_a(f,g)=\sum_{\tlP \in \fc^{\tlG}_{P_0}(\tlM_1)} \eps_{\tilde{P}}^{\tilde{G}} \sum_{\delta\in P_X(F)\back G(F)} \hat{\sigma}_{\tilde{P}}(H_{\tilde{P}}(\delta g)-T_{\tlP})    \\
  \sum_{Y\in \nc_{\tmgo_{\tlP^-}}(F)} \, \int_{\tilde{\ngo}_{\tlP^-}(\AAA)} f((\delta g)^{-1}\cdot \iota_X(X^+,X^-+Y+U))dU.
\end{eqnarray*}
Pour conclure, il suffit de rassembler les $\tlP$ dans la somme suivant la valeur de $\tlP^-$ et de tenir compte de l'égalité $P_X=\tlP^-\cap H^-$.
\end{preuve}
\end{paragr}

\subsection{Démonstration du théorème \ref{thm:desc}}

\begin{paragr}
  On continue avec les notations de la section précédente. Soit $\tlH^+=GL_F(V^+\oplus F e_0)$. On fixe un sous-groupe de Borel $\tlB^+$ de $\tlH^+$ tel que $\tlB^+\cap GL_F(V^+)$ soit un sous-groupe de Borel de $H^+$. Le groupe $\tlB^+\times \tlP_0$ est alors un sous-groupe parabolique défini sur $F$ minimal de $\tlH=\tlH^+\times \tlH^-$. Soit $T'$ un point de $\ago_{\tlB^+}^{\tlH^+}\oplus \ago_{ \tlP_0   }^{\tlH^-}$ dans la chambre de Weyl positive associée à $\tlB^+\times \tlP_0$. Par action du groupe de Weyl, on déduit de $T'$ des points $T'_{\tlS}$ pour tout sous-groupe parabolique $\tlS\in \fc(\tlT^+\times \tlM_0)$ où $\tlT^+$ est un sous-tore maximal de $\tlB^+$. Dans la suite, on ne s'intéressera qu'aux sous-groupes paraboliques de $\tlH$ qui appartiennent à $\fc^{\tlH}(\tlH^+\times \tlM_0)$. Cet ensemble est en bijection évidente avec $\fc^{\tlH^-}(\tlM_0)$. Rappelons qu'on a  $\zgo_{\tlM_{0}} =  \zgo_{\tlM_{1}} \subset  \ago_{\tlM_{1}}$.  Pour tout $ \tlR \in \fc^{\tlH^{-}}(\tlM_{0}) $, soit  
$$ \hat T'_{\tlR} $$
 la projection orthogonale de $ T'_{\tlH^+\times \tlR} $ sur $ \zgo_{\tlM_{0}} $. On a donc $ \hat T'_{\tlR} \in \zgo_{\tlR} \subset \ago_{\tlM_{1}}$. De plus,  si $ \tlR  \subset \tlS\subset \tlH^-$ la projection orthogonale de $ \hat T'_{\tlR} $ sur $ \zgo_{\tlS} $ est égale à $ \hat T'_{\tlS}  $.

\end{paragr}

\begin{paragr}   Soit $K_{\tlH}\subset \tlH(\AAA)$, resp. $K_{H}\subset H(\AAA)$, un sous-groupe compact maximal de $H(\AAA)$, resp. $\tlH(\AAA)$,  \og en bonne position \fg{} par rapport au sous-groupe de Levi   $\tlT^+ \times \tlM_0$, resp. $(\tlT^+ \times \tlM_0)\cap H$. Pour tout $x\in H(\AAA)$ et tout sous-groupe parabolique  $R$ défini sur $F$ de $H$, soit $k_R(x)\in K_H$ un élément  tel que $x\in R(\AAA) k_H(x)$.
 
Soit $\tlR\in \fc^{\tlH^-}(\tlM_0)$ et $\tlP\in \fc_{\tlR}(\tlM_1)$. On a $\zgo_{\tlH^+}=0$ et donc on a 
$$\zgo_{\tlH^+\times \tlR}=\zgo_{\tlH^+}\oplus \zgo_{\tlR}=\zgo_{\tlR}.$$

Le choix de $K_{\tlH}$ fait qu'on dispose d'une application de  Harish-Chandra 
$$ H_{\tlH^+\times \tlR} : \tlH(\AAA) \rightarrow \ago_{ \tlH^+\times \tlR}.$$
  Soit $ \hat H_{\tlR} $ l'application composée de  $H_{\tlH^+\times \tlR}$ avec la projection orthogonale $ \ago_{\tlH^+\times \tlR} \rightarrow \zgo_{\tlR}$.

Soit $x\in H(\AAA)$ et $y\in G(\AAA)$. Les projections orthogonales sur $\zgo_{\tlP}$ de   $\hat H_{\tlR}(x) $  et 
$$H_{\tlP}( xy)-H_{\tlP}(k_R(x)y)
$$
sont donc égales. Il s'ensuit qu'on a 
\begin{equation}
  \label{eq:desc-sig}
  \hat \sigma_{\tlP}(H_{\tlP}( xy) - T_{\tlP}) = 
\hat \sigma_{\tlP}(\hat H_{\tlR}(x) - \hat T'_{\tlR} - Y^{T,T'}_{\tlP}(x,y)).
\end{equation}
où l'on introduit
$$Y^{T,T'}_{\tlP}(x,y)= -H_{\tlP}(k_{R}(x)y) + T_{\tlP} - \hat T'_{\tlR} \in \ago_{\tlP}.
$$
On obtient en particulier une famille 
$$\yc_{\tlR}^{T, T'}(x,y) = (Y_{\tlP}^{T, T'}(x,y))_{\tlP \in \fc^{0}_{\tlR}(\tlM_{1})}$$
 qui est orthogonale-positive. Soit $\tlS\in \fc^{\tlH^-}(\tlR)$ et $\tlQ\in  \fc_{\tlS}^{\tlG}(\tlM_1)$ tel que $\tlP\subset \tlQ$. La projection orthogonale de $Y^{T,T'}_{\tlP}(x,y)$ sur $\ago_{\tlQ}$ est égale à $Y^{T,T'}_{\tlQ}(x,y)$.

\begin{lemme}
  \label{lem:Upsilon} Avec les notations ci-dessus, on a l'égalité suivante où le membre de droite est défini en \eqref{eq:Upsilon}
\begin{equation*}
  \label{eq:Upsilon2}
   \Upsilon^T_{\tlR}(xy) = \sum_{\tlS \in \fc^{\tlH^{-}}(\tlR)}
\varepsilon_{\tlR}^{\tlS}
\hat \sigma_{\tlR}^{\tlS}(\hat H_{\tlR}( x) - \hat T'_{\tlR}  )
\mathrm{B}_{\tlS}^{\tlG}(\hat  H_{\tlS}( x) - \hat T'_{\tlS}, \yc_{\tlS}^{T,T'}(x,y)).
\end{equation*}
\end{lemme}

\begin{preuve}
  En utilisant \eqref{eq:desc-sig}, puis  \eqref{eq:42star}, on obtient
\begin{eqnarray*}
  \Upsilon^T_{\tlR}(xy)&=&\sum_{ \tlP \in \fc_{\tlR}(\tlM_{1})}
 \varepsilon_{\tlP}^{\tlG}
 \hat \sigma_{\tlP}(\hat H_{\tlR}( x) - \hat T'_{\tlR} - Y^{T, T'}_{\tlP}(x,y))\\
&=& \sum_{\tlS \in \fc^{\tlH^{-}}(\tlR)}
\varepsilon_{\tlR}^{\tlS}
\hat \sigma_{\tlR}^{\tlS}( \hat H_{\tlR}( x) - \hat T'_{\tlR}  )
\mathrm{B}_{\tlS}^{\tlG}(\hat  H_{\tlR}( x) - \hat T'_{\tlR}, \yc_{\tlS}^{T,T'}(x,y)).
\end{eqnarray*}
En utilisant le fait que $\mathrm{B}_{\tlS}^{\tlG}(  H, \yc_{\tlS}^{T,T'}(x,y))$ ne dépend que de la projection de $H$ sur $\zgo_{\tlS}$, on aboutit au lemme.
\end{preuve}

\end{paragr}

\begin{paragr}[Cas de fonctions $\al$ particulières.] --- \label{S:al-compact} On reprend les hypothèses et les notations du § \ref{S:Omega}. Soit 
  \begin{equation}
    \label{eq:al-compact}
    \al= \be_S \otimes \phi_S
  \end{equation}
avec $\be_S \in  \Cc(G(\AAA_S))$ et $\phi_S\in \Cc(\Om_H)$. On définit alors les fonctions $f_\al$ et $f_\al^{H,\eta}$ (cf. formules intégrales \eqref{eq:falpha} et \eqref{eq:fHalpha}).  Soit $f=f_\al \otimes \mathbf{1}_{\tggo(\oc^S)}$ et $f^H\in \Cc(\thgo(\AAA))$ donnée par
$$f^H= \frac{\vol(G(\oc^S))}{\vol(H(\oc^S))}  f_\al^{H,\eta}\otimes \mathbf{1}_{\thgo(\oc^S)}.
$$
Soit  $v\notin S$ une place de $F$.  Il résulte du lemme  \ref{lem:integrite} que,  pour tout $y\in G(F_v)$ et  tout $Z\in \thgo_{a_H}(F_v)$, on a 
\[ 
\mathbf{1}_{\tggo(\oc_v)}(y^{-1}\cdot \iota_H(Z)) = 
\vol (H(\oc_v))^{-1}\int_{H(F_v)} 
\mathbf{1}_{G(\oc_v)}(y^{-1}h)
\mathbf{1}_{\thgo(\oc_v)}( h^{-1} \cdot Z) dh.
\]

Soit $\beta=\vol(H(\oc^S))^{-1} \al_S\otimes \mathbf{1}_{G(\oc^S)}\in \Cc(G(\AAA))$ et $\phi= \phi_S\otimes \mathbf{1}_{\thgo(\oc^S)}\in \Cc(\thgo(\AAA))$. On a alors  pour tous $y\in G(\AAA)$ et  $Z\in \thgo_{a_H}(\AAA)$
\begin{equation}
  \label{eq:f-beta}
  f(y^{-1}\cdot \iota_H(Z)) = 
\int_{H(\AAA)} 
\be(y^{-1}h)
\phi( h^{-1} \cdot Z) dh.
\end{equation}

\begin{proposition}\label{prop:dvpT-T'}
Pour tout $T$ dans un certain translaté de $\ago_{\tl}^+$, l'intégrale 
  $$
\int_{[G]} \tlj_a^T(f,g) \, \eta(g) \, dg 
$$
est égale à la somme sur $\tlR_1 \in \fc^{\tlH^{-}}(P_{0})$ de 
\begin{multline}\label{eq:descFin3}
  \int\limits_{G(\AAA)}
\int\limits_{R_1(\rmF) \backslash H(\AAA)} 
\mathrm{B}_{\tlR_1}^{\tlG}(\hat H_{\tlR_1}(x) - \hat T'_{\tlR_1}, \yc_{\tlR_1}^{T,T'}(x,y))
\sum_{\tlR \in \fc^{\tlR_1}(P_{0}) }\varepsilon_{\tlR}^{\tlR_1}
\sum_{\delta \in R(\rmF) \backslash R_1(\rmF)}
\hat \sigma_{\tlR}^{\tlR_1}(\hat H_{\tlR}(\delta x) - \hat T_{\tlR}') \\
\sum_{Y\in \nc_{\tmgo_{\tlR}}(F)} \, \int_{\tilde{\ngo}_{\tlR}(\AAA)}
 \beta(y^{-1}) \phi((\delta x)^{-1} \cdot \iota_H(X^+, X^-+Y+U))dU$$
\eta(xy)dxdy.
\end{multline}

\end{proposition}

\begin{preuve} À l'aide du lemme \ref{lem:formj}, l'intégrale
$$
\int_{[G]} \tlj_a^T(f,g) \, \eta(g) \, dg 
$$
s'écrit encore
\begin{multline}\label{eq:descFin}
\int\limits_{H(\AAA)\back G(\AAA)}\int\limits_{H(F)\back H(\AAA)}
\sum_{\tlR \in \fc^{\tlH^-}(P_0)} \sum_{\delta\in R(F)\back H(F)} \Upsilon^T_{\tlR}(\delta xy) \\\sum_{Y\in \nc_{\tmgo_{\tlR}}(F)} \, \int_{\tilde{\ngo}_{\tlR}(\AAA)} f((\delta xy)^{-1}\cdot \iota_H(X^+,X^-+Y+U))dU\, 
\eta(xy)dxdy.
\end{multline}
En utilisant l'expression \eqref{eq:f-beta}, on voit qu'il suffit de démontrer  la convergence suivante :
\begin{multline}\label{eq:descFin-valabs}
\int\limits_{G(\AAA)}\int\limits_{H(F)\back H(\AAA)}
|\sum_{\tlR \in \fc^{\tlH^-}(P_0)} \sum_{\delta\in R(F)\back H(F)} \Upsilon^T_{\tlR}(\delta xy) \\\sum_{Y\in \nc_{\tmgo_{\tlR}}(F)} \, \int_{\tilde{\ngo}_{\tlR}(\AAA)} \beta(y^{-1}) \phi(h^{-1}(\delta x)^{-1}\cdot (X^+,X^-+Y+U))dU| \,dxdy<\infty.
\end{multline}
Mais celle-ci est une conséquence des propositions \ref{prop:nonH} et  \ref{prop:nonH2}.
\end{preuve}
\end{paragr}

\begin{paragr}[Terme \eqref{eq:descFin3} pour $\tlR_1  \not=\tlH^{-}$.] --- Il ne s'agit pas de l'expliciter, la proposition suivante va suffire.

  \begin{proposition}
   \label{prop:nonH}
Soit $ \tlR_1 \in \fc^{\tlH^{-}}(P_{0})$ tel que $\tlR_1 \not=\tlH^{-}$.  L'expression \eqref{eq:descFin3} pour un tel  $ \tlR_1 $ converge absolument et c'est la restriction d'un polynôme-exponentielle en $T$ et $T'$ dont la partie purement polynomiale est nulle.
  \end{proposition}

  \begin{preuve}  En utilisant la décomposition d'Iwasawa, on décompose l'intégrale $ R_1(F) \backslash H(\AAA) $ en une intégrale sur les variables $(n, H, m, k)$ parcourant 
\[
[N_{R_1}] \times \zgo_{\tlR_1} \times (H^+(\AAA)\times  [G_{\tlR_1}] \times \MM_{\tlR_1}(F) \backslash  \MM_{\tlR_1}(\AAA)^{1}) 
\times K_{H},
\]
où l'on décompose $ M_{\tlR_1} = \MM_{\tlR_1} \times G_{\tlR_1} $ comme dans le paragraphe \ref{S:Levi}. Au niveau des mesures, cette décomposition donne $dx = e^{-2\rho_{R_1}(H+H_{R}(m))}dndH dmdk$. Observons que dans \eqref{eq:descFin3}, le facteur $\hat \sigma_{\tlR}^{\tlR_1}(\hat H_{\tlR}(\delta x) - \hat T_{\tlR}') $ est égal à $\hat \sigma_{\tlH^+\times \tlR}^{\tlH^+\times\tlR_1}( H_{\tlH^+\times \tlR}(\delta x) - T_{\tlH^+\times \tlR}') $. Après quelques changements de variables, on voit que l'intégrale \eqref{eq:descFin3} est égale à 
\begin{multline}
  \label{eq:newj}
\int_{G(\AAA)}  \int_{ [G_{\tlR_1}] \times [\MM_{\tlR_1}]^1} \kappa^{T'}_{a_H}(\Phi_y^{T,T'},h) \eta(h)\, dh \, \beta(y)\eta(y)dy
\end{multline}
où
\begin{itemize}
\item la fonction $\Phi_y^{T,T'}$ appartient à  $\Cc(\thgo^+(\AAA)\times \tmgo_{\tlR_1}(\AAA))$ et est donnée par la formule intégrale
\[
\Phi_{y}^{T,T'}(Y^+,Y^-) = \int_{K_H}\int_{\tngo_{\tlR_1}(\AAA)} \phi(k^{-1}\cdot  (Y^{+}, Y + U ) ) u^{T,T'}(k,y)\eta(k) \, dkdU
\]
pour $Y^+\in \thgo^+(\AAA)$ et $Y^- \in \tmgo_{\tlR_1}(\AAA)$ ; la fonction $u^{T,T'}(k,y)$ étant donnée par la formule intégrale
\begin{equation}
  \label{eq:poids}
u^{T,T'}(k,y)=  \int_{\zgo_{\tlR_1}} e^{\upla_{\tlR_1}(H)}\mathrm{B}_{\tlR_1}^{\tlG}(H - \hat{T}'_{\tlR_1}, \yc_{\tlR_1}^{T,T'}(k,y))\, dH
\end{equation}
\item $\kappa^{T'}_{a_H}$ est défini comme le noyau $\kappa_{a}^{T}$ au § \ref{par:noyTronq} mais  dans le contexte de $H^+\times M_{R_1}$ agissant sur $\thgo^+\times \tmgo_{\tlR_1}$. 
\end{itemize}
Le lemme suivant permet de séparer dans $u^{T,T'}(k,y)$ la dépendance en $(k,y)$ d'une part et en $(T,T')$ d'autre part.

\begin{lemme}\label{lem:lastExp}
  \begin{enumerate}
  \item L'intégrale \eqref{eq:poids} qui définit $u^{T,T'}(k,y)$ converge absolument. 
  \item De plus, comme fonction de $T,T',k$ et $y$, elle appartient à l'espace engendrée par les fonctions

$$p_{\tlQ, \tlP}(-H_{\tlQ}(ky)) p_{\tlQ, \tlP}'(T_{\tlQ} - \hat T'_{\tlR_1})  e^{(\upla_{\tlQ} - \upla_{\tlP})(\hat T'_{\tlR_1})} 
e^{\upla_{\tlP}(T_{\tlQ})} e^{ - \upla_{\tlP}(H_{\tlQ}(ky))}
$$
où
\begin{enumerate}
\item $\tlQ \in \fc_{\tlR_1}^{0}(\tlM_1)$ et $\tlQ\subset \tlP$ ;
\item $p_{\tlQ, \tlP}$ et $p_{\tlQ, \tlP}'$ sont des polynômes ;
\item les expressions $(\upla_{\tlQ} - \upla_{\tlP})(\hat T'_{\tlR_1})$ et $\upla_{\tlP}(T_{\tlQ}) $, vues respectivement comme forme linéaire en $T'$ et $T$, ne sont pas génériquement nulles simultanément.
\end{enumerate}
\end{enumerate}
\end{lemme}

\begin{preuve}
  En vertu du lemme \ref{lem:gamSumGam0}, l'intégrale en question égale
\[
\sum_{\tlQ \in \fc_{\tlR_1}^{0}(\tlM_1)} 
\int_{\zgo_{\tlQ}}
e^{\upla_{\tlQ}(H)}
\mathrm{B}_{\tlQ}(H - \hat T'_{\tlR_1}, -H_{\tlQ}(ky)  + T_{\tlQ} - \hat T'_{\tlR_1})dH
\]
Soit $\tlQ \in \fc_{\tlR_1}^{0}(\tlM_1)$. On a alors $\zgo_{\tlQ}=\zgo_{\tlR_1}$. Comme   $\hat T'_{\tlR_1}\in \zgo_{\tlR_1}$, on peut faire un changement de variables ; dans la somme ci-dessus, le terme indexé par $\tlQ$ est donc égal à 
$$
e^{\upla_{\tlQ}(\hat T'_{\tlR_1})} \int_{\zgo_{\tlQ}} e^{\upla_{\tlQ}(H)}
\mathrm{B}_{\tlQ}(H , -H_{\tlQ}(ky)  + T_{\tlQ} - \hat T'_{\tlR_1})dH.
$$
D'après le lemme \ref{lem:upGamma} assertion 3, cette dernière intégrale est égale

$$
\sum_{\tlQ\subset \tlP}p_{\tlQ, \tlP}(-H_{\tlQ}(ky)  + T_{\tlQ} - \hat T'_{\tlR_1})  e^{(\upla_{\tlQ} - \upla_{\tlP})(\hat T'_{\tlR_1})} 
e^{\upla_{\tlP}(T_{\tlQ})} e^{ - \upla_{\tlP}(H_{\tlQ}(ky))}.
$$
On en déduit les conclusions cherchées hormis 2.c. Vérifions donc 2.c. La  forme linéaire 
$$T\mapsto  \upla_{\tlP}(T_{\tlQ})=  \upla_{\tlP}(T_{\tlP})
$$ 
se factorise par $\ago_{\tlP}$ mais n'est pas nulle sauf si $\tlP=\tlG$  (cf.  lemme 4.2 \textit{ii)} de \cite{leMoi2}). Supposons $\tlP=\tlG$. On a alors
$$(\upla_{\tlQ} - \upla_{\tlP})(\hat T'_{\tlR_1})= (\upla_{\tlQ} )(\hat T'_{\tlR_1}).
$$
Mais $\tlR_1\not=\tlH^-$. Il s'ensuit que $\tlQ\not=\tlG$. Donc la forme linéaire  $\upla_{\tlQ}$ est non nulle sur $\zgo_{\tlQ}=\zgo_{\tlR_1}$.  Mais la projection $T'\to  T'_{\tlR_1}$  est une surjection sur $\zgo_{\tlR_1}$. Donc l'expression ci-dessus vue comme forme linéaire en $T'$ n'est pas identiquement nulle. 
\end{preuve}

Soit $\tlQ \in \fc_{\tlR_1}^{0}(\tlM_1)$ et $\tlQ\subset \tlP$. Soit   $\Phi_y$ ue fonction dépendant continûment de $y$. 
Il suffit pour conclure de montrer le produit de
\begin{equation}
  \label{eq:exponentielle}
  e^{(\upla_{\tlQ} - \upla_{\tlP})(\hat T'_{\tlR_1})} \cdot e^{\upla_{\tlP}(T_{\tlQ})}
\end{equation}
avec l'expression
\begin{equation}
  \label{eq:kappaT'}
  \int_{G(\AAA)}  \int_{ [G_{\tlR_1}] \times [\MM_{\tlR_1}]^1} \kappa^{T'}_{a_H}(\Phi_y,h) \eta(h)\, dh \, \beta(y)\eta(y)dy
\end{equation}
est  un polynôme-exponentielle en $T'$ et $T$ dont les exposants sont non nuls. C'est déjà le cas pour le facteur \eqref{eq:exponentielle}. Observons que comme fonction de $T'$ l'expression    \eqref{eq:exponentielle} ne dépend que de la projection orthogonale de $T'$ sur le facteur $a_{\tlR_1}$. L'expression \eqref{eq:kappaT'} est quant à elle un polynôme-exponentielle en la projection de $T'$ sur le facteur $\ago_{\tlP_0}^{\tlR_1}$.  Cela se démontre comme dans \cite{leMoi2} section 4.3.
\end{preuve}
\end{paragr}

\begin{paragr}[Terme \eqref{eq:descFin3} pour $\tlR_1 =\tlH^{-}$.] --- La proposition suivante en donne la partie polynomiale.

  \begin{proposition}
       \label{prop:nonH2}
  L'expression  \eqref{eq:descFin3} pour $ \tlR_1 =\tlH^{-}$ converge absolument et ne dépend pas de $T$. Comme fonction de $T'$ c'est la restriction d'un polynôme-exponentielle en $T'$ dont la partie purement polynomiale est égale à 
$$
\frac{\vol(G(\oc^S))}{\vol(H(\oc^S))}\cdot I^{H,\eta}_{a_H}(f^{H,\eta}_\al\otimes \mathbf{1}_{\thgo(\oc^S)}).
$$

 \end{proposition}

 \begin{preuve}
   Pour $\tlR_1=\tlH^-$, le terme   \eqref{eq:descFin3} se réduit à 
$$
\int_{G(\AAA)}  \beta(y^{-1})\eta(y) \,dy \cdot \int_{[H]} \kappa^{H,T'}_{a_H}(\phi,h) \eta(h) \, dh.
$$
Il ne dépend donc pas de $T$. Comme $\eta$ est trivial sur $G(\oc^S)$ (cf. §\ref{S:hypsurS}), on a 
$$\int_{G(\AAA)}  \beta(y^{-1})\eta(y) \,dy = \frac{\vol(G(\oc^S))}{\vol(H(\oc^S))}\cdot \int_{G(\AAA_S)}  \al_S(y)\eta(y)^{-1} \,dy_S.
$$ 
D'après la proposition \ref{prop:newTroncKer1} (ou plutôt sa variante à $H$), l'intégrale 
 $$\int_{[H]} \kappa^{H,T'}_{a_H}(\phi,h) \eta(h) \, dh$$
est la restriction d'un polynôme-exponentielle en $T'$ dont le terme purement polynomial est constant égal à $I^{H,\eta}_{a_H}(\phi)$. La proposition s'ensuit si l'on observe que d'après les définitions on a 
$$(\int_{G(\AAA_S)}  \al_S(y)\eta(y)^{-1} \,dy_S)\phi=f^{H,\eta}_\al\otimes \mathbf{1}_{\thgo(\oc^S)}.
$$
 \end{preuve}
\end{paragr}

\begin{paragr}[Démonstration du théorème \ref{thm:desc}.] ---\label{S:demodethmdesc}
On commence par le démontrer pour $\al$ comme au §\ref{S:al-compact}. D'après la proposition \ref{prop:newTroncKer2}, l'expression $I^{G,\eta}_a(f_\al\otimes \mathbf{1}_{\tggo(\oc^S)})$ est la partie polynomiale (en fait constante) de l'expression 
$$
\int_{[G]} \tlj_a^T(f,g) \, \eta(g) \, dg. 
$$
En utilisant les propositions \ref{prop:dvpT-T'}, \ref{prop:nonH} et \ref{prop:nonH2}, on voit que  $I^{G,\eta}_a(f_\al\otimes \mathbf{1}_{\tggo(\oc^S)})$ est la partie polynomiale du terme  \eqref{eq:descFin3} pour $ \tlR_1 =\tlH^{-}$ c'est-à-dire 
$$
\frac{\vol(G(\oc^S))}{\vol(H(\oc^S))}\cdot I^{H,\eta}_{a_H}(f^{H,\eta}_\al\otimes \mathbf{1}_{\thgo(\oc^S)}).
$$

Cela prouve donc le théorème \ref{thm:desc} pour un tel $\al$. En général, on peut supposer $\al=\al_{S_\infty}\otimes \al_{S^\infty}$ où  $S=S_\infty \cup S^{\infty}$ est la partition de $S$ en places archimédienne et non-archimédiennes. On peut également supposer que  $\al_{S^\infty}$ est aussi du type considéré au  §\ref{S:al-compact} (c'est-à-dire un produit tensoriel de fonctions à support compact). En général,  $\al_{S_\infty}$ est au moins une limite de fonctions du type considéré au  §\ref{S:al-compact}. On conclut alors en invoquant les deux faits suivants :
\begin{itemize}
\item la continuité des  distributions  $I_{a}^{G, \eta}$ et $I_{a_{H}}^{H, \eta}$ (cf. assertion 4 du théorème \ref{thm:I}) ;
\item les applications 
$ \al_{S_\infty} \mapsto f_{ \al_{S_\infty}\otimes \al_{S^\infty}   } $ et $ \al \mapsto f_{\al_{S_\infty}\otimes \al_{S^\infty}}^{H, \eta} $
sont continues. 
\end{itemize}
\end{paragr}

\part{Le cas infinitésimal hermitien}\label{partie:2}

\section{Préliminaires algébriques} \label{sec:PrelimU}

\subsection{Stratification} \label{ssec:stratU}

\begin{paragr}\label{S:cadreu}
Soit $E/F$ une extension quadratique de  corps de caractéristique $0$. Soit $\sigma$ le générateur du groupe de Galois $\Gal(E/F)$. Soit $(V,\Phi)$ un couple formé d'un $V$ un $E$-espace vectoriel de dimension $n$ et d'une forme $\sigma$-hermitienne  $\Phi$  non dégénérée (notre convention sera que $\Phi$ est $\sigma$-linéaire en la première variable et linéaire en la seconde).  Soit $U(V,\Phi)$ le groupe unitaire, défini sur $F$, des automorphismes de $V$ qui préservent la forme $\Phi$. Pour alléger les notations, lorsque $(V,\Phi)$ ou $V$ est sous-entendu, on note 
$$U=U_\Phi=U(V,\Phi).$$
Soit $\tugo(V,\Phi)$ le $F$-espace vectoriel formé des couples $(A,b)$ où $b\in V$ et $A$ est  un endomorphisme de $V$  auto-adjoint pour la forme $\Phi$. On pose alors
$$\tugo=\tugo_{\Phi}=\tugo(V,\Phi).$$
Cet espace est muni d'une action linéaire à gauche et définie sur $F$ du groupe $U$ donnée par 
$$g\cdot (A,b)= (gAg^{-1},gb)
$$
pour $g\in U$ et $(A,b)\in \tugo$.
\end{paragr}

\begin{paragr}\label{S:quotientu}
On suit les notations du § \ref{S:quotient-categor}.  Soit
$$\Ac(V,\Phi)=\tugo //U$$
 le quotient catégorique et  
\begin{equation}
  \label{eq:aU}
  a:\tugo\to \Ac
\end{equation}
le morphisme canonique. Pour $a\in \Ac(V,\Phi)$, on rappelle que  $\tugo_a$ est alors la fibre en $a$  du morphisme \eqref{eq:aU}.

Pour alléger, on sous-entend parfois les données $(V,\Phi)$ et on note
$$
 \Ac=\Ac_\Phi= \Ac(V,\Phi).
$$

Soit $X=(A,b)\in \tugo$. Soit $t^n+a_1t^{n-1}+\ldots+a_n\in F[t]$ le polynôme caractéristique de $A$. Pour $1\leq i\leq n$, les fonctions $a_i$ et $b_i=\Phi(b,A^{i-1}b)$ sont des générateurs de l'algèbre $F[\tugo]^{U}$ algébriquement indépendants. Par le biais de ces fonctions, on identifie $\Ac$ à l'espace affine $\AAA_{2n}$ de dimension $2n$. De cette identification, on déduit des constructions du §\ref{S:Areg} des parties localement fermées $\Ac^{(r)}$, $\Ac^{\geq r}$ etc. construites à l'aide du morphisme $d_r$. Pour $X=(A,b)\in \tugo$, on a explicitement $d_r(X)=\det(\Delta_r(X))$ où
$$\Delta_r(X)=(\Phi(A^i b,A^{j}b))_{0\leq i,j\leq r-1}.
$$
On définit de manière évidente les parties localement fermées $\tugo^{(r)}$,  $\tugo^{(\geq r)}$ etc. L''ouvert dense $\tugo^{(\rs)}=\tugo^{(n)}$ est formé des éléments semi-simples et réguliers. On a le lemme suivant (évident vu le lemme \ref{lem:A0}) :

\begin{lemme}
  \label{lem:A0u}
Le fermé formé des $(A,b)\in \tugo$ tels que $\Phi(b,A^ib)=0$ pour tout $i$ est exactement $\tugo^{(0)}$.
\end{lemme}
 
\end{paragr}

\subsection{Sommes directes}
\begin{paragr}[Somme directe canonique $V^+\oplus V^-$ associé à $X$.]\label{S:Xru} --- 
Soit $0\leq r\leq n$ et $X=(A,b)\in \tugo^{(r)}$. Soit  $V^+$  le sous-$E$-espace de $V$ engendré par $A^ib$ pour $0\leq i\leq r-1$. Ce sous-espace $V^+$ est non-dégénéré pour la forme $\Phi$. Soit $V^-$ son orthogonal. On a donc 
 \begin{equation}
    \label{eq:direct2u}
    V=V^+\oplus V^-.
  \end{equation}
C'est la \emph{somme directe canonique} associée à $X$. 
\end{paragr}

\begin{paragr}[Morphisme associé à une somme directe.] --- \label{S:somme-directeu} Soit $V^+$ un sous-$E$-espace \emph{non nul} de $V$  non-dégénéré pour la forme $\Phi$. Soit $r=\dim(V^+)$ et $\tsgo(V^+)\subset \tugo$ la partie localement fermée formée des   $X=(A,b)\in \tugo$ tels que la famille $(A^ib)_{0\leq i\leq r-1}$ soit une base de $V^+$. Soit $V^-$ l'othogonal de $V^+$. On a donc
\begin{equation}
    \label{eq:direct3u}
    (V,\Phi)=(V^+,\Phi^+)\oplus (V^-,\Phi^-)
  \end{equation}
et les formes $\Phi^\pm$ sont non dégénérées.
Soit $\tugo_{\Phi^\pm}$ l'espace attaché à $ (V^\pm,\Phi^\pm)$.
Soit $(A,b)\in \tsgo(V^+)$ et $b'\in V^+$ défini par la condition
\begin{equation}
  \label{eq:bprimeu}
  \Phi(b',A^ib)=\left\lbrace
  \begin{array}{l}
    0 \text{  si  } 0\leq i<r-1 \\
1 \text{  si  } i=r-1.
  \end{array}\right.
\end{equation}

Suivant la décomposition \eqref{eq:bprimeu}, on a des projections $p_\pm: V\to V^\pm$ et des injections $i_\pm : V^\pm \to V$. On pose alors $A^\pm= p_\pm \circ A \circ i_\pm$.  On a $(A^+,b) \in \tugo_+^{(r)}$. Soit $L=p_-\circ A \circ i_+$. Pour $0\leq i \leq r-2$, on a $L A^ib=0$. Par conséquent, on a  $L=\Phi(b',i_+(\cdot)) b^- $ pour un certain $b^-\in V^-$. On vérifie qu'on a alors $p_+\circ A \circ i_-= \Phi(b^-,i_-(\cdot)) b'$.
  On définit alors un isomorphisme 
  \begin{equation}
    \label{eq:iota-su}
    \iota=\iota_{V^+}: \tugo^{\rs}_{\Phi^+} \times \tugo_{\Phi^-} \to \tsgo(V^+)
  \end{equation}
  
par
$$
\iota((A,b),(A^-,b^-))= (  \begin{pmatrix}
  A & \Phi(b^-,i_-(\cdot)) b' \\ \Phi(b',i_+(\cdot)) b^-   & A^-
\end{pmatrix}, b)
$$
où $b'\in V^+$ est le vecteur associé à $(A,b)$ par la condition \eqref{eq:bprimeu}.

Par commodité, pour $V^+=(0)$, on définit $\iota_{(0)}$ comme l'identité et $\tsgo((0))=\tugo$. Pour tout  $X\in \tugo^{(r)}$, on a une décomposition associée $V=V^+\oplus V^-$, cf. \eqref{eq:direct2u}. On pose alors
$$\iota_X=\iota_{V^+}.$$

Soit $U_{\Phi^\pm}$ le groupe unitaire associé à $(V^\pm,\Phi^\pm)$. Le morphisme $\iota$ est $U_{\Phi^+}\times U_{\Phi^-}$-équivariant si l'on identifie naturellement $U_{\Phi^+}\times U_{\Phi^-}$ à un sous-groupe de $U$. Le morphisme induit au niveau des quotients catégoriques n'est autre  que le morphisme \eqref{eq:iota} (via les identifications du §\ref{S:quotientu}).

\begin{lemme}\label{lem:entierku} Soit $(X,Y)\in \tugo_{\Phi^+}^{(r)} \times \tugo_{\Phi^-}$ et $k\in \NN$. On a $Y\in  \tugo_{\Phi^-}^{(k)}$ si et seulement si $\iota(X,Y)\in  \tugo^{(r+k)}$.
\end{lemme}

\begin{preuve}
Le lemme est une conséquence immédiate du lemme   \ref{lem:entierk}.
\end{preuve}

\end{paragr}

\subsection{Décomposition de Jordan}\label{ssec:Jordanu}

\begin{paragr}
  Cette section est parallèle à la section \ref{ssec:Jordan}. On construit une décomposition de Jordan pour tout $X\in \tugo$. Un tel $X$  est dit  nilpotent si son invariant $a$ est nul. On a le lemme suivant. 

  \begin{lemme}\label{lem:ssu}(Rallis-Schiffmann, cf. \cite{rallSchiff} théorème 17.2)
    Un élément $X=(A,b)\in \tugo$ est semi-simple si et seulement si dans la somme directe canonique  \eqref{eq:direct2u} de $X$ les espaces $V^{\pm}$ sont stables par $A$ et si l'endomorphisme de $V^-$ induit par $A$ est semi-simple au sens usuel.
\end{lemme}
\end{paragr}

\begin{paragr}[Décomposition de Jordan.] --- \label{S:Jor1u} Soit $X\in \tugo$. On a $X \in \tugo^{\rs}$ si et seulement si $X$ est semi-simple et régulier   (cf.,  à la terminologie près, \cite{rallSchiff} théorème 17.1). On pose dans ce cas $X_s=X$ et $X_n=0$.  Supposons à l'opposé que $X=(A,b)\in \tgl(V)^{(0)}$.  Soit $A=A_s+A_n$ la décomposition de Jordan usuelle de $A$. On  pose alors 
$$X_s=(A_s,0)$$
et
$$X_n=(A_n,b).
$$

Soit $1\leq r <n$. Supposons $X=(A,b)\in \tugo^{(r)}$. Avec les notations de §\ref{S:somme-directeu}, on a $X=\iota_X(X^+,X^-)$ avec $X^+\in  \tugo_{\Phi^+}^{(r)}$ et $X^-\in  \tugo_{\Phi^-}$. D'après le lemme \ref{lem:entierku}, on a $X^-\in  \tugo_{\Phi^-}^{(0)}$. D'après ce qui précède, on a $X^+_s=X^+$ et $X^-_s$. On pose alors
  \begin{equation}
    \label{eq:ssu}
  X_s=\iota_X(X^+_s,X^-_s)
\end{equation}
et
 \begin{equation}
    \label{eq:nilpu}
X_n=X-X_s.
\end{equation}

\begin{remarque}
  \label{rq:rationalite}
On est quelque peu imprécis dans la notation $X\in \tugo$. On peut interpréter celle-ci comme $X\in \tugo(F)$ ou $X\in \tugo(\bar{F})$ où $\bar{F}$ est une clôture algébrique de $F$. Les constructions valent sur ces deux ensembles. Notons bien que, par construction, on a $X_s,X_n\in \tugo(F)$ si $X\in \tugo(F)$. 
\end{remarque}

\begin{lemme}\label{lem:J1u}
  On a $a(X)=a(X_s)$ et les éléments $X_s$ et $X_n$ sont respectivement semi-simples et nilpotents. De plus, pour tout $\delta\in U$, on a $(\delta\cdot X)_s=\delta\cdot X_s  $ et $(\delta\cdot X)_n=\delta\cdot X_n$.
\end{lemme}

\begin{preuve}
cf. démonstrations des lemmes \ref{lem:J1}, \ref{lem:J2} et \ref{lem:action}.
\end{preuve}
\end{paragr}

\subsection{Orbites semi-simples associés à un invariant}\label{ssec:orbss}

\begin{paragr}[Décomposition de $a\in \Ac(F)$.] ---  \label{S:dec-a}Soit $E/F$ une extension quadratique de  corps de caractéristique $0$ et $\sigma$ le générateur du groupe de Galois $\Gal(E/F)$. Soit $V$ un $F$-espace vectoriel de dimension $n$. Soit $\Ac=\Ac_V$ (cf. §\ref{S:AV}) et $a\in \Ac(F)$. 

Il existe un unique entier $r$ tel que $a\in \Ac^{(r)}(F)$. Soit $V=V^+\oplus V^-$ une décomposition de $V$ en sous-espaces telle que $\dim(V^+)=r$. Alors $a$ correspond à un unique couple $(a_+,a_-)\in \Ac_{V^+}^{\rs}(F)\times \Ac^{(0)}_{V^-}(F)$. La donnée de $a_-$ est en fait la donnée d'un polynôme $P$ unitaire de degré $n-r$. Soit
$$P=\prod_{i\in I} P_i^{n_i} \prod_{j\in J} P_j^{n_j}$$
sa décomposition en polynômes irréductibles où les facteurs $P_i$ et $P_j$ sont irréductibles, unitaires et deux à deux distincts. On suppose que $P_i$ reste irréductible sur $E$ alors que $P_j$ lui se décompose sur ce corps. On a donc  $P_j=Q_j Q_j^\sigma$ où $Q_j$ est un polynôme unitaire à coefficients dans $E$ et irréductible sur ce corps.

Pour tout $i\in I\cup J$, soit $F_i=F[t]/(P_i)$ et $E_i=F_i\otimes_F E$. Alors $F_i$ est une extension de $F$ et $E_i$ est une algèbre étale de dimension $2$ sur $F_i$. L'algèbre $E_i$ est un corps si et seulement  $i\in I$. Soit $\sigma_i$ l'involution $1\otimes \sigma$ de $E_i$.  Soit $\al_i\in F_i$ la classe du monôme $t$.
\end{paragr}

\begin{paragr}[Relèvement de $a$ dans $\tgl(V)$.] --- D'après le lemme \ref{lem:fibress} et sa preuve, il existe un triplet $(A^+\oplus A^-,b,c)\in \tgl(V)$ tel que $(A^+,b,c)\in \tgl(V^+)$ est d'invariant $a_+$ et $A^-$ est un endomorphisme semi-simple de $V^-$ de polynôme caractéristique $P$. En particulier, le triplet $X=(A^+\oplus A^-,b,c)$ est un élément semi-simple de $\tgl(V)$ d'invariant $a$. On a donc une décomposition de $V^-$ en espaces propres généralisés de $A^-$
$$V^-=(\oplus_{i\in I\cup J}V_i )$$
Chaque sous-espace $V_i$ pour $i\in I\cup J$ est stable par $A^-$ ; soit $A_i$ l'endomorphisme de $V_i$ induit par $A^-$. Le polynôme minimal de $A_i$ est $P_i$. L'endomorphisme  $A_i$ fait donc de $V_i$ un $F_i$-espace vectoriel, l'action de $A_i$ correspondant à la multiplication par le scalaire $\al_i\in F_i$.
\end{paragr}

\begin{paragr}[Décomposition sur $E$.] --- On affuble un $F$-espace vectoriel d'un indice $E$ pour indiquer le $E$-espace vectoriel obtenu par extension des scalaires à $E$. Ainsi $V_E= V\otimes_F E$. La décomposition précédente induit donc une décomposition en $E$-espaces vectoriels
  \begin{equation}
  \label{eq:VE}
  V_E= V^+_E \oplus( \oplus_{i\in I\cup J}V_{i,E}).
\end{equation}
D'après ce qui précède, chaque facteur $V_{i,E}$ est naturellement un $E_i$-module libre de type fini. 
\end{paragr}

\begin{paragr}[Formes hermitiennes sur $V_{i,E}$.] --- Pour tout $i\in I\cup J$, soit $\mathcal{X}_i$ l'ensemble des formes $\sigma_i$-hermitiennes non-dégénérées $\Phi_i$ sur le $E_i$-module $V_{i,E}$. Le groupe $G_i=GL_{E_i}(V_{i,E})$ agit naturellement sur $\mathcal{X}_i$.

 \begin{lemme}\label{lem:transi}
Pour $j\in J$, le groupe $G_j$ agit transitivement sur $\mathcal{X}_j$. 
  \end{lemme}

  \begin{preuve}
    Soit $j\in J$. Soit $\tau_1$ et $\tau_2$ les  deux morphismes de $E$ dans $F_j$ comme $F$-extensions. On a bien sûr $\tau_2=\tau_1\circ \sigma$.  On a un isomorphisme déterminé par $x\otimes y \mapsto (x\tau_1(y),x\tau_2(y))$ 
$$E_j\simeq F_j\times F_j$$
qui échange l'involution $\sigma_i$ en l'involution qui permute les coordonnées.
On a donc un isomorphisme
$$V_{j,E}=V_j \otimes_{F_j} E_j\simeq V_j\otimes_{F_j} (F_j\times F_j)\simeq  V_j\times V_j$$
via lequel le groupe $G_j$ s'identifie à $GL_{F_j}(V_j)\times GL_{F_j}(V_j)$. Pour toute forme $\Phi_j\in \mathcal{X}_j$, on note encore la forme qui s'en déduit sur  $V_j\times V_j$.  Les sous-espaces   $V_j\times (0)$ et $(0)\times V_j$ sont nécessairement totalement isotropes. On voit que la donnée de $\Phi_j$ équivaut à la donnée d'une forme $F_j$-bilinéaire  non-dégénérée 
$$\bg \cdot,\cdot\bd :V_j \times V_j \to F_j,$$
par la relation $\Phi_j((v_1,v_2),(w_1,w_2))=(\bg w_1,v_2\bd,\bg v_1, w_2\bd )$.  L'action de  $G_j$ sur les formes $\sigma_j$-hermitiennes non dégénérées $\Phi_j$ se traduit en l'action évidente $GL_{F_j}(V_j)\times GL_{F_j}(V_j)$ sur les formes $F_j$-bilinéaires  non-dégénérées et cette dernière est évidemment transitive. 
  \end{preuve}

Pour tout $\Phi_i\in \mathcal{X}_i$, soit $\Phi_{i,E}$ la forme $\sigma$-hermitienne non-dégénérée sur le $E$-espace $V_{i,E}$  donnée par 
$$\Phi_{i,E}(v,w)=\trace_{E_i/E}(\Phi_i(v,w)).
$$
Pour cette forme, l'endomorphisme $A_i\otimes 1$ est auto-adjoint. On va utiliser le lemme suivant dont la vérification est laissé au lecteur.

\begin{lemme}
  \label{lem:bij-forme}
Pour tout $i\in I\cup J$, l'application $\Phi_i\in \mathcal{X}_i \mapsto   \Phi_{i,E}$  est une bijection de  $\mathcal{X}_i$ sur l'ensemble des formes  $\sigma$-hermitiennes non-dégénérées sur le $E$-espace $V_{i,E}$ pour lesquelles l'endomorphisme $A_i\otimes 1$ est auto-adjoint.
\end{lemme}

Soit $\mathcal{X}_{I,J}=\prod_{i\in I\cup J}\mathcal{X}_i$ muni de l'action de $G_{I,J}=\prod_{i\in I\cup J} G_i$. Pour tout $\Phi=(\Phi_{i})_{i\in I\cup J}\in \mathcal{X}_{I,J}$, soit $\Phi_{E}$ la forme 
$$\oplus^{\perp}_{i\in I\cup J} \Phi_{i,E}$$
sur $\oplus_{i\in I\cup J}V_i$.
\end{paragr}

\begin{paragr}[Formes hermitiennes sur $V^+_{E}$.] --- Il existe une forme $\Phi^+$ (unique à équivalence près) et un couple $(A_0^+,b_0)$ tel que $(A^+_0,b_0)\in \tugo(V^+\otimes_F E,\Phi^+)$ soit d'invariant $a_+$ (cela se démontre comme le lemme 2.3 de \cite{AFL}).  On fixe de tels éléments.
\end{paragr}

\begin{paragr}
  Soit $\mathcal{Y}_a$ l'ensemble des couples $(\Phi, X)$ où $\Phi$ est une forme hermitienne non-dégénérée sur $V_E$ et $X\in \tugo(V_E,\Phi)_a(F)$ est semi-simple. Le groupe $G=GL_E(V_E)$ agit sur  $\mathcal{Y}_a$ par $g\cdot \Phi=(\Phi\circ g^{-1},g\cdot X)$.

Pour tout $\Phi\in\mathcal{X}_{I,J} $, soit $A_0=  A_0^+ \oplus (A^-\otimes 1)$ et 
$$X_{\Phi}=(\Phi^+\oplus^\perp  \Phi_{E},  (A_0,b_0)).$$
C'est un élément de $\mathcal{Y}_a$.

On regarde $G_{I,J}$ comme le sous-groupe évident de $GL_E(V_E)$ qui préserve chaque facteur de la décomposition \eqref{eq:VE} et agit trivialement sur le acteurs $V^+_E$. L'application $X:\Phi  \mapsto  X_{\Phi}$ est alors $G_{I,J}$-équivariante. Elle induit donc un foncteur entre groupoïdes quotients
$$[X]:[\mathcal{X}_{I,J}/G_{I,J}]\to [\mathcal{Y}_a/ GL_E(V_E)].$$

\begin{proposition}\label{prop:equiv}
Le foncteur  $[X]$ est une équivalence de catégories.
\end{proposition}

\begin{preuve}  Soit $\Phi \in  \mathcal{X}_{I,J}$. Écrivons $X_{\Phi}=(\Phi_0,A_0,b_0)$. Prouvons que le foncteur est essentiellement surjectif. Soit $(\Phi_1, A_1,b_1)\in \mathcal{Y}_a$. On lui associe le triplet $(A_1,b_1,\Phi_1(b_1,\cdot))\in \tgl_E(V_E)$ où l'on voit $ \Phi_1(b_1,\cdot)$ comme une forme linéaire sur $V_E$. D'après les lemmes \ref{lem:ss} et \ref{lem:ssu}, c'est un élément semi-simple et l'invariant de ce triplet est encore $a$. De même, on associe à $ (\Phi_0,A_0,b_0)$ le triplet   $(A_0,b_0,\Phi_0(b_0,\cdot))\in \tgl_E(V_E)$ semi-simple d'invariant $a$. D'après le lemme \ref{lem:fibress}, il existe $g\in  GL_E(V_E)$ qui conjugue ces deux triplets. Quitte à utiliser cet élément, on peut et on va supposer que $A_1=A_0$, $b_1=b_0$ et $\Phi_1(b_0,\cdot)=\Phi_0(b_0,\cdot)$. Il s'ensuit que l'orthogonal de $V^+_E$ pour $\Phi_1$ est égal à $V^-_E$. Comme $(A^i_0 b_0)_{0\leq i\leq r-1}$ est une base de  $V^+_E$, les formes $\Phi_1$ et $\Phi_0$ coïncident sur cette base (par définition de l'invariant $a$). On a donc $\Phi_{|V^+_E}=(\Phi_0)_{|V^+_E}$. L'endomorphisme $A_0$ étant auto-adjoint pour la forme $\Phi_1$, il est alors clair que la décomposition \ref{eq:VE} est orthogonale pour $\Phi_1$ puis que $\Phi_1=\Phi'_{E}$ pour un élément $\Phi'\in \mathcal{X}_{I,J}$. Ainsi le triplet de départ est dans l'image essentielle du foncteur $[X]$.

Soit  $\Phi,\Phi' \in  \mathcal{X}_{I,J}$. Soit $g\in G$ qui conjugue $X_\Phi$ en $X_{\Phi'}$. Alors $g$ centralise $(A_0,b_0)$ mais le centralisateur de cet élément est exactement  $ G_{I,J}$.  Il s'ensuit que $\Phi_{E}$ et $\Phi'_E$ sont conjugués sous $G_{I,J}$ donc les formes  $\Phi$ et $\Phi'$ sont dans la même orbite sous $G_{I,J}$.
 \end{preuve}

 \begin{corollaire}\label{cor:orbite-ss}
Le foncteur $[X]$ induit une bijection de l'ensemble quotient $\prod_{i\in I}\mathcal{X}_i/ G_i$ sur l'ensemble des $U_\Phi(F)$-orbites des éléments semi-simples dans $\tugo_{\Phi,a}(F)$ lorsque $\Phi$ décrit un système de représentants de formes hermitiennes sur $V_E$.
 \end{corollaire}

\begin{preuve}
  Le quotient $\mathcal{Y}_a/G$ s'identifie à l'ensemble des $U_\Phi(F)$-orbites des éléments semi-simples dans $\tugo_{\Phi,a}(F)$ pour $\Phi$ décrivant un système de représentants de formes hermitiennes sur $V_E$. La proposition donne une bijection de $\mathcal{X}_{I,J}/G_{I,J}$ sur $\mathcal{Y}_a/G$. Comme $X_j/G_j$ est réduit à un point, le corollaire s'ensuit. 
\end{preuve}

\end{paragr}

\subsection{Sections transverses}\label{ssec:transvu}

\begin{paragr}\label{S:V+U}
Soit $E/F$ une extension quadratique de  corps de caractéristique $0$. Soit $\sigma$ le générateur du groupe de Galois de $\Gal(E/F)$. Soit $(V^+,\Phi^+)$ un espace $\sigma$-hermitien non dégénéré. Soit $\tugo_ +=\tugo_{(V^+,\Phi^+)}$,  $U_+=U(V^+,\Phi^+)$ et $\Ac_+=\Ac(V^+,\Phi^+)$ (cf. § \ref{S:cadreu}).
\end{paragr}

\begin{paragr}\label{S:VijU}
Soit $I$ et $J$ deux ensembles d'indices, finis et disjoints (éventuellement vides). Pour tout $i\in I\cup J$, soit
\begin{itemize}
\item  $F_i$  une extension finie de $F$ ;
\item $E_i=F_i\otimes_F E$ et $\sigma_i$ l'involution $\Id\otimes \sigma$ ;
\item $W_i$ un $F_i$-espace vectoriel de dimension finie et $V_i=W_i\otimes_{F_i}E_i$
\item $\Phi_i$ une forme $\sigma_i$-hermitienne non dégénérée sur $V_i$ ;
\item $\Phi_{i,E}=\trace_{E_i/E}\circ \Phi_i$ ; c'est une forme $\sigma$-hermitienne non dégénérée sur $V_i=W_i\otimes_{F}E$ vu comme $E$-espace vectoriel.
\end{itemize}
On suppose que $i$ appartient au sous-ensemble $I$ si et seulement si $E_i$ est un corps.
\end{paragr}

\begin{paragr}
Supposons $i\in I$.  Soit $\tugo_ i=\tugo(V_i,\Phi_i)$ l'espace construit au §  \ref{S:cadreu} relativement au $E_i$-espace hermitien $(V_i,\Phi_i)$. 
 \end{paragr}
 
 \begin{paragr}  Soit $j\in J$. Comme dans la preuve du lemme \ref{lem:transi}, on identifie $E_j$ à $F_j\times F_j$ et $V_j$ à $W_j\times W_j$. La donnée de la forme $\Phi_j$ est équivalente à la donnée d'une forme $F_j$-bilinéaire non dégénérée 
   \begin{equation}
     \label{eq:laforme}
   \bg \cdot,\cdot\bd  \, : \, W_j \times W_j \to F_j.
   \end{equation}
 Soit $\tggo_j=\tgl_{F_j}(W_j)$ (cf.  §\ref{S:tggo}).  
\end{paragr}

\begin{paragr}
Soit le $E$-espace vectoriel
$$V^-=\oplus_{i\in I\cup J}V_i 
$$
On le munit de la forme hermitienne non-dégénérée (qui fait de la somme ci-dessus une somme orthogonale)
$$\Phi^-= \oplus^\perp_{i\in I\cup J},\Phi_{i,E} .
$$
Soit
$$(V,\Phi)=(V^+,\Phi^+)\oplus  (V^-,\Phi^-)
$$
\end{paragr}

\begin{paragr} Soit
$$\tugo_-= \big(\oplus_{i\in I} \tugo(V_i,\Phi_{i}) \big) \oplus    \big(\oplus_{j\in J} \tggo_j\big) $$
 muni de l'action du groupe 
$$U_-= \prod_{i\in I} U(V_i,\Phi_i) \times  \prod_{j\in J} GL_{F_j}(W_j).
$$
Ce dernier est naturellement un sous-groupe du groupe unitaire $U(V^-,\Phi^-)$ (pour $j\in J$, le facteur $GL_{F_j}(W_j)$ agit naturellement sur le premier facteur et par l'inverse de son adjoint par rapport à \eqref{eq:laforme} sur le second facteur de $V_j=W_j\times W_j$). Soit $\Ac_-$ le quotient catégorique de $\tugo_-$ par $U_-$ et soit $a: \tugo_-\to \Ac_{-}$ le morphisme canonique. Par restriction des scalaires, on voit tous ces objets comme des objets définis sur $F$.
\end{paragr}

\begin{paragr}
  Soit
\begin{equation}
  \label{eq:iota-}
\iota_-:   \tugo_-  \to \tugo(V^-,\Phi^- )
\end{equation}
le $F$-morphisme qui induit
\begin{itemize}
\item pour $i\in I$, le morphisme évident  $\tugo(V_i,\Phi_{i}) \to   \tugo(V_i,\Phi_{i,E})  $ d'oubli de la $E_i$-structure ;
\item pour $j\in J$, le morphisme $\tggo_j \to \tugo(V_j ,\Phi_{j,E})$  donné par 
$$ (A,b,\bg \cdot, c\bd ) \mapsto ( A\oplus A^\vee  , (b,c))
$$
l'on identifie $V_j$ à $W_j\times W_j$, $\End_{E_j}(V_j)$ à $\End_{F_j}(V_j) \oplus \End_{F_j}(V_j)$,  $V_j$ à $V_j^*$ par $c\mapsto \bg \cdot ,c\bd$ (où l'accouplement est celui considéré en \eqref{eq:laforme}, et $A^\vee$ est l'adjoint à droite de $A$ pour la forme  \eqref{eq:laforme}.

\end{itemize}
Le morphisme $\iota_-$ est équivariant pour l'action du groupe $U_{-}$ : il se descend donc en un morphisme homonyme 
\begin{equation}
  \label{eq:iota-quot}
\iota_- :   \Ac_{-} \to \Ac(V^-,\Phi^-).
\end{equation}
\end{paragr}

\begin{paragr} Soit  $D$ la fonction sur $\Ac_{-}$ qui, à un élément $((A_i,b_i)_{ i\in I},(A_j,b_j,c_j)_{j\in J})\in \tugo_-$, associe le produit 
$$\prod_{i,j} \mathrm{Res}(P_i,P_j)$$
 où
 \begin{itemize}
 \item les couples $(i,j)$ parcourent l'ensemble $(I\coprod (J\times \Gal(E/F)))^2$ privé de sa diagonale,
 \item $P_i$ est  le polynôme caractéristique de $A_i$ vu comme $E$-endomorphisme de $V_i$ pour $ i \in I\coprod J$
\item $P_{(j,\tau)}=\tau(P_j)$ pour tous $(j,\tau)\in J\times  \Gal(E/F)$
 \item $\mathrm{Res}$ désigne le résultant.
 \end{itemize}
Soit  $\Ac_{-}'$ l'ouvert de  $\Ac_{-}$ défini par la condition $D\not=0$. Le morphisme induit par \eqref{eq:iota-quot} sur $\Ac_{-}'$ est  étale. Soit $\tugo_-'$  l'image inverse de  $\Ac_{-}'$ par le morphisme canonique $a$. 
\end{paragr}

\begin{paragr}
  Soit $\tugo^\flat=\tugo_+\oplus \tugo_-$ muni de l'action du groupe $U^\flat=U_+\times U_-$. Soit  $\Ac^\flat=\Ac_{+}\times \Ac_{-}$ le quotient catégorique.  
 En composant \eqref{eq:iota-su} avec \eqref{eq:iota-}, on  obtient un morphisme    $U^\flat$ -équivariant 
\begin{equation}
  \label{eq:iotaU}
 \iota: \tugo_+^{\rs}\times \tugo_- \to \tugo
\end{equation}
ce qui donne un morphisme homonyme  sur les quotients catégoriques
$$\iota: \Ac_+^{\rs}\times \Ac_-  \to \Ac_U
$$
Soit $\tugo^{\flat,\Urs}$ (resp. $\Ac^{\flat,\Urs}$) l'image réciproque par $\iota$ de l'ouvert  $\tugo^{\rs}$ (resp.  $\Ac_U^{\rs}$).

Soit $(\Ac^\flat)'=\Ac_{+}^{\rs}\times \Ac_{-}'$.  Soit $(\tugo^{\flat})'$ l'ouvert de $\tugo^{\flat}$ obtenu par image inverse de $(\Ac^\flat)'$ par le morphisme canonique. Notons qu'on a $\tugo^{\flat,\Urs}\subset (\tugo^{\flat})' $ et $\Ac^{\flat,\Urs}\subset  (\Ac^\flat)'$.
Le morphisme $\iota$ induit un morphisme étale sur l'ouvert $(\Ac^\flat)'$ (cf. \cite{Z1} appendice B). On a enfin un isomorphisme (\emph{loc. cit.})
\begin{equation}
  \label{eq:isocrucialU}
   U \times^{U^\flat} (\tugo^\flat)' \to \tugo \times_{\Ac_\Phi} (\Ac^\flat)'
 \end{equation}
qui est donné par $(g,Y) \mapsto (g\cdot \iota (Y), a(Y))$ et  où la source désigne le quotient de  $U \times  (\tugo^\flat)'$ par l'action libre à droite de $U^\flat$ donnée par $(g,Y)\cdot h=(gh,h^{-1}\cdot Y)$.
\end{paragr}

\begin{paragr}[Situation sur un corps de nombres.] ---\label{S:UsurA}
On suppose de plus que $E/F$ est une extension quadratique de  corps de nombres. On va étendre les constructions et les résultats ci-dessus sur une base un peu plus générale. Soit   $S$ est un ensemble fini de places de $F$ contenant les places archimédiennes. Soit  $A$ l'anneau des entiers de $F$ \og hors $S$ \fg{} et $A_E$ la clôture intégrale de $A$ dans $E$. On suppose que les anneaux $A$ et  $A_E$ sont principaux et que le morphisme $A\to A_E$ est étale. On se donne alors les objets suivants  pour $i\in I\cup J$ :
\begin{itemize}
\item $A_i$ et $B_i$ les clôtures intégrales de $A$ dans $F_i$ et $E_i$ (on suppose que le morphisme $A\to A_i$ est étale) ; 
\item  $W_i$ un $A_i$-module libre et $V_i=W_i\otimes_{A_i}B_i$ ;
\item  $\Phi_i$ une forme $\sigma_i$-hermitienne non dégénérée sur le $B_i$-module $V_i$.
\end{itemize}
Soit $V^+$ un $A_E$-module libre muni d'une forme  $\sigma$-hermitienne non dégénérée $\Phi^+$. Les contructions précédentes donnent alors des schémas en groupes réductifs sur $A$ notés encore $U^\flat$ et $U$ agissant sur les $A$-modules $\tugo^\flat$ et $\tugo$. On dispose encore de morphismes vers  les espaces affines $\Ac^\flat$ et $\Ac$ sur $A$. Le lemme suivant se démontre comme le lemme \ref{lem:integrite}.

   \begin{lemme}\label{lem:integriteu}
    Soit $v\notin S$ et $\oc_v$ le complété de $A$ en $v$. Soit $F_v$ le corps des fractions. Soit $a\in (\Ac^\flat)'(\oc_v)$. 
    \begin{enumerate}
    \item Pour tous $Y\in \tugo^\flat_a(F_v)$ et $g\in U(F_v)$, les assertions suivantes sont équivalentes :
      \begin{enumerate}
      \item $g^{-1}\cdot \iota(Y)\in \tugo(\oc_v)$.
      \item Il existe $h\in U^\flat(F_v)$ tel que $g\in hU(\oc_v)$ et $h^{-1}\cdot Y\in \tugo^\flat(\oc_v)$
      \end{enumerate}
    \item Pour tout $Y\in \tugo^\flat_a(F_v)$, les assertions suivantes sont équivalentes
      \begin{enumerate}
      \item $\iota(Y)\in \tugo(\oc_v)$.
      \item $Y\in \tugo^\flat(\oc_v)$.
      \end{enumerate}
    \item Pour tout $g\in U(F_v)$ et tout $X\in \tugo^{\rs}(\oc_v)$ les assertions suivantes sont équivalentes
      \begin{enumerate}
      \item $g^{-1}\cdot X\in \tugo^{\rs}(\oc_v)$.
      \item $g\in U(\oc_v)$.
      \end{enumerate}
          \end{enumerate}
  \end{lemme}

\end{paragr}

\section{Combinatoire des cônes}\label{sec:combiU}

\subsection{Sous-espaces paraboliques}

\begin{paragr}\label{S:parabu}
  Soit 
$$(0)=W_0\subsetneq W_1 \subsetneq \ldots \subsetneq W_r$$
un drapeau de sous-espaces totalement isotropes de $V$. On lui associe les objets suivants :

\begin{itemize}
\item le sous-groupe parabolique $P\subset U$ qui stabilise le drapeau ;
\item le sous-espace $V_P=W_r$ ;
\item le sous-espace $\tpgo $ de $\tugo$ formé des couples $(A,b)$ tels que  $b\in V_P^\perp$ (orthogonal de $V_P$ pour $\Phi$) et $A$ est un   endomorphisme auto-adjoint de $V$ qui vérifie  $AW_i\subset W_i$ pour $1\leq i \leq r$ ;
\item le sous-espace $\tngo=\tngo_{P}$ de $\tugo$ formé des couples $(A,b)$ tels que $b\in V_P$ et $A$ est un   endomorphisme auto-adjoint de $V$ qui vérifie  $AW_i\subset W_{i-1}$ et $AW_r^{\perp}\subset W_r$ pour $1\leq i \leq r$ ;
\end{itemize}
Les sous-espaces $\tpgo$ et $\tngo$ sont  stables sous l'action de $P$. Notons enfin que la donnée d'un drapeau est équivalente à la donnée d'un sous-groupe parabolique $P$ de $U$. 
\end{paragr}

\begin{paragr}[Facteur de Levi.] ---\label{S:LeviU}
  Le choix d'un facteur de Levi $M$ du sous-groupe parabolique  $P$ correspond au choix d'une somme directe  orthogonale
$$V=(V_1\oplus V_{-1})  \bigoplus^{\perp} \ldots\bigoplus^{\perp} (V_r  \oplus V_{-r}) \bigoplus^{\perp} V_0 $$
qui vérifie les propriétés suivantes
\begin{itemize}
\item  pour $1\leq i\leq r$, les espaces $V_i$ et $V_{-i}$ sont totalement isotropes ;
\item Pour $1\leq i \leq r$, on a 
$$W_i=V_1\oplus\ldots\oplus V_i$$
et
$$W_{i}^{\perp}=W_r \oplus V_0 \oplus V_{-r}\oplus \ldots \oplus V_{-(i+1)}
$$
(avec la convention $V_{-(r+1)}=(0)$.
\end{itemize}

En particulier, on a $V_P^\perp=V_P\oplus V_0$. Soit $\tmgo$ le sous-espace de $\tugo$ formé des couples $(A,b)$ tels que $AV_i\subset V_i$ pour $i\in \{-r,\ldots,0,\ldots,r\}$ et $b\in V_0$.  
\end{paragr}

\begin{paragr}[Les fonctions $\upla$.]  ---
Soit $ \upla_{P} \in \ago_{P}^{*} $ le déterminant de l'action du tore $ A_{P} $ 
sur $ V_{P} $. Alors, si on note $ \rho_{\tlP} $ (resp. $ \rho_{P} $)
la demi-somme des poids pour l'action de $ A_{P} $ sur $ \tngo_{P} $ 
(resp. sur $ \ngo_{P} $) on a
\[ 
\upla_{P} = 2\rho_{\tlP} - 2\rho_{P}.
\]
\end{paragr}

\begin{paragr}[Lemmes géométriques.] --- On utilise les notations du paragraphe \ref{S:tau}. 
Pour tout sous groupe parabolique $ P $ de $ U $ et tous $H,X\in \ago_{0}$
on rappelle la fonction $ \Gamma' $ d'Arthur \cite{arthur2} section 2:
\begin{equation*}
  \Gamma'_{P}(H,X) = 
\sum_{R\in \fc(P) }\varepsilon_{P}^{R}
\hat{\tau}_{R}(H - X) \tau_{P}^{R}(H).
\end{equation*}

On a alors:

\begin{lemme}\label{lem:upGammaU}
\begin{enumerate}
\item On a
\[
\hat{\tau}_{P}(H-X) = \sum_{R\in \fc(P) } \varepsilon_{R}^{U}
\hat{\tau}_{P}^{R}(H)
\Gamma'_{R}(H,X), \quad X,H \in \ago_{0}.
\]
\item Pour tout $X \in \ago_{0}$, la fonction $H \in  \ago_{P} \mapsto \Gamma_{P}'(H,X)$ est à support compact. 
\item Pour tous $P \subset R$, il existe un  polynôme $ p_{P,R} $ sur $ \ago_{R} $ tel que pour tout $X\in \ago_{0}$ on ait 
\[ 
\int_{\ago_{P}}e^{\upla_{P}(H)}\Gamma'_{P}(H,X)dH = 
\sum_{P \subset R} e^{\upla_{R}(X_{R})}p_{P,R}(X_{R})
\]
où $ X_{R} $ désigne la projection orthogonale de $ X $ sur $ \ago_{R} $. 
\end{enumerate}
\end{lemme}
\begin{preuve}
  Le point 1 est expliqué dans \cite{arthur2}, §2 
  et le deuxième est démontré dans le lemme 2.1 de \emph{loc. cit}.
  Le dernier c'est le lemme 4.3 de \cite{leMoi}. 
\end{preuve}
\end{paragr}

\subsection{Descente et  combinatoire des cônes}\label{ssec:desc-combiU}

\begin{paragr}
 On reprend les notations du §\ref{ssec:transvu}.
\end{paragr} 
 
 \begin{paragr}[Les groupes $U_{i}$] --- 
 Soit $ i \in I $. 
 On écrit $ U_{i} = U(V_i,\Phi_{i}) $. 
 En utilisant le théorème de décomposition de Witt, on décompose 
 l'espace hermitien en une somme orthogonale 
 \[
V_{i} = \bigoplus_{l=1}^{n_{i}}(D_{i, l} \oplus D_{i,l}^{*, \sigma_{i}})
\oplus V_{i}^{\flat} 
 \]
 où les $D_{i,l}$ sont des $E_{i}$-droites, 
 les $D_{i,l}^{*, \sigma_{i}}$ s'identifient aux $\sigma_{i}$-duales de $D_{i,l}$ 
et $V_{i}^{\flat}$  est anisotrope. 
 Ces données définissent uniquement un $ F_{i} $-sous-groupe de Levi minimal, noté $ L_{i} $, 
 de $ U_{i} $. 
 \end{paragr}

 \begin{paragr}[Les groupes $G_{j}$ et $\tlG_{j}$] --- 
 Soit $ j \in J $.
 On fixe l'isomorphisme $ V_{j,E} \cong W_{j} \times W_{j}$ 
 comme dans la preuve du lemme \ref{lem:transi}. 
 Soit $\oplus_{l=1}^{n_{j}} D_{j,l}$ 
une décomposition de $ W_{j} \times (0) \subset V_{j,E}$
 en $ F_{j} $-droites. 
 Puisque $ W_{j} \times (0)  $ est anisotrope 
 et $ (0) \times W_{j}  $ est en dualité avec lui, on obtient 
 la décomposition duale  $ \oplus_{l=1}^{n_{j}} D_{j,l}^{*, \sigma_{j}}$ de ce dernier.

 La décomposition de $ W_{j} $ 
 définit un sous-groupe de Levi minimal $ L_{j} $ 
 de $ G_{j} := GL_{F_{j}}(W_{j}) $.
 
 On pose aussi $ \tlG_{j} := GL_{E_{j}}(W_{j} \oplus F_{j}e_{j})  $ 
 où $ e_{j} $ est un vecteur fixé. 
 On identifie $G_{j}$ 
  avec le sous-groupe de $\tlG_{j}$ qui agit trivialement sur $F_{j} e_{j}$.
On note  $ \tlL_{j} $ le sous-groupe de Levi minimal 
 de $ \tlG_{j} $ 
  défini comme le centralisateur de 
 $ L_{j} $ dans $  \tlG_{j}  $. 
 On se trouve alors dans la situation du §\ref{ssec:parab}. 
 Les éléments de $\fc^{\tlG_{j}}(\tlL_{j})$ sont donc les sous-groupes 
 paraboliques relativement standards, relativement au sous-groupe de Levi $L_{j}$ de $G_{j}$.
\end{paragr}

\begin{paragr}[Sous-groupes $\tlU_-$, $\tlM_{0}^{-}$, $M_{1}$.] ---  \label{S:M1U} 
Soient 
$$
\tlU_{-}= \prod_{i\in I} U_i \times \prod_{j \in J} \tlG_{j}, \quad
\tlM_0^{-} =  \prod_{i\in I} L_i \times \prod_{j\in J} \tlL_j, \quad
M_1= U_{+} \times \prod_{i\in I \cup J}L_{i}.
 $$
 Le groupe $\tlM_{0}^{-}$ est un $F$-sous-groupe de Levi minimal de $\tlU_{-}$ 
 (si on considère ces groupes sur $F$, par restriction de scalaires).
 On identifie aussi $M_{1}$ à un sous-groupe de Levi de $U$. 
\end{paragr}

\begin{paragr}[L'application $ Q\mapsto \tlQ^-$.] --- \label{S:QQ-U}Soit $Q\in \fc^{U}(M_1)$
on va lui associer le groupe $\tlQ^-$, élément de $\fc^{\tlU_{-}}(\tlM_0^{-})$.

Pour tout $ i \in I$, on définit le groupe 
$ Q_{i} \in \fc^{U_{i}}(L_{i}) $ 
par le drapeau des espaces isotropes obtenu 
par l'intersection du drapeau 
des sous-espaces isotropes de $ V $ dont $ Q $ est le stabilisateur
avec l'espace $ V_{i} $. 

Pour tout $ j \in J$ on définit le sous-groupe parabolique 
$ \tlQ_{j} \in \fc^{\tlG_{j}}(\tlL_{j}) $ (nécessairement relativement standard) 
comme associé aux données suivantes(cf. §\ref{S:parab}):
\begin{itemize}
\item le drapeau  de $F_j$-espaces vectoriels 
\begin{equation}
  \label{eq:drapeau3}
 (0) = W_{j,0} \subsetneq W_{j,1} \subsetneq \cdots \subsetneq W_{j,s_{j}} = W_{j} \times (0)
\end{equation}
où les $W_{j,l}$ sont obtenus par intersection du drapeau complet dans $V$ dont
$ Q $ est le stabilisateur avec $ W_{j} \times (0) \subset V_{j} $ et élimination des doublons ;
\item le couple  d'entiers $(i_0',j_0')$ définis par la condition 
$ W_{j,i_{0}'} =  V_{Q} \cap (W_{j} \times (0)) $ et  $W_{j,j_{0}'} =  
V_{Q}^{\perp} \cap (W_{j} \times (0))$.
\end{itemize}

On pose alors 
$$\tlQ^-=  \prod_{i\in I} Q_i \times \prod_{j \in J} \tlQ_{j}.
$$
\end{paragr}

\begin{paragr}[Les espaces $\zgo_{\tlQ^{-}}$.] --- \label{S:zgoQ}
Fixons $ \tlR = \prod_{i \in I}R_{i} \times \prod_{j \in J}\tlR_{j} \in 
\fc^{\tlU_{-}}(\tlM_{0}^{-}) $. Soit $ M_{\tlR}  = \prod_{i \in I}M_{R_{i}} \times \prod_{j \in J} M_{\tlR_{j}} $ 
son facteur de Levi contenant contenant $ \tlM_0^{-} $, où 
$ R_{i} \in \pc^{U_{i}}(M_{R_{i}})$ 
et $\tlR_{j} \in \pc^{\tlG_{j}}(\tlL_{j})$. 
On pose 
$ \MM_{\tlR} := \prod_{i \in I} M_{R_{i}} \times \prod_{j \in J}\MM_{\tlR_{j}}$ 
et $ \tlG_{\tlR} := \prod_{j \in J}\tlG_{\tlR_{j}} $ 
de sorte que $ M_{\tlR} = \MM_{\tlR} \times \tlG_{\tlR} $. 
Notons que
pour tout $ Q\in \fc^{U}(M_1)$  tel que $ \tlQ^{-} = \tlR $
on a $ \MM_{\tlR} \subset M_{Q} $.

On pose $ \zgo_{\tlR} = \prod_{i} \ago_{R_{i}} \times \prod_{j} \zgo_{\tlR_{j}}$ 
où $\zgo_{\tlR_{j}} = \ago_{\MM_{\tlR}}$ (cf. §\ref{S:zgoP}).
L'espace $ \zgo_{\tlR} $ est naturellement plongé dans $\ago_{M_{1}}$. 
La définition ne dépend que de $M_{\tlR}$, on écrit donc parfois $\zgo_{M_{\tlR}} = \zgo_{\tlR}$. 
Si $ \tlR = \tlQ^{-} $ pour un $ Q \in \fc^{U}(M_{1}) $ on a 
$ \ago_{Q}  \subset \zgo_{\tlR}$ 
avec égalité si $ Q $ est minimal avec la propriété $  \tlQ^{-}  = \tlR$.
De plus, on a 
$ \ago_{M_{1}}  = \zgo_{\tlM_{0}^{-}}$. 
\end{paragr}

\begin{paragr}[Ensembles $\fc$.] --- \label{S:combiU}  Pour tout $\tlR \in \fc^{\tlU_{-}}(\tlM_0^{-})$, 
on définit des ensembles  $\fc_{\tlR}^{0}(M_1)\subset \fc_{\tlR}^{}(M_1)\subset\overline{\fc}_{\tlR}^{}(M_1) \subset \fc^{U}(M_1)$ par les conditions
\begin{eqnarray*}
\overline{\fc}_{\tlR}^{}(M_1) & = & \{P \in \fc^{U}(M_1) | \,    \tlR \subset \tlP^{-} \}\\
\fc_{\tlR}^{}(M_1) & = & \{P \in \fc^{U}(M_1) | \, \tlP^{-} = \tlR\}, \\ 
\fc_{\tlR}^{0}(M_1) & = &\{P \in \fc_{\tlR}^{}(M_1) | \, \ago_{P} = \zgo_{\tlR}\}.
\end{eqnarray*}

\begin{lemme}
  \label{lem:convexiteU}
Pour tout $\tlR\in  \fc^{\tlU_{-}}(\tlM_0)$, il existe un unique sous-groupe de Levi $M$ de 
$U$ contenant $M_1$ tel que
$$ \fc_{\tlR}^{0}(M_1) \subset\pc^{U}(M). 
$$
En outre,  l'ensemble $ \fc_{\tlR}^{0}(M_1) $ est une famille convexe dans 
$ \pc^{U}(M) $   au sens  de l'appendice \ref{App:domConvP}.
\end{lemme}

\begin{preuve}
  Soient $ \tlR = \prod_{i \in I}R_{i} \times \prod_{j \in J}\tlR_{j}$ 
  et $ M_{\tlR} = \prod_{i \in I}M_{R_{i}} \times \prod_{j \in J} M_{\tlR_{j}}  $ 
  comme dans §\ref{S:zgoQ}. 
  
  Pour tout $ i \in I $, le groupe $M_{R_{i}} $ 
  est le stabilisateur dans $ U_{i} $ des sous-espaces qui apparaissent dans la décomposition
  orthogonale 
  \[ 
  V_{i} = V_{i,1} \oplus \cdots V_{i,s_{i}} \oplus 
  V_{i,s_{i}}^{*, \sigma_{i}} \oplus \cdots \oplus  V_{i,1}^{*, \sigma_{i}} 
  \oplus V_{i}^{\flat}
  \]
  où les espaces $ V_{i,l} $ sont isotropes et en 
  $ \sigma_{i} $-dualité avec $ V_{i,l}^{*, \sigma_{i}} $ 
  et la restriction de $ \Phi_{i} $ à $ V_{i}^{\flat} $ est non-dégénérée. 
 De plus, $R_i$ est le stabilisateur du drapeau isotrope
\[
(0)\subsetneq V_{i,1}\subsetneq (V_{i,1} \oplus V_{i,2}) \subsetneq \ldots 
\subsetneq (V_{i,1} \oplus \cdots \oplus  V_{i,s_{i}}). 
\]

Pour tout $ j \in J $ le groupe $ M_{\tlR_{j}} $ 
est le stabilisateur 
dans $ \tlG_{j} $ des sous-espaces qui apparaissent dans la décomposition
\[ 
W_{j} \oplus F_{j}e_{j} = W_{j,1} \oplus \cdots \oplus W_{j,s_{j}}
\]
et il existe un unique indice $l_j$ tel que $e_j \in W_{i,l_j}$. 
De plus, $\tlR_j$ est le stabilisateur du drapeau
\[
(0)\subsetneq W_{j,1}\subsetneq (W_{j,1} \oplus W_{j,2}) \subsetneq \ldots 
\subsetneq (W_{j,1} \oplus \cdots \oplus  W_{j,s_{j}}) = W_{j} \oplus F_{j}e_{j}. 
\]

Si $ l \neq l_{j} $, l'espace $ W_{j,l} $, vu comme un sous-espace de $ W_{j} \times (0) $ 
est une somme des certaines droites $ D_{j,k} $. On définit alors $ W_{j,l}^{*, \sigma_{j}} 
\subset (0) \times W_{j}$ 
comme la somme des $ D_{j,k}^{*, \sigma_{j}} $ pour tout $ D_{j,k} \subset W_{j,l} $. 
On voit alors que le sous-groupe de Levi $M$ est nécessairement le stabilisateur des sous-espaces dans la décomposition
\begin{equation}
  \label{eq:decompoU}
V  =\oplus_{(i,k)} (V_{i,k} \oplus V_{i,k}^{*, \sigma_{i}}) 
\oplus_{(j,l)} (W_{j,l} \oplus W_{j,l}^{*, \sigma_{j}}) \oplus W
\end{equation}
où l'on somme sur les couples $(i,k)$ et 
$(j,l)$
tels que $i \in I$, $1 \leq k \leq s_i$ et 
$ j \in J $, $1\leq l \leq s_j$ et 
$l \not= l_j$ et où l'on pose 
$W= V^{+} \oplus_{i \in I }V_{i}^{\flat} \oplus_{j \in J }(W_{j,l_j} \cap V)$.

Les racines réduites de $A_{M}$ dans $U$ sont alors naturellement indexées par les couples $(V',V'')$ formés de deux éléments distincts parmi les sous-espaces qui apparaissent dans la décomposition \eqref{eq:decompoU}. Soit $\Sigma(\tlR)$ le sous-ensemble des racines associées à l'un des couples suivants pour 
$i\in I$ et $ j \in J $ :
\begin{itemize}
\item  $(V_{i,l}, V_{i, l'})$ pour $l < l'$;
\item  $(V_{i,l}, W)$;
\item   $(W_{j,l},W_{j,l'})$ pour $l<l'$ et $l \not=l_j$, $l'\not=l_j$;
\item  $(W_{j,l}, W)$ pour  $l < l_{j}$ ;
\item $(W,W_{j,l})$ pour $l > l_j$.
\end{itemize}
On a alors 
\[
 \fc_{\tlR}^{0}(M_1) = \bigcap_{\al \in \Sigma(\tlR)}H(\al)^+ ,
 \]
avec les notations de l'appendice \ref{App:domConvP}. Le lemme résulte alors du lemme \ref{lem:hfSpCnv}.
\end{preuve}
\end{paragr}

\begin{paragr}[Des ensembles $ \hat \Pi_{\tlR} $.] --- On reprend les notations du §\ref{S:widePi}. On pose pour $\tlR\in  \fc^{\tlU_{-}}(\tlM_0^{-})$:
\[ 
 \widehat \Pi_{\tlR} :=  \bigsqcup_{i \in I} \hat \Delta_{R_{i}}  \bigsqcup_{j \in J} \widehat \Pi_{\tlR_{j}}
\]
vu comme un sous-ensemble de $ \zgo_{\tlR}^{*} \subset \ago_{M_{1}}^{*} $.
\end{paragr}

\begin{paragr}[Des lemmes combinatoires.] --- Soit $\tlR\in  \fc^{\tlH^-}(\tlM_0)$.

\begin{lemme}\label{lem:artLem51U}
La somme
$$
\sum_{  P \in \overline{\fc}_{\tlR}^{}(M_1)}\eps_{P}^{U} \, \hat{\tau}_{P},
$$
vue comme fonction sur $\ago_{M_1}$ est égale à la fonction caractéristique du cône fermé
\[
\{H \in \ago_{M_1} | \gamma(H) \leq 0, \ \forall \ \gamma \in \widehat{\Pi}_{\tlR}\}.
\]
\end{lemme}

\begin{preuve}
En raisonnant comme dans la preuve du lemme \ref{lem:artLem51}, 
en utilisant les résultats de l'appendice \ref{App:domConvP} 
et le lemme \ref{lem:convexiteU} (analogue du lemme \ref{lem:convexite}) 
on voit que la somme en question est la fonction caractéristique 
de
\[
\{H \in \ago_{M_1}| \varpi(H) \le 0, \ 
\forall \ P \in \fc_{\tlR}^{0}(M_1) \ \forall \ \varpi \in \hat \Delta_{P}\}.
\]

Il est facile de voir que 
\[
\widehat \Pi_{\tlR} \subset \bigcup_{P \in \fc_{\tlR}^{0}(M_1)}\hat \Delta_{P}.
\]
D'autre part, tout $\delta \in \bigcup_{P \in \fc_{\tlR}^{0}(M_1)}\hat \Delta_{P}$ s'écrit comme 
une somme d'éléments de $\widehat \Pi_{\tlR}$ à coefficients positifs d'où le résultat. 
 
\end{preuve}
\end{paragr}

\begin{paragr}[Fonctions $\hat  \sigma_{\tlR}^{\tlS} $.] ---

Pour $ \tlS = \prod_{i \in I}S_{i} \times \prod_{j \in J}\tlS_{j} \in \fc^{\tlU_{-}}(\tlM_{0}^{-}) $ 
contenant $ \tlR $  on définit les fonctions caractéristiques 
$ \sigma_{\tlR}^{\tlS} $ 
et $ \hat \sigma_{\tlR}^{\tlS} $
définies sur $ \ago_{M_{1}} $
comme les produits 
\[ 
\prod_{i \in I} \tau_{R_{i}}^{S_{i}} \cdot \prod_{j \in J} \sigma_{\tlR_{j}}^{\tlS_{j}}, \quad 
\prod_{i \in I} \hat \tau_{R_{i}}^{S_{i}} \cdot \prod_{j \in J} \hat \sigma_{\tlR_{j}}^{\tlS_{j}}
\]
respectivement (voir §\ref{S:sigma} pour la définition des fonctions $\sigma$ et $\hat \sigma$).
\end{paragr}

\begin{paragr}
Soient $\tlR \in \fc^{\tlU_{-}}(\tlM_{0}^{-})$ et 
$M$ le sous-groupe de Levi de $U$ tel que 
$\fc_{\tlR}^{0}(M_{1}) \subset \pc^{U}(M)$ (cf. lemme \ref{lem:convexiteU}). 
Soit
$\yc_{\tlR} = (Y_{P})_{P \in _{\fc_{\tlR}^{0}(M_{1})}}$ 
une famille des vecteurs dans $\ago_{M}$ qui est
 $A_{M}$-orthogonale positive au sens de la définition \eqref{eq:orthPosit}. 
Pour tout $Q \in \overline{\fc}_{\tlR}(M_{1})$ soit $P \in \fc_{\tlR}^{0}(M_{1})$ contenu dans $Q$ et $Y_{Q} \in \ago_{Q}$ la projection orthogonale de $Y_{P}$ sur $\ago_{Q}$. Cette  définition ne dépend du choix de $P \subset Q$. De cette manière, on obtient pour tout $\tlS \in \fc^{\tlU_{-}}(\tlR)$  une  famille positive $\yc_{\tlS} := (Y_{Q})_{Q \in \fc_{\tlS}^{0}(M_{1})}$ .

Pour tout $ H \in \zgo_{\tlR}$, soit
\[
\mathrm{B}_{\tlR}^{U}(H, \yc_{\tlR}) = \sum_{\tlS \in \fc^{\tlU_{-}}(\tlR)}
\sigma_{\tlR}^{\tlS}(H) 
\big(
\sum_{Q \in \fc_{\tlS}(M_{1})}
\varepsilon_{Q}^{U}
\hat \tau_{Q}(H - Y_{Q})
\big).
\]

En vertu du lemme \ref{lem:langlands} et du son analogue classique, on a alors:
\begin{equation}\label{eq:42starU}
\sum_{ P \in \fc_{\tlR}(M_{1})}
\varepsilon_{P}^{U}
\hat \tau_{P}(H - Y_{P}) = 
\sum_{\tlS \in  \fc^{\tlU_{-}}(\tlR)}
\varepsilon_{\tlR}^{\tlS}
\hat \sigma_{\tlR}^{\tlS}(H)
\mathrm{B}_{\tlS}^{U}(H, \yc_{\tlS}).
\end{equation}

Les résultats de la section \ref{S:combi} se généralisent formellement 
au cas du groupe unitaire, le seul résultat non-évident étant 
le lemme \ref{lem:artLem51U} ci-dessus.
On obtient alors l'analogue du lemme \ref{lem:gamSumGam0} suivant. 

\begin{lemme}\label{lem:gamSumGam0U} À un ensemble de mesure $0$ près, on a l'égalité des fonctions sur $\zgo_{\tlR}$ suivante:
\[
\mathrm{B}_{\tlR}^{U}(\cdot, \yc_{\tlR}) = 
 \sum_{P \in \fc_{\tlR}^{0}(M_{1})} 
 \Gamma_{P}'(\cdot, Y_{P}).
\]
\end{lemme}
\end{paragr}

\section{Une formule des traces infinitésimale}\label{sec:RTFinfU}

\subsection{Classes de conjugaison semi-simples}

\begin{paragr}[Notations.] --- On adopte les notations de \ref{ssec:stratU}. On suppose de plus que $F$ est un corps de nombres.  
\end{paragr}

\begin{paragr}[Lien avec le point de vue de \cite{leMoi}]--- Dans cette section, on va utiliser grandement les résultats de \cite{leMoi}. Comme nous adoptons un point de vue légèrement différent, expliquons brièvement comment faire le lien. On considère ici l'action de $U(V)$ sur $\tugo$. Dans \cite{leMoi}, on considère l'espace hermitien  $W = V \oplus Ee_{0}$ muni de la forme $\Phi\oplus^\perp 1$. Le groupe $U(V)$ est le sous-groupe de $U(W)$ qui fixe $e_0$. Il agit donc sur l'algèbre de Lie de $U(W)$. Le choix d'un élément non nul de $E$ de trace $0$ permet d'identifier de manière $U(V)$-équivariante $\tugo$ à un hyperplan de  $\Lie(U(W))$. De plus, l'élément $e_0$ détermine une droite dans  $\Lie(U(W))$ sur laquelle $U(V)$ agit trivialement et supplémentaire de $\tugo$ dans  $\Lie(U(W))$. On se contentera donc dans la suite de citer \cite{leMoi}.
\end{paragr}

\begin{paragr}[Compatibilité des classes et sous-groupes paraboliques] ---
On continue avec les notations du paragraphe précédent. 

\begin{lemme}\label{lem:LeviSSU} Soit $ M $ un facteur  de Levi d'un sous-groupe parabolique $ P  $ de $ U $. 
Pour tout $ X \in \tmgo(F) $,  on a $ X_{s} \in \tmgo(F) $. 
\end{lemme}

\begin{preuve} 
Soit $X=(A,b)$.  Avec les notations du paragraphe \ref{S:LeviU},  à  $ M $ est associée une décomposition 
orthogonale

\[ 
V= V_1^{\sharp}  \bigoplus^{\perp} 
\ldots \bigoplus^{\perp} V_r^{\sharp} \bigoplus^{\perp} V_0 
\]
où $ V_{i}^{\sharp} = V_{i} \oplus V_{-i} $ et $ b \in V_{0} $. 
On peut écrire $ A = \sum_{i = 0}^{r} A_{i} $ où  $ A_{i}$ est un endomorphisme auto-adjoint de $V_{i}^{\sharp}$ et $(A_{0}, b)  $ est un élément de $ \tugo_{V_{0}}(F) $.  Pour $1\leq i\leq r$, soit $A_i'$ la partie semi-simple de $A_i$ (au sens usuel). Soit $(A_0',b') \in  \tugo_{V_{0}}(F) $ la partie semi-simple de $(A_0,b)$. Alors $  (\sum_{i=0}^{r}A_{i}', b') $ est la partie semi-simple de $X$ donc appartient aussi à $\tmgo(F)$.
\end{preuve}
\end{paragr}

\begin{paragr}[L'ensemble $ \oc(\Phi) $.] --- 
Soit $ \oc = \oc(\Phi) $ l'ensemble de classes de $ U(F) $-conjugaison d'éléments semi-simples dans $ \tugo(F) $. 

En vertu du lemme \ref{lem:J1u}, l'application canonique $a: \tugo(F)\to \Ac(F)$ se factorise à travers l'application  $ U(F) $-invariante $ \ogo : \tugo(F) \rightarrow \oc $ qui, à $ X \in \tugo(F) $, associe la classe de conjugaison de sa partie semi-simple. On note encore $a:\oc \to \Ac(F)$ l'application qui s'en déduit.

Pour tout $ \ogo \in \oc $ et tout  sous-$F$-espace $ \thgo \subset\tugo $, soit 
$$ \thgo(F)_{\ogo} =\{X\in \thgo(F) \mid \of(X)=\of\}.$$ 
\end{paragr}

\begin{paragr}[Compatibilité des classes et sous-groupes paraboliques] ---\label{S:ogoLevi}

\begin{lemme}\label{lem:ogoLevi}
Soient $ P $ un sous-groupe parabolique 
de $ U $ de décomposition de Levi $ P = MN $ 
et $ X \in \tmgo(F) $. 
Pour tout $ Y \in \tngo(F) $ on a $ \ogo(X) = \ogo(X + Y) $.
\end{lemme}

\begin{preuve} On a une décomposition orthogonale $ V_{P}^{\perp} = V_{P} \oplus V_{0} $ (cf. § \ref{S:LeviU}) de sorte que $ \tmgo = \mgo \oplus V_{0} $. Soit $ X = (A,b) \in \tmgo(F)$ et $Y=(W,w)\in \tngo =\ngo(F) \oplus V_{P}$. On a $a(X)=a(X+Y)$. Soit $r\geq 0$ tel que $a(X)\in \Ac^{(r)}$. Alors $V^+$, le sous-espace de dimension $r$ engendré par $b,Ab,\ldots,A^{r-1}b$ est inclus dans $V_0$. Avec les notations du §\ref{S:somme-directeu}, on a $X=\iota_{V^+}(X^+,X^-)$. La famille $(b+w), (A+W)(b+w),\ldots, (A + W)^{r-1}(b + w)$ est une famille libre de $V_P\oplus V_0$.  Il existe donc $ n \in N(F) $ tel que $ n((A + W)^{i}(b + w)) = A^{i}b $ pour $ i = 0, \ldots, r-1 $. Il suffit clairement de prouver qu'on a $\ogo(X)=\ogo(n\cdot(X+Y))$. On a $a( n\cdot(X+Y))=a(X)$ (cf. \cite{leMoi} proposition 2.5) et   $n\cdot(X+Y)=\iota_{V^+}(Z^+,Z^-)$. D'après le lemme \ref{lem:iso-r}, $X^+$ et $Z^+$ sont des éléments semi-simples réguliers de $\tugo_{V^+}$ tels que $a(X^+)=a(Z^+)$ : ils sont donc conjugués sous $U(V^+)(F)$ (cf. corollaire \ref{cor:orbite-ss}). Pour conclure, il reste à voir que les parties semi-simples  $X^-_s$ et $Z^-_s$ sont $U(V^-)(F)$-conjugués. Soit $ P^{-} = U(V^{-}) \cap P $. Par construction, on a $X^-=(A^-, b^- )$ avec $A^-\in \mgo_{P^-}$ et $Z^--=(A_1^-,b_1^-)$ avec $A_1^-\in A^-+\ngo_{P^-}$. On a alors $X^-_s=(A^-_s,0)$ et $Z_s^-=(A_{1,s}^-,0)$ ; le fait que les parties semi-simples $A^-_s$ et $A_{1,s}^-$ sont conjugués résultent du corollaire 2.6 de \cite{PH1}. Cela conclut.
\end{preuve}

On a alors le corollaire du lemme \ref{lem:ogoLevi} 
immédiat suivant. 

\begin{corollaire}
On a:
\[ 
\tmgo(F)_{\ogo} + \tngo(F) = \tpgo(F)_{\ogo}.
\]
\end{corollaire}

\end{paragr}

\subsection{Distributions globales}\label{ssec:distGlobale}

\begin{paragr}[Choix auxiliaires.] ---\label{S:choixAuxU}
Soit $P_0$ un sous-groupe parabolique minimal  de $ U $ de décomposition de Levi $ M_{0}N_{0} $. 
On note $ \ago_{0}^{+}= \ago_{P_{0}}^{+} $. Soit  $K = \prod_{v} K _{v}$ de $U(\AAA)$ un sous-groupe compact maximal adapté à $M_{0}$ (cf. § \ref{S:HP}). On utilisera parfois un modèle de $U$ sur un anneau d'entiers \og hors $S$\fg{} pour un ensemble $S\subset \vc$ fini. On procèdera de la façon suivante : on considère une $E$-base de $V$. Presque pour tout $v\in V$, le  $\oc_{E\otimes_F F_v}$-réseau $V(\oc_v)$ dans $V\otimes_F F_v$ engendré par cette base est auto-dual et le groupe $U(\oc_v)$ est le stabilisateur de $V(\oc_v)$. Presque partout, on a $K_v=U(\oc_v)$. De même, $\tugo(\oc_v)$ est alors formé des couples $(A,b)$ où $A$ stabilise $V(\oc_v)$ et $b\in V(\oc_v)$.

 Pour presque toute place $v$ non-archimédienne, l'espace hermitien $V \otimes_{F} F_{v}$  admet un réseau autodual dont $K_{v}$ est le stabilisateur.  Pour $P\in \fc(P_0)$ on en déduit une application $H_P:U(\AAA)\to \ago_P^*$ (cf. §\ref{S:HP}).
\end{paragr}

\begin{paragr}[Noyaux paraboliques.] ---\label{S:noyauU}
Soit $ f \in \Sc(\tugo(\AAA)) $ et $ \ogo \in \oc $.

Soit $ P \in \fc(P_{0}) $ et $ MN $ sa décomposition de 
Levi standard. Pour tout $ g \in U(\AAA) $ soit:
\begin{equation}\label{eq:kU}
 k_{P,\ogo}(g)=  k_{P,\ogo}(f,g) = 
  \sum_{X\in \tmgo(F)_{\ogo}}  
  \int_{\tngo(\AAA)} f(g^{-1}\cdot (X+U))dU.
\end{equation}
\end{paragr}

\begin{paragr}[Noyau tronqué.] --- \label{S:noy-tronqueU}
Pour $ T \in \ago_{0} $ on définit le noyau tronqué:
\begin{equation}\label{eq:kTU}
k^T_{\ogo}(g) = k^T_{\ogo}(f,g)= 
\sum_{ P\in \fc(P_{0})} 
\eps_{P}^{U} 
\sum_{\delta\in P(F)\back U(F)} 
\hat{\tau}_{P}(H_{P}(\delta g)-T_{P}) \,  k_{P,\ogo}(f,\delta g).
\end{equation}
\end{paragr}

\begin{paragr}[Convergence d'intégrales.] ---\label{S:cvU}

\begin{theoreme}
  \label{thm:cvU} (\cite{leMoi} théorème 3.1, théorème 4.5)
  \begin{enumerate}
  \item Il existe un point $T_+ \in  \ago_{0}^+$ tel que pour tout $T \in T_++ \ago_{0}^+$, l'intégrale 
$$I^{U, T}_{\ogo}(f)= \int_{[U]} k^T_{\ogo}(f,g) \, dg
$$
converge absolument.
\item L'application $T \mapsto I^{U, T}_{\ogo}(f)$ est la restriction d'une fonction exponentielle-polynôme en $T$.
\item Le terme purement polynomial de cette exponentielle-polynôme est constant.
  \end{enumerate}
\end{theoreme}
\end{paragr}

\begin{paragr}[Distributions $I^{U}_{\ogo}$.] --- \label{S:IaU}On continue avec les notations des sections précédentes. 
 On note $I_{\ogo}(f) = I^U_{\ogo}(f)$ le terme constant de l'intégrale $I^{U,T}_{\ogo}(f)$ dans le théorème \ref{thm:cvU}. On obtient ainsi une distribution $I^{U}_{\ogo}$ qui vérifie les propriétés suivantes.

  \begin{theoreme} (\cite{leMoi}  § 4.5, théorème 4.7, théorème 5.1)
\label{thm:IU}
\begin{enumerate}
\item La distribution $I^{U}_{\ogo}$ ne dépend que du choix de la mesure de Haar sur $U(\AAA)$.
\item La distribution $I^{U}_{\ogo}$ est invariante au sens où pour tous 
$f\in\Sc(\tugo(\AAA))$ et $g\in U(\AAA)$, on a 
$$I^{U}_{\ogo}(f^g)= I^{U}_{\ogo}(f).
$$
\item Le support de la distribution $I_{\ogo}^{U}$ 
est contenu dans l'ensemble des $X = (X_{v})_{v}$ dans $\tugo(\AAA)$  
tels que pour tout $v$ la partie semi-simple de $ X_{v}$ est $U(F_{v})$-conjuguée à un élément de $\of$. 
\item Soit $ a \in \Ac(F) $ et $\oc_a=\{\of\in \oc \mid a(\of)=a\}$. Pour tout $f\in \Sc(\tugo(\AAA))$, la somme 
$$ I_a^U(f)=\sum_{\of\in \oc_a} I^{U}_{\ogo}(f)$$
converge absolument et le support de la  distribution $ I^{U}_{a} $ est inclus dans $\tugo_{a}(\AAA)$. 
En particulier $I_{a}^{U} = 0$ si $\tugo_{a}(F) = \emptyset$. 
\item Soit $S \subset \vc_\infty$ un ensemble de places archimédiennes de $ F $. Soit $f^S\in \Sc(\tugo(\AAA^S))$. La forme linaire sur $\Sc(\tugo(\AAA_S))$ donnée par 
$$f_S \mapsto I^{U}_{\ogo}(f_S\otimes f^S)
$$
est continue pour la topologie usuelle sur $\Sc(\tugo(\AAA_S))$.
\end{enumerate}
      \end{theoreme}

      \begin{remarque}
        Le lecteur attentif aura remarqué que le symbole $\of$ n'a pas la même signification que dans \cite{leMoi}. En fait, la preuve de \emph{loc. cit.} marche sans changement pour la définition qu'on utilise ici, le point clef étant l'utilisation du lemme \ref{lem:ogoLevi} en lieu et place de la proposition 2.5 de \emph{loc. cit.}. 
Les assertions concernant le support sont évidentes par construction. Enfin l'assertion de continuité découle de la démonstration   des théorèmes 3.1 et 3.2 de  \emph{loc. cit.}.
\end{remarque}
\end{paragr}

\subsection{Formule des traces infinitésimale}
\begin{paragr}[Transformation de Fourier partielle.] --- \label{S:TFPU} Considérons la forme bilinéaire symétrique, non dégénérée et $U$-invariante sur $\tugo$ donnée pour $X=(A,b)\in \tugo$
\begin{equation}
  \label{eq:bilinU}
  \bg X,X\bd= \trace(A^2)+2\Phi(b,b).
\end{equation}
Suivant les constructions du § \ref{S:TFPV}, pour tout sous-espace $\tugo_1\subset \tugo$ qui est $U$-invariant et non dégénéré pour $\bg \cdot,\cdot \bd$, on dispose de la transformée de Fourier partielle
\begin{equation}\label{eq:TFPU}
f\mapsto \hat{f}_{\tugo_1}
\end{equation}
qui est un automorphisme de $\Sc(\tugo(\AAA))$.  Les seuls espaces $\tugo_1$ possibles sont les suivants et leurs sommes directes 
\begin{itemize}
\item $F$ vu comme le sous-espace des homothéties auto-adjointes de $V$ ;
\item $\mathfrak{su}_F(V)$ le sous-espace des endomorphismes autoadjoints de $(V, \Phi)$ de trace nulle ;
\item $V$.
\end{itemize}

\begin{remarque}
  \label{rq:TFPU}
Dans la suite, les sous-espaces $\tugo_1$ les plus utiles seront $\mathfrak{su}_F(V)$, $V$ et $\mathfrak{su}_F(V) \oplus V$.
\end{remarque}
\end{paragr}

\begin{paragr}[Formule des traces infinitésimale.] --- \label{S:RTFinfU}Pour toute  transformation de Fourier partielle $f\mapsto \hat{f}$ (cf. §\ref{S:TFPU}), soit $D\mapsto \hat{D}$ la transformation duale au niveau des distributions. Notons que celle-ci préserve le fait d'être équivariante. 

Pour tout $\al \in F$, soit $\oc^{\al}$ l'ensemble des classes de conjugaison d'éléments semi-simples $(A,b)\in \tugo(F)$ tels que la trace de l'endomorphisme $A$ soit égale à $\al$.

  \begin{theoreme}(cf. \cite{leMoi} théorème 3.1,  théorème 5.1)
    \label{thm:RTFinfU}
    \begin{enumerate}
    \item Pour tout $f\in \Sc(\tugo(\AAA))$, la somme $\sum_{\ogo \in \oc} I_{\ogo}^{U}(f)$ converge absolument ce qui définit une distribution $\sum_{\ogo \in \oc} I_{\ogo}^{U}$.
    \item On a 
$$\sum_{\ogo \in \oc} I_{\ogo}^{U} =\sum_{\ogo \in \oc} \hat I_{\ogo}^{U}
$$
pour toute transformation de Fourier partielle au sens du §\ref{S:TFPU}.
\item Si la transformation de Fourier partielle est associée à l'un des trois espaces décrits dans la remarque \ref{rq:TFPU}, alors pour tout $\al \in F$, on a 
$$\sum_{\ogo \in \oc^{\al}} I_{\ogo}^{U}= \sum_{\ogo \in \oc^{\al}} \hat I_{\ogo}^{U}.
$$
    \end{enumerate}
  \end{theoreme}
  
  \begin{remarque}
     L'assertion $3$ n'apparaît pas explicitement dans \cite{leMoi} mais c'est  juste une variante de l'assertion $2$ que l'on obtient en considérant une formule sommatoire de Poisson limitée à un des  sous-espaces invariants de  $\mathfrak{su}_F(V)\oplus V$.
  \end{remarque}
\end{paragr}

\subsection{Généralisation}\label{ssec:produitU}

\begin{paragr}
On reprend les notations de la section \ref{ssec:transvu} ($F$ étant toujours un corps de nombres). On considère le  cas de  $\tugo^{\flat}$ muni de l'action du groupe  $U^{\flat}$. Ce cas est mixte en ce sens qu'il comprend des facteurs du type précédent (un groupe unitaire $U$ agissant sur un espace $\tugo$) mais aussi des facteurs linéaires ($GL(n)$ agissant sur $\tgl(n)$) qui ont été étudiés à la section \ref{sec:RTFinf}.
\end{paragr}

\begin{paragr}\label{S:produitU}
  Les constructions et les théorèmes des sections  ci-dessus,   ainsi que celles de la section \ref{sec:RTFinf} avec le caractère quadratique trivial,   se généralisent au cas de $U^{\flat}$ agissant sur $\tugo^{\flat}$. 
  Soit  $ \oc^{\flat} $ l'ensemble de classes de 
  $ U^{\flat}(F) $-conjugaison d'éléments semi-simples de  $ \tugo^{\flat}(F) $ (la notion de semi-simple est l'analogue de celle utilisée au sens des paragraphes \ref{ssec:Jordan}  et \ref{ssec:Jordanu}  ; elle correspond d'ailleurs à la notion d'être semi-simple composante par composante).
  Pour tout $ \ogo \in \oc^{\flat} $   on dispose donc d'une distribution $I^{U^{\flat}}_{\ogo}$ sur $\Sc(\tugo^{\flat}(\AAA))$ qui est invariante et à support dans   l'ensemble des $(X_{v})_{v} \in \tugo^{\flat}(\AAA)$   tels que, pour toute place  $v$, la partie semi-simple de $ X_{v} $ est $ U^{\flat}(F_{v}) $ conjuguée à un élément de $ \ogo $. On laisse au lecteur le soin de formuler  les analogues des théorèmes \ref{thm:IU} et \ref{thm:RTFinfU} dans ce cadre. 
\end{paragr}

\section{Le théorème de densité}\label{sec:densiteU} 


\subsection{Intégrales orbitales locales}\label{ssec:IOlocU}

\begin{paragr}
  On continue avec les notations de la section \ref{sec:RTFinfU}.
\end{paragr}

\begin{paragr}[Mesure de Haar.] --- \label{S:ensSU}  Soit $S$  un ensemble fini  de places de $F$. Soit $dg_S$ et $dg^S$ des mesures de Haar sur $U(\AAA_S)$ et $U(\AAA^S)$ de sorte que $dg=dg_S\otimes dg^S$.
\end{paragr}

\begin{paragr}[Intégrales orbitales locales.] --- Soit $f \in \Sc(\tugo(\AAA_S))$ et $a\in \Ac^{\rs}(\AAA_S)$. 
On introduit alors \emph{l'intégrale orbitale semi-simple régulière locale}
  \begin{equation}\label{eq:IOPloc2U}
    I_a(f) = 
    \begin{cases}
    0 \quad \text{ si } \tugo_{a}(\AAA_{S}) = \emptyset, \\
   \displaystyle \int_{U(\AAA_S)} f(g^{-1}\cdot X) \, dg_S, \quad \text{ pour } \, X \in \tugo_{a}(\AAA_{S}) \quad \text{ sinon}.
     \end{cases}
  \end{equation}

\begin{remarque}
  Cette construction dépend implicitement du choix de la mesure de Haar $dg_S$. Cependant, elle ne dépend pas du choix de $X\in \tugo_{a}(\AAA_{S}) $ car, lorsqu'il est non-vide, cet ensemble est une orbite pour l'action de $U(\AAA_S)$ (cf. corollaire \ref{cor:orbite-ss}).
\end{remarque}

\begin{remarque}
  \label{rq:decompositionU}
  Soit $S'\subset \vc$ un ensemble fini de places disjoint de $S$ et $S''=S\cup S'$. On suppose les mesures de Haar vérifient l'égalité $dg_S\otimes dg_{S'}=dg_{S''}$. Soit $a=(a_S,a_{S'})\in \Ac^{\rs}(\AAA_{S''})$.
   Pour toutes fonctions  $f _{S}\in \Sc(\tugo(\AAA_S))$ et $f _{S'}\in \Sc(\tugo(\AAA_{S'}))$, on a 
 $$I_a(f_S\otimes f_{S'})=I_{a_{S}}(f_{S})\cdot I_{a_{S'}}(f_{S'}).
$$
\end{remarque}
\end{paragr}

\subsection{Distributions stables}

\begin{paragr}[Fonctions instables.] --- \label{S:instU}Soit 
$$\Sc(\tugo(\AAA_S))_{0} \subset \Sc(\tugo(\AAA_S))$$
 le sous-espace des fonctions instables c'est-à-dire des fonctions $f\in  \Sc(\tugo(\AAA_S))$ telles que $I_a(f)=0$ pour tout $a\in \Ac^{\rs}(\AAA_S)$. Ce sous-espace ne dépend pas du choix de la mesure de Haar qui intervient dans la construction \eqref{eq:IOPloc2U}.
 \end{paragr}

\begin{paragr}
  Une distribution, c'est-à-dire une forme linéaire sur $\Sc(\tugo(\AAA_S))$, est dite stable si elle s'annule sur le sous-espace des fonctions instables.  Les distributions stables sont évidemment invariantes et on peut raisonnablement s'attendre à ce que l'inverse soit vrai. Faute de le savoir (cf. néanmoins \cite{Zhrang3} pour le cas lorsque $\dim(V)\leq 2$), on a du moins le théorème suivant.

  \begin{theoreme}
    \label{thm:densiteU}
Soit $\ogo \in \oc(\Phi, F)$ et $f^S\in \Sc(\tugo(\AAA^S))$. La distribution
$$f \in \Sc(\tugo(\AAA_S)) \mapsto I_{\ogo}(f\otimes f^S)$$
est stable.
  \end{theoreme}

La démonstration de ce théorème se trouve à la section \ref{ssec:demo-densiteU} et va utiliser les résultats des  deux sections qui suivent.
\end{paragr}

\subsection{Stabilité et transformée de Fourier}

\begin{paragr}[Stabilité et transformation de Fourier.] --- Comme dans le cas global, on définit des transformées de Fourier partielles relatives à des sous-espaces non dégénérées de $\tugo$ (cf. §§ \ref{S:TFPV} et  \ref{S:TFPU}). On en fixe une notée $f\mapsto \hat{f}$. On a le théorème dû à W. Zhang \cite{Z1} et H. Xue \cite{xue} dans le cas archimédien, dont la preuve est analogue à celle du théorème \ref{thm:TFstable}.

  \begin{theoreme}\label{thm:TFstableU}
    L'espace $\Sc(\tugo(\AAA_S))_0$ est invariant par la transformée de Fourier  $f\mapsto \hat{f}$.
  \end{theoreme}

\begin{corollaire}
  \label{cor:TFstableU}
Si $D$ est stable, il en est de même de $\hat{D}$.
\end{corollaire}
\end{paragr}

\subsection{Descente}\label{ssec:descenteU}

\begin{paragr}[Situation.] --- \label{S:situation-descU}

Soit $a\in \Ac(F)$. Le paragraphe \ref{S:dec-a} attache à $a$ un ensemble d'indices $I\cup J$ et pour $i\in I\cup J$ une extension $F_i$ de $F$, une $F_i$-algèbre quadratique $E_i$, un élément $\al_i\in F_i$.  Complétons ces données en des données comme aux §§\ref{S:V+U} et \ref{S:VijU}. On suppose qu'on a une décomposition
$$V=V^+\oplus (\oplus_{i\in I\cup J} V_i)
$$
et que la forme $\Phi$ égale à la forme $\Phi^+\oplus^\perp (\oplus_{i\in I\cup J}^\perp \Phi_{i,E})$. Suivant la section \ref{ssec:transvu}, on dispose alors des objets suivants :
\begin{itemize}
\item le sous-groupe $U^{\flat}$  de $U$ ;
\item l'espace $\tugo^{\flat}$ muni de l'action de $U^{\flat}$ dont le quotient catégorique est $\Ac^{\flat}$ ;
\item les ouverts denses $(\Ac^{\flat})' \subset \Ac^{\flat}$ et  $(\tugo^{\flat})' \subset \tugo^{\flat}$  ;
\item le  morphisme $U^{\flat}$-équivariant $\iota : (\tugo^{\flat})' \to \tugo$ qui induit un  morphisme étale $\iota : (\Ac^{\flat})' \to \Ac$ ;
\item un isomorphisme
\begin{equation}
  \label{eq:isocrucialUII}
   U \times^{U^\flat} (\tugo^\flat)' \to \tugo \times_{\Ac_\Phi} (\Ac^\flat)'
 \end{equation}
\end{itemize}
Soit  $X^{\flat} = (X^{+}, (X_{i})_{i \in I\cup J})\in \ugo^\flat(F)$ de caractéristique notée  $a^\flat\in (\Ac^\flat)'(F)$ tel que
\begin{itemize}
\item pour $ i \in I $, on a $ X_{i} = (\al_{i}, 0) $ vu comme élément de $\tugo_{i}$ ;
\item Pour $ j \in J $, on a $ X_{j} = (\al_{j}, 0, 0) $ vu comme élément de $\tggo_{j}$ ;
\item $X^+\in \tugo_+^{\rs}(F)$ est tel que $\iota(a^\flat)=a$.
\end{itemize}
Un tel $X^+$ existe. Soit $X=\iota(X^\flat)$. Alors $X$ est un élément semi-simple de $\tugo_a(F)$. Soit $\of$ la $U(F)$-orbite de $X$. Le corollaire \ref{cor:orbite-ss} implique que toute orbite semi-simple de  $\tugo_a(F)$ est atteinte par une telle construction. Soit $\of^\flat$ la $U^\flat(F)$-orbite de $X^\flat$.
\end{paragr}

\begin{paragr}\label{S:hypsurSU}
  On a fixé au §\ref{S:ensSU} un ensemble $S$ fini de places. Quitte à l'agrandir, on suppose vérifiées les propriétés supplémentaires du §\ref{S:UsurA}. On dispose donc d'un anneau $A\subset F$ et tous les objets viennent avec une structure sur $A$. Quitte à agrandir encore $S$, on va supposer  qu'on a  $a^{\flat}\in(\Ac^{\flat})'(A)$.
\end{paragr}

\begin{paragr}
  Soit $dg=dg_S\otimes dg^S$ une mesure de Haar sur $U^{\flat}(\AAA)=U^{\flat}(\AAA_S)\times U^{\flat}(\AAA^S)$ qui se décompose en mesures de Haar sur les facteurs.
\end{paragr}

\begin{paragr}\label{S:OmegaU} Soit $v\in S$. Le morphisme $\iota : (\Ac^{\flat})' \to \Ac$ s'identifie au morphisme $\iota_H: \Ac_H'\to \Ac$  utilisé au §\ref{S:Omega} (pour un $H$ évident qui dépend bien sûr des données fixées au §\ref{S:situation-descU}). Via cette identification, on dispose donc d'ouverts $\om_v \subset \Ac(F_v)$ et $\om^{\flat}_{v} \subset (\Ac^{\flat})'(F_v)$ qui correspondent respectivement aux ouverts notés $\om_v$ et $\om_{H,v}$ au §\ref{S:Omega}.

Soit $ \Omega^{\flat}_{v} \subset (\tugo^{\flat})'(F_{v})$ et $\Om_v\subset \tugo(F_v)$ les images réciproques respectives de $\om^{\flat}_{v} $ par le morphisme $ a $.
L'isomorphisme \eqref{eq:isocrucialUII} induit un difféomorphisme 
$$(U \times^{U^\flat} (\tugo^\flat)')(F_v) \to (\tugo \times_{\Ac_\Phi} (\Ac^\flat)')(F_v).$$
Soit $\Ker(H^1(F_v,U^\flat)\to H^1(F_v,U))$ l'ensemble des classes de $1$-cocyles galoisiens à valeurs dans $U^\flat$ dont l'image $U$ est la classe du cocyle trivial. L'espace $(U \times^{U^\flat} (\tugo^\flat)')(F_v)$ est alors une réunion disjointe d'ouverts fermés (cf. corollaire A.1.6 de \cite{AGhChandra}) indexés par l'ensemble fini $\Ker(H^1(F_v,U^\flat)\to H^1(F_v,U))$. Il en est donc de même de  $(\tugo \times_{\Ac_\Phi} (\Ac^\flat)')(F_v)$ et de son ouvert $\Om_v$. On écrit donc
\begin{equation}
  \label{eq:dec-Omphi}
  \Om_v=\bigcup_{ \varphi \in \Ker(H^1(F_v,U^\flat)\to H^1(F_v,U)) } \Om_{v}^{\varphi}.
\end{equation}
Les   $ \Om_{v}^{\varphi}$ ainsi obtenus  sont ouverts et $U(F_v)$-invariants. Dans cette décomposition, il existe un unique ouvert qui contient l'orbite $\of$ : c'est l'ouvert  $\Om_{v}^1$ correspondant à la classe $1$ du cocycle trivial. L'isomorphisme \eqref{eq:isocrucialUII} induit  un difféomorphisme
\begin{equation}\label{eq:fibTriv}
U(F_v) \times^{U^\flat(F_v)} \Om_v^\flat \to \Om_v^1
\end{equation}
induit par la submersion
\begin{equation}\label{eq:fibTrivII}
U(F_v) \times  \Om_v^\flat \to \Om_v^1
\end{equation}
donnée par $ (g, Y) \mapsto g \cdot \iota(Y) $.

\begin{lemme} \label{lem:dansOm1v}
  \begin{enumerate}
  \item Pour tout $a'\in \om_v^\flat$, on a 
$$\Om^\flat_v \cap \tugo^\flat_{a'}(F_v)\not=\emptyset$$
si et seulement si
$$\Om^1_v \cap \tugo_{\iota(a')}(F_v)\not=\emptyset.$$
\item  Pour tout $a\in \om_v\cap \Ac^{\rs}(F_v)$, la condition 
$$\Om^1_v \cap \tugo_{a}(F_v)\not=\emptyset$$
implique que pour tout $\varphi\not=1\in \Ker(H^1(F_v,U^\flat)\to H^1(F_v,U))$, on a 
$$\Om^\varphi_v \cap \tugo_{a}(F_v)=\emptyset.$$
\item Soit $Z\in \tugo(F_v)$. L'élément $Z$ appartient à $\Om_v^1$ si et seulement si sa partie semi-simple $Z_s$ appartient à   $\Om_v^1$.
\item Tout $Z\in \tugo(F_v)$ dont la partie semi-simple est $U(F_v)$-conjugué à un élément de $\of$ appartient en fait à $ \Omega^1_{v}$.
  \end{enumerate}

\end{lemme}

\begin{preuve} 
L'assertion $1$ découle de l'isomorphisme \eqref{eq:fibTriv} et du fait que $\iota$ induit une bijection de  $\om_v^\flat$ sur $\om_v$.
L'assertion $2$ résulte du fait que,  d'un part, pour  $a$ est régulier semi-simple tel que  $\tugo_{a}(F_v)$ soit non-vide,  $\tugo_{a}(F_v)$ est une $U(F_v)$-orbite (cf. corollaire \ref{cor:orbite-ss}) et  d'autre part,  les ouverts $\Om^\varphi_v $ sont $U(F_v)$-invariants et disjoints.

Prouvons l'assertion 3. Pour tout $Z\in \Om^1_v$ la partie semi-simple $Z_s$ appartient aussi à $\Om^1_v$. En effet, un tel $Z$ s'écrit $g\cdot\iota(Y)$ pour $g\in U(F_v)$ et $Y\in \Om_v^\flat$. D'après les définitions, on voit que $Z_s=g\cdot\iota(Y_s)$. D'après le lemme  \ref{lem:J1u}, on a  $a(Y_s)=a(Y)$ donc $a(Y_s)\in \om_v^\flat$ et  $Y_s\in \Om_v^\flat$. Il s'ensuit que $Z_s\in \Om^1_v$.
Mais le même énoncé vaut pour  $\Om^\varphi_v$  pour tout $\varphi\in H^1(F_v,U^\flat)$ (il suffit de remplacer $\tugo^\flat$ par une variante tordue par un cocycle représentant  $\varphi$). 
Donc si $Z_s\in \Omega^1_v$, on a $a(Z)=a(Z_s)\in \om_v$ donc $Z\in \Omega_v$. Il existe un unique $\varphi\in  H^1(F_v,U^\flat)$  tel que $Z\in \Omega^\varphi_v$. On a aussi $Z_s\in \Omega^\varphi_v$ donc $\varphi=1$.

Prouvons l'assertion 4. L'hypothèse implique que $Z_s\in  \Omega^1_v$. La conclusion résulte alors de l'assertion 3.
\end{preuve}
\end{paragr}

\begin{paragr}[Ouverts $\Omega^1$.] --- \label{S:Om1} Soit $\om=\prod_{v\in S } \om_{v} \subset \Ac(\AAA_S)$,   $\om^{\flat}=\prod_{v\in S } \om^{\flat}_{v} \subset (\Ac^{\flat})'(\AAA_S)$, 
$\Om^{\flat}=\prod_{v\in S} \Om^{\flat}_{v} \subset (\tugo^{\flat})'(\AAA_S)$ et $\Om^1=\prod_{v\in S} \Om_v^1 \subset \tugo(\AAA_S)$.
Posons aussi $\tugo(\AAA_{S})_{\ogo} = \prod_{v \in S}\tugo(F_{v})_{\ogo}$.

Comme au §\ref{S:Omega}, on dispose alors de deux applications linéaires surjectives,  continues sur les composantes archimédiennes, entre espaces de Schwartz (désignés par le symbole $\Sc$) :

\begin{itemize}
\item 
$$
\begin{array}{lll}
  \Sc(U(\AAA_S) \times \Om^{\flat})&\to &\Sc(\Om^1)\\
\al &\mapsto& f_\al
\end{array}
$$
où $f_\al$ est déterminé par l'égalité
\begin{equation}
  \label{eq:falphaU}
  f_\al(g\cdot \iota(Y))=\int_{U^{\flat}(\AAA_S)} \al(gh,h^{-1}\cdot Y) \, dh
\end{equation}
pour tout $g\in U(\AAA_S)$ et $Y\in \Om^{\flat}$ ;
\item 
$$
\begin{array}{lll}
\Sc(U(\AAA_S) \times \Om^{\flat})& \to &\Sc(\Om^{\flat})\\
\al&\mapsto& f^{\flat}_\al
\end{array}
$$
où $f^{\flat}_\al$ est définie par
\begin{equation}
  \label{eq:fHalphaU}
f^{\flat}_\al(Y)=\int_{U(\AAA_S)} \al(g, Y) \, dg
\end{equation}
pour tout $Y\in \Om^{\flat}$.
\end{itemize}
\end{paragr}

\begin{paragr}  On définit les intégrales orbitales locales $I_{a}^{U^{\flat}}$  pour tout $a \in \Ac^{\flat, \rs}(\AAA_S)$ comme dans \eqref{eq:IOPloc2U}. On a alors l'analogue du lemme \ref{lem:desc-locale}.

\begin{lemme}
   \label{lem:desc-localeU}
Soit  $\al \in \Sc(U(\AAA_S) \times \Om^{\flat})$, $a_1\in  \Ac^{\rs}(\AAA_S)$ et  $a_2\in \Ac^{\flat,\rs} (\AAA_S)$. 
\begin{enumerate}
\item Si $a_1\notin \om$,  l'intégrale orbitale locale de $f_\al$ est nulle
$$I_{a_1}^{U}(f_\al)=0.
$$
\item Si $a_2\notin \om^\flat$, l'intégrale orbitale locale de $f_\al^{\flat}$ est nulle
$$I_{a_2}^{U^\flat}(f_\al^{\flat})=0.
$$
\item Si $a_2\in  \om^\flat$ et $a_1=\iota(a_2)$, on a
  \begin{equation}
    \label{eq:desc-locU}
    I_{a_1}^{U}(f_\al)= I_{a_2}^{U^\flat}(f_\al^{\flat}).
  \end{equation}
\end{enumerate}
\end{lemme}

\begin{preuve}
  Les assertions $1$ et $2$ découlent des propriétés de support des fonctions. Attardons-nous sur l'assertion $3$. Le membre de gauche de \eqref{eq:desc-locU} est nul sauf si $\tugo_{a_1}(\AAA_S)\cap \Om^1\not=\emptyset$.  Cette dernière condition équivaut à $\tugo_{a_2}^\flat(\AAA_S)\cap \Om^\flat\not=\emptyset$ d'après le lemme \ref{lem:dansOm1v}. Et si cette dernière n'est pas satisfaite, alors le membre de droite de \eqref{eq:desc-locU} est nul. On peut alors supposer que les deux intersections sont non vides : l'assertion $3$ résulte alors de manipulations élémentaires.
\end{preuve}

\end{paragr}

\begin{paragr}[Descente des distributions globales.] --- Rappelons qu'on a défini des éléments $\of\in \oc$ et $\of^\flat\in \oc^{\flat}$ au §\ref{S:situation-descU}. On considère les distributions globales $I^{U}_{\ogo}$ et $I^{U^{\flat}}_{\ogo^{\flat}}$  relatives respectivement à l'action de $U$ sur $\tugo$ et $U^{\flat}$ sur $\tugo^{\flat}$ et définies aux §§\ref{S:IaU} et \ref{S:produitU}.

  \begin{theoreme}
    \label{thm:descU}
Pour tout $\al \in \Sc(U(\AAA_S) \times \Om^{\flat})$, on a
$$I^{U}_{\ogo}(f_\al\otimes \mathbf{1}_{\tugo(\oc^S)})= 
\frac{\vol(U(\oc^S))}{\vol(U^{\flat}(\oc^S))}\cdot I^{U^{\flat}}_{\ogo^{\flat}}(f^{\flat}_\al\otimes \mathbf{1}_{\tugo^\flat(\oc^S)}).$$
  \end{theoreme}
Ce théorème sera démontré à la fin de la section \ref{sec:descU}. 
\end{paragr}

\subsection{Démonstration du théorème \ref{thm:densiteU}}\label{ssec:demo-densiteU}

\begin{paragr}[Hypothèse de récurrence.] --- En utilisant le cas linéaire (théorème \ref{thm:densite}) et une hypothèse de récurrence sur les facteurs unitaires, on peut et on va supposer que le théorème  \ref{thm:densiteU} vaut pour l'action de $U^\flat$ sur $\tugo^\flat$ (notations du §\ref{S:situation-descU}) lorsque  $U^\flat$ est un sous-groupe propre de $U$.
  \end{paragr}

  \begin{paragr}
     Soit $ \ogo \in \oc$ et $a=a(\of)\in \Ac(F)$.
  \end{paragr}

\begin{paragr}[Ensemble $S$.] --- 
 Soit $f\in \Sc(\tugo(\AAA_S))_{0}$ et $f^S\in \Sc(\tugo(\AAA^S))$. En vertu de la remarque 
\ref{rq:decompositionU},   
  il suffit de prouver le théorème \ref{thm:densiteU} pour un ensemble fini $S$ de places aussi grand que l'on veut. Quitte à agrandir $S$, on peut et on va supposer que $S$ satisfait toutes les propriétés requises par la section \ref{ssec:descenteU}. On peut également supposer que $f^S$ égale   $\mathbf{1}_{\tugo(\oc^S)}$ la fonction caractéristique de $\tugo(\oc^S)$.
\end{paragr}

\begin{paragr}[Cas régulier semi-simple.] --- \label{S:casrssU} Supposons $\of$ est l'orbite d'un élément semi-simple régulier $X\in \tugo(F)$. On a donc  $a=a(\of) \in \Ac^{\rs}(F)$. D'après le corollaire \ref{cor:orbite-ss}, on a $\tugo_{a}(F)=\of$. Dans ce cas, on a en utilisant le lemme \ref{lem:integriteu} assertion 3
  \begin{eqnarray*}
    I_\of^U(f\otimes \mathbf{1}_{\tugo(\oc^S)})&=&\int_{U(\AAA)}   (f\otimes \mathbf{1}_{\tugo(\oc^S)})(g^{-1}X)\, dg \\
&=& \vol(U(\oc^S),dg^S)\cdot I_{a_S}^U(f).
  \end{eqnarray*}
Le théorème \ref{thm:densiteU} est donc évident pour de telles distributions.
\end{paragr}

\begin{paragr}[Cas de descentes.] ---\label{S:cas-descU} On suppose qu'on a $a\notin  \Ac^{\rs}(F)$ et que l'orbite $\of$ est obtenue à partir des données considérées au §\ref{S:situation-descU} (comme c'est toujours possible, cf. corollaire \ref{cor:orbite-ss}). On dispose donc du groupe $U^{\flat}$, de l'espace  $\tugo^{\flat}$ etc. On suppose que $U^\flat$ est  un sous-groupe propre de $U$. On suit les notations de la section \ref{ssec:descenteU}. On dispose d'un ouvert $\om \subset \Ac(\AAA_S)$ et  d'un ouvert $U(\AAA_s)$-invariant $\Om^1\subset \tugo(\AAA_S)$ (cf. §\ref{S:Om1}).

Soit $\theta \in \Cc(\om)$ qui vaut $1$ dans un voisinage de $a_S$. Soit $\mathbf{1}_{\Om^1}$ la fonction caractéristique de $\Om^1$. Rappelons que $\Om^1$ est un ouvert-fermé de $a^{-1}(\om)$. La fonction
$$f'= \mathbf{1}_{\Om^1} \cdot (\theta\circ a)\cdot  f$$
appartient à l'espace  $\Sc(\Om^1)\cap  \Sc(\tugo(\AAA_S))_{0}$. 

Tout élément de $\tugo(\AAA_S)$ de partie semi-simple qui, en chaque place, est conjugué à un élément de $\of$ doit appartenir à $\Om^1$ (cf. assertion 4 du lemme \ref{lem:dansOm1v}). Il résulte alors de  la propriété de support de $I^U_\of$ (assertion 3 du théorème \ref{thm:IU}) qu'on a l'égalité
$$I_{\ogo}^U(f \otimes  \mathbf{1}_{\tugo(\oc^S)}  )= I_{\ogo}^U(f' \otimes  \mathbf{1}_{\tugo(\oc^S)}  ).$$
On est donc ramené à prouver la nullité du membre de droite. Mais $f'=f_\al$ pour $\al\in \Sc(U(\AAA_S)\times \Omega^{\flat})$. Par ailleurs, on a par le  théorème \ref{thm:descU}
\begin{equation}
  \label{eq:ofofb}
  I^{U}_{\ogo}(f_\al\otimes \mathbf{1}_{\tugo(\oc^S)})= 
\frac{\vol(U(\oc^S))}{\vol(U^{\flat}(\oc^S))}\cdot I^{U^{\flat}}_{\ogo^{\flat}}(f^{\flat}_\al\otimes \mathbf{1}_{\tugo^\flat(\oc^S)}).
\end{equation}
Le lemme \ref{lem:desc-localeU} implique facilement qu'on a $f^{\flat}_\al\in \Sc(\tugo^\flat(\AAA_S))_{0}$ (le raisonnement est le même qu'au §\ref{S:cas-desc}). On conclut par l'hypothèse de récurrence que le membre de droite dans l'égalité \eqref{eq:ofofb}est nul.
\end{paragr}

\begin{paragr}[Fin de la démonstration.] --- Soit $\ogo \in \oc$ qui échappe aux deux cas précédents. Alors $\of$ est la classe de conjugaison d'un élément $(\al,0)$ où $\al\in F$ est vu comme une  homothétie de $V$. Quitte à effectuer une translation par $\al$ sur la fonction $f$, on peut supposer que $\al=0$. Notons $0$ l'image de l'élément nul dans $\Ac(F)$. Par le corollaire \ref{cor:orbite-ss}, on a 
$$I^U_{\of}=I^U_0.$$
La fin de la preuve du théorème \ref{thm:densiteU} est alors analogue à celle du théorème \ref{thm:densite}. Esquissons les points cruciaux
\begin{itemize}
\item En utilisant la formule des traces donnée 
par le théorème \ref{thm:RTFinfU} 3., l'hypothèse de récurrence et le corollaire 
\ref{cor:TFstableU} 
on arrive à l'égalité 
\begin{equation}
  \label{eq:RTF0U}
I_{0}^{U}(f\otimes f^S) = \hat{I}_{0}^{U}(f\otimes f^S).
\end{equation} 
pour tout $f^S\in \Cc(\tugo (\AAA^S))$ et $ \hat I $ est l'une des transformées de Fourier de $I$ considérées au §\ref{S:TFPU}.
\item 
Soit $u \notin S$ une  place non-archimédienne, décomposée dans $ E $ et  $f^{S'}\in \Sc(\tugo(\AAA^{S'}))$ où $S'=S\cup \{u\}$. En utilisant l'isomorphisme $\tugo(F_u)\simeq \tgl(n)(F_u)$, on dispose d'une distribution sur $\Cc( \tgl(n)(F_u))$
$$f_u \mapsto I_{0}^{U}(f\otimes f_u \otimes f^{S'})$$
qui vérifie les trois assertions du théorème \ref{thm:incertitude}  et qui est donc nulle ce qui permet de conclure ; les propriétés requises de support résulte des propriétés de support de $I_0^U$, cf. théorème \ref{thm:IU} assertion 4, et de l'égalité \eqref{eq:RTF0U} (pour plus de détails, on pourra se reporter à la fin du §\ref{S:findemodensite}).
\end{itemize}
\end{paragr}

\section{Descente globale} \label{sec:descU}

\subsection{Un noyau auxiliaire}

\begin{paragr}
 On reprend les notations de §\ref{ssec:distGlobale}. Dans cette section, on donne une première étape dans la preuve du  théorème \ref{thm:descU} : on montre qu'on peut remplacer dans la construction des distributions globales (cf. section \ref{sec:RTFinfU}) le noyau $k^T_{\ogo}$   (cf. \eqref{eq:kU}) par un noyau plus adapté à la descente.
\end{paragr}

\begin{paragr}
Soit $P \in \fc(M_{0})$. Soit $P = MN$ sa décomposition de Levi associée. 
Soit $\tmgo$ l'espace associé (cf. §\ref{S:LeviU}).  Soit $X \in \tmgo(F)$ et $X_s$ sa partie semi-simple (cf. §\ref{S:Jor1u}). On a alors $X_s\in \tmgo(F)$, en vertu du lemme \ref{lem:LeviSSU}.
\end{paragr}

\begin{paragr}[Centralisateurs.] ---\label{S:centraU}  On écrit ensuite $X=\iota_X(X^+,X^-)$ avec $X^+\in \tugo_{+}^{(r)}(F)$ et $X^-\in \tugo_{-}^{(0)}(F)$. La partie semi-simple $X_s^-$ est donc de la forme $(Y,0) \in  \tugo_{-}(F)$. Soit $\ugo_{-}(V^-,X_s^-)$ la sous-algèbre des $Z\in \ugo_{-}$ qui vérifient $[Z,Y]=0$. On identifie d'ailleurs $\ugo_{-}$ à la sous-algèbre de $\ugo(V)$ des endomorphismes qui laissent stable $V^-$ et sont nuls sur $V^+$. 

Soit $N_{X_s}$ le centralisateur de $X_s$ dans $N$ et $\tngo(X_s)$ le sous-espace de $\tugo_{-}= \ugo_{-} \oplus  V^- $ défini par
$$\tngo(X_s)=[\ngo\cap \ugo_{-}(V^-,X_s^-)]\oplus V_{P}.$$

Observons que l'analogue du lemme \ref{lem:ss-U} est vrai dans notre contexte (et s'énonce mot par mot de même façon). 
\end{paragr}

\begin{paragr}[Un nouveau noyau tronqué.] --- Soit $f\in \Sc(\tugo(\AAA))$,  
$\ogo \in \oc(\Phi)$ et $P \in \fc(M_{0})$. On reprend les notations des §\ref{S:noyauU} et §\ref{S:noy-tronqueU}. On introduit la variante suivante de \eqref{eq:kU}
  \begin{equation}
    \label{eq:j-U}
    j_{P,\ogo}(f,g)=\sum_{X \in \tmgo_{P}(F)_{\ogo}} \sum_{\nu \in N_{X_s}(F)\back N(F)} \int_{\tngo(X_s,\AAA)} f((\nu g)^{-1}\cdot \iota_X(X^+,X^-+U))dU. 
  \end{equation}

On a également la variante suivante de \eqref{eq:kTU} pour $T\in \ago_{0}$
\begin{equation}
  \label{eq:jT-U}
  j_{\ogo}^T(f,g)=\sum_{P \in \fc(P_{0})} \eps_{P}^{U} \sum_{\delta\in P(F)\back U(F)} \hat{\tau}_{P}(H_{P}(\delta g)-T_{P})j_{P,\ogo}(f,\delta g).
\end{equation}

\end{paragr}

\begin{paragr}[Théorème de comparaison.] --- On a l'analogue suivant 
du théorème \ref{thm:cvj}.

  \begin{theoreme}\label{thm:cvjU}
Soit $\ogo \in \oc$ et $f\in \Cc(\tugo(\AAA))$. Il existe $T_f$ tel que pour tout $T\in T_f+\ago^+$,  on a 
\begin{equation*}
  \int_{[U]} |j^T_{\ogo}(f,g) |\, dg <\infty, \quad 
  \text{ et }  \quad  \int_{[U]} j^T_{\ogo}(f,g)\, dg = \int_{[U]} k^T_{\ogo}(f,g) \, dg.
\end{equation*}

\end{theoreme}

\begin{preuve}
La preuve suit les pas de la preuve du théorème \ref{thm:cvj}. 
On démontre alors d'abord l'analogue évident du lemme \ref{lem:resommation} ensuite 
on utilise les méthodes de \cite{leMoi}, théorème 3.1 
pour effectuer les réductions comme dans la preuve de \ref{thm:cvj}.
Puisque on utilise la troncature d'Arthur, la théorie de réduction 
est classique, en particulier 
la fonction $ \chi_{P_{1}, P_{2}}^{T}(\cdot) $ c'est juste 
$ F^{P_{1}}(\cdot, T)\sigma_{P_{1}}^{P_{2}}(H_{P_{0}}(x) - T) $ 
dans la notation de \cite{ar1}, preuve du théorème 7.1. 
Ensuite on procède comme dans la preuve du théorème \ref{thm:cvj}.
\end{preuve}
\end{paragr}

\subsection{Combinatoire de la descente} 

\begin{paragr}[Sous-groupes $M_1$ et $P_1$.] --- On reprend les notations et les hypothèses de la section \ref{ssec:transvu} 
ainsi que celles du §\ref{ssec:desc-combiU}. En particulier cf. §\ref{S:M1U}, 
on dispose d'un sous-groupe de Levi $M_1$ de $U$. Quitte à faire agir $U(F)$, on peut et on a va supposer que $M_1$ contient le sous-groupe de Levi minimal fixé $M_{0}$ et qu'il est le facteur de Levi semi-standard d'un sous-groupe parabolique $P_1$ contenant $P_{0}$ (qui est un sous-groupe parabolique minimal de $ U $ fixé auparavant).
\end{paragr}

\begin{paragr}[Sous-groupe $M_0^{-}$ et $P_0^{-}$.] --- \label{S:parU} D'après les notations et la construction du §\ref{S:QQ-U}, on dispose d'une application $Q \mapsto \tlQ^-$ de $\fc^{U}(M_1)$ dans $\fc^{\tlU_{-}}(\tlM_0^{-})$. Le groupe $U_- = \prod_{i \in I}U_{i} \times \prod_{j \in J}G_{j}$ 
est naturellement un sous-groupe de $\tlU_-$. 
Soient alors
$M_0^{-} := \tlM_{0}^{-} \cap U_-$, 
$ \tlP_{0}^{-} \in \fc^{\tlU_{-}}(\tlM_0^{-})$ le groupe associé à 
$ P_{1} \in  \fc^{U}(M_1)$ et 
$P_0^{-}=\tlP_0^{-} \cap U_- \in \fc^{U_{-}}(M_{0}^{-})$. Alors $P_0^{-}$ est un $F$-sous-groupe parabolique  minimal de  $U_-$.
Soit $$\fc^{U}_{P_0^{-}}(M_1)$$
 l'ensemble des $Q\in \fc^{U}(M_1)$ tels que $P_0^{-}\subset \tlQ^-$.
\end{paragr}

\begin{paragr}[Un élément semi-simple $X$.] ---  
On reprend les notations du §\ref{S:situation-descU} 
et on suppose que les données du paragraphe 
\ref{S:parU} sont celles de ce paragraphe. 
On fixe alors un $X \in \tugo(F)_{\ogo}$ semi-simple et 
$X^{\flat} = (X^{+}, X^{-})$, où $X^{-} = ((X_{i})_{i \in I}, (X_{j})_{j \in J}) \in \tugo_{-}(F)$ 
et $\iota^{\flat}(X^{\flat}) = X$. 

On a aussi $a^{\flat} = a(X^{\flat})$ ainsi que $\ogo^{\flat} \in \oc^{\flat}$ - l'unique 
classe de $U^{\flat}(F)$-conjugaison semi-simple dans $\tugo^{\flat}(F)$ telle que $a(\ogo^{\flat}) = a^{\flat}$. 
\end{paragr}

\begin{paragr} Soit $T\in \ago_{0}^+$.  Pour tout $g\in U(\AAA)$ et $\tlR \in \fc^{\tlU_-}(R_0)$, soit
  \begin{equation}
    \label{eq:UpsilonU}
     \Upsilon^T_{\tlR}(g)= \sum_{P \in \fc^{U}_{\tlR}(M_1)}   \eps_{P}^{U} \,  \hat{\tau}_{P}(H_{P}(g)-T_{P}),
   \end{equation}
   où  $\fc^{U}_{\tlR}(M_1)$ est défini au §\ref{S:combiU}, et soit $\nc_{\tmgo_{\tlR}}$ le cône nilpotent de $\tmgo_{\tlR}$ c'est-à-dire le fermé de $\tmgo_{\tlR}$ obtenu par intersection avec la fibre en $0$ de l'application $\tugo_- \to \Ac_{-}$. Soit  $f\in \Sc(\tugo(\AAA))$.  On a introduit en \eqref{eq:tljT} un noyau  $\tlj_{\ogo}^{T}(f,g)$ pour $g\in U(\AAA)$. Le lemme suivant en fournit une expression alternative.

\begin{lemme}\label{lem:formjU} 
Pour tout $g\in U(\AAA)$,  on a 
  \begin{eqnarray*}
\tlj^T_{\ogo}(f,g)=\sum_{\tlR \in \fc^{\tlU_-}(P_0^{-})} \sum_{\delta\in R(F)\back U(F)} \Upsilon^T_{\tlR}(\delta g) \sum_{Y\in \nc_{\tmgo_{\tlR}}(F)} \, \int_{\tilde{\ngo}_{\tlR}(\AAA)} f((\delta g)^{-1}\cdot \iota^{\flat}(X^+,X^-+Y+U))dU
\end{eqnarray*}
où $R=\tlR \cap U_-$. 
\end{lemme}

\begin{preuve}
La preuve est analogue à la preuve du lemme \ref{lem:formj}. 
\end{preuve}
\end{paragr}

\subsection{Démonstration du théorème \ref{thm:descU}}

\begin{paragr}
  On continue avec les notations de la section précédente. 
On pose $\tlU^{\flat} = U_{+} \times \tlU_{-}$.   
  On fixe un sous-groupe 
parabolique minimal   
   $P_{0}^+$ de $U_+$ 
ainsi que sa partie de Levi $M_{0}^{+}$. 
Ceci défini la chambre positive dans $\ago_{P_{0}^+\times \tlP_{0}^{-}}^{\tlU^{\flat}}$    
et pour tout $T'$ dans ce chambre
on en déduit des points $T'_{\tlS}$ pour tout sous-groupe parabolique $\tlS \in \fc(M_{0}^{+} \times \tlM_{0}^{-})$ par action du groupe de Weyl associé. 
Pour $\tlR \in \fc^{\tlU_{-}}(\tlM_{0}^{-})$ on note $\hat T'_{\tlR}$ 
la projection orthogonale de $T_{U_{+} \times \tlR}$ à $\zgo_{\tlM_{0}^{-}} = \ago_{M_{1}}$. 
\end{paragr}

\begin{paragr}   Soit $K_{\tlU^{\flat}}\subset \tlU^{\flat}(\AAA)$, resp. $K_{U^{\flat}}\subset U^{\flat}(\AAA)$, un sous-groupe compact maximal de 
$\tlU^{\flat}(\AAA)$, resp. 
$U^{\flat}(\AAA)$,  \og en bonne position \fg{} par rapport au sous-groupe de Levi   $M_{0}^+ \times \tlM_0^{-}$, resp. $M_{0}^+ \times M_0^{-}$. 
Pour tout $x\in U^{\flat}(\AAA)$ et tout sous-groupe parabolique  $R$ défini sur $F$ de $U^{\flat}$, soit 
$k_R(x)\in K_{U^{\flat}}$ un élément  tel que $x\in R(\AAA) k_{R}(x)$.

Le choix de $K_{\tlU^{\flat}}$ fait qu'on dispose d'une application de  Harish-Chandra 
$$ H_{\tlU^{\flat} \times \tlR} : \tlU^{\flat}(\AAA) \rightarrow \ago_{ \tlU^{\flat} \times \tlR}.$$
  Soit $ \hat H_{\tlR} $ l'application composée de  $H_{\tlU^{\flat} \times \tlR}$ avec la projection orthogonale 
  $ \ago_{\tlU^{\flat} \times \tlR} \rightarrow \zgo_{\tlR}$.

Soit $\tlR \in \fc^{\tlU_{-}}(\tlM_{0}^{-})$ et notons $R = U_{-} \cap \tlR \in 
\fc^{U_{-}}(M_{0}^{-})$.
On introduit 
la famille orthogonale-positive
$\yc_{\tlR}^{T, T'}(x,y) = (Y_{P}^{T, T'}(x,y))_{P \in \fc^{0}_{\tlR}(M_{1})}$
où 
$$
Y^{T,T'}_{P}(x,y)= -H_{P}(k_{R}(x)y) + T_{P} - \hat T'_{\tlR} \in \ago_{P} = \zgo_{\tlR}.
$$

On a alors l'analogue suivant du lemme \ref{lem:Upsilon} qui se démontre de même façon que celui-ci. 
\begin{lemme}
  \label{lem:UpsilonU} Avec les notations ci-dessus, on a l'égalité suivante où le membre de droite est défini en \eqref{eq:UpsilonU}
\begin{equation*}
   \Upsilon^T_{\tlR}(xy) = \sum_{\tlS \in \fc^{\tlU_{-}}(\tlR)}
\varepsilon_{\tlR}^{\tlS}
\hat \sigma_{\tlR}^{\tlS}(\hat H_{\tlR}( x) - \hat T'_{\tlR}  )
\mathrm{B}_{\tlS}^{\tlG}(\hat  H_{\tlS}( x) - \hat T'_{\tlS}, \yc_{\tlS}^{T,T'}(x,y)).
\end{equation*}
\end{lemme}
\end{paragr}

On a défini tout pour finir la preuve du théorème \ref{thm:descU}, qui a partir de maintenant va être la même que celle du théorème  \ref{thm:desc}. 
On se limite donc de résumer son schéma:
\begin{enumerate}
\item On suppose d'abord que, avec les notations du §\ref{S:OmegaU}, on a 
$\al = \beta_{S} \otimes \phi_{S}$ avec $\beta_{S} \in C_{c}^{\infty}(U(\AAA_{S}))$ 
et $\phi_{S} \in C_{c}^{\infty}(\Omega^{\flat})$. On lui 
associe les fonctions $f_{\al}$ et  $f_{\al}^{\flat}$ par les formules \ref{eq:falphaU} 
et \ref{eq:fHalphaU}. 
\item On pose 
$$
f = f_{\al} \otimes \mathbf{1}_{\tugo(\oc^{S})} \in C_{c}^{\infty}(\tugo(\AAA)), \quad
\beta = \beta_{S} \otimes \mathbf{1}_{U(\oc^{S})} \in C_{c}^{\infty}(U(\AAA)), \quad
\phi = \phi_{S} \otimes \mathbf{1}_{\tugo^{\flat}(\oc^{S})} \in C_{c}^{\infty}(\tugo^{\flat}(\AAA)).
$$
\item En utilisant les lemmes \ref{lem:integriteu}, 
\ref{lem:formjU} et \ref{lem:UpsilonU}, on montre l'analogue de la proposition \ref{prop:dvpT-T'} suivant. 
Pour tout $T \in \ago_{0}^{+}$ assez positif par rapport à support de $f$ on a que l'intégrale
  $$
\int_{[U]} \tlj_{\ogo}^T(f,g) \, dg 
$$
est égale à la somme sur $\tlR_1 \in \fc^{\tlU_{-}}(P_{0}^{-})$ de 
\begin{multline}\label{eq:descFin3U}
  \int\limits_{U(\AAA)}
\int\limits_{R_1(\rmF) \backslash U^{\flat}(\AAA)} 
\mathrm{B}_{\tlR_1}^{U}(\hat H_{\tlR_1}(x) - \hat T'_{\tlR_1}, \yc_{\tlR_1}^{T,T'}(x,y))
\sum_{\tlR \in \fc^{\tlR_1}(P_{0}^{-}) }\varepsilon_{\tlR}^{\tlR_1}
\sum_{\delta \in R(\rmF) \backslash R_1(\rmF)}
\hat \sigma_{\tlR}^{\tlR_1}(\hat H_{\tlR}(\delta x) - \hat T_{\tlR}') \\
\sum_{Y\in \nc_{\tmgo_{\tlR}}(F)} \, \int_{\tngo_{\tlR}(\AAA)}
 \beta(y^{-1}) \phi((\delta x)^{-1} \cdot \iota^{\flat}(X^+, X^-+Y+U))dU dxdy.
\end{multline}
\item On démontre que l'expression \eqref{eq:descFin3U} est absolument convergente et 
est la restriction d'un polynôme-exponentielle 
en $(T, T')$ dont la partie purement-polynômiale est
\begin{itemize}
\item nulle, si $\tlR_{1} \neq \tlU_{-}$,
\item égale à $\frac{\vol(U(\oc^S))}{\vol(U^{\flat}(\oc^S))}\cdot I^{U^{\flat}}_{\ogo^{\flat}}(f^{\flat}_\al\otimes \mathbf{1}_{\tugo(\oc^S)})$
si $\tlR_{1} = \tlU_{-}$.
\end{itemize}
\item Le point précédant démontre le théorème \ref{thm:descU} dans le cas où $\al$ est comme ci-dessus. 
\item On conclut la preuve par un argument de densité comme dans §\ref{S:demodethmdesc}.
\end{enumerate}

\part{Le transfert singulier global}\label{partie:3}

\section{Le cas infinitésimal}\label{sec:lecasinf}

\subsection{Factorisation de distribution : le cas linéaire}\label{ssec:thlin}

\begin{paragr}\label{S:fac-not}
  Soit $E/F$ une extension quadratique de corps de nombres. Soit $\eta$ le caractère quadratique de $\AAA_F^\times$ associé à $E/F$ par la théorie du corps de classes. Soit $n\geq 1$ et $V=F^n$. 
   Soit $\tggo=\tgl_F(V)$ muni de l'action du groupe $G=GL_F(V)$. Soit $\Ac$ le quotient catégorique. On reprend les notations et les constructions des sections \ref{sec:RTFinf} et \ref{sec:densite}.
\end{paragr}

\begin{paragr}[Mesures de Haar.] ---\label{S:Haarlin}
  On fixe pour tout $v\in V$ une mesure de Haar sur $G(F_v)$. On suppose que pour tout $v$ fini en dehors d'un ensemble fini de places, la mesure donne le volume $1$ au sous-groupe compact ouvert $G(\oc_v)=GL_n(\oc_v)$. On en déduit une mesure de Haar sur $G(\AAA)$.
\end{paragr}

\begin{paragr}
  Pour toute place $v$ de $F$, on dispose d'un espace de Bruhat-Schwartz $\Sc(\tggo(F_v))$ (cf. §\ref{S:Schw}) et d'un sous-espace  $\Sc(\tggo(F_v))_{\eta}$ de fonctions instables (cf. \ref{S:inst}). Soit $\eta_v$ la composante en $v$ du caractère $\eta$ et 
$$\Ic(\eta_v)=  \Sc(\tggo(F_v))/ \Sc(\tggo(F_v))_{\eta_v}$$
le quotient. Pour tout ensemble fini $S$ de places de $F$, soit 
$$\Ic(\eta)_S=\otimes_{v\in S} \Ic(\eta_v)$$
le produit tensoriel. Pour tous  ensembles finis $S\subset S'$ de places de $F$, contenant les places archimédiennes, on a une application 
$\Ic(\eta)_S \to \Ic(\eta)_{S'}$ donnée par $f \mapsto f\otimes (\otimes_{v\in S'\setminus S} \mathbf{1}_v)$ où, pour toute place non-archimédienne de $F$,  on note $\mathbf{1}_{v}$ l'image de la fonction caractéristique de $\tggo(\oc_v)=\tgl_{\oc_v}(\oc_v^n)$.

Soit 
$$\Ic(\eta)= \varinjlim_{S}  \Ic(\eta)_S$$
\end{paragr}

\begin{paragr} En utilisant la base canonique de $V=F^n$, on construit un facteur $\tilde{\eta}$   (cf. §\ref{S:eta}). Le choix de ce facteur et des mesures de Haar permet de définir  des intégrales orbitales locales semi-simples régulières (cf. § \ref{S:IOPloc}). Pour tout $v\in \vc$, on a  une application linéaire injective
  \begin{equation}
    \label{eq:Ieta}
    I^{\tilde{\eta}}_v: \Ic(\eta_v) \to  \mathcal{C}^\infty( \Ac^{\rs}(F_v))
  \end{equation}
  donnée par 
$$I^{\tilde{\eta}}(f)_v:  a \mapsto I_{a,v}^{\tilde{\eta}}(f).$$
 Pour tout ensemble fini $S\subset \vc$,  on définit, plus généralement, une application linéaire injective
 \begin{equation}
    \label{eq:IetaS}
    I^{\tilde{\eta}}_S: \Ic(\eta)_S \to  \mathcal{C}^\infty( \Ac^{\rs}(\AAA_S)).
  \end{equation}
Soit
 $$\mathcal{C}^\infty(\Ac^{\rs})=\varinjlim_S  \mathcal{C}^\infty( \Ac^{\rs}(\AAA_S))$$
où la limite est prise sur les ensembles finis $S\subset\vc$ et les applications de transition sont données pour $S\subset S'$ (et $S$ assez grand contenant les places archimédiennes) par
$$\phi \mapsto \phi  \prod_{v\in S' \setminus S}   I_v^{\tilde{\eta}}(\mathbf{1}_{v}).$$
Les applications $ I^{\tilde{\eta}}_S$ induisent alors une application injective
$$\Orb^{\eta}: \Ic(\eta) \to \mathcal{C}^\infty(\Ac^{\rs}).$$
Celle-ci est construite à l'aide du facteur $\tilde{\eta}$ (qui dépend du choix d'une base de $V$) mais il est facile de voir que l'application $\Orb^{\eta}$ ne dépend que de $\eta$.
\end{paragr}

\begin{paragr}[Distributions $I_a^\eta$.]  ---  On a construit à la section \ref{sec:RTFinf} des distributions $I_a^\eta$ pour $a\in \Ac(F)$. D'après le théorème \ref{thm:densite}, on peut voir ces distributions comme des formes linéaires sur $\Ic(\eta)$ qu'on note encore $I_a^\eta$.
\end{paragr}

\subsection{Factorisation de distributions : le cas hermitien}\label{ssec:thuni}

\begin{paragr}
  Pour tout groupoïde $\xc$, soit $|\xc|$ l'ensemble des classes d'isomorphismes.
\end{paragr}

\begin{paragr}\label{S:hc}
On continue avec les notations du §\ref{S:fac-not}.   Soit $\sigma$ le générateur de $\Gal(E/F)$. Soit $\hc=\hc_{n,E/F}$ le groupoïde quotient de l'ensemble des formes $\sigma$-hermitiennes non-dégénérées sur $V_E=V\otimes_F E$ par l'action de $GL_E(V_E)$. Les classes d'isomorphismes de $\hc$ sont les classes d'équivalences de  formes $\sigma$-hermitiennes non-dégénérées sur $V_E$. Pour tout $\Phi\in \hc$, on dispose du groupe unitaire $U_\Phi$,  de l'espace $\tugo_\Phi$ sur lequel agit $U_\Phi$, d'un morphisme  $\tugo_\Phi\to \Ac$  qui identifie $\Ac$ au quotient catégorique de $\tugo_\Phi$ par  $U_\Phi$  (cf. §§\ref{S:cadreu}, \ref{S:quotientu})
\end{paragr}

\begin{paragr}\label{S:reseau}
  Pour toute place $v$ de $F$, soit $F_v$ le complété de $F$ et $E_v=F_v\otimes_F E$. Soit  $\sigma_v=1\otimes \sigma$. Soit $\hc_v$ le groupoïde quotient de l'ensemble des forme $\sigma_v$-hermitiennes non-dégénérées sur $V_{E_v}=V\otimes_F E_v$ par l'action de $GL_{E_v}(V_{E_v})$. Comme précédemment, on associe à $\Phi\in \hc_v$, les objets $U_\Phi, \tugo_{\Phi}$ définis sur $F_v$.
On fixe sur $U_\Phi(F_v)$ une mesure de Haar. On exige les conditions de compatibilité suivantes :

\begin{itemize}
\item Pour tout isomorphisme entre  $\Phi$ et $\Phi'$, on demande que les mesures de Haar se correspondent via l'isomorphisme
$$U_\Phi \to U_{\Phi'}$$
qui s'en déduit. Un tel isomorphisme est bien défini à conjugaison près par un élément de  $U_\Phi(F_v)$ ;
\item si $v$ est fini, non ramifié dans $E$ et que  $V_{E_v}$ possède un $\oc_{E_v}$-réseau autodual pour $\Phi$, on normalise la mesure de Haar en imposant qu'elle donne le volume $1$ au sous-groupe qui stabilise ce réseau.
\end{itemize}
Il est possible de satisfaire ces conditions : les automorphismes de $\Phi$ agissent par automorphismes intérieurs sur $U_\Phi(F_v)$ et préservent donc toute mesure de Haar ; dans le second cas, ils agissent  transitivement sur l'ensemble des réseaux autoduaux. Toujours dans le second cas, on fixe et on note $V(\oc_{E_v})$ un tel réseau. On dispose alors d'un $\oc_{E_v}$-réseau $\tugo(\oc_F) \subset \tugo(F)$ formé des couples $(A,b)$ où $A$ est un endomorphisme qui stabilise le réseau $V(\oc_{E_v})$ et $b\in V(\oc_{E_v})$. 

Si $v$ est décomposé dans $E$,  l'ensemble $|\hc_v|$ est réduit à un singleton. Pour $\Phi\in \hc_v$,  on identifie l'action de $U_{\Phi}$ sur $\tugo_{\Phi}$ à celle de $GL_{F_v}(n)$ sur $\tgl_{F_v}(n)$.

Supposons que  $v$ n'est pas décomposé dans $E$.  L'ensemble  $|\hc_v|$  est  fini et il possède $n+1$ éléments ou $2$ éléments selon que $v$ est archimédienne ou non. Supposons de plus que $v$ est non-archimédienne et inerte dans $E$. Les deux objets de $\hc_v$, disons $\Phi_0$ et $\Phi_1$ (à isomorphisme près), sont  distingués par le fait que le discriminant $\disc(\Phi_0)$ est une norme alors que $\disc(\Phi_1)$ n'en est pas une, ou encore qu'il existe un $\oc_{E_v}$-réseau auto-dual dans $V_{E_v}$ pour $\Phi_0$  alors qu'il n'y en a pas pour $\Phi_1$.

Notons enfin que l'extension des scalaires induit un morphisme de localisation $\hc \to \hc_v$. 
\end{paragr}

\begin{paragr}[Espace $\Ic(\hc_v)$.] ---
Soit $v\in V$ et $\Phi\in \hc_v$. On note $U$, $\ugo,\ldots $ au lieu de $U_\Phi$, $\tugo_\Phi$ etc. On a défini (cf. §\ref{S:Schw})  l'espace de Schwartz-Bruhat $\Sc(\tugo(F_v))$. Puisqu'on dispose d'une mesure de Haar sur $U(F_v)$, on a une application linéaire
$$\Sc(\tugo(F_v)) \to \mathcal{C}^\infty(\Ac^{\rs})$$
donnée par les intégrales orbitales locales semi-simples régulières (cf. section \ref{ssec:IOlocU}). Elle se factorise par le quotient
$$\Ic(\Phi)=\Sc(\tugo)/\Sc(\tugo)_0.
$$
(cf. §\ref{S:instU}) en une application injective
\begin{equation}
  \label{eq:IOlocPhi}
  I^{\Phi}_a : \Ic(\Phi) \to  \mathcal{C}^\infty(\Ac^{\rs}).
\end{equation}

\begin{lemme}\label{lem:12.2.4.1}
Soit $\Phi$ et $\Phi'$ deux objets de $\hc_v$ isomorphes. Les espaces $\Ic(\Phi)$ et $\Ic(\Phi')$ sont canoniquement isomorphes.
\end{lemme}

\begin{preuve}
  Il suffit d'observer que $U_\Phi$ agit trivialement sur  $\Ic(\Phi)$.
\end{preuve}

Soit
$$\Ic(\hc_v)=  \bigoplus_{\Phi\in |\hc_v|} \Ic(\Phi).
$$
Cet espace est muni de projections $f\mapsto f_\Phi\in \Ic(\Phi)$ pour tout $\Phi\in \hc_v$.  Lorsque $v\in \vc$ est finie et non-ramifié dans $E$, on définit un élément $\mathbf{1}_v\in \Ic(\hc_v)$ de la manière suivante.

\begin{enumerate}
\item Soit $v$ est inerte dans $E$. On a alors $\Ic(\hc_v)=  \Ic(\Phi_0)\oplus \Ic(\Phi_1)$ (cf. §\ref{S:reseau}).  Soit $\mathbf{1}_v$ l'image dans  $\Ic(\hc_v)$ du couple $(\mathbf{1}_{\tugo(\oc_F)},0)$. Cet élément est indépendant du choix du réseau $V(\oc_{E_v})$.
\item Soit $v$ est décomposé dans $E$, $|\hc_v|$ est réduit à un singleton $\{\Phi\}$ et $\Ic(\hc_v)=  \Ic(\Phi)$. Soit  $\mathbf{1}_v$ l'image dans  $\Ic(\hc_v)$ de $\mathbf{1}_{\tugo(\oc_F)}$. Cet élément est indépendant du choix du réseau $V(\oc_{E_v})$.
\end{enumerate}
\end{paragr}

\begin{paragr}
On définit  une application linéaire
  \begin{equation}
    \label{eq:Iv}
    I_v:\Ic(\hc_v) \to \mathcal{C}^\infty(\Ac^{\rs}(F_v))
  \end{equation}
  induite par 
$$f  \mapsto  \left(a\mapsto \sum_{ \Phi\in |\hc_v|  } I^{\Phi}_a(f_\Phi)\right),$$
où  $I^{\Phi}_a$ est l'application \og intégrale orbitale \fg{} définie en \eqref{eq:IOlocPhi}.
\end{paragr}

\begin{paragr}\label{S:IS}
Soit $S$ un ensemble fini de places de $F$, suffisamment grand pour contenir les places archimédiennes et les places non-archimédiennes ramifiées dans $F$. Soit $\hc_S$ le produit des $\hc_v$ sur $v\in S$. On pose alors
$$\Ic(\hc_S)=\otimes_{v\in S}\Ic(\hc_v).
$$
Il est muni d'une projection $f\mapsto f_\Phi$
\begin{equation}
  \label{eq:proj-Ic}
   \Ic(\hc_S) \to \otimes_{v\in S} \Ic(\Phi_v).
 \end{equation}
pour tout $\Phi=(\Phi_v)_{v\in S}\in \hc_S$. 

L'application \eqref{eq:Iv} s'étend en une application linéaire injective
$$I_S : \Ic(\hc_S) \to  \mathcal{C}^\infty( \Ac^{\rs}(\AAA_S)).$$
On définit alors
$$\Ic(\hc)=\varinjlim_{S} \Ic(\hc_S)
$$
où les applications de transition $\Ic_S \to \Ic_{S'}$ sont données par $f \mapsto f\otimes (\otimes_{v\in S'\setminus S} \mathbf{1}_v)$.
\end{paragr}

\begin{paragr}[Distributions $I^{\hc}$.] ---\label{S:Ihc}
  Soit $\Phi \in \hc$ et $\tugo, U$, etc. les objets qui lui sont attachés. Soit $$S_\Phi\subset \vc$$ l'ensemble fini des places archimédiennes, des places ramifiées dans $E$ et des places tels que le discriminant $\disc(\Phi_v)$ ne soit pas une norme de $E_v/F_v$.  Le groupe $U(\AAA)$ est muni de la mesure de Haar produit des mesures fixées sur les $U(F_v)$.
Soit $\of\in \oc(\Phi)$ une classe de conjugaison semi-simple dans $\tugo(F)$. On dispose d'une distribution (cf. §\ref{S:IaU})
$$I_\of^\Phi= I_\of^{U_\Phi}$$ 
sur l'espace de Schwartz $\Sc(\tugo(\AAA))$.

 Supposons que $S$ contienne $S_\Phi$ et que $U$ et $\tugo$ soient définis sur l'anneau des entiers hors $S$. D'après le théorème \ref{thm:densiteU}, la  forme linéaire linéaire 
$$f\in \otimes_{v\in S}\Sc(\tugo(F_v))\mapsto I_\of^{\Phi}(f\otimes \mathbf{1}_{\tugo(\oc^S)})   $$
se factorise par  $\otimes_{v\in S} \Ic(\Phi_v)$. Par composition avec \eqref{eq:proj-Ic}, on obtient une forme linéaire sur  $\Ic(\hc_S)$. Cette forme ne dépend pas du choix du modèle entier (cf. §\ref{S:choixAuxU}) toujours par le théorème \ref{thm:densiteU}. Dans la suite, on note simplement 
$$\mathbf{1}^S= \mathbf{1}_{\tugo(\oc^S)}$$
sans préciser le modèle entier choisi. 

 Par une nouvelle application du théorème \ref{thm:densiteU}, on obtient une forme linéaire
$$I_\of^{\Phi}: \Ic(\hc) \to \CC.$$
Celle-ci ne dépend que de la classe de $\Phi$ (plus exactement une équivalence entre des  formes $\Phi$ et $\Phi'$ permet d'identifier l'orbite $\of$ à une orbite $\of'$ associée à $\Phi'$ et les formes linéaires $I_\of^{\Phi}$ et $I_{\of'}^{\Phi'}$ ainsi obtenues sont égales).

\begin{lemme}\label{lem:Sf}
  Soit $f\in \Ic(\hc)$. Il existe un ensemble fini $S_f \subset \vc$ tel que si $S_\Phi\not\subset S_f$, on a 
$$I_\of^{\Phi}(f)=0$$
pour tout $\of\in \oc(\Phi)$.
\end{lemme}

\begin{preuve}
  Soit $f\in \Ic(\hc)$. Soit $S\subset\vc$ un ensemble fini tel que 
  \begin{enumerate}
  \item $S$ contienne les places archimédiennes et les places de $F$ ramifiées dans $E$ ;
  \item $f\in \Ic(\hc)$ est représentée par $f_S\in \Ic(\hc_S)$.
  \end{enumerate}

Supposons que $S_\Phi\not\subset S$. Soit $S'=S\cup S_\Phi$ et $f_{S'}$ l'image de $f_S$ dans $\Ic(\hc_{S'})$. L'image de $f_{S'}$ par l'application (cf. \eqref{eq:proj-Ic}) 
$$
  \Ic(\hc_{S'}) \to \otimes_{s\in S'} \Ic_{(V_v,\Phi_v)}
$$
est nulle. On a donc pour tout $\of\in \oc(\Phi)$
\begin{eqnarray*}
  I_\of^\Phi(f)&=& I_\of^\Phi(f_S \otimes \mathbf{1}^S )\\ &=& I_\of^\Phi(f_{S'} \otimes \mathbf{1}^{S'} ) \\
&=&0
\end{eqnarray*}
\end{preuve}

On définit pour tout $a\in \Ac(F)$ et $f\in \Ic(\hc)$
$$I_a^\Phi(f)=\sum_{\of\in \oc(\Phi)_a} I_\of^\Phi(f),$$
cf. assertion 4 du théorème  \ref{thm:IU}, la somme ci-dessus étant \emph{a priori} infinie mais du moins absolument convergente.

On définit alors pour tout $a\in \Ac(F)$ et $f\in \Ic(\hc)$
\begin{equation}
  \label{eq:Iahc}
  I_a^{\hc}(f)= \sum_{\Phi\in |\hc|  } I_a^{\Phi}(f).
\end{equation}
Cela fait sens par le lemme suivant.

\begin{lemme}\label{lem:finitude}
Pour tout $f$, il existe un ensemble fini $H_f\subset |\hc|$ tel que pour tous $a\in \Ac(F)$ et $\Phi\in |\hc|$, on a 
  $$I_a^{\Phi}(f)=0$$
sauf si $\Phi\in H_f$. En particulier, dans la somme \eqref{eq:Iahc}, il n'y a qu'un nombre fini de termes non nuls.
\end{lemme}

\begin{preuve}
  D'après le lemme \ref{lem:Sf}, il existe un ensemble fini $S_f$ contenant les places non-archimédiennes tel que seules les formes hermitiennes $\Phi$ vérifiant $S_\Phi\subset S_f$ contribuent effectivement. Par conséquent, la classe de telles formes  $\Phi$ est déterminée hors $S_f$. Mais il n'y a alors qu'un nombre fini $H_f$ de telles formes $\Phi$ à équivalence près.
\end{preuve}
\end{paragr}

\subsection{Le théorème de transfert}

\begin{paragr}
  On poursuit avec les notations des deux sections précédentes.
\end{paragr}

\begin{paragr}[Correspondance locale.] --- \label{S:corresp}On dit que des éléments $f\in \Ic(\eta_v)$ et $f'\in \Ic(\hc_v)$ se correspondent si on l'égalité
\begin{equation}\label{eq:corr-I}
  I^{\tilde{\eta}}_v(f)    = I^{}_v(f').
\end{equation}

Lorsque $v$ est fini ou scindé dans $E$, cette correspondance induit en fait un isomorphisme 
\begin{equation}
  \label{eq:isom-Ic}
\Ic(\eta_v) \simeq  \Ic(\hc_v).
\end{equation}
C'est tautologique si $v$ est scindé dans $E$ et c'est la traduction de l'énoncé du transfert non-archimédien pour $v$ fini non scindé (cf. \cite{Z1} théorème 2.6). Il en serait  de même pour  une place archimédienne mais non-scindée si le transfert archimédien était connu (pour des résultats partiels, cf. \cite{xue}). 
\end{paragr}

\begin{paragr}[Le lemme fondamental.] ---\label{S:LF}
  Pour $v\in \vc$ fini,  non-ramifié dans $E$ et également en dehors d'un ensemble fini fixé de \og mauvaises\fg{} places inertes, l'isomorphisme \eqref{eq:isom-Ic} envoie l'élément $\mathbf{1}_v\in \Ic(\eta_v)$ sur l'élément  $\mathbf{1}_v\in  \Ic(\hc_v)$. Lorsque $v$ est décomposé, cet énoncé est tautologique alors que pour $v$ inerte c'est le lemme fondamental conjecturé par Jacquet-Rallis \cite{jacqrall} et démontré par Yun et Gordon \cite{yun} lorsque la caractéristique résiduelle est assez grande.

  \begin{remarque}
    Le lecteur attentif ne manquera pas d'observer  que le facteur de transfert utilisé par Yun ne correspond pas à celui que nous utilisons. En fait, les deux définitions diffèrent par $\eta_v(d_n(X))$ (notations du §\ref{S:Areg}). Mais lorsque $d_n(X)$ est de valuation impaire, l'intégrale orbitale sur le groupe linéaire est nulle et ce quelque soit le facteur de transfert.
  \end{remarque}
\end{paragr}

\begin{paragr}
  Le lemme fondamental implique que les applications $I_S$ du §\ref{S:IS} induisent une application linéaire injective
$$\Orb: \Ic(\hc) \to \mathcal{C}^\infty(\Ac^{\rs}).$$
  On dit alors que deux éléments $f\in \Ic(\eta)$ et $f'\in \Ic(\hc)$ se correspondent si 
$$\Orb^{\eta}(f)=\Orb(f').$$

On peut alors énoncer le théorème principal de cette section. 

\begin{theoreme}\label{thm:transfert}
  Soit $f\in \Ic(\eta)$ et $f'\in \Ic(\hc)$ qui se correspondent. On a alors
$$I_a^\eta(f)=  I_a^{\hc}(f')
$$
pour tout $a\in \Ac(F)$.
\end{theoreme}

On écarte le cas immédiat de $a\in \Ac^{\rs}(F)$ au paragraphe suivant. La démonstration des autres cas du théorème se trouve aux section \ref{ssec:tr-desc} et \ref{ssec:fin-tr}.
\end{paragr}

\begin{paragr}[Cas $a\in \Ac^{\rs}(F)$.] --- \label{S:areg} Dans ce cas, il existe une unique forme  $\Phi\in \hc$ (à isomorphisme près) telle que $\tugo_{\Phi,_a}(\AAA)\not=\emptyset$ et dans ce cas on a aussi $\of=\tugo_{\Phi,a}(F)\not=\emptyset$ est une orbite semi-simple. Soit $U=U_\Phi$. On a donc
$$I_a^{\hc}(f')=  I_\of^{\Phi}(f').
$$
Soit $S$ un ensemble fini de places contenant les places archimédiennes et les place ramifiées tel que
\begin{itemize}
\item Pour tout $v\notin S$, le discriminant de $\Phi_v$ est une norme ;
\item $f$ et $f'$ sont représentés par respectivement par $f_S\in  \Ic(\eta_S)$ et $f'_S \in \Ic(\hc_S)$ ;
\item $a\in \Ac^{\rs}(\oc_v)$ pour $v\notin S$.
\end{itemize}
On a alors
$$I_a^{\eta}(f)= I^{\tilde{\eta}}_{a,S}(f_S) $$
et
$$  I_\of^{\Phi}(f')=  I_{a,S}(f'_S) $$
et l'égalité à démontrer résulte directement du fait que $f$ et $f'$ se correspondent.
\end{paragr}

\subsection{Cas de descente}\label{ssec:tr-desc}

\begin{paragr} \label{S:a_0} Soit $a_0\in \Ac^{(r)}(F)$ avec $0\leq r<n$. On reprend les notations de la section \ref{ssec:orbss}. On dispose donc d'une décomposition associée à $a_0$
$$V=V^+ \oplus (\oplus_{i\in I} V_i) \oplus (\oplus_{j\in J}V_j)$$
où $V^+$ est un $F$-espace et $V_i$ est un $F_i$-espace pour tout $i\in I\cup J$. 
L'invariant $a_0$ définit une collection $(a^+, (a_i)_{i\in I\cup J})$ où $a^+\in \Ac_{V^+}^{\rs}(F)$ et $a_i\in \Ac_{V_i}^{(0)}(F)$  pour tout $i\in I\cup J$.  On dispose d'un élément  $\al_i\in F_i$ pour tout $i\in I\cup J$.
\end{paragr}

\begin{paragr}
  Comme à la section \ref{ssec:transv}, on dispose du groupe 
$$H=GL_F(V^+)\times \prod_{i\in I\cup J} GL_{F_i}(V_i)$$
agissant sur $\thgo=\tgl_F(V^+)\oplus (\oplus \tgl_{F_i}(V_i))$. Soit $\Ac_H$ le quotient catégorique. Par restriction des scalaires, on regarde $H$ comme un sous-$F$-groupe de $G$. Dans toute cette section, on suppose que $H$ est un sous-groupe propre de $G$. Le caractère $\eta\circ \det$ induit donc un caractère sur $H(\AAA)$ qui est trivial sur chaque facteur    $GL_{F_i}(V_i)(\AAA)$ pour $i\in J$. 
On fixe sur $V^+$, resp. sur $V_i$ pour tout $i\in I\cup J$, une $F$-base, resp. une $F_i$-base. Cela permet de définir un facteur $\tilde{\eta}_H$ (cf. §\ref{S:desc-locale}) et d'identifier chaque facteur de $H$ à un certain $GL(m)$ sur un corps $F$ ou $F_i$. On suit alors les préconisations du §\ref{S:Haarlin} pour le choix d'une mesure de Haar sur $H(\AAA)$ et ses composantes locales.

Comme dans la section \ref{ssec:thlin}, on définit dans ce contexte les objets suivants pour tout $S\subset \vc$ fini :
\begin{itemize}
\item un espace $\Ic^H(\eta_v)=  \Sc(\thgo(F_v))/ \Sc(\thgo(F_v))_{\eta_v}$ ;
\item des espaces $\Ic^H(\eta)_S=\otimes_{v\in S} \Ic^H(\eta_v)$ et $\Ic^H(\eta)= \varinjlim_{S}  \Ic^H(\eta)_S$ ;
\item une application intégrale orbitale $I_S^{\tilde{\eta}_H}: \Ic^H(\eta)_S \to \mathcal{C}^\infty( \Ac_H^{\rs}(\AAA_S))$
\item une application injective 
$$\Orb^\eta: \Ic^H(\eta) \to \mathcal{C}^\infty(\Ac_H^{\rs})=\varinjlim_S  \mathcal{C}^\infty( \Ac^{\rs}_H(\AAA_S)).$$
\item des distributions pour tout $a\in \Ac_H(F)$ (cf. section \ref{ssec:produit})
$$I_a^{H,\eta} : \Ic^H(\eta) \to \CC.$$
\end{itemize}
\end{paragr}

\begin{paragr} On note par un indice $E$ les $E$-espaces vectoriels obtenus par extension des scalaires d'un $F$-espace vectoriel. Par exemple, pour $i\in I\cup J$, on  a $V_{i,E}=V_i \otimes_F E$ qui  est aussi un $E_i$-module pour $E_i=F_i\otimes_F E$.

Soit $\hc^\flat$ la catégorie dont les objets sont les familles $(\Phi^+,(\Phi_i)_{i\in I})$ où $\Phi^+$, resp.   $\Phi_i$, est une forme $\sigma$-hermitienne, resp.   $\sigma_i$-hermitienne,  non-dégénérée sur $V^+_E$, resp sur le $E_i$-espace vectoriel $V_{i,E}$. Les morphismes sont donnés par l'équivalence des formes, composante par composante.  

Pour tout $\Phi^\flat =(\Phi^+,(\Phi_i)_{i\in I})    \in \hc^\flat$, on dispose d'un groupe
$$U_{\Phi^\flat}= U_{\Phi^+} \times \prod_{i\in I} U_{\Phi_i} \times \prod_{j\in J} GL_{F_j} (V_j).
$$
qui agit sur l'espace
$$\tugo_{\Phi^\flat}= \tugo_{\Phi^+} \oplus  (\oplus_{i\in I} \tugo_{\Phi_i} )\oplus  \oplus_{j\in J} \tgl_{F_j} (V_j)$$
qu'on voit, par restriction des scalaires, comme des objets  définis sur $F$. Le quotient catégorique est l'espace $\Ac_H$ introduit plus haut.
Fixons également pour $j\in J$ une  forme  $\sigma_j$-hermitienne  non-dégénérée, $\Phi_j$, sur le $E_j$-module $V_{j,E}$. On dispose alors d'un morphisme 
$$\hc^\flat \to \hc$$
qui consiste à envoyer $\Phi^\flat =(\Phi^+,(\Phi_i)_{i\in I})    \in \hc^\flat$ sur la forme $\Phi$ sur $V_E$ donnée par (cf. notations du §\ref{S:VijU})
\begin{equation}
  \label{eq:laformePhi}
  \Phi=\Phi^+\oplus^\perp (\oplus^\perp_{i\in I\cup J} \Phi_{i,E}).
\end{equation}
Au moyen de ce choix supplémentaire, on identifie $U_{\Phi^\flat}$ à un sous-groupe de $U_\Phi$. Suivant les préconisations des §§\ref{S:Haarlin} et \ref{S:reseau}, on dispose de mesures de Haar sur $U_{\Phi^\flat}(\AAA)$ et ses composantes locales.

Bien sûr, pour tout $v\in \vc$, on définit de  manière évidente un pendant local $\hc_v^\flat$ à $\hc^\flat$, des objets $U_{\Phi^\flat}$ et $\tugo_{\Phi^\flat}$ pour $\Phi^\flat\in\hc_v^\flat$, un morphisme de localisation etc.

Pour tout  $\Phi^\flat\in \hc^\flat$, soit $\oc(\Phi^\flat)$ l'ensemble des orbites semi-simples de $U_{\Phi^\flat}(F)$ dans $\tugo_{\Phi^\flat}(F)$. Pour tout $a\in \Ac_H(F)$, soit  $\oc(\Phi^\flat)_a\subset \oc(\Phi^\flat)$ le sous-ensemble des orbites d'invariant $a$.

Comme précédemment, on introduit les objets suivants pour   tout $v\in \vc$ et tout $S\subset \vc$ fini :
\begin{itemize}
\item pour tout $\Phi^\flat \in \hc^\flat_v$,  un espace $\Ic(\Phi^\flat)=  \Sc(\tugo_{\Phi^\flat}(F_v))/ \Sc(\tugo_{\Phi^\flat}(F_v))_{0}$ ;
\item pour tout $v\in \vc$ un espace $$\Ic(\hc_v^\flat)=\oplus_{ \Phi ^\flat\in \hc_v ^\flat}\Ic(\Phi^\flat)$$
\item des espaces $\Ic(\hc^\flat)_S=\otimes_{v\in S} \Ic(\hc^\flat_v)$ et $\Ic(\hc^\flat)= \varinjlim_{S}  \Ic(\hc^\flat)_S$ ;
\item des applications \og intégrale orbitale \fg{} $I_S: \Ic(\hc^\flat)_S \to \mathcal{C}^\infty( \Ac_H^{\rs}(\AAA_S))$
\item d'une application injective 
$$\Orb: \Ic( \hc^\flat) \to \mathcal{C}^\infty(\Ac_H^{\rs})$$
\item pour tout $\Phi^\flat\in \hc^\flat$ et tout $\of\in \oc(\Phi^\flat)$ une  distribution (cf. §\ref{S:produitU})
$$I^{\Phi^\flat}_\of :  \Ic( \hc^\flat)  \to \CC$$
\item pour tout $a\in \Ac_H(F)$ une distribution
  $$I_a^{\flat} : \Ic( \hc^\flat)  \to \CC$$
définie par 
$$I_a^\flat= \sum_{\Phi^\flat\in |\hc^\flat|} \sum_{  \of\in \oc(\Phi^\flat)_a  }  I^{\Phi^\flat}_\of .$$
\end{itemize}
\end{paragr}

\begin{paragr}[Hypothèse de récurrence.] --- En raisonnant par récurrence sur la dimension de $G$, on peut et on va supposer que le théorème \ref{thm:transfert} vaut pour le sous-groupe propre $H$. Par conséquent, pour tous  $f\in \Ic^H(\eta)$ et $f'\in \Ic(\hc^\flat)$ telles que 
$$\Orb^\eta(f)=\Orb(f')$$
on a 
$$I_a^{H,\eta}(f)=  I_a^{\flat}(f')
$$
pour tout $a\in \Ac_H(F)$.
  
\end{paragr}

\begin{paragr}[Réécriture de $I_{a_0}^{\hc}$.] ---\label{S:reecriture}
L'invariant $a_0$ fixé au §\ref{S:a_0} correspond à une famille 
$$(a^+, (a_i)_{i\in I\cup J})\in \Ac_H.$$
Comme $a^+\in \Ac_{V^+}(F)$, il existe une unique forme hermitienne non-dégénérée $\Phi^+$ sur $V^+_E$ (à équivalence près) telle que $\tugo_{\Phi^+,a^+}(F)\not=\emptyset$. On fixe une telle forme $\Phi^+$ ainsi qu'un élément $X^+=(A^+,b^+)\in \tugo_{\Phi^+,a^+}(F)$. Soit $\hc_I$ le groupoïde quotient de l'ensemble des familles  $(\Phi_i)_{i\in I}$ de  formes hermitiennes non-dégénérées $\Phi_i$  sur $V_i$ par l'action de $\prod_{i\in I} GL_{E_i}(V_{i,E})$. Si $I=\emptyset$, le groupoïde $\hc_I$ par convention est le quotient d'un singleton par l'action du groupe trivial. On dispose de deux morphismes $\hc_I\mapsto \hc^\flat$ et $\hc_I\mapsto \hc$. Le premier est donné par l'application qui à $\Phi=(\Phi_i)_{i\in I} \in\hc_I$ associe
  $$\Phi^\flat=(\Phi^+,(\Phi_i)_{i \in I})\in \hc^\flat.$$ 
Le second est induit par l'application qui au même $\Phi$ associe la forme $\tilde{\Phi}$  sur $V_E$ obtenue par l'égalité (cf. notations du §\ref{S:VijU})
\begin{equation}
  \label{eq:laformePhi2}
  \tilde{\Phi}=\Phi^+\oplus^\perp (\oplus^\perp_{i\in I\cup J} \Phi_{i,E}).
\end{equation}
Introduisons les éléments $X_0=(A^+\oplus (\oplus_{i\in I\cup J} \al_i),b^+)$ et $Y_0=(X^+,(\al_i,0)_{i\in I\cup J})$. Pour tout $\Phi\in \hc_I$, ce sont des éléments semi-simples respectivement de  $\tugo_{\tilde{\Phi},a_0}(F)$ et $\tugo_{\Phi^\flat,a_0'}(F)$ pour un certain $a_0'\in \Ac_H(F)$ indépendant de $\Phi^\flat$.

Il résulte du corollaire \ref{cor:orbite-ss} que, lorsqu'on fait varier $\Phi\in |\hc_I|$,  les classes de $U_{\tilde{\Phi}}(F)$-conjugaison de $X_0$ parcourent de manière biunivoque l'ensemble  $\cup_{\Phi_0 \in \hc} \oc(\Phi_0)_{a_0}$. De même,  quand $\Phi$ décrit $ |\hc_I|$, les classes de $U_{\Phi^\flat}(F)$-conjugaison de $Y_0$ sont exactement les éléments de $\cup_{\Phi_1 \in \hc^\flat} \oc(\Phi_1)_{a_0'}$. 

Par conséquent,  si, pour tout $\Phi\in \hc_I$, on note  $I_{0}^{\tilde{\Phi}}$, resp. $I_{0}^{\Phi^\flat}$, la distribution sur $\Ic(\hc)$, resp. $\Ic(\hc^\flat)$, associée à la forme $\tilde{\Phi}$ et l'orbite de $X_0$, resp. la forme $\Phi^\flat$ et l'orbite de $Y_0$, on obtient  le lemme suivant.

\begin{remarque}
  On prendra garde que $I_{0}^{\tilde{\Phi}}$ dépend vraiment de $\Phi$ et pas seulement de $\tilde{\Phi}$ comme la notation pourrait le laisser croire. On espère que cela ne crée pas de confusion.
\end{remarque}

\begin{lemme}
  \label{lem:Ihc0}
On a  les deux égalités :
  \begin{enumerate}
  \item 
$$I_{a_0}^{\hc}=\sum_{\Phi \in |\hc_I|}   I_{0}^{\tilde{\Phi}}.
$$
\item $$I_{a_0'}^{\flat}=\sum_{\Phi \in |\hc_I|}   I_{0}^{\Phi^\flat}.
$$
\end{enumerate}
\end{lemme}
\end{paragr}

\begin{paragr}
Soit $f\in \Ic(\eta)$ et $f'\in \Ic(\hc)$ qui se correspondent. Par définition,  pour tout $a\in \Ac(F)$, on a 
\begin{equation}
  \label{eq:IaPhi}
  I_a^{\hc}(f')=\sum_{\Phi\in |\hc|} \sum_{\of\in \oc(\Phi)}I_a^\Phi(f')
\end{equation}

Pour tout ensemble $S\subset \vc$ fini de places contenant les places archimédiennes, soit $\hc^S\subset \hc$ et $\hc_I^S\subset \hc_I^S$ les sous-catégories des formes dont le discriminant est une norme hors $S$. Les ensembles  $|\hc^S|$ et  $|\hc_I^S|$ sont finis.

D'après le lemme \ref{lem:Sf}, il existe un ensemble $S_1$ fini contenant $\vc_\infty$ tel que $I^\Phi_\of(f)=0$ pour tout $\Phi$ tel que $S_\Phi\not\subset S_1$ et $\of\in \oc(\Phi)$. Autrement dit, dans\eqref{eq:IaPhi}, on peut limiter la somme aux $\Phi\in \hc^{S_1}$.

Quitte à agrandir $S_1$, on suppose qu'il contient toutes les places de $F$ qui se ramifient dans l'un des corps $E,F_i$ ou $E_i$ pour $i\in I\cup J$. On travaille avec l'ensemble fini $|\hc_I^{S_1}|$ identifié à un système de représentants. 

Soit $S$ contenant $S_1$ tel que $f$ et $f'$ sont représentés respectivement par $f_S\otimes \mathbf{1}^S$ et $f_S'\otimes \mathbf{1}^S$ et tel que les conditions du §\ref{S:UsurA} soient satisfaites. On suppose également que pour $v\notin S$, on a $a_0'\in \Ac_H'(\oc_v)$.

\begin{proposition}\label{prop:f-fH}
  Il existe $f^H_S\in \Ic^H(\eta)_S$ et $f^\flat_S\in \Ic(\hc^\flat)_S$  telles que si l'on  note $f^H \in \Ic^H(\eta)$ et $f^\flat\in \Ic(\hc^\flat)$ les éléments représentés respectivement par $f^H_S\otimes \mathbf{1}^S\in \Ic^H(\eta)$ et  $f^\flat_S\otimes \mathbf{1}^{S}$,  les trois conditions  sont satisfaites :
\begin{enumerate}
\item les fonctions $f^H$ et  $f^\flat$ se correspondent.
\item $I_{a_0}^\eta(f)=I_{a_0'}^{\eta}(f^H)$
\item Pour toute forme $\Phi\in \hc_I$, on a 
$$I_{0}^{\tilde{\Phi}}(f')=I_{0}^{\Phi^\flat}(f^\flat).
$$
\end{enumerate}
\end{proposition}

La démonstration occupe les §§\ref{S:cst-f} à \ref{S:ass3}. Notons tout de suite le corollaire suivant de la proposition \ref{prop:f-fH} ci-dessus. Ce corollaire démontre donc le théorème \ref{thm:transfert} pour l'élément $a=a_0$.

\begin{corollaire}
On a 
  $$ I_{a_0}^\eta(f) =I^\hc_{a_0}(f').
$$
\end{corollaire}

\begin{preuve}  En utilisant successivement l'assertion $2$ de la proposition \ref{prop:f-fH}, l'hypothèse de récurrence qu'on applique aux éléments $f^H$ et  $f^\flat$ qui se correspondent par l'assertion 1 de de la proposition \ref{prop:f-fH},  le lemme \ref{lem:Ihc0} assertion 2,  l'assertion 3  de  proposition  \ref{prop:f-fH} et enfin  l'assertion 1 du  lemme \ref{lem:Ihc0}, on obtient
\begin{eqnarray*}
   I_{a_0}^\eta(f) &=& I_{a_0'}^{\eta}(f^H)\\
 &=& I_{a_0'}^\flat(f^\flat)\\
&=& \sum_{\Phi \in |\hc_I|   } I^{\Phi^\flat }_{0}(f^\flat)\\
 &=&  \sum_{\Phi \in |\hc_I|   } I^{ \tilde{\Phi} }_{0}(f')\\
&=& I^\hc_{a_0}(f')
\end{eqnarray*}
\end{preuve}

\end{paragr}

\begin{paragr}[Construction des fonctions $f^H_S$ et $f_S^\flat$ dans la proposition \ref{prop:f-fH}.]---  \label{S:cst-f}Comme au §\ref{S:Omega}, on dispose d'ouverts $\om_H\subset \Ac_H'(\AAA_S)$ et  $\om\subset \Ac(\AAA_S)$ voisinages respectifs de $a_0'$ et $a_0$. De plus, le morphisme $\iota_H$ induit un isomorphisme entre ces ouverts. On pourra restreindre si besoin est ces ouverts.

  Pour $v\in S$, soit $\hc_{I,v}$ le groupoïde quotient de l'ensemble des familles de formes hermitiennes $(\Phi_{i,v})_{i\in I}$ sur les $V_{i,E_v}$ par le groupe $\prod_{i\in I} GL_{E_v}(V_{i,E_v})$. Rappelons qu'on dispose d'une forme $\Phi^+$ et pour $j\in J$ d'une forme $\Phi_j$. Cela permet de définir comme dans le cas global (cf. \ref{S:reecriture}) des applications $\Phi\in \hc_{I,v} \mapsto \Phi^\flat\in \hc^\flat_v$ et $\Phi\in \hc_{I,v} \mapsto \tilde{\Phi}\in \hc_v$ . 

Soit  $\hc_{I,S}$ le produit en un sens évident des $\hc_{I,v}$ pour $v\in S$. De même, on définit  $\hc^\flat_S$. Pour tout $\Phi \in \hc_{I,S}$, on sait construire alors deux actions de groupes au-dessus de $\AAA_S$ : celle de $U_{\tilde{\Phi}}$ sur $\tugo_{\tilde{\Phi}}$ et celle de $U_{\Phi^\flat}$ sur $\tugo_{\Phi^\flat}$. On prendra garde que ces situations ne viennent pas toujours d'un situation globale ; c'est le cas néanmoins si  $\tilde{\Phi} $ vient par extension des scalaires d'un élément de $\hc$ dans le premier cas et si $\Phi$ vient  par extension des scalaires  d'un élément de $\hc_I$.

Dans la suite, on travaille avec les ensembles de classes d'isomorphismes qu'on identifie souvent à un système de représentants. 

Pour tout $\Phi_0\in |\hc_S|$ et $\Phi_1\in |\hc_S^\flat|$, on dispose d'ouverts $\Omega^{\Phi_1}=a^{-1}(\om_H)\subset \tugo_{\Phi_1}(\AAA_S)$ et $\Omega^{\Phi_0}=a^{-1}(\om)\subset \tugo_{\Phi_0}(\AAA_S)$. 

D'après la discussion du §\ref{S:OmegaU}, pour tout $\Phi\in |\hc_{I,S}|$, on a une submersion 
\begin{equation}
  \label{eq:unesubmersion}
  U_{\tilde{\Phi}}(\AAA_S) \times \tugo_{\Phi^\flat}'(\AAA_S) \to \tugo_{\tilde{\Phi}}(\AAA_S) \ ;
\end{equation}
soit $\Omega^{\Phi}$ l'image de $U_{\tilde{\Phi}}(\AAA_S) \times \Omega^{\Phi^\flat}$ par celle-ci. Alors $\Omega^{\Phi}$ est un ouvert de  $\tugo_{\tilde{\Phi}}(\AAA_S)$. 

\begin{lemme}\label{lem:dec-OmPhi0}
Quitte à restreindre $\om_H$ (et $\om$ en conséquence),  pour tout $\Phi_0\in  |\hc_S|$,  on a une réunion disjointe 
$$\Om^{\Phi_0}=\bigcup_{ \{\Phi\in |\hc_{I,S}| \, \mid \tilde{\Phi}=\Phi_0 \}} \Om^{\Phi}.
$$
\end{lemme}

\begin{preuve}
Si  $\Phi_0$ appartient à l'image de $|\hc_{I,S}|\to |\hc_S|$, l'égalité résulte de \eqref{eq:dec-Omphi}.    Si $\Phi_0$ n'appartient pas à l'image de $|\hc_{I,S}|\to |\hc_S|$,  il s'agit de voir que l'ouvert $\Om^{\Phi_0}$ est vide . Il suffit en fait de voir que $\tugo_{\Phi_0}^{\rs}(\AAA_S)\cap  \Om^{\Phi_0}=\emptyset$. Il suffit donc de voir pour tout  $a\in \Ac^{\rs}(\AAA_S)\cap \om$, on a   $\tugo_{\Phi_0,a}(\AAA_S)=\emptyset$. Un tel  $a$ est l'image d'un élément $a' \in \om_H\cap \Ac_H^{\Grs}(\AAA_S)$. D'après le corollaire  \ref{cor:orbite-ss}, il existe donc $\Phi_1\in |\hc^\flat_S|$ et $Y\in \tugo_{\Phi_1}(\AAA_S)$ tel que $a(Y)=a'$. La première composante de $\Phi_1$ (c'est-à-dire celle associée à $V^+_E$) est déterminée par la composante correspondante de $a'$. Quitte à restreindre l'ouvert $\om_H$, on peut et  on va supposer que la première composante de $\Phi_1$ est égale à la forme obtenue à partir de $\Phi^+$. Il s'ensuit que $\Phi_1=\Phi^\flat$ pour $\Phi\in \hc_{I,S}$. Par hypothèse, les formes  $\Phi_0$ et $\tilde{\Phi}$ ne peuvent être égales. Il est clair alors que   le morphisme  $\iota : \tugo'_{\Phi^\flat}(\AAA_S) \to \tugo_{\tilde{\Phi}}(\AAA_S)$ (cf. \eqref{eq:iotaU})  envoie $Y$ sur un élément $Z$ régulier semi-simple tel que $a(Z)=a$. Donc $\tugo_{\tilde{\Phi},a}(\AAA_S)\not=\emptyset$ et forcément  $\tugo_{\Phi_0,a}(\AAA_S)=\emptyset$ (par le corollaire \ref{cor:orbite-ss} $a$ régulier semi-simple possède une fibre non vide pour une seule forme).

Désormais, on suppose que $\om_H$ et $\om$ sont tels que les conclusions du lemme \ref{lem:dec-OmPhi0} sont valables.

L'élément $f_S'$ est en fait représenté par une collection $(f'_{S,\Phi_0})_{\Phi_0\in \hc_S}$ où l'on voit chaque composante $f'_{S,\Phi_0}$ comme un élément de $\otimes_{v\in S} \Sc( \tugo_{\Phi_0}(F_v))$. Soit $\theta \in \Cc(\om)$ qui vaut $1$ au voisinage de $a_0$. On peut et on suppose que $\theta$ est un tenseur pur.

\begin{lemme}\label{lem:dec-sur-Phi}
Pour tout $\Phi_0\in |\hc_S|$, on a 
$$  (\theta \circ a) f'_{S,\Phi_0}= \sum_{  \{\Phi\in |\hc_{I,S}|\,  \mid \tilde{\Phi}=\Phi_0 \}   }  f'_{S,\Phi}
$$
où l'on définit pour $\Phi\in |\hc_{I,S}|$ l'élément $f'_{S,\Phi}\in \Sc(\Omega^\Phi)$ par 
 $$f'_{S,\Phi}=  (\theta \circ a) \mathbf{1}_{\Omega^\Phi}  f'_{S,\tilde{\Phi}}.
$$
\end{lemme}

\begin{preuve}
Le lemme  \ref{lem:dec-OmPhi0}  implique que  la restriction de $ \theta \circ a$ à $\tugo_{\Phi_0}(\AAA_S)$ est égale à 
$$\sum_{  \{\Phi\in |\hc_{I,S}|\,  \mid \tilde{\Phi}=\Phi_0 \}   }  (\theta \circ a) \mathbf{1}_{\Omega^\Phi} .$$
Par ailleurs $\Omega^\Phi$ est ouvert et fermé dans $\Om^{ \tilde{\Phi}  }$. La fonction $f'_{S,\Phi}$ est encore de Schwartz.
\end{preuve}

 Soit  $\Phi\in |\hc_{I,S}|$. D'après les considérations du §\ref{S:Om1}, il existe $\al\in \Sc(U_{\tilde{\Phi}}(\AAA_S)\times \Om^{\Phi^\flat})$ tel que $f'_{S,\Phi}= f_\al$ ; on pose alors 
$f^\flat_{S,\Phi}=f_\al^\flat$.

On obtient ainsi une collection de fonctions $(f^\flat_{S,\Phi})_{\Phi\in \hc_{I,S}}$ qu'on interprète comme un élément
$$f_S^\flat \in \Ic(\hc^\flat)_S.$$ 

Du côté linéaire, la construction est plus simple. La fonction $(\theta\circ a)f_S$ est de la forme $f_\al$ pour $\al\in \Sc(G(\AAA_S)\times \Om_H)$ où $\Om_H$ est l'image inverse de $\om_H$ dans $\thgo(\AAA_S)$ (avec les notations du  §\ref{S:Omega}). On pose alors $f_S^H=f_{\al}^{H,\eta}$ ce qui définit un élément de $\Ic(\eta)_S$. 
\end{preuve}
\end{paragr}

\begin{paragr}[Vérification de l'assertion 1. ] ---Par hypothèse, les fonctions $f_S$ et $f_S'$ ont les mêmes intégrales orbitales semi-simples régulières. Il en est donc de même pour  $(\theta\circ a) f_S$ et $(\theta\circ a)f_S'$. En tenant compte du lemme \ref{lem:dec-sur-Phi}, on a donc pour tout $a\in \Ac^{\rs}(\AAA_S)$
\begin{eqnarray*}
  I_{S,a}^{\tilde{\eta}}((\theta\circ a)f_S)&=&\sum_{\Phi \in |\hc_{I,S}|} I_{S,a}^{\tilde{\Phi}}(f'_{S,\Phi})
\end{eqnarray*}
(NB: l'intégrale orbitale semi-locale $I_{S,a}^{\tilde{\Phi}}$ se définit comme en \eqref{eq:IOPloc2U}). Essentiellement par construction des fonctions $f^H_S$ et  $f_S^\flat$ (on utilise les  lemmes \ref{lem:desc-locale} et  \ref{lem:desc-localeU} et même une généralisation évidente de ce dernier), on en déduit que pour tout $a\in \om^H \cap \Ac_H^{\Grs}(\AAA_S)$ on a 
\begin{equation}\label{eq:IafSH}
  I_{S,a}^{\tilde{\eta}}(f_S^H)=\sum_{\Phi \in |\hc_{I,S}|} I_{S,a}^{\Phi^\flat}(f^\flat_{S,\Phi}).
\end{equation}
Par lissité des intégrales orbitales sur  $\Ac_H^{\rs}(\AAA_S)$, cette égalité se prolonge à $ \om^H \cap \Ac_H^{\rs}(\AAA_S)$ puis à tout  $\Ac_H^{\rs}(\AAA_S)$ puisque les deux membres \eqref{eq:IafSH} sont nuls pour $a\in \Ac_H^{\rs}(\AAA_S)$  mais $a\notin \om_H$. Par conséquent, $f^H_S$ et $f^\flat_S$ se correspondent. Le lemme fondamental permet de conclure.
\end{paragr}

\begin{paragr}[Vérification de l'assertion 2 de la proposition \ref{prop:f-fH}.] --- Par la propriété de support de $I_{a_0}^{\eta}$ (cf. théorème \ref{thm:I}), on a 
$$I_{a_0}^\eta(f_S\otimes\mathbf{1}^S)=I_{a_0}^\eta(f \otimes\mathbf{1}^S)( (\theta\circ a)f_S \otimes\mathbf{1}^S).$$
Le théorème \ref{thm:desc} donne 
$$
  I_{a_0}^\eta(f \otimes\mathbf{1}^S)( (\theta\circ a)f_S \otimes\mathbf{1}^S)=I_{a_0'}^{\eta}(f^H_S  \otimes\mathbf{1}^S).
$$
La combinaison de ces deux égalités donne l'assertion de la proposition \ref{prop:f-fH}.
\end{paragr}

\begin{paragr}[Vérification de l'assertion 3 de la proposition \ref{prop:f-fH}.] --- \label{S:ass3}Soit $\Phi\in \hc_I$. Cette forme \og globale\fg{}  donne 
\begin{itemize}
\item des formes globales $\tilde{\Phi}\in \hc$ et $\Phi^\flat \in \hc^\flat$ ;
\item par extension des scalaires une forme locale $\Phi_S\in \hc_{I,S}$ ;
\item des formes locales $\tilde{\Phi}_S \in \hc_S$ et $\Phi_S^\flat \in \hc^\flat_S$ (levons l'ambiguïté sur la notation ; les deux interprétations possibles sont équivalentes  : ces formes  sont obtenues à partir de la forme locale $\Phi_S$ ou sont des  localisées des formes globales  $\tilde{\Phi}$ et $\Phi^\flat$ déduites de $\Phi$).
\end{itemize}

Si $\tilde{\Phi}\notin \hc^S$ alors $\Phi\notin \hc_I^S$ et les deux membres de l'assertion 3 de la proposition \ref{prop:f-fH} sont nuls. On suppose désormais qu'on a $\tilde{\Phi}\in \hc^S$.

\begin{lemme}\label{lem:etape}
  On a 
$$I_0^{\tilde{\Phi}}(f')=I^{\tilde{\Phi}}_0(f'_{S,\Phi_S}\otimes\mathbf{1}_{\tugo_{\tilde{\Phi}}(\oc^S)}).
$$
\end{lemme}

\begin{preuve}
Par définition de $I_0^{\tilde{\Phi}}(f')$, on a 
$$I_0^{\tilde{\Phi}}(f')=I^{\tilde{\Phi}}_0(f'_{S,\tilde{\Phi}_S}\otimes\mathbf{1}_{\tugo_{\tilde{\Phi}}(\oc^S)}).$$
  Par les propriétés de support de $I^{\tilde{\Phi}}_0$ (cf. théorème \ref{thm:IU} assertion 3), on peut remplacer dans le membre de droite, $f '_{S,\tilde{\Phi}_S}$ par $(\theta\circ a) f '_{S,\tilde{\Phi}_S}$. En utilisant le lemme \ref{lem:dec-sur-Phi}, il suffit de prouver que pour $\Phi_1\in \hc_{I,S}$ tel que $\Phi_1\not=\Phi_S$ et $\tilde{\Phi}_1=\tilde{\Phi}_S$ on a 
$$I^{\tilde{\Phi}}_0(f'_{S,\Phi_1}\otimes\mathbf{1}_{\tugo_{\tilde{\Phi}}(\oc^S)})=0.$$
Or, d'après le lemme \ref{lem:dansOm1v} assertion 4, l'ouvert $\Omega^{\Phi_1}$ ne rencontre pas  la $U_{\tilde{\Phi}}(\AAA_S)$-orbite de $X_0$ puisque celle-ci est incluse dans $\Omega^{\Phi_S}$. Pour des raisons de support  (cf. théorème \ref{thm:IU} assertion 3), on obtient la nullité annoncée.
\end{preuve}

Supposons $\Phi\in \hc_I^S$. Alors, on a,  d'après le  théorème \ref{thm:descU} de descente,
\begin{equation}
  \label{eq:avoir}
  I^{\tilde{\Phi}}_0(f'_{S,\tilde{\Phi}_S}\otimes\mathbf{1}_{\tugo_{\tilde{\Phi}}(\oc^S)})=I_0^{\Phi^\flat}( f^\flat_{S,\Phi_S }\otimes \mathbf{1}_{\tugo_{\Phi^\flat}(\oc^S)}).
\end{equation}
D'après  le lemme \ref{lem:etape} ci-dessus, le membre de gauche est $I_0^{\tilde{\Phi}}(f')$ alors que le membre de droite est par définition $I_0^{\Phi^\flat}( f^\flat)$. On obtient bien l'assertion 3  de la proposition \ref{prop:f-fH}.

Supposons $\Phi\in \hc_I \setminus \hc_I^S$. Le membre de droite de l'assertion 3  de la proposition \ref{prop:f-fH} est alors nul. Il s'agit donc de prouver

\begin{equation}
  \label{eq:nullite-finale}
  I_0^{\tilde{\Phi}}(f')=0.
\end{equation}

Soit $v\notin S$ tel que $\Phi_v$ ait une composante dont le discriminant ne soit pas une norme. On fixe $\om_{H,v}$ un voisinage de $a_0'$dans $\Ac_H'(\oc_v)$ et $\om_{v}$ un voisinage de $a_0$ dans $\Ac(\oc_v)$ de sorte que le morphisme $\iota_H$ (cf. \eqref{eq:iotaHquotient}) induise un isomorphisme de $\om_{H,v}$ sur $\om_v$. Soit $\Om_v$ l'image inverse de $\om_v$ dans $\tugo_{\tilde{\Phi}_v}(F_v)$. Soit $\Phi_1 \in |\hc_{I,v}|$ tel que $\tilde{\Phi}_1= \tilde{\Phi}_v$. Soit  $\Om^{\Phi_1^\flat}$ l'image inverse de $\om_{H,v}$ dans $\tugo_{\Phi^\flat_v}(F_v)$. Soit $\Om^{\Phi_1}$ l'image de $U_{\Phi_1}(F_v)\times \Om^{\Phi_1^\flat}$ par la submersion (cf. \eqref{eq:isocrucialUII})
$$U_{\Phi_1}(F_v)\times \tugo_{\Phi_1^\flat}' (F_v)\to \tugo_{\tilde{\Phi}}(F_v).
$$

Comme au lemme \ref{lem:dec-OmPhi0}, quitte à restreindre les ouverts $\om_{H,v}$ et $\om_v$,  on a une réunion disjointe 
$$\Om_v=\bigcup_{\{\Phi_1 \in |\hc_{I,v}|\, \mid \tilde{\Phi}_1=\Phi_v\}} \Om^{\Phi_1}.$$

Soit $\of_v$ la $U_{\tilde{\Phi}_v}(F_v)$-orbite de $X_0$. Par construction, on a $\of\subset \Om^{\Phi_v}$. Par les propriétés de support de la distribution $I_0^{\tilde{\Phi}}$ (cf. théorème \ref{thm:IU}) et le théorème de densité (cf. théorème \ref{thm:densiteU}), pour prouver \eqref{eq:nullite-finale}, il suffit de voir que pour tout $a\in \Om_v\cap \Ac^{\rs}(F_v)$, l'intégrale orbitale suivante  s'annule
$$I_a^{ \tilde{\Phi}_v }(\mathbf{1}_{\Om^{\Phi_v}} \mathbf{1}_{\tugo_{\tilde{\Phi}_v}(\oc_v)})=0.$$

Si ce n'est pas le cas, il existe $X \in  \tugo_{\tilde{\Phi}_v}(\oc_v)\cap \Om^{\Phi_v}$ tel que $a(X)=a$. 
Soit $\Phi_0\in \hc_{I,v}$ dont toutes les composantes sont de discriminant une norme et tel que $\tilde{\Phi}_0=\tilde{\Phi}_v$. On a un isomorphisme de $\oc_v$-schéma (cf. \eqref{eq:isocrucialUII})
$$U_{\tilde{\Phi}_v}\times^{U_{\Phi_0^\flat}} \tugo_{\Phi_0^\flat}' \to \tugo_{\tilde{\Phi}_v} \times_{\Ac }\Ac_H'
$$
Soit $a'$ est l'image inverse de $a$ par $\iota_H$. Le couple $(X,a')$ définit un $\oc_v$-point de  $\tugo_{\tilde{\Phi}_v} \times_{\Ac }\Ac_H'$ donc un  $\oc_v$-point de $U_{\tilde{\Phi}_v}\times^{U_{\Phi_0^\flat}} \tugo_{\Phi_0^\flat}' $. Le morphisme  
$$
U_{\tilde{\Phi}_v}\times \tugo_{\Phi_0^\flat}' \to U_{\tilde{\Phi}_v}\times^{U_{\Phi_0^\flat}} \tugo_{\Phi_0^\flat}' 
$$
est lisse et induit une application surjective au niveau des points sur le corps résiduel (lemme de Lang) : il s'ensuit que l'application
$$
U_{\tilde{\Phi}_v}(\oc_v)\times \tugo_{\Phi_0^\flat}' (\oc_v) \to (U_{\tilde{\Phi}_v}\times^{U_{\Phi_0^\flat}} \tugo_{\Phi_0^\flat}' )(\oc_v)
$$
est surjective. On en déduit que $X\in \Om^{\Phi_0}$. Or par hypothèse $\Phi_v\not=\Phi_0$ donc $\Om^{\Phi_0}\cap \Om^{\Phi_v}=\emptyset$. Contradiction.
\end{paragr}

\subsection{Fin de la démonstration du théorème \ref{thm:transfert}}\label{ssec:fin-tr}

\begin{paragr}
  On continue avec les notations des section précédentes. Soit 
$$\psi: F\back  \AAA \to \CC^\times$$
un caractère additif, continu et non trivial. Pour tout $v\in \vc$, on en déduit un caractère local $\psi_v: F_v\to \CC^\times$.
\end{paragr}

\begin{paragr}[Transformation de Fourier : cas hermitien.] ---
  Pour tout $\Phi\in \hc_v$, l'espace $\tugo_{\Phi}$ est muni de la forme bilinéaire donnée en \eqref{eq:bilinU}. L'action de $GL_E(V_E)$ fait que les  sous-$F$-espaces $\tugo_1$ de $\tugo_{\Phi}$, qui sont $U_\Phi$-invariants  et non-dégénérés pour la forme précédente,  se correspondent naturellement lorsque $\Phi$ décrit $\hc$. La donnée d'un tel sous-espace $\tugo_1$ définit  donc une transformation de Fourier partielle $\mathfrak{F}_{\tugo_1}$ sur $\Sc(\tugo_{\Phi}(F_v))$ pour tout La même remarque vaut pour $\Phi\in \hc$ selon l'analogue local de la formule \eqref{eq:TFPU} ; la mesure sur $\tugo_{\Phi}(F_v)$ est l'unique mesure auto-duale. 

D'après le théorème \ref{thm:TFstableU}, la transformation de Fourier  $\mathfrak{F}_{\tugo_1}$ induit un automorphisme de $\Ic(\Phi_v)$ encore noté $\Fgo_{\tugo_1}$ et donc un automorphisme de $\Ic(\hc_v)$.

Soit $S_0\subset \vc$  fini, qui contient les places archimédiennes et les places ramifiées dans $E$ et tel que $\psi_v$ soit de conducteur $\oc_v$. Alors pour tout $v\notin S$, on a $\Fgo_{\tugo_1}\mathbf{1}_v=\mathbf{1}_v $. Il s'ensuit que  $\Fgo_{\tugo_1}$ induit un automorphisme de $\Ic(\hc)$.

\end{paragr}

\begin{paragr}[Transformation de Fourier : cas linéaire.] --- On munit $\tggo$ de la  la forme bilinéaire donnée en \eqref{eq:fbs}. Soit $v\in \vc$. Comme dans le cas hermitien, le caractère additif $\psi_v$ et tout sous-espace $\tggo_1$ de $\tggo$ non-dégénéré et $G$-invariant donnent une transformation de Fourier $\Fgo_{\tggo_1}$ sur $\Sc(\tggo(F_v))$ donnée par l'analogue local de la formule intégrale \eqref{eq:TFPg1} pour la mesure auto-duale sur  $\tggo_1(F_v)$. Par le théorème \ref{thm:TFstable},  on en déduit un automorphisme $\Fgo_{\tggo_1}$  de $\Ic(\eta_v)$, puis un automorphisme encore noté $\Fgo_{\tggo_1}$ de $\Ic(\eta)$.  
\end{paragr}

\begin{paragr}[Correspondances et transformation de Fourier.] --- Soit $\tggo_1$ et $\tugo_1$ des sous-espaces invariants et non dégénérés comme ci-dessus. On suppose qu'ils se correspondent au sens où on les identifie après extension des scalaires à $E$. On a alors le lemme suivant.

  \begin{lemme}\label{lem:commu-tr-cor}
    Soit $f\in \Ic(\eta)$ et $f'\in \Ic(\hc)$ qui se correspondent. Alors il en est de même de $\Fgo_{\tggo_1}(f)$ et de $\Fgo_{\tugo_1}(f')$.
  \end{lemme}

  \begin{preuve}
    On utilise \cite{Z1} théorème 4.17 et \cite{xue} lemmes 9.2 à 9.4. D'après ces auteurs, il existe un entier $m$ attaché à $\tggo_1$ tel que pour tout $v\in \vc$, si les éléments $f_v\in \Ic(\eta_v)$ et $f'_v\in \Ic(\hc_v)$ se correspondent alors il en est de même des éléments 
$$\eps(\eta_v,\frac12,\psi_v)^m  \Fgo_{\tggo_1}(f_v) \ \  \text{    et   } \ \  \Fgo_{\tugo_1}(f'_v).$$
Le lemme résulte alors du fait que pour toute partie $S\subset\vc$ finie contenant les places archimédiennes, les places ramifiées dans $E$ et les places pour lesquels $\psi_v$ est ramifiée on a 
$$\prod_{v\in S} \eps(\eta_v,\frac12,\psi_v)=1.
$$
  \end{preuve}
\end{paragr}

\begin{paragr}[Factorisation de la formule des traces infinitésimale.] --- \label{S:factoRTF}On continue avec les hypothèses précédentes. On suppose de plus que $\tggo_1$ est l'un des sous-espaces de la remarque \ref{rq:TFP}. Par conséquent, $\tugo_1$ est l'un des sous-espaces de la remarque \ref{rq:TFPU}.

  \begin{theoreme} Soit $\al\in F$. 
    \label{thm:RTF-fact}

    \begin{enumerate}
    \item  Pour tout $f\in \Ic(\eta)$, on  a 
$$\sum_{a\in \Ac_\al(F)} I_a^\eta(f)= \sum_{a\in \Ac_\al(F)} I_a^\eta(\Fgo_{\tggo_1}(f)).$$
\item Pour tout $f\in \Ic(\hc)$, on  a 
$$\sum_{a\in \Ac_\al(F)} I_a^\hc(f)= \sum_{a\in \Ac_\al(F)} I_a^\hc(\Fgo_{\tugo_1}(f)).$$
    \end{enumerate}
  \end{theoreme}

  \begin{preuve}
    L'assertion 1 se déduit directement du  théorème \ref{thm:RTFinf}. Pour l'assertion 2, il faut se rappeler qu'à $f$ fixé, seul un nombre fini de formes hermitiennes $\Phi$ intervient dans la définition de  $I_a^\hc(f)$ (cf. lemmes \ref{lem:Sf} et \ref{lem:finitude}).  L'assertion 2 découle alors du théorème \ref{thm:RTFinfU} appliqué à un nombre fini de formes $\Phi$.
  \end{preuve}
\end{paragr}

\begin{paragr}[Fin de la démonstration du théorème \ref{thm:transfert}.] --- Soit $f\in \Ic(\eta)$ et $f'\in \Ic(\hc)$ qui se correspondent. Soit $S\subset \vc$ fini, assez grand pour que $f$ et $f'$ soient respectivement représentées par $f_S\otimes \mathbf{1}^S$ et $f_S'\otimes \mathbf{1}^S$. 

Fixons une place auxiliaire $v\notin S$ non-archimidienne et pour simplifier décomposée dans $E$. On a alors $|\hc_v|=1$, le caractère $\eta_v$ est trivial, $\tggo_v\simeq \tugo_{\Phi}$ pour $\Phi\in \hc_v$. Il en résulte qu'on a une identification $\Ic(\eta_v)\simeq \Ic(\hc_v)$ (ce fait est général, du moins pour $v$ non-archimédien, mais tout-à-fait non trivial si $v$ n'est pas décomposé dans $E$, cf. §\ref{S:corresp}) 

Soit $S'=S\cup\{v\}$. En utilisant l'identification $\Ic(\eta_v)\simeq \Ic(\hc_v)$, on voit que les classes des fonctions $f_S\otimes f_v \otimes  \mathbf{1}^{S'}$ et $f_S'\otimes f_v \otimes  \mathbf{1}^{S'}$ se correspondent. 

Soit $f_v\in \Cc(\tggo(F_v))$. Pour tout $a\in \Ac(F)$, on définit
$$T_a(f_v)=I_a^\eta(  f_S\otimes f_v \otimes  \mathbf{1}^{S'})- I_a^\hc(  f_S'\otimes f_v \otimes  \mathbf{1}^{S'}  )$$
où on confond dans les notations $f_v$ avec sa classe dans $\Ic(\hc_v)$.
On va montrer que $T_a=0$ ce qui  implique bien sûr le théorème \ref{thm:transfert}.  On sait que $T_a(f_v)=0$ pour tout les $a$ qui sont réguliers (cf. §\ref{S:areg}) ou qui donnent lieu à une descente (on utilise alors une récurrence, cf. section \ref{ssec:tr-desc}). Il reste donc à traiter le cas de $a\in \Ac^{(0)}$ qui ne donnent pas lieu à une descente. Ces $a$-là correspondent à la donnée d'un polynôme scindé sur $F$ qui n'a qu'une racine $\al\in F$. L'homothétie de rapport $\al$ s'interprète comme un endomorphisme de $V$ et par extension des scalaires comme un endomorphisme de $V_E$ auto-adjoint pour toute forme hermitienne. Quitte à utiliser des translations par cet endomorphisme sur les fonctions $f$ et $f'$ (ce qui ne change pas le fait qu'elles se correspondent), on est ramené à traiter le seul cas $a=0$.

On note aussi pour tout sous-espace $\tggo_1$ (et l'espace $\tugo_1$ qui lui est attaché, cf. §\ref{S:factoRTF}) 
$$T_a^\vee(f_v)= I_a^\eta(  \Fgo_{\tggo_1}(f_S\otimes f_v \otimes  \mathbf{1}^{S'}))- I_a^\hc(   \Fgo_{\tugo_1}( f_S'\otimes f_v \otimes  \mathbf{1}^{S'}  )).$$

Le théorème \ref{thm:RTF-fact} implique qu'on a (pour $\al=0$)
\begin{equation}
  \label{eq:RTF-T}
  \sum_{a\in \Ac_0(F)} T_a(f_v)= \sum_{a\in \Ac_0(F)} T_a^\vee(f_v).
\end{equation}
Comme on l'a déjà dit, les termes associés à $a\not=0$ sont nuls dans le membre de gauche : celui-ci se simplifie en $T_0(f_v)$.  Le même argument couplé au lemme \ref{lem:commu-tr-cor} montre que le membre de droite se simplifie en $T_0^\vee(f_v)$.
Soit $\hat{f}_v$ la transformation de Fourier partielle de $f_v$ attachée au sous-espace $\tggo_1$. On note $\hat{T}_0$ la transformation de Fourier (duale) de $T_0$. L'égalité \eqref{eq:RTF-T}implique donc qu'on a 
$$\hat{T}_0(f_v)=I_0^\eta(  \hat{ \hat {f _v}}\otimes  \Fgo_{\tggo_1}(f_S \otimes  \mathbf{1}^{S'}))- I_a^\hc( \hat{ \hat {f _v}}\otimes    \Fgo_{\tugo_1}( f_S' \otimes  \mathbf{1}^{S'}  )).
$$
Il résulte de cette égalité et des propriétés de support des distributions $I_0^\eta$ et $I_0^{\Phi}$ qu'on a les propriétés suivantes pour $T_0$ : 
\begin{enumerate}
\item $T_0$ est une distribution $G(F_v)$-invariante de support inclus dans le cône nilpotent $\tggo_{a=0}(F_v)$ ;
\item $\hat{T}_0$ est une distribution $G(F_v)$-invariante de support inclus dans le cône nilpotent $\tggo_{a=0}(F_v)$  pour toute transformation de Fourier associée  à un sous-espace  $\tggo_1$ comme dans le théorème \ref{thm:RTF-fact}.
\end{enumerate}

Il résulte alors du théorème d'Aizenbud (cf. théorème \ref{thm:incertitude}) qu'on a $T_0=0$.
\end{paragr}

\section{Distributions géométriques pour $ GL_n \times GL_{n+1} $}\label{sec:densiteS}

\subsection{Actions de groupes considérées}\label{ssec:prelimS}

\begin{paragr} 
Soit $E/F$ une extension quadratique de corps de nombres. Soit $\AAA$ et $\AAA_{E}$ les anneaux des adèles de $F$ et $E$. 
On invoque les conventions du début de la section \ref{sec:RTFinf}.
\end{paragr}

\begin{paragr}\label{S:prelimGps} --- 
Pour toute $ F $-variété $Y$, soit $ Y_{E} $ la restriction de scalaires de $ E $ à $ F $ de $ Y \times_{F} E $. On note encore $\sigma$ l'automorphisme de  $Y_{E}$ associé au générateur $\sigma$ du groupe de Galois $\Gal(E/F)$.
\end{paragr}

\begin{paragr}\label{S:baseFix} ---
  Soit $n \in \NN$. Soit $V$ un espace vectoriel de dimension $n$ qu'on munit d'une base. Soit $W = V \oplus Fe_{0}$. Soit $ e_0^{*} $ la forme linéaire sur $W$ de noyau $V$ qui vérifie  $ e_0^{*}(e_0) = 1 $. Soit $V_{E} = V \otimes_{F} E$ et $W_{E} = W \otimes_{F} E$. 

On considère les $F$-groupes suivants
\begin{itemize}
 \item $G = GL(V)$, 
\item $\tlG = GL(W)$.
\item $G' = G \times \tlG$.
\end{itemize}
 On identifie $G$ au sous-groupe de $\tlG$  qui fixe $e_{0}$ et stabilise $V$. On a alors un plongement diagonal $G\subset G'$. Par changement et restriction de base, on a une inclusion $G_E\subset G'_E$. On utilise aussi l'inclusion naturelle $G'\subset G'_E$.
 Les sous-groupes $G_{E}$ et $G'$ agissent sur $G'_E$ par translation  respectivement à gauche et à droite : pour $ (g, h) \in G_E \times G' $ et $ g' \in G'_E $ on a $ g\cdot g' \cdot  h = g g' h $. L'application $(g,\tilde{g})\in G'_E\mapsto g^{-1}\tilde{g} \in \tilde{G}_E$ identifie le quotient $G_E\back G_E'$ à $\tilde{G}_E$. Dans cette identification, l'action à droite de $G'=G\times \tilde{G}$ devient l'action à droite de $G'$ sur $\tilde{G}_E$ donnée par $\tilde{g}\cdot (g,h)= g^{-1}\tilde{g}h$ pour $\tilde{g}\in \tilde{G}_E$ et  $(g,h)\in G'=G\times\tlG$.
\end{paragr} 

\begin{paragr}[La variété $X$.] --- \label{S:varX}
Soit 
\[
X  = \{g \in \tlG_{E} \mid g \sigma(g) = 1\}.
\]
Le groupe $G$ agit sur $X$ par conjugaison. On note $\tilde{\Ac}=X // G$ le quotient catégorique. 

Après extension des scalaires à $E$, on a 
$$X \times_F E=\{ (g_1,g_2)\in GL_E(W_E)^2 \mid   g_1g_2=1\} $$
et donc un isomorphisme  $  GL_E(W_E)\simeq X \times_F E$ donné par $g\mapsto (g,g^{-1})$. Dans cet isomorphisme l'action par conjugaison de $GL_E(V_E)$ sur $X \times_F E$ devient l'action par conjugaison de $GL_E(V_E)$ sur   $  GL_E(W_E)$. On en déduit une identification des quotients catégoriques :
$$\tilde{\Ac}\times_F E \simeq  GL_E(V_E)\back\back GL_E(W_E).
$$
Le  quotient de droite est un ouvert du quotient catégorique $GL_E(V_E)\back\back \End_E(W_E)$ qui  s'identifie à l'espace affine sur $E$ de dimension $2n+1$ : un système de générateurs indépendants de l'algèbre des fonctions invariantes sur $ \End_E(W_E)$ est donnée par les fonctions suivantes de la variable $Y\in \End_E(W_E)$ : ce sont d'une part  les $n+1$ coefficients du polynôme caractéristique de $Y$ et d'autre part les fonctions $Y\mapsto e_0^*Y^ie_0$ pour $1\leq i \leq n$ (cf. §\ref{S:A}). L'ouvert de  $ \End_E(W_E)$  où le déterminant 
$$\det( (e_0^*Y^{i+j}e_0)_{1\leq i,j\leq n+1})$$
ne s'annule pas est formé des éléments réguliers semi-simples. Ce même déterminant qui est invariant définit un ouvert du quotient catégorique et donc un ouvert de $\tilde{\Ac}\times_F E $. Cet ouvert provient par extension des scalaires d'un ouvert $\tilde{\Ac}^{\rs}$. Pour tout $a\in \tilde{\Ac}$, soit $X_a$ la fibre au-dessus de $a$ du morphisme canonique
$$ X\to \tilde{\Ac}.$$
Soit $X^{\rs}$ l'image inverse de $\tilde{\Ac}^{\rs}$. C'est l'ouvert formé  des éléments réguliers et semi-simples. Pour tout  $a\in \tilde{\Ac}^{\rs}$, la fibre $X_a$ est une $G$-orbite géométrique. Il résulte de l'annulation du premier ensemble de cohomologie pour $G$, que pour tout $a\in \tilde{\Ac}^{\rs}(F)$, la fibre $X_a(F)$ est une $G(F)$-orbite (donc non vide).
\end{paragr}

\begin{paragr}[Identification de quotients.] ---
  Par le théorème 90 de Hilbert, l'application $\tilde{g}\mapsto \tilde{g}\tilde{g}^{-\sigma}$ identifie le quotient  $\tilde{G}_E/\tilde{G}$ à $X$. Dans cette identification, l'action à droite de $G$ sur  $\tilde{G}_E/\tilde{G}$  devient l'action par conjugaison de $G$ sur $X$.
Soit
\begin{equation}\label{eq:nu}
\nu : G'_E \to X, \quad \nu((h, \tlh)) = (g^{-1} \tlg) (g^{-1} \tlh)^{-\sigma}.
\end{equation}
L'application est surjective, induit un isomorphisme du quotient géométrique  $G_E\back G_E'/ \tlG$ sur $X$ et un isomorphisme entre quotients catégoriques
$$G_E\back \back G_E'// \tlG \simeq \tilde{\Ac}=X//G.$$
Certaines propriétés de $X$ se transposent à $G'_E$. Énonçons-en quelques-unes. Pour tout $a\in \tilde{\Ac}$, soit $G'_{E,a}$ la fibre au-dessus de $a$ du morphisme canonique
$$ G_E'\to \tilde{\Ac}.$$
Soit $(G_E')^{\rs}$ l'image inverse de $\tilde{\Ac}^{\rs}$. C'est l'ouvert formé  des éléments réguliers et semi-simples. Pour tout  $a\in \tilde{\Ac}^{\rs}$, la fibre $G'_{E,a}$ est une $G_E\times G'$-orbite géométrique. Si $a\in \tilde{\Ac}^{\rs}(F)$, la fibre  $G'_{E,a}(F)$ est une $(G_E\times G')(F)$-orbite (non vide).
\end{paragr}

\subsection{Les distributions géométriques}

 \begin{paragr}[Sous-groupes de Levi, paraboliques et compacts.] --- 
 On reprend les notations et les conventions de la section \ref{sec:RTFinf}.  Pour tout $\tlP \in \fc^{\tlG}(T_{\tl})$, soit 
$$P=\tlP\cap G\subset G$$
et
$$P'=P\times \tlP\subset G'.$$
On dispose aussi du sous-groupe parabolique $P'_E$ de $G'_E$. Le groupe $G'_E$ possède $T_0 \times T_{\tl}$ comme sous-tore déployé maximal. Alors $P'_E$ contient  $T_0 \times T_{\tl}$. On dispose d'une décomposition de Levi $P'_E=M_{P'_E} N_{P'_E}$ où le facteur de Levi $M_{P'_E}$ contient $T_0 \times T_{\tl}$. Pour tout $a\in \tilde{\Ac}(F)$ soit
$$M_{P'_E,a}=M_{P'_E}\cap G'_{E,a}.$$
On munit  $ N_{P'_E}(\AAA)  $ de  la mesure de Haar qui donne  le volume $1$ au quotient  $ N_{P'_E}(F) \back N_{P'_E}(\AAA) $.
 \end{paragr}
 
 \begin{paragr}[Les noyaux $k_{\tlP,a}$.] --- 
 Soit $f \in C_{c}^{\infty}(G'_{E}(\AAA))$.
Pour tous $a \in \tilde \Ac(F)$, $\tlP \in \fc^{\tlG}(B)$ et  $x,y \in N_{P'_E }(\AAA) M_{P'_E}(F) \back G'_{E}(\AAA)$ on pose
 \[
 k_{\tlP, a}(f,x,y) = 
 \sum_{\gamma \in M_{P'_E,a}(F)}
 \int_{N_{P'_E}(\AAA)}f(x^{-1} \gamma n y)dn.
 \]
Notons que la somme $ \sum_{a \in \tilde \Ac(F)}  k_{\tlP, a}(f)$ est le noyau usuel de l'opérateur de convolution à droite associé à $f$ agissant sur l'espace $L^{2}(N_{P'_E}(\AAA) M_{P'_E}(F) \back G'_E(\AAA))$.
 \end{paragr}
 
 \begin{paragr}[Les distributions $I_{a}^{\eta}$.] --- Soit  $\eta$ le caractère quadratique de $\AAA_{F}^{\times}$ associé à $E/F$. Pour $\tlP\in \fc^{\tlG}(T_{\tl})$, on dispose de l'application $H_{\tlP}:\tlG(\AAA) \to \ago_{\tlP}$ (cf. §\ref{S:HP} et §\ref{S:auxil}). 
Pour tout $g'=(g, \tlg) \in G'(\AAA) = G(\AAA) \times \tlG(\AAA)$,  on note encore, par abus de notations, 
\begin{eqnarray*}
\eta(g') &=& \eta(\det g)^{n+1}\eta(\det \tlg)^{n}
\end{eqnarray*}
et
$$H_{\tlP}(g')=H_{\tlP}(g)$$
où l'on voit $g$ comme un élément du sous-groupe $G(\AAA)$ de $\tlG(\AAA)$.

Soit $T \in \ago_{\tl}$.  On pose alors pour $x \in [G_{E}]$ et $y\in [G']$ 
 \[
 k_{a}^{T}(f, x, y) = \sum_{\tlP \in \fc^{\tlG}(B)}
 \varepsilon_{\tlP}^{\tlG}
 \sum_{
 \begin{subarray}{c}
 \delta_{1} \in P_E(F)\back G_E(F) \\
 \delta_{2} \in P'(F)\back G'(F)
  \end{subarray}}
 \hat \tau_{\tlP}(H_{\tlP}(\delta_{2} y) - T_{\tlP})
 k_{\tlP, a}(f, \delta_{1}x, \delta_{2}y).
 \]

 \begin{theoreme}\label{thm:cvgHe} Soit   $dx$ et $dy$ des mesures de Haar sur $G_E(\AAA)$ et $G'(\AAA)$. 
    \begin{enumerate}
   \item Il existe un point $T_+\in  \ago_{\tl}^+$ tel que pour tout $T\in T_++ \ago_{\tl}^+$, l'intégrale 
 $$ I_{a}^{\eta,T}(f) = 
 \int_{[G_{E}]} \int_{[G']}
 k_{a}^T(f,x, y) \, \eta(y)\, dx dy
 $$
 converge absolument.
 \item L'application $T\mapsto I_{a}^{ \eta,T}(f)$ est la restriction d'une fonction exponentielle-polynôme en $T$.
 \item Le terme purement polynomial de cette exponentielle-polynôme, noté $I_{a}^{\eta}(f)$, 
 est constant et est indépendant des choix auxiliaires autres que celui des mesures de Haar sur $G_E(\AAA)$ et $G'(\AAA)$. 
 \item La distribution $ I_{a}^{\eta} $ est 
invariante par translation  à gauche de  $G_{E}(\AAA)$ et $\eta$-équivariante par translation à droite de $G'(\AAA)$.
 \item Le support de la distribution $I_{a}^{\eta}$ est contenu dans $G'_{E, a}(\AAA)$.
   \end{enumerate}
 \end{theoreme}
 
 \begin{preuve}
Le théorème est démontré dans \cite{leMoi3}.
 \end{preuve}
 \end{paragr}

 \subsection{Densité}
 
 \begin{paragr}[Facteur $\Omega$.]\label{S:factOmega} ---  On suit à quelques variantes mineures près les constructions de \cite{Z1} section 2.4. 
Soit $\eta'$ un caractère de $E^{\times} \back \AAA_{E}^{\times}$   dont la restriction à $\AAA_{F}^{\times}$ est égale à  $\eta$.
 Pour tout $x \in X^{\rs}(\AAA)$ on pose:
 \[
 \Omega(x) := \eta'(\det(x)^{-[(n+1)/2]} \det (e_{0}, xe_{0}, \ldots, x^{n}e_{0})),
 \]
où $ [(n+1)/2] $ désigne la partie entière de $(n+1)/2$.  On a pour tous $x \in X^{\rs}(\AAA)$ et $ g \in G(\AAA) $ 
$$  \Omega(gxg^{-1}) =  \eta(g)\Omega(x) $$
et 
$$\Omega(x)=1 \quad \forall x \in  X^{\rs}(F).$$
 
 Pour tout $h = (g, \tlg) \in (G'_{E})^{\rs}(\AAA)$ 
 on pose (cf. \ref{eq:nu})
 \[
 \Omega(h) := 
 \begin{cases} \Omega(\nu(h)) \quad & \text{ si }n\text{ est pair }; \\
 \eta'(\det(g^{-1} \tlg))\, \Omega(\nu(h)) \quad &\text{ si } n\text{ est impair}.
 \end{cases}
 \]
 On a alors pour tous $h  \in (G'_{E})^{\rs}(\AAA)$, $ x\in G_E(\AAA) $ et $ y \in G'(\AAA) $   
$$ \Omega(x^{-1} h y) = \eta(y) \Omega(h)$$
et
$$\Omega(h)=1 \quad \forall h \in  (G'_{E})^{\rs}(F).$$

 Soit $S \subset \vc$ fini.  Les constructions et les propriétés ci-dessus  s'étendent à des éléments de  $X^{\rs}(\AAA_S)$ ou $(G'_{E})^{\rs}(\AAA_S)$. Il suffit dans les définitions de remplacer  $\eta'$ par sa restriction $\eta'_S$ à $(\AAA_S \otimes_F E) ^\times$.

 \end{paragr} 

 \begin{paragr}[L'ensemble $S_\natural$ de \og mauvaises\fg{} places.]\label{S:mauv} ---  On fixe un ensemble fini 
 $$
 S_{\natural} \subset \vc
 $$
qui contient  les places archimidiennes, celles qui se ramifient dans $E$ et tel que les caractères $\eta$ et $\eta'$ soient non ramifiés hors de $S_\natural$. On suppose également que les variétés considérées $G'$, $G'_E$, $G_{E}$, $ X $ etc. sont naturellement définies sur  l'anneau des entiers de $ F $ hors $ S_{\natural} $.
 \end{paragr}

\begin{paragr}[Décomposition de mesures de Haar.] ---
 Pour tout groupe $H \in \{G,\tlG, G_{E}\}$  et toute place $v \in \vc$ on fixe une mesure de Haar sur $H(F_{v})$  de sorte que pour tout $v \notin S_{\natural}$ le volume de $H'(\oc_{v})$ soit $1$. En prenant le produit, on en déduit une mesure de Haar sur $H(\AAA)$ et $H(\AAA_S)$. On met alors la mesure produit sur $G'(F_v)=G(F_v)\times \tlG(F_v)$. De même, on en déduit une mesure de Haar sur $G'(\AAA)$ et $G'(\AAA_S)$.
\end{paragr}

 \begin{paragr}[Intégrales orbitales semi-simples régulières locales.]\label{S:IOlocH} --- 
 Soit $a \in \tilde \Ac^{\rs}(\AAA_{S})$ et $f \in C_{c}^{\infty}(G'_{E}(\AAA_{S}))$.
 On introduit l'intégrale orbitale
 \begin{equation*}
I^{\Omega}_{a}(f)  = \int_{G_{E}(\AAA_{S})}\int_{G'(\AAA_{S})}
 f(g^{-1} \gamma h)\, \Omega(g^{-1}\gamma h)
 dh dg 
 \end{equation*}
où $\gamma$ est un élément quelconque  de  $G'_{E,a}(\AAA_{S})$ (qui est en fait une $G_E(\AAA_S)\times G'(\AAA_S)$-orbite).
 \end{paragr}

 \begin{paragr}[Le théorème de densité.]\label{S:instH} --- 
 Soit 
 $$C_{c}^{\infty}(G'_{E}(\AAA_S))_\eta\subset C_{c}^{\infty}(G'_{E}(\AAA_S))$$
 le sous-espace des fonctions $\eta$-instables au sens où leurs intégrales $I^{\Omega}_{a}$ s'annulent pour tout $a\in \tilde \Ac^{\rs}(\AAA_S)$. Cette définition ne dépend en fait que du caractère $\eta$ (le facteur $\Omega$ sera utile plus tard pour le transfert). Toute distribution qui s'annule sur $C_{c}^{\infty}(G'_{E}(\AAA_S))_\eta$ est dite  $\eta$-stable.
 
   \begin{theoreme}
     \label{thm:densiteH}
 Soit $a\in \tilde \Ac(F)$ et $f^S \in C_{c}^{\infty}(G'_{E}(\AAA^S))$ une fonction auxiliaire. La distribution
 $$f \in C_{c}^{\infty}(G'_{E}(\AAA_S)) \mapsto I_a^{\eta}(f\otimes f^S)$$
 est $\eta$-stable.
   \end{theoreme}
   \end{paragr}
   
   On démontre ce théorème au §\ref{ssec:demDens}.    La méthode consiste à se ramener au cas infinitésimal (cf. théorème \ref{thm:densite}). 
   
\subsection{Réduction au cas infinitésimal}

\begin{paragr}[Une légère extension des résultats de la partie \ref{partie:1}.]\label{S:extToEnd} ---  La décomposition $W=V\oplus Fe_0$ induit une écriture matricielle
$$Y=
\begin{pmatrix}
  A & b\\ 
c & d
\end{pmatrix}
$$
pour tout endomorphisme $Y\in \End_F(W)$. On en déduit, avec les notations de la section \ref{sec:prelim}, une décomposition en somme directe:
\begin{equation}\label{eq:endDec}
\End_F(W)= \tggo \oplus F(e_0^{*}\otimes e_0).
\end{equation}
Le groupe $G$ agit sur $\End_F(W)$ par conjugaison en laissant la droite $F(e_0^{*}\otimes e_0)$ fixe et en agissant sur $\tggo$ comme dans §\ref{S:tggo}. L'ouvert régulier semi-simple s'écrit alors
$$\End_F(W)^{\rs}=\tggo^{\rs} \oplus F(e_0^{*}\otimes e_0).
$$
Le quotient catégorique de $\End_F(W)$ par $G$ s'identifie avec $\Ac \times \AAA $, où $ \AAA $ est la droite affine sur $F$.  Il est évident que les résultats de la partie \ref{partie:1} s'étendent au cas de $G$ agissant sur $\End_F(W)$. Plus précisément, pour tout $\mu \in F$ et $f \in \Sc(\End_F(W)(\AAA))$, soit $f_{\mu} \in \Sc(\tggo(\AAA))$ défini par $f_{\mu}(X) = f(X + \mu)$ selon la décomposition \eqref{eq:endDec}. Soit $\tlP \in \fc^{\tlG}(B)$ et $(a, \mu) \in \Ac(F) \times F$. 
On définit alors 
\[
k_{\tlP, (a, \mu)} (f):= k_{ \tlP, a}(f_\mu)
\]
où le membre de droite est défini au §\ref{S:noyau}. 
On définit de même façon $k_{(a, \mu)}^{T}(f)$, pour $T \in \ago_{\tl}$ et la distribution
$I_{(a, \mu)}^{\eta, T}$. Les théorèmes \ref{thm:cv} et  \ref{thm:I} valent dans ce contexte et définissent des distributions $\eta$-équivariantes $I_{(a, \mu)}^{\eta}$ sur $\Sc(\End(W)(\AAA))$. Les définitions du paragraphe \ref{ssec:IOloc} se transportent aussi  sans grand changement ; on utilise le facteur de transfert défini pour tout $Y \in \tggo^{\rs}(\AAA_{S})$ par 
\[
\tilde \eta(Y) := \eta_{S}((-1)^{n}\det (e_{0}, Ye_{0}, \ldots, Y^{n}e_{0})),
\] 
où $\eta_S$ est la restriction de $\eta$ à $\AAA_S^\times$ et le déterminant est pris dans la base de $W$ obtenue à partir de la base de $V$ et $e_0$. Le signe n'est là que pour satisfaire l'identité  $ \tilde \eta(Y) = \tilde \eta(Y') $ si l'on écrit 
$ Y = (Y' ,d) $ selon la décomposition \eqref{eq:endDec} et  où le membre de droite est défini en §\ref{S:eta}. L'analogue du théorème \ref{thm:densite} est également vrai. 

En vertu de ces remarques, dans la suite on réserve la notation $\tggo$ pour $\End(W)$ (ou encore $\tggo$ est désormais l'algèbre de Lie de $\tlG$) et $\Ac$ pour le quotient catégorique de $\End(W)$ par $G$. On transposera la plupart des notations de la partie \ref{partie:1} en remplaçant simplement l'ancien $\tggo $  par le nouveau.
\end{paragr}

\begin{paragr}[L'isomorphisme de Cayley.]\label{S:cayley} --- Soit $\tau \in E^{\times}$ tel que $\sigma(\tau) = - \tau$. Soit $\tggo_\tau\subset \tggo$ et $\Ac_\tau\subset \Ac$ les ouverts complémentaires des fermés définis par $\det(Y^2-\tau^2)=0$. On note aussi 

$$\tggo^{\rs}_\tau=\tggo^{\rs}\cap \tggo_{\tau}.
$$
Soit 
$$\xi\in E^{1} = \{ \xi \in E\mid  \, \xi \sigma( \xi) = 1\}.$$
Soit $X_\xi$ et $\tilde{\Ac}_\xi$ les ouverts complémentaires des fermés définis par la condition $N_{E/F}(\det(g-\xi))=0$.

\begin{remarque}\label{rq:recouv}
Quand $\xi$ décrit $E^1$,   les ouverts   $\tilde{\Ac}_\xi$ recouvrent $\tilde{\Ac}$ (cf. lemme 3.4 de \cite{Z1}).
\end{remarque}

Pour tout $Y\in \tggo_\tau$, soit
\begin{equation}\label{eq:cayleyDef}
\kappa(Y)=\kappa_{\xi}(Y) = -\xi (1+ \tau^{-1}Y)(1-\tau^{-1}Y)^{-1}.
\end{equation}
Le lemme suivant est une conséquence du lemme 3.5 de \cite{Z1}.

\begin{lemme}\label{lem:kappaIsoG} 
Le morphisme $\kappa$ induit un $F$-isomorphisme $G$-équivariant de $\tggo_{\tau}$ sur $X_{\xi}$ et un isomorphisme (encore noté $\kappa$)
 \[ 
 \kappa  : \Ac_{\tau} \simeq\tilde \Ac_{\xi}.
 \]
\end{lemme}

Soit $ G'_{E, \xi} \subset G'_E$ l'ouvert obtenu comme image inverse de $\tilde \Ac_{\xi}$ par le morphisme canonique.
\end{paragr}

\begin{paragr}[Compatibilité des facteurs $ \Omega $ et $ \tilde \eta $.] --- 
 Avec les notations du §\ref{S:eta} on a:
 
 \begin{lemme}\label{lem:factCayl} Il existe $\mu_n \in E^{\times}$ qui dépend de $\tau$ et $\xi$  tel que pour 
 tout $ Y \in \tggo_{\tau}^{\rs}(\AAA_{S}) $ on a
 \[ 
 \Omega(\kappa(Y)) = \begin{cases}
 \eta_{S}'(\mu_n) \, \tilde \eta(Y) \quad &\text{ si }n\text{ est pair}, \\
  \eta_{S}'(\mu_n) \, \eta'_{S}(\det (Y - \tau))  \,\tilde \eta(Y) \quad &\text{ si }n\text{ est impair}.
 \end{cases}
 \]
 \end{lemme}
 
 \begin{preuve}
C'est un  calcul direct. 
 \end{preuve}
\end{paragr}

\begin{paragr}[Correspondance des fonctions.]\label{S:fonHggo} --- On suit \cite{Z1} section 2.1. Soit $ S \subset \vc_{F} $ fini et $f \in C_{c}^{\infty}(G'_{E}(\AAA_{S}))$.  Soit $x \in X(\AAA_{S})$ et $(h,\tlh)\in G'_E(\AAA_S)$ tel que $x=\nu((h,\tlh))$ (cf. \eqref{eq:nu}). On  définit (on vérifie que l'intégrale de droite est indépendante du choix de  $(h,\tlh)$)
 \[
f^{X}(x) = 
\begin{cases}
\displaystyle \int_{G_{E}(\AAA_{S})}\int_{\tlG(\AAA_{S})}f(g(h, \tlh \tlg)) \, \eta'(\det(\tlh \tlg)) \, dgd\tlg 
\quad &\text{ si }n\text{ est impair}, \\
\displaystyle \int_{G_{E}(\AAA_{S})}\int_{\tlG(\AAA_{S})}f(g (h, \tlh \tlg))\, dg d\tlg 
\quad &\text{ si }n\text{ est pair}.
\end{cases}
 \]
Alors $f\mapsto f^X$ définit une application linéaire surjective $  C_{c}^{\infty}(G'_{E}(\AAA_{S}))\to C_{c}^{\infty}(X(\AAA_{S}))$. 
 
 On note le résultat suivant, conséquence du lemme 3.1 de \cite{AFL}.
 
 \begin{lemme}\label{lem:charHtoX}
 Soit $ v $ une place en dehors de l'ensemble $ S_{\natural} $. 
 Alors $ \mathbf{1}_{G'_{E}(\oc_v)}^X = \mathbf{1}_{X(\oc_v)} $.
 \end{lemme}

 Supposons de plus que $f \in C_{c}^{\infty}(G'_{E, \xi}(\AAA_{S}))$. Dans ce cas, on a $ f^{X} \in  C_{c}^{\infty}(X_{\xi}(\AAA_{S}))$.
 On définit alors $\phi_{f} \in C_{c}^{\infty}(\tggo_{\tau}(\AAA_{S}))$ par

 \begin{equation}
   \label{eq:phi-f}
   \phi_{f}(Y) =
\begin{cases}
  \dfrac{ \Omega(\kappa(Y))}{\tilde \eta(Y)}
f^{X}(\kappa(Y)) \text{ si }  Y \in \tggo_{\tau}(\AAA_{S}) \ ;\\
0  \text{ sinon}.
\end{cases}
\end{equation}

\begin{lemme}\label{lem:comp-kappa}
Pour tout $a \in \Ac^{\rs}(\AAA_{S})$
   \begin{equation}\label{eq:IOHggo}
 I_{a}^{\tilde \eta}(\phi_{f}) = 
 \begin{cases}
 I_{\kappa(a)}^{\Omega}(f) \quad & \text{si }a \in \Ac_{\tau}(\AAA_{S}), \\
 0 \quad &\text{ sinon},
  \end{cases}
\end{equation}
où  la distribution à gauche est définie par \eqref{eq:IOPloc2}. 
\end{lemme}

\begin{preuve}
  La vérification ne pose aucune difficulté (cf. aussi \cite{Z1} section 2.1).
\end{preuve}
  
Le lemme précédent a un pendant global que voici.
  
 \begin{lemme}\label{lem:IaHggo} Soit $a \in \Ac_\tau(F)$. Soit $S\subset\vc$ fini contenant l'ensemble $S_{\natural}$ de §\ref{S:mauv} tel que 
   \begin{enumerate}
   \item Les éléments $\tau$, $\xi$ et $\mu_n$ (cf. lemme \ref{lem:factCayl}) soient des unités hors $S$ ;
   \item $a\in  \Ac_\tau(\oc_v)$ pour tout $v\notin S$. 
   \end{enumerate}

 Alors, $\tla=\kappa(a)\in  \tilde \Ac_{\xi}(\oc_v)$  pour tout $v\notin S$. De plus,   pour tout $f \in C_{c}^{\infty}(G'_{E, \xi}(\AAA_{S}))$, on a
 \[
 I_{\tla}^{\eta}(f \otimes \mathbf{1}_{G'_{E}(\oc^{S})}) =  I_{a}^{\eta}(\phi_{f}  \otimes \mathbf{1}_{\tggo(\oc^{S})})
 \]
 où la distribution $ I_{a}^{\eta} $ est essentiellement celle définie au §\ref{S:Ia} (cf. §\ref{S:extToEnd}).
  \end{lemme}
  
  \begin{preuve} Pour alléger les notations, posons pour la démonstration $f'= f \otimes \mathbf{1}_{G'_{E}(\oc^{S})}$ et $\phi'= \phi_{f}  \otimes \mathbf{1}_{\tggo(\oc^{S})}$. 
  Il résulte des théorèmes \ref{thm:cvgHe} et \ref{thm:cv} 
  et de la définition \ref{eq:kT} qu'il 
  suffit de montrer que pour tout 
   $\tlP \in \fc^{\tlG}(B)$
   et tout $g \in [G]$ qu'on a:
   \[
   \int_{[G_{E}]} \int_{[\tlG]}
   \sum_{
   \begin{subarray}{c}
   \delta_{1} \in P_E(F)\back G_E(F)
   \\
   \delta_{2} \in \tlP(F)\back \tlG(F)
   \end{subarray}}
   k_{\tlP, \tla}(f',\delta_{1}x ,(g,\delta_{2}\tlg))
   \, \eta(g, \tlg)
   \,d\tlg dx = 
   k_{\tlP, a}(\phi',g)
   \eta(g)
   ,
   \]
   où $k_{\tlP, \cdot }(\phi')$ est défini au §\ref{S:noyau} compte tenu du § \ref{S:extToEnd}. Cette égalité 
   découle du lemme \ref{lem:charHtoX} et des conditions imposées sur l'ensemble $ S $. 
   Les détails se trouvent dans \cite{leMoi3}, section 5.
  \end{preuve}
  \end{paragr}

    \begin{paragr}[Fonction $\theta_{S}$.]\label{S:thetaSaxi} --- 
    Soit $S \subset \vc$ un ensemble fini de places de $F$. Soit $\tla \in \tilde \Ac_{\xi}(\AAA_{S})$. Pour tout $v\in S$, soit $\theta_{v} \in C_{c}^{\infty}(\tilde \Ac_{\xi}(F_v))$ une fonction  qui vaut $1$ dans un voisinage de $\tla_v$. Soit
$$\theta_{S}=  \otimes_{v \in S}\theta_v \in C_{c}^{\infty}(\tilde \Ac_{\xi}(\AAA_{S})).$$
    \end{paragr}

    \subsection{Preuve du théorème \ref{thm:densiteH}}\label{ssec:demDens}

    \begin{paragr} --- On se place sous les hypothèse du  théorème \ref{thm:densiteH}. Soit $S\subset \vc$ fini et $ f \in C_{c}^{\infty}(G'_{E}(\AAA_{S}))_{\eta} $. Soit $\tla\in \tilde \Ac(F) $. Quitte à changer $\xi  \in E^{1} $ on peut et on va supposer que $\tla\in \tilde{\Ac}_\xi(F) $ (cf. remarque \ref{rq:recouv}). Soit $a\in \Ac_\tau(F)$ tel que $\kappa(a)=\tla$. On peut élargir $S$ en ajoutant des composantes à $f$ sans que cela altère l'instabilité de $f$. Ainsi on peut et on va  supposer que $S$ vérifie les conditions du lemme \ref{lem:IaHggo}. Quitte à encore élargir $S$, on voit qu'il suffit de montrer qu'on a 
$$ I_{\tla}^{ \eta}(f \otimes \mathbf{1}_{G'_{E}(\oc^{S})}) = 0.$$ 
Soit   $\theta_{S}$ comme au §\ref{S:thetaSaxi}. Par composition avec le morphisme canonique, on voit $\theta_S$ comme une fonction invariante sur $G'_{E}(\AAA_S)$. Soit
     \[ 
    \bar f = f \cdot \theta_{S}.
     \]
Cette fonction vérifie les propriétés suivantes, certaines héritées de $f$  :
\begin{itemize}
\item $ \bar f \in  C_{c}^{\infty}(G'_{E,\xi}(\AAA_S)) $ ;
\item $\bar f$ est $\eta$-instable ;
\item $f-\bar f$ est nulle au voisinage de $G'_{E,\tla}(\AAA_S)$
\item  $I_{\tla}^{ \eta}(f \otimes \mathbf{1}_{G'_{E}(\oc^{S})}) = I_{\tla}^{ \eta}(f \otimes \mathbf{1}_{G'_{E}(\oc^{S})}) $ (qui résulte de la propriété de support de   $I_{\tla}^{ \eta}$ (cf. théorème \ref{thm:cvgHe}).
\end{itemize}
Quitte à remplacer $f$ par $\bar f$, on peut et on va supposer que $f\in   C_{c}^{\infty}(G'_{E,\xi}(\AAA_S)) $.     Soit $ \phi_{f} \in C_{c}^{\infty}(\tggo_{\tau}(\AAA_{S})) $ défini par \eqref{eq:phi-f}. Il résulte du lemme \ref{lem:comp-kappa} que les intégrales orbitales semi-simples régulières de $\phi_f$ s'annulent.  On a alors en combinant le lemme \ref{lem:IaHggo} et le  théorème \ref{thm:densite}
     \[
     I_{\tla}^{\eta}(f \otimes \mathbf{1}_{G'_E(\oc^{S})}) =      I_{a}^{\eta}(\phi_{ f} \otimes \mathbf{1}_{\tggo(\oc^{S}})=0,
     \]
ce qu'il fallait démontrer.
      \end{paragr}

\section{Distributions géométriques pour $ U_n \times U_{n+1} $}\label{sec:densiteGpU}
 \subsection{Préliminaires}

 \begin{paragr} 
 On reprend les notations du §\ref{ssec:prelimS}.
  \end{paragr}
 
 \begin{paragr}\label{S:tlU}
Pour toute  forme $\sigma$-hermitienne non dégénérée $\Phi$ sur $V_{E}$,  on définit la forme $\sigma$-hermitienne non-dégénérée 
$\tilde \Phi$ sur $W_{E}$ de sorte que la restriction de $\tilde \Phi$ à $V_{E}$ soit égale à  $\Phi$ et qu'on ait $\tilde \Phi(e_{0}, e_{0}) = 1$  et $\tilde \Phi(V, e_{0}) = 0$. On dispose d'une inclusion de groupes unitaires $U=U_\Phi \subset \tlU=U_{\tilde{\Phi}}$ où $U$ s'identifie au sous-groupe de $\tlU$ qui agit trivialement sur $e_{0}$ et stabilise $V_{E}$.  On introduit aussi 
$$U'=U'_\phi= U_\Phi\times U_{\tilde{\Phi}}
$$
et $U$ se plonge dans $U'$ diagonalement.
\end{paragr}
 
 \begin{paragr}[Actions et quotients.]\label{S:acQuotUU} --- 
Le groupe $U$ agit sur $\tlU$ par conjugaison. On note $\tlU // U$ le quotient catégorique.
  
On a aussi l'action du groupe $U \times U$ sur $U'$ à droite. 
Si $ (g_{1}, g_{2}) \in U \times U$ et $ (h, \tlh) \in U' $, alors $  (h,\tlh) \cdot(g_{1}, g_{2}) = (g_{1}^{-1}hg_{2},g^{-1}_1 \tlh g_{2})$. 
On note:
\begin{equation}\label{eq:nuU}
\nu : U' \mapsto \tlU, \quad \nu((h , \tlh)) = h^{-1}\tlh.
\end{equation}
L'application $ \nu $ induit un isomorphisme du  quotient géométrique $ U' / (U \times \{1\}) $ sur $\tlU$. L'action du second facteur  $\{1\}\times U$ à la source donne l'action de $U$ par conjugaison sur le but. On obtient ainsi une  identification des quotients catégoriques 
$$ 
U' // (U \times U) \cong \tlU // U.
$$
 \end{paragr}
 
\begin{paragr}[Identification de quotients catégoriques.]\label{S:quottlU} --- 
On utilisera les notations de la section \ref{sec:densiteS}.

\begin{lemme}\label{lem:Atilde}
  On a un isomorphisme canonique $\tilde \Ac =  X //G \simeq \tlU // U$. De plus, deux
  éléments $Y_1\in X(F)$ et $ Y_2\in \tlU(F)$ ont même image dans le quotient catégorique si
  et seulement s'ils ont le même polynôme caractéristique et si on a 
$$e_0^*(Y_1^i e_0)= \tilde \Phi(e_0,Y_2^ie_0)$$
pour tout $i$.
\end{lemme}

\begin{preuve}
On a (cf. §\ref{S:varX})
$$X \times_F E\simeq GL_E(W_E)$$
et
$$\tlU \times_F E=\{(g_1,g_2) \in GL_E(W_E)^2  \mid \,^tg_2 \tilde{\Phi} g_1=\tilde{\Phi}\}.
$$
On a donc un $E$-isomorphismes $ X \times_F E \simeq\tlU \times_F E$ qui est induit par $g\mapsto (g, \,^t(\tilde{\Phi} g^{-1} \tilde{\Phi}^{-1}))$ et qui est $GL_E(V_E)$-équivariant. On a donc un isomorphisme au niveau des quotients catégoriques après extension des scalaires à $E$. Il suffit de vérifier que cet isomorphisme se descend à $F$. En fait, les variétés $X$ et $\tilde{U}$ sont des formes l'une de l'autre et que  la torsion galoisienne est donnée, à une conjugaison près par un élément de $GL_E(V_E)$, par  la transposition. Cette conjugaison agit évidemment trivialement sur les quotients catégoriques. Pour conclure, il suffit de montrer que la transposition agit elle-aussi trivialement sur les quotients catégoriques. C'est évident au vu de la description du quotient rappelé au §\ref{S:varX}. Cette même description et le fait que $\Phi(e_0,\cdot)=e_0^*$ donne la dernière assertion.
\end{preuve}

En utilisant le lemme ci-dessus et la discussion du §\ref{S:acQuotUU}, on a donc des morphismes
$$U'\to \tilde\Ac$$
et
$$\tlU \to \tilde \Ac $$
qui sont invariants respectivement sous l'action de $U\times U$ et $U$. Pour tout $\tla\in \tilde{\Ac}$, soit $U'_{\tla}$ et $\tlU_{\tla}$ les fibres au-dessus du $a$ des morphismes ci-dessus. Lorsque de plus $\tla\in \tilde{\Ac}^{\rs}(F)$, les ensembles  $U'_{\tla}(F)$ et $\tlU_{\tla}(F)$ sont soit simultanément vides soit ce sont respectivement une $U(F)\times U(F)$-orbite et une $U(F)$-orbite (cf. \cite{AFL} lemme 2.3).
\end{paragr}

\subsection{Les distributions géométriques}
 
\begin{paragr}[Sous-groupes paraboliques.] --- Soit $P_0$ un sous-groupe parabolique minimal  de $ U $ de décomposition de Levi $ M_{0}N_{0} $. On note $ \ago_{0}^{+} := \ago_{P_{0}}^{+} $.  Pour tout $P \in \fc^{U}(P_{0})$, soit $\tlP$ le sous-groupe parabolique de $\tlU$ qui est le stabilisateur dans $\tlU$ du drapeau isotrope qui définit $P$ et 
$$P'=P\times\tlP.$$
On a une décomposition de Levi $P'=M_{P'}N_{P'}$ où $M_{P'}$ est le facteur de Levi qui contient $M_0\times M_0$. On munit $N_{P'}(\AAA)$ de la mesure de Haar qui donne le volume $1$ à $N_{P'}(F)\back N_{P'}(\AAA)$.
\end{paragr}

\begin{paragr}[Le noyau $k_{P}$.] --- 
Soit $ a \in \tilde \Ac(F) $, $P \in \fc^{U}(P_{0})$ et $ f \in C_{c}^{\infty}(U'(\AAA)) $. Pour tous $x,y \in  N_{P'}(\AAA)M_{P'}(F) \back U'(\AAA)$
Soit
 \[
k_{P, a}(f,x,y) := 
\sum_{\gamma \in M_{P',a} (F)}\int_{N_{P'}(\AAA)}f(x^{-1} \gamma n y)dn, 
 \]
où  $M_{P',a}=M_{P'}\cap U'_a$. Notons qu'on a  $ k_{P ,a} =0$ si $ U'_{a}(F) = \emptyset$.
\end{paragr}

\begin{paragr}[Les distributions $I_{a}^{U'}$]\label{S:IUGp} ---  Avec les choix et les notations de la section \ref{ssec:distGlobale}, pour $T \in \ago_{0}$ on pose
 pour tous $ x,y \in [U]$
\begin{equation*}\label{eq:katDefU}
 k_{a}^{T}(x, y) = 
 \sum_{P \in \fc^{U}(P_{0})}
 \varepsilon_{P}^{U}
 \sum_{
 \delta_{1}, \delta_{2} \in P(F)\back U(F)}
 \hat \tau_{P}(H_{P}(\delta_{2} y) - T_{P})
 k_{P, a}(f,\delta_{1} x, \delta_{2}y).
  \end{equation*}
 
 Le théorème suivant est démontré dans \cite{leMoi3}, section 6. 
 
 \begin{theoreme}\label{thm:cvgU'} Fixons une mesure de Haar sur $U(\AAA)$.
   \begin{enumerate}
   \item Il existe un point $T_+\in  \ago_{0}^+$ tel que pour tout $T\in T_++ \ago_{0}^+$, l'intégrale 
 $$
 I^{\Phi,T}_a(f) =\int_{[U]}\int_{[U]} k^T_a(f,x,y) \, dx dy
 $$
 converge absolument.
 \item L'application $T\mapsto I^{\Phi,T}_a(f)$ est la restriction d'une fonction exponentielle-polynôme en $T$.
 \item Le terme purement polynomial de cette exponentielle-polynôme, 
 noté $ I^{\Phi}_a(f) $, 
 est constant et ne dépend pas des choix auxiliaires autres que celui de la mesure de Haar sur $U(\AAA)$.
 \item La distribution $ I^{\Phi}_a$ est $ U(\AAA) \times U(\AAA) $-invariante.
 \item Le support de $ I_{a}^{\Phi} $ est inclus dans $ U'_{a}(\AAA) $.
   \end{enumerate}
 \end{theoreme}
 \end{paragr}

 \subsection{Densité}

\begin{paragr}[Mesures de Haar.]\label{S:reseauUtlU} --- 
On peut naturellement définir des modèles pour $U$, $\tlU$ et donc $U'$ sur un anneau d'entiers hors un ensemble. Pour presque tout $v\in \vc$, le groupe $U(\oc_v)$ (resp. $\tlU(\oc_v)$) est alors le stabilisateur du réseau autodual engendré par la base de $V$ (resp. et $e_0$) qu'on a fixée (cf. §\ref{S:baseFix}). On a $U'(\oc_v)=U(\oc_v)\times \tlU(\oc_v)$.
Pour toute place $ v \in \vc $ on fixe une mesure de Haar sur $ U(F_{v}) $ telle que pour presque tout $ v \in \vc $ le volume de $U(\oc_{v}) $ soit $ 1 $. On en déduit la mesure de Haar sur $ U(\AAA) $.
\end{paragr}

  \begin{paragr}[Intégrales orbitales semi-simples régulières locales $I_a^\Phi$.]\label{S:IOlocU} --- 
Soit $S \subset \vc$ fini et $a \in \tilde \Ac^{\rs}(\AAA_{S})$. On a l'alternative suivante : soit la fibre $U'_a(\AAA_S)$ est vide auquel cas on pose $I_{a}^{\Phi}=0$ soit cette fibre est une $U(\AAA_{S})\times U(\AAA_{S})$-orbite d'un élément $\gamma$. On pose alors pour tout $f \in C_{c}^{\infty}(U'(\AAA_{S}))$ 
\begin{equation*}
I_{a}(f) = I_{a}^{\Phi}(f) = \int_{U(\AAA_{S})}\int_{U(\AAA_{S})}
f(h_1^{-1} \gamma h_2)\,dh_1dh_2 
\end{equation*}
et cette définition est évidemment indépendante du choix de $ \gamma \in U'_{a}(\AAA_{S})$.
\end{paragr}

 \begin{paragr}[Le théorème de densité.]\label{S:fonInstU} --- 
Soit l'espace des fonctions instables 
 $$
 C_{c}^{\infty}(U'(\AAA_S))_0 = 
 \{
 f \in C_{c}^{\infty}(U'(\AAA_S)) \mid
 I_{a}(f)=0 \ \forall a\in \tilde \Ac^{\rs}(\AAA_S) \}.
$$
Il est indépendant du choix de la mesure de Haar.

   \begin{theoreme}
     \label{thm:densiteUU}
 Soit $a\in \tilde \Ac(F)$ et $f^S\in C_{c}^{\infty}(U'(\AAA^S))$. La distribution
 $$f \in C_{c}^{\infty}(U'(\AAA_S)) \mapsto I_a^{\Phi}(f\otimes f^S)$$
est stable au sens où elle  s'annule sur  le sous-espace $ C_{c}^{\infty}(U'(\AAA_S))_0$.
   \end{theoreme}
   \end{paragr}
   
   La preuve de ce théorème, essentiellement identique à celle du théorème    \ref{thm:densiteH}, procède par réduction au cas infinitésimal (cf. théorème \ref{thm:densite}). Pour s'en convaincre, et aussi  pour les besoins de la section \ref{sec:factorisationGpU},    on va introduire les constructions qui interviennent dans cette réduction. 
   
\subsection{Réduction au cas infinitésimal}
\begin{paragr}[Une légère extension des résultats de la partie \ref{partie:2}.] ---\label{S:ext-endII} Comme au §\ref{S:extToEnd}, il est aisé d'étendre les résultats de la partie  \ref{partie:2} au cas de $U$ agissant sur l'espace des endomorphismes de $W_E$ qui sont auto-adjoints, cet espace se décomposant en $\tugo\oplus Fe_0\otimes e_0^*$ (où $\tugo$ est défini  à la section \ref{sec:PrelimU}). 
Désormais on réserve la notation 
\[
\tugo = \tugo_{\Phi} 
\]
pour cet espace d'endomorphismes auto-adjoints (qu'on ne confondra pas avec l'algèbre de Lie de $\tlU$, les deux étant néanmoins $F$-isomorphes comme représentation linéaire de $U$). On utilisera les constructions et les résultats  de la partie \ref{partie:2} dans cette situation sans plus de  commentaire. Par exemple, la notation $\Ac$ désigne le quotient catégorique $\tugo//U$ qui s'identifie d'ailleurs au quotient $\tggo//G$ du  §\ref{S:extToEnd}.
\end{paragr}

\begin{paragr}[Mauvaises places.]\label{S:mauvU} --- 
Soit 
$$
S_{\natural}^{U} \subset \vc
$$
fini, contenant les places 
infinies, celles ramifiées dans $E/F$ et tel pour toute place $v \in \vc \smallsetminus S_{\natural}^{U}$ les groupes $U(\oc_v)$ et $\tlU(\oc_v)$ sont les stabilisateurs des réseaux de $V_E\otimes_F F_v$ et $W_E\otimes_F F_v$ considérés au §\ref{S:reseauUtlU}. On note $\tugo(\oc_v)$ le stabilisateur dans $\tugo(F_v)$ du réseau considéré dans $W_E\otimes_F F_v$.
\end{paragr}

\begin{paragr}[Les variétés $\tlU_{\xi}$ et $ U'_{\xi} $.] --- Reprenons les notations du §\ref{S:cayley}. Soit $\tlU_{\xi}\subset \tlU$ et $ U'_{\xi} \subset U'$ les ouverts, images inverses de $ \tilde \Ac_{\xi} $ par le morphisme  canonique. Quand $\xi\in E^1$, ceux-ci recouvrent $\tlU$ et $U'$. Soit $\tugo_\tau\subset\tugo$  l'image inverse de $\Ac_\tau$. Par la formule \eqref{eq:cayleyDef}, on obtient un isomorphisme $U$-équivariant, noté encore $\kappa$, de $\tugo_\tau$ sur $\tlU_\xi$ (cf.  lemme 3.5 de \cite{Z1}). 
\end{paragr}

  \begin{paragr}[Correspondance des fonctions.]\label{S:fonUugo}  --- 
   Soit $ S \subset \vc $ fini. On définit un morphisme surjectif
$$f \in C_{c}^{\infty}(U'(\AAA_S)) \mapsto f^{\tlU}\in C_{c}^{\infty}(\tlU(\AAA_S))
$$
par la formule    
   \[
  f^{\tlU}(\tlx) = \int_{U(\AAA_S)}f(g^{-1}(h, \tlh))dg, \quad \tlx \in \tlU(\AAA),
   \]
   où, avec la notation de \eqref{eq:nuU}, $ (h, \tlh) \in U'(\AAA_S)$
   est tel que $ \nu(h,\tlh) = \tlx$. 
 
 Pour tout $f \in C_{c}^{\infty}(U'(\AAA_S)) $, soit  $\phi_{f} \in C_{c}^{\infty}(\tugo_{\tau}(\AAA_{S}))$ la fonction déterminée par 
 \begin{equation}
   \label{eq:phi-fU}
      \phi_{f}(Y) = f^{\tlU}(\kappa(Y)). 
    \end{equation}
pour tout  $Y \in \tugo_{\tau}(\AAA_{S})$.
   
 En utilisant la définition \eqref{eq:IOPloc2U}, 
 on a pour tout $a \in \Ac^{\rs}(\AAA_{S})$ 
 \begin{equation}\label{eq:IOUugo}
 I_{a}(\phi_{f}) = 
 \begin{cases}
 I_{\kappa(a)}(f), \quad & \text{si }a \in \Ac_{\tau}(\AAA_{S}), \\
 0 \quad &\text{ sinon}.
  \end{cases}
 \end{equation}
  
  On a aussi l'analogue du lemme \ref{lem:IaHggo} suivant.
  
 \begin{lemme}\label{lem:IaUugo} Soit  $a \in \Ac_\tau(F)$ et $S\subset\vc$ fini contenant l'ensemble $S_{\natural}^U$  du §\ref{S:mauvU} qui vérifient les hypothèses  1 et 2 du lemme \ref{lem:IaHggo}.
Pour tout $f \in C_{c}^{\infty}(U'_{\xi}(\AAA_{S}))$,  on a
 \[
 I_{\kappa(a)}^{\Phi}(f  \otimes \mathbf{1}_{U'(\oc^{S})}) = 
 I_a^{U}(\phi_{f} \otimes \mathbf{1}_{\tugo(\oc^{S})}),
 \]
 où la distribution à droite est celle définie au théorème \ref{thm:IU}.
  \end{lemme}
  
  \begin{preuve}
  La preuve est similaire à celle du lemme \ref{lem:IaHggo} est résulte 
  essentiellement des résultats de \cite{leMoi3}. 
  \end{preuve}
 \end{paragr}

  \begin{paragr}[Preuve du théorème \ref{thm:densiteUU} ] ---   Ces définitions et propriétés posées, la preuve du théorème \ref{thm:densiteUU} consiste à se ramener au théorème \ref{thm:densiteU} et est alors en tout point semblable à celle   du théorème \ref{thm:densiteH} donnée au §\ref{ssec:demDens}. On laisse les détails au lecteur.   
  \end{paragr}

\section{Le théorème de transfert}\label{sec:factorisationGpU}

Dans cette section finale, on énonce et démontre le résultat principal de l'article (cf. théorème \ref{thm:transfertGp}). On reprend les notations des sections \ref{sec:densiteS} et \ref{sec:densiteGpU}. On va suivre de près le formalisme développé dans les sections \ref{ssec:thlin} et \ref{ssec:thuni}.

\subsection{Factorisation de distributions}\label{ssec:thlinS}

\begin{paragr}
Soit $ v \in \vc $. Avec les notations du §\ref{S:instH}
on pose 
$$\Ic(\eta_v) =  
C_{c}^{\infty}(G'_{E}(F_v))/ C_{c}^{\infty}(G'_{E}(F_v))_{\eta}.
$$
Pour tout ensemble fini $S$ de places de $F$, soit 
$$\Ic(\eta)_S=\otimes_{v\in S} \Ic(\eta_v)$$
le produit tensoriel. Pour tous  ensembles finis $S\subset S'$ de places de $F$, contenant les places archimédiennes, on a une application 
$\Ic(\eta)_S \to \Ic(\eta)_{S'}$ donnée par 
$f \mapsto f\otimes (\otimes_{v\in S'\setminus S} \mathbf{1}_v)$ où, pour toute place non-archimédienne de $F$, 
on note $\mathbf{1}_{v}$ l'image de la fonction caractéristique de $G'_{E}(\oc_v)$.

Soit 
$$\Ic(\eta)= \varinjlim_{S}  \Ic(\eta)_S.$$
\end{paragr}

\begin{paragr} 
Le choix de facteur $\Omega$ au §\ref{S:factOmega} et des mesures de Haar permet de définir des intégrales orbitales locales semi-simples régulières 
(cf. §\ref{S:IOlocH}). Pour tout $S \subset \vc$ fini, 
on a  une application linéaire injective
 \begin{equation}
    \label{eq:IetaSH}
    I^{\Omega}_S: \Ic(\eta)_S \to  \mathcal{C}^\infty(\tilde \Ac^{\rs}(\AAA_S))
  \end{equation}
  donnée par 
$$I^{\Omega}_S(\otimes_{v \in S}f_{v}):  (a_{v})_{v \in S} \mapsto \prod_{v \in S}I_{a_{v}}^{\Omega}(f_{v}).$$
Soit
 $$\mathcal{C}^\infty(\tilde \Ac^{\rs})=\varinjlim_S  \mathcal{C}^\infty(\tilde \Ac^{\rs}(\AAA_S))$$
où la limite est prise sur les ensembles finis $S\subset\vc$ et les applications de transition sont données pour $S\subset S'$ (et $S$ assez grand contenant les places archimédiennes) par
$$
\phi \mapsto \phi  \cdot  
I_{S' \smallsetminus S}^{\Omega}( \otimes_{v\in S' \smallsetminus S}\mathbf{1}_{v}).$$
On en déduit une application injective:
$$\Orb^{\Omega}: \Ic(\eta) \to \mathcal{C}^\infty(\tilde \Ac^{\rs}).$$
Comme au §\ref{S:Ihc}, le théorème de densité \ref{thm:densiteH} implique que pour tout $a\in \Ac(F)$  les distributions globales $I_a^{\eta}$ de la section \ref{sec:densiteS} définissent des formes linéaires sur $\Ic(\eta)$ qu'on note encore $I_a^{\eta}$.
\end{paragr}

\begin{paragr}[Groupes unitaires.]\label{S:reseauGp} ---
On adapte maintenant dans le cas des groupes  les constructions du §\ref{ssec:thuni} et §\ref{S:reseau}. On note $\hc$ le groupoïde des fomes hermitiennes sur $V_E$ (cf. §\ref{S:hc}) ; on dispose de  son pendant local $\hc_v$ pour $v\in S$. Pour tout $\Phi\in \hc$,  (resp. $\Phi\in \hc_v$), on dispose des $F$-groupes  (resp. $F_v$-groupes $U_\Phi\subset \tlU_{\Phi}$ et  $ U_\Phi\subset  U'_{\Phi} = U_\Phi \times \tlU_{\Phi}$ (cf.  §\ref{S:tlU}). 

Pour tout $\Phi\in \hc_v$, le groupe $U_\Phi(F_v)$ est muni d'une mesure de Haar qui vérifie les conditions de compatibilité du §\ref{S:reseau}. Cela munit $U_\Phi(\AAA)$ d'une mesure de Haar pour tout $\Phi\in \hc$. 

Soit $v$ une place non-archimédienne, non ramifiée dans $E$. Si $ \Phi_{v} \in \hc_v$ admet un réseau auto-dual dans $ V_{E_{v}} $, on en fixe  un noté $V(\oc_{E_v})$. Soit  $W(\oc_{E_v})= V(\oc_{E_v}) \oplus \oc_{v}e_{0} $ ; c'est un réseau auto-dual dans $ W_{E_{v}} $. Soit $U(\oc_v)$ et $ \tlU(\oc_{v}) $ les centralisateurs respectifs de ces réseaux dans $U(F_v)$ et $ \tlU(F_{v}) $. Soit  $ U'(\oc_v) = U(\oc_v) \times \tlU(\oc_v) $.
\end{paragr}

\begin{paragr}[Espace $\Ic(\hc_v)$.] ---\label{S:IvdefU}
Soit $v\in \vc$ et $\Phi\in \hc_v$. On note $U$, $\tlU,\ldots $ au lieu de $U_\Phi$, $\tlU_\Phi$ etc. 
Avec les notations du §\ref{S:fonInstU} soit
$$
\Ic(\Phi) :=C_{c}^{\infty}(U'(F_{v}))/C_{c}^{\infty}(U'(F_{v}))_0,
$$
qui ne dépend que de la classe d'isomorphisme de $ \Phi $.

Soit 
$$\Ic(\hc_v)=  \bigoplus_{\Phi\in |\hc_v|} \Ic(\Phi).
$$

Soit $v\in \vc$ une place finie et non-ramifiée dans $E$. On distingue deux cas :
\begin{enumerate}
\item La place $v$ est inerte dans $E$. Dans ce cas $\Ic(\hc_v)=  \Ic(\Phi_0)\oplus \Ic(\Phi_1)$ (cf. §\ref{S:reseauGp}).  
On pose $\mathbf{1}_v$ l'image dans  $\Ic(\hc_v)$ du couple $(\mathbf{1}_{U'(\oc_v)},0)$.
\item La place $v$ est décomposée dans $E$. Alors 
$|\hc_v|$ est réduit à un singleton $\{\Phi\}$ et $\Ic(\hc_v)=  \Ic(\Phi)$. Soit  $\mathbf{1}_v$ l'image dans  $\Ic(\hc_v)$ de $\mathbf{1}_{U'(\oc_v)}$.
\end{enumerate}

Pour tout $f\in \Ic(\hc_v)$ et $\Phi\in \hc_v$, soit $f_\Phi\in \Ic(\Phi)$ la composante de $f$ sur  $\Ic(\Phi)$. Par le choix de la  mesure de Haar sur $U(F_{v})$ (cf. §\ref{S:reseauGp}), on dispose d'intégrales orbitales locales semi-simples régulières (cf. §\ref{S:IOlocU}) et donc d'une application linéaire
  \begin{equation}
    \label{eq:IvGp}
    I_v:\Ic(\hc_v) \to \mathcal{C}^\infty(\tilde \Ac^{\rs}(F_v))
  \end{equation}
qui à $f\in  \Ic(\hc_v)$ associe l'application 
$$a \mapsto 
I_{a}^{\hc_{v}}(f) =  \sum_{ \Phi\in |\hc_v|  } I^{\Phi}_a(f_\Phi).
 $$
\end{paragr}

\begin{paragr}\label{S:ISGp}
Soit $S$ un ensemble fini de places de $F$, suffisamment grand pour contenir les places archimédiennes et les places non-archimédiennes ramifiées dans $F$. On pose alors
$$\Ic(\hc_S)=\otimes_{v\in S}\Ic(\hc_v).
$$
On a une projection $f\mapsto f_\Phi \in \Ic(\Phi)$ pour tout $\Phi=(\Phi_v)_{v\in S}$
On étend l'application \eqref{eq:IvGp} en une application linéaire injective
$$I_S : \Ic(\hc_S) \to  \mathcal{C}^\infty(\tilde \Ac^{\rs}(\AAA_S)),
\quad I_{S}(f) = \left ( a \mapsto 
I_{a}^{\hc_{S}}(f)=\prod_{v} I_{a_v}^{\hc_v}(f_v)
\right ).
$$

On définit alors
$$\Ic(\hc)=\varinjlim_{S} \Ic(\hc_S)
$$
où les applications de transition $\Ic_S \to \Ic_{S'}$ sont données par $f \mapsto f\otimes (\otimes_{v\in S'\setminus S} \mathbf{1}_v)$.
\end{paragr}

\begin{paragr}[Distributions $I^{\hc}$.] --- On procède comme au §\ref{S:Ihc}.   Soit $\Phi \in \hc$ et $U$, $U'$ etc. les objets qui lui sont attachés. Soit $S_\Phi\subset \vc$ l'ensemble fini des places archimédiennes, des places ramifiées dans $E$ et des places tels que le discriminant de $(\Phi_v$ ne soit pas une norme de $E_v/F_v$.

 Soit   $a\in \tilde \Ac(F)$. D'après le théorème \ref{thm:densiteUU}, la  forme linéaire linéaire 
$$f\in \otimes_{v\in S}C_{c}^{\infty}(U'(F_v))\mapsto I_a^{\Phi}(f\otimes \mathbf{1}_{U'(\oc^S)} ),  $$
induite par la distribution $I_a^{\Phi}$ définie au  §\ref{S:IUGp},  se factorise par  $\otimes_{s\in S} \Ic(\Phi_v)$. 

Comme au §\ref{S:Ihc}, on en déduit  une forme linéaire
$$I_a^{\Phi}: \Ic(\hc) \to \CC.$$

\begin{lemme}\label{lem:SfGp}
  Soit $f\in \Ic(\hc)$. Il existe un ensemble fini $S_f \subset \vc$ tel que si $S_\Phi\not\subset S_f$, on a 
$$I_a^{\Phi}(f)=0$$
pour tout $a \in \tilde \Ac(F)$.
\end{lemme}

\begin{preuve}
La preuve est analogue à celle du lemme \ref{lem:Sf}
\end{preuve}

On définit pour tout $a\in \Ac(F)$ et $f\in \Ic(\hc)$
\begin{equation}
  \label{eq:IahcGp}
  I_a^{\hc}(f)= \sum_{\Phi\in |\hc|  } I_a^{\Phi}(f).
\end{equation}
À $f$ fixée, il n'y a qu'un nombre fini de termes non nuls (la preuve de ce fait est analogue à celle du  lemme \ref{lem:finitude}).
\end{paragr}

\subsection{Le théorème de transfert}

\begin{paragr}
  On poursuit avec les notations des deux sections  précédentes.
\end{paragr}

\begin{paragr}[Correspondance.] --- \label{S:correspGp}On dit que des éléments $f\in \Ic(\eta_v)$ et $f'\in \Ic(\hc_v)$ se correspondent si on l'égalité

\begin{equation}\label{eq:corr-IGp}
  I^{\Omega}_v(f)    = I^{}_v(f').
\end{equation}

Lorsque $v$ est fini ou scindé dans $E$, cette correspondance induit en fait un isomorphisme 
\begin{equation}
  \label{eq:isom-IcGp}
\Ic(\eta_v) \simeq  \Ic(\hc_v).
\end{equation}
Lorsque $v$ est scindé dans $E$, cet énoncé est aisé (il repose sur le théorème de factorisation de Dixmier-Malliavin, cf. \cite{DM}, si la place $v$ est archimédienne). C'est essentiellement l'énoncé du transfert non-archimédien si la place $v$ est finie et  non scindée (cf. \cite{Z1} théorème 2.6). Cet isomorphisme vaudrait également pour $v$  archimédienne mais non-scindée si le transfert archimédien était connu (pour des résultats partiels, cf. \cite{xue}). 
\end{paragr}

\begin{paragr}[Le lemme fondamental.] ---\label{S:LFGp}
  Pour $v\in \vc$ fini,  non-ramifié dans $E$ et également en dehors d'un ensemble fini fixé
  de \og mauvaises\fg{} places inertes, l'isomorphisme \eqref{eq:isom-IcGp} envoie l'élément
  $\mathbf{1}_v\in \Ic(\eta_v)$ sur l'élément  $\mathbf{1}_v\in  \Ic(\hc_v)$. Lorsque $v$ est décomposé, cet énoncé est élémentaire et pour $v$ inerte,   via le lemme \ref{lem:charHtoX}, il résulte des travaux de Yun et Gordon \cite{yun} lorsque la caractéristique résiduelle est assez grande.
\end{paragr}

\begin{paragr}
  Le lemme fondamental implique que les applications $I_S$ du §\ref{S:ISGp} induisent une application linéaire injective
$$\Orb: \Ic(\hc) \to \mathcal{C}^\infty(\Ac^{\rs}).$$
  On dit alors que deux éléments $f\in \Ic(\eta)$ et $f'\in \Ic(\hc)$ se correspondent si 
$$\Orb^{\Omega}(f)=\Orb(f').$$

On peut alors énoncer le théorème principal de cette section. 

\begin{theoreme}\label{thm:transfertGp}
  Soit $f\in \Ic(\eta)$ et $f'\in \Ic(\hc)$ qui se correspondent. On a alors
$$I_a^\eta(f)=  I_a^{\hc}(f')
$$
pour tout $a\in \tilde \Ac(F)$.
\end{theoreme}

\begin{preuve} On démontre le théorème par réduction au cas infinitésimal (cf. théorème \ref{thm:transfert} ou plutôt de sa légère généralisation aux situations des §§\ref{S:extToEnd} et \ref{S:ext-endII}).
Soit $ {\tla} \in \tilde \Ac(F)$ et $ \xi \in E^{1} $ tel que $ {\tla} \in \tilde \Ac_{\xi}(F) $. Soit $S\subset \vc$ assez grand. Pour tout $v\in S$ et tout $\Phi \in |\hc_v|$ soit $f_v\in \Cc(G_E'(F_v))$ et $f'_{\Phi_v} \in \Cc(U'(F_v))$. Soit $f_S=\otimes_{v\in S} f_v$, $f'_\Phi= \otimes_{v\in S} f'_{\Phi_v}$ pour tout $\Phi\in \hc_S$ et $f_S'=(f'_\Phi)_{\Phi\in \hc_S}$. Quitte à agrandir $S$, on peut et on va supposer que $f$ et $f'$ sont représentées par $ f_S\otimes \mathbf{1}^S$ et $   f'_S\otimes \mathbf{1}^S$. On a alors
$$I_{\tla}^\eta(f)=I_{\tla}^\eta(f_S\otimes \mathbf{1}_{G'_E(\oc^S)})$$
et
$$I_{\tla}^{\hc}(f')=\sum_{\Phi \in \hc^S  } I_{\tla}^{\Phi}( f'_{\Phi_S}\otimes \mathbf{1}_{U'_\Phi(\oc^S)}),$$
où
\begin{itemize}
\item $\hc^S$ désigne le groupoïde des formes hermitiennes de discriminant une norme hors $S$ 
\item $U'_\Phi(\oc^S)=\prod_{v\notin S} U'_\Phi(\oc_v)$ est associé au choix (sans importance) d'une base de $V_E$ et pour tout $v\notin S$ d'un réseau auto-dual $R_v$ de $V_E\otimes_F F_v$  tel que pour presque tout $v\notin S$, $R_v$ soit le réseau engendré par cette base ;  le  groupe $U'_\Phi(\oc_v)$ est le stabilisateur dans $U_\Phi'(F_v)$ du réseau $R_v\times (R_v\oplus \oc_{E\otimes_F F_v}e_0)$.
\end{itemize}

Quitte à multiplier chaque fonction par $\theta_S$ (vue comme fonction sur le groupe, cf. §\ref{S:thetaSaxi}), on peut supposer que $f_v\in \Cc(G_{E,\xi}'(F_v))$ et $f'_{\Phi_v} \in \Cc(U'_\xi(F_v))$ (par les propriétés de support des distributions cela ne change aucun des termes et cela n'affecte pas que ces fonctions se correspondent). Les constructions de \eqref{eq:phi-f} et \eqref{eq:phi-fU} (prises pour $S=\{v\}$) donnent des fonctions $\phi_{v}=\phi_{f_v}\in \Cc(\tggo(F_v))$ et $\phi_{\Phi_v}= \phi_{f'_{\Phi_v} }\in \Cc(\tugo_{\Phi_v}(F_v))$.
 Pour tout $\Phi\in \hc_S$ soit $\phi_\Phi=\otimes_{v\in S} \phi_{\Phi_v}$ et $\phi_S=(\phi_\Phi)_{\Phi\in \hc_S}$.

On met un exposant $\inf$ pour distinguer les constructions de la section \ref{sec:lecasinf} de celles utilisées pour les groupes. L'élément $(\otimes_{v\in S}\phi_v) \otimes \mathbf{1}^S$ définit un élément $\phi$ de $\Ic^{\inf}(\eta)$. De même,  $\phi_S\otimes  \mathbf{1}^S$ définit un élément $\phi'$ de $\Ic^{\inf}(\hc)$. Comme $f$ et $f'$ se correspondent, il résulte de  \eqref{eq:IOHggo} et \eqref{eq:IOUugo} que $\phi$ et $\phi'$ se correspondent. Soit $a\in \Ac_\tau(F)$ tel que $\kappa(a)=\tla$ (cf. lemme \ref{lem:kappaIsoG}). Quitte à élargir encore $S$, on a,  d'après les lemmes \ref{lem:IaHggo} et \ref{lem:IaUugo},  
$$I_a^{\eta,\inf}(\phi)=I_{\tla}^\eta(f)$$
et
$$I_a^{\hc,\inf}(\phi')=I_{\tla}^{\hc}(f').$$
L'énoncé cherché se déduit alors  du théorème \ref{thm:transfert}.
\end{preuve}
\end{paragr}
\appendix
\section{Les domaines convexes} \label{App:domConvP}

\begin{paragr}
Soit $G$ un $F$-groupe réductif, 
$M_{0}$ un sous-groupe de Levi de $G$ et 
$A_{0}$ son $F$-tore maximal central 
et déployé. Contrairement à la section \ref{ssec:notations-combi}, on ne suppose pas que $M_{0}$ est minimal. 
\end{paragr}

\begin{paragr}   
On appelle \textit{facette} de 
$\mathfrak{a}_{0} := \ago_{M_{0}}$ 
 tout ensemble de type:
\[
\mathfrak{a}_{Q}^{+} = 
\{H \in \mathfrak{a}_{Q} | \al(H) >0 \ \forall \al \in \Delta_{Q}\}, \quad Q \in \fc(M_{0}).
\]
Par dimension d'une facette on entend 
la dimension sur $\RR$ de plus petit sous-espace de 
$\mathfrak{a}_{0}$ la contenant. 
Les facettes de dimension maximale 
correspondent aux sous-groupes paraboliques 
$P \in \mathcal{P}(M_{0})$ et sont appelées \textit{chambres}. 
On identifie alors $\mathcal{P}(M_{0})$ avec l'ensemble de chambres de 
$\mathfrak{a}_{0}$.
On dit que $P,P' \in \mathcal{P}(M_{0})$ sont adjacents s'ils ont 
une facette de dimension $\dim_{\RR}\mathfrak{a}_{0}-1$ 
en commun, on désigne cette situation 
par $P \sim P'$. Une suite 
$\mathcal{G} = (P_{1},\ldots, P_{n})$ de groupes $P_{i} \in 
\mathcal{P}(M_{0})$ est appelée \textit{galerie} 
si $P_{i} \sim P_{i+1}$ pour 
$1 \le i \le n-1$. On dit alors que $\mathcal{G}$
relie $P_{1}$ et $P_{n}$ et on appelle $n$
sa longueur. 
Pour $P_{1},P_{2}\in \mathcal{P}(M_{0})$
on définit alors 
$d(P_{1},P_{2})$ - la distance entre $P_{1}$ et $P_{2}$
- comme le minimum de longueur de galeries 
reliant $P_{1}$ et $P_{2}$. 
On appelle galerie minimale entre $P_{1}$ et $P_{2}$
toute galerie reliant $P_{1}$ et $P_{2}$ 
de longueur $d(P_{1},P_{2})$. On dit aussi qu'une galerie 
$ (P_{1},\ldots, P_{n})$ est minimale si 
$n = d(P_{1},P_{n})$.

En suivant \cite{wallach}, on donne une autre 
définition de $d(\cdot,\cdot)$. 
Soit $P \in \mathcal{P}(M_{0})$, on note $\bar P$ l'unique 
élément de $\pc(M_{0})$ dont l'ensemble des 
racines simples égale $-\Delta_{P}$.
Notons aussi $\Sigma(P)$ (resp. $\Sigma(A_{0})$) 
l'ensemble des racines réduites de $A_{0}$ sur $N_{P}$ (resp. sur $G$).
On a alors 
$\Sigma(A_{0}) = \Sigma(P) \bigsqcup \Sigma(\bar P)$.
Pour $P_{1},P_{2} \in 
\mathcal{P}(M_{0})$ posons $\Sigma(P_{2}|P_{1}) = 
\Sigma(\bar P_{2})\cap \Sigma(P_{1})$.
 On a alors $d(P_{2},P_{1}) = |\Sigma(P_{2}|P_{1})|$. 

On va utiliser les deux résultats bien connus.

\begin{lemme}\label{lem:Wallach}
 Soient $P, P' \in \mathcal{P}(M_{0})$ et 
soit $(P_{1},\ldots, P_{d})$
une galerie minimale telle que 
$P_{1} = P$ et $P_{d} = P'$.
Dans ce cas, si $\Sigma(P_{i+1}|P_{i}) = \{\al_{i}\}$
pour $i = 1,\ldots, d-1$ alors 
$\Sigma(P'|P) = \{\al_{1},\ldots, \al_{d}\}$. 
\end{lemme}

\begin{lemme}\label{lem:keyCoxeter}
 Soient $P,P_{1},P_{2} \in \mathcal{P}(M_{0})$. 
Soit $d = d(P_{1},P)$ et supposons que 
$\Sigma(P_{2}|P_{1}) = \{\al\}$. Dans ce cas,
si $\al \in \Sigma(P|P_{1})$ alors $d(P_{2},P) = d-1$, 
sinon $d(P_{2},P)=d+1$.
\end{lemme}
\end{paragr}

\begin{paragr}[Famille convexe.] --- On dit que $\mathcal{S} \subseteq \mathcal{P}(M_{0})$ est une \textit{famille convexe} 
(de chambres) si pour tout $P,P' \in \Sc$ 
et pour toute galerie minimale 
$(P_{1},\ldots, P_{d})$ où $P_{1} = P$ et $P_{d} = P'$,
on a $P_{i} \in \mathcal{S}$ pour $i=1,\ldots, d$. 
Il est claire qu'une intersection des familles convexes 
est une famille convexe. L'ensemble vide, 
aussi bien que $\mathcal{P}(M_{0})$, est un exemple d'une 
famille convexe. Le plus simple exemple non-trivial 
d'une famille convexe est donné par l'ensemble $H(\al)^{+}$ 
qu'on va définir maintenant. 
Soit $\al \in \Sigma(A_{0})$ et notons 
$H(\al)^{+}$ 
l'ensemble de $P\in \mathcal{P}(M_{0})$ 
tels que $\al \in \Sigma(P)$. 

\begin{lemme}\label{lem:hfSpCnv}
 Soit $\al \in \Sigma(A_{0})$, 
alors $H(\al)^{+}$ est convexe. 
\end{lemme} 

\begin{preuve} Soient $P,P' \in H^{+}(\al)$ et 
soit $\mathcal{G} = (P_{1},\ldots, P_{d})$ 
une galerie minimale telle que $P_{1}=P, P_{d}=P'$. 
Raisonnons par absurde. 
Supposons qu'il 
existe $1 < j < d$ tel que $P_{j} \notin H(\al)^{+}$ et $P_{j+1} \in H(\al)^{+}$. 
Donc $-\al \in \Sigma(\overline P_{j+1}) \cap \Sigma(P_{j})$ 
ce qui implique que $\Sigma(P_{j+1}|P_{j}) = \{-\al\}$. 
En utilisant le lemme \ref{lem:Wallach} on voit que 
$-\al \in \Sigma(P'|P)$, d'où $\{-\al\} \in \Sigma(P)$. 
Mais $P \in H(\al)^{+}$ donc 
$-\al \notin \Sigma(P)$ ce qui n'est pas possible. 
\end{preuve}

\begin{lemme}\label{lem:langlType}
 Soit $\mathcal{S} \subseteq \mathcal{P}(M_{0})$ 
une famille convexe et $P \in \mathcal{S}$. 
Soit $Q \in \fc(M_{0})$ contenant un élément de 
$\mathcal{S}$. 
Alors, il existe un unique $P_{1} \in \mathcal{S}$ 
contenu dans $Q$ 
tel que $\Delta_{P_{1}}^{Q} \subseteq \Sigma(P)$.
\end{lemme}

\begin{preuve} Notons d'abord que l'unicité vient du fait qu'il existe 
un unique $P_{1} \in \mathcal{P}(M_{0})$ contenant $Q$ 
tel que $\Delta_{P_{1}}^{Q} \subseteq \Sigma(P)$, 
il s'agit alors de montrer que $P_{1} \in \mathcal{S}$. 
Soit $Q$ comme dans l'énoncé du lemme. 
Prenons un $P_{1} \in \mathcal{S}$ contenant $Q$ 
pour lequel $d(P_{1},P)$ soit minimal. 
On prétend que 
$\Delta_{P_{1}}^{Q} \subseteq \Sigma(P)$. Sinon soit 
$\al \in \Delta_{P_{1}}^{Q}$ tel que 
$\al \notin \Sigma(P)$ 
et soit $P_{2} \in \mathcal{P}(M_{0})$ adjacent 
à $P_{1}$ tel que $-\al \in \Delta_{P_{2}}$. 
On a alors $\al \in \Sigma(P|P_{1})$, 
donc $d(P_{2},P) < d(P_{1},P)$ d'après 
le lemme \ref{lem:keyCoxeter} ce qui veut 
dire que $P_{2}$ appartient \`a une galerie minimale 
liant $P$ et $P_{1}$. Par convexité de $\mathcal{S}$ 
on a $P_{2} \in \mathcal{S}$. 
Mais comme $\mathfrak{a}_{Q}^{+} \subseteq 
\overline{\mathfrak{a}_{P_{1}}^{+}} \cap 
\overline{\mathfrak{a}_{P_{2}}^{+}}$  (la barre désigne l'adhérence)
on a bien $P_{2} \subseteq Q$ ce qui 
contredit la minimalité de $P_{1}$. 
\end{preuve} 
\end{paragr}

\begin{paragr} Soit $\Lambda \in \mathfrak{a}_{0}^{*}$ 
et $P \in \mathcal{P}(M_{0})$. 
On définit 
$\epsilon_{P}(\Lambda)$ comme $(-1)$ \`a la 
puissance cardinal de l'ensemble 
$\{\al^{\vee} \in \Delta_{P}^{\vee}|\  \Lambda(\al^{\vee}) \le 0\}$.
Ensuite, soit $\phi_{P}(\Lambda,H)$ la 
fonction caractéristique 
de l'ensemble:
\[
\{H \in \mathfrak{a}_{0}| \varpi_{\al}(H) > 0 \text{ si } 
\Lambda(\al^{\vee}) \le 0, \ 
\varpi_{\al}(H) \le 0 \text{ si }\Lambda(\al^{\vee}) > 0, \ \forall 
\al \in \Delta_{P}\}
\]
où l'on note $\omega_{\al}$ l'unique 
élément de $\hat \Delta_{P}$ qui ne s'annule pas sur 
$\al^{\vee} \in \Delta_{P}^{\vee}$.

En suivant \cite{arthur4}, paragraphe 3, on dit qu'un ensemble 
$\mathcal{Y} = \{Y_{P}\}_{P \in \mathcal{P}(M_{0})}
\subseteq \mathfrak{a}_{0}$ 
de points indexes par $P \in \mathcal{P}(M_{0})$ est 
$A_{0}$-\textit{orthogonal positif} si 
pour tout $P_{1},P_{2}\in \mathcal{P}(M_{0})$ adjacents 
on a:
\begin{equation}\label{eq:orthPosit}
Y_{P_{1}} - Y_{P_{2}} = r\al^{\vee}, \quad r \ge 0, \ 
\Sigma(P_{2}|P_{1}) = \{\al\}.
\end{equation}
Si $\mathcal{Y} = \{Y_{P}\}_{P \in \mathcal{P}(M_{0})}$  
est un ensemble $A_{0}$-orthogonal positif et 
si $Q \supseteq P$ o\`u $P \in \mathcal{P}(M_{0})$ on notera 
$Y_{Q} \in \mathfrak{a}_{Q}$ 
pour projection de $Y_{P}$ sur $\mathfrak{a}_{Q}$. 
Alors $Y_{Q}$ ne dépend pas du choix de $P$ 
contenu dans $Q$.
Pour $P \in \mathcal{P}(M_{0})$ soit:
\[
(\mathfrak{a}_{P}^{*})^{+} = \{\Lambda \in (\mathfrak{a}_{0})^{*} | 
\Lambda(\al^{\vee}) > 0,  \ \forall \ \al \in \Delta_{P}\} 
\]
et pour $\mathcal{S} \subseteq \mathcal{P}(M_{0})$ notons 
$(\mathfrak{a}_{\mathcal{S}}^{*})^{+} = \bigcup_{P \in 
\mathcal{S}}(\mathfrak{a}_{P}^{*})^{+}$.
\end{paragr}

\begin{paragr}
  Soit $\mathcal{S} \subseteq \mathcal{P}(M_{0})$ et 
$\mathcal{Y} = \{Y_{P}\}_{P \in \mathcal{P}(M_{0})}$ 
un ensemble $A_{0}$-orthogonal positif, on définit alors:
\begin{gather}
\psi_{\mathcal{S}}(H,\mathcal{Y}) = 
\sum_{\begin{subarray}{c}
Q \in \fc(M_{0})\\
\exists P \in \Sc, \  P \subseteq Q
\end{subarray}
}
\varepsilon_{Q}^{G}\widehat{\tau}_{Q}(H-Y_{Q}), \ 
H \in \mathfrak{a}_{0}, \label{eq:psiAsSumHtau}
\\
\psi_{\mathcal{S}}(\Lambda, H, \mathcal{Y}) = 
\sum_{P \in \mathcal{S}}\epsilon_{P}(\Lambda)\phi_{P}(\Lambda, H-Y_{P}), \ 
\Lambda \in (\mathfrak{a}_{\mathcal{S}}^{*})^{+}, H \in \mathfrak{a}_{0}.
\notag
\end{gather}

\begin{lemme}\label{lem:noLanglands}
 Soient $\mathcal{S} \subseteq \mathcal{P}(M_{0})$ une 
famille convexe et 
$\mathcal{Y}$ un ensemble $A_{0}$-orthogonal positif. 
Alors, pour tout 
$\Lambda \in (\mathfrak{a}_{\mathcal{S}}^{*})^{+}$ 
et tout $H \in \mathfrak{a}_{0}$ l'on a:
\[
\psi_{\mathcal{S}}(H,\mathcal{Y}) = 
\psi_{\mathcal{S}}(\Lambda, H, \mathcal{Y}).
\]
\end{lemme}

\begin{preuve} 
Fixons un $P_{1} \in \mathcal{S}$ et un
$\Lambda \in (\mathfrak{a}_{P_{1}}^{+})^{*}$. 
Pour tout $P \in \mathcal{S}$ soit 
$Q_{P}^{1} \supseteq P$ 
défini par $\Delta_{P} \cap \Sigma(P_{1})$. En 
utilisant le lemme \ref{lem:langlType}, on 
s'aperçoit alors que $\psi(H,\mathcal{Y})$ 
égale:
\[
\sum_{P \in \mathcal{S}}
\sum_{P \subseteq Q \subseteq Q_{P}^{1}}
\varepsilon_{Q}^{G}
\widehat{\tau}_{Q}(H-Y_{P}).
\]

Fixons $P \in \mathcal{S}$, on montrera l'égalité:
\begin{equation}\label{eq:noLanglands}
\sum_{P \subseteq Q \subseteq Q_{P}^{1}}
\varepsilon_{Q}^{G}
\widehat{\tau}_{Q}(H-Y_{P}) = \epsilon_{P}(\Lambda)\phi_{P}(\Lambda, H-Y_{P}).
\end{equation}
En effet, soit $H$ tel que $\phi_{P}(\Lambda, H-Y_{P}) =1$, 
on a alors:
\[
\widehat{\tau}_{Q}(H-Y_{P}) = 
\begin{cases}
1 & \text{ si }Q=Q_{P}^{0} \\
0 & \text{ sinon} \\
\end{cases}, \ 
P \subseteq Q \subseteq Q_{P}^{0}.
\]
De plus \[ \varepsilon_{Q_{P}^{1}}^{G} = \epsilon_{P}(\Lambda).\]
Supposons par contre que $\phi_{P}(\Lambda, H-Y_{P}) =0$. 
Soit $Q' \supseteq P$ maximal tel que 
$\varpi_{\al}(H-Y_{P}) \le 0$ pour tout $\al \in \Delta_{P}^{Q'}$.  

Si $\Delta_{P}^{Q'} \cap (\Delta_{P} \smallsetminus
\Delta_{P}^{Q_{P}^{1}}) \neq \varnothing$ 
on a $\widehat{\tau}_{Q}(H-Y_{P}) =0$
pour tout $P \subseteq Q \subseteq Q_{P}^{1}$. 
Sinon $Q' \subsetneq Q_{P}^{1}$ et donc:
\[
\sum_{P \subseteq Q \subseteq Q_{P}^{1}}
\varepsilon_{Q}^{G}
\widehat{\tau}_{Q}(H-Y_{P}) = 
\sum_{Q'  \subseteq Q \subseteq Q_{P}^{1}}
\varepsilon_{Q}^{G} = 0
\]
ce qui démontre l'égalité \eqref{eq:noLanglands} et 
par conséquent le lemme \ref{lem:noLanglands}. 
\end{preuve}

\begin{corollaire} Avec les notations comme dans le lemme \ref{lem:noLanglands}, 
la fonction $\psi_{\mathcal{S}}(\Lambda, H, \mathcal{Y})$ 
est indépendante de $\Lambda \in (\mathfrak{a}_{\mathcal{S}}^{+})^{*}$.
\end{corollaire}

Le résultat suivant
est une généralisation du lemme 
3.2 de \cite{arthur4}. 

\begin{lemme}\label{lem:arthGenrlz} 
Soient $\mathcal{S} \subseteq \mathcal{P}(M_{0})$ une 
famille convexe et 
$\mathcal{Y}$ un ensemble $A_{0}$-orthogonal positif. 
Soit $H \in \mathfrak{a}_{0}$, alors 
les conditions suivantes sont équivalentes:
\begin{enumerate}
\item $\psi_{\Sc}(H,\mathcal{Y}) \neq 0$.
\item $\varpi(H-Y_{P}) \le 0, \ \forall P \in \mathcal{S}, 
\ \forall
\varpi \in \hat \Delta_{P}$.
\item $\psi_{\Sc}(H,\mathcal{Y})=1$.
\end{enumerate}
\end{lemme}

\begin{preuve} On suit de tout près la preuve du lemme 3.2
dans \textit{loc. cit. }
Le lemme étant trivial pour $\mathcal{S} = \varnothing$
on suppose $\mathcal{S}$ non-vide. 

\textit{i) $\Rightarrow $ ii)}. 
Supposons $\psi_{\Sc}(H,\mathcal{Y}) \neq 0$. 
Fixons un $P \in \mathcal{S}$ et 
soit $\Lambda \in (\mathfrak{a}_{P}^{*})^{+}$. 
D'après le lemme \ref{lem:noLanglands} il 
existe un $P' \in \mathcal{S}$ 
tel que 
\[
\Lambda(\al^{\vee})\varpi_{\al}(H-Y_{P'}) \le 0, \ \forall \al \in 
\Delta_{P'}.
\]
En sommant sur tout $\al \in \Delta_{P'}$ on obtient:
\[
\Lambda(H-Y_{P'}) \le 0.
\]
Il est facile de voir que $\Lambda(Y_{P'}-Y_{P}) \le 0$, 
on a alors:
\[
\Lambda(H-Y_{P}) = 
\Lambda(H-Y_{P'}) + \Lambda(Y_{P}-Y_{P'}) \le 0, \ \forall
\Lambda \in (\mathfrak{a}_{P}^{*})^{+}.
\]
La fonction 
$(\mathfrak{a}_{P}^{*})^{+} \ni \Lambda \mapsto \Lambda(H-Y_{P})$ étant continue, 
on obtient:
\[
\Lambda(H-Y_{P})  \le 0, \ 
\forall
\Lambda \in \overline{(\mathfrak{a}_{P}^{*})^{+}}.
\]
Donc, en particulier $\varpi(H-Y_{P}) \le 0$ 
pour tout $\varpi \in \hat \Delta_{P}$. 

\textit{ii) $ \Rightarrow $ iii)}. 
Fixons un $P \in \mathcal{S}$ et $\Lambda \in (\mathfrak{a}_{P}^{*})^{+}$. 
Si $H$ vérifie les conditions du point \textit{ii)} 
on a bien 
$\phi_{P}(\Lambda,H- Y_{P})=1$ ainsi que 
$\epsilon_{P}(\Lambda)=1$. 
Par contre, si 
$P'\in \mathcal{S} \smallsetminus \{P\}$ on a $\Lambda \notin 
(\mathfrak{a}_{P'}^{*})^{+}$ et donc 
$\phi_{P'}(\Lambda,H-Y_{P}) = 0$, d'o\`u \textit{iii)}. 

L'implication \textit{iii) $\Rightarrow$ i)} est évidente. 
\end{preuve} 
\end{paragr}

\bibliography{bib}
\bibliographystyle{alpha}

\begin{flushleft}
Pierre-Henri Chaudouard \\
Université Paris Diderot (Paris 7) et Institut Universitaire de France\\
 Institut de Mathématiques de Jussieu-Paris Rive Gauche \\
 UMR 7586 \\
 Bâtiment Sophie Germain \\
 Case 7012 \\
 F-75205 PARIS Cedex 13 \\
 France
\medskip

Adresse électronique :\\
Pierre-Henri.Chaudouard@imj-prg.fr \\

\medskip

Micha\l\  Zydor \\
Department of Mathematics \\
The Weizmann Institute of Science \\
Rehovot 76100 \\
Israel \\
\medskip

Adresse électronique :\\
michalz@weizmann.ac.il 

\end{flushleft}

\end{document}